\documentclass{amsart}
\usepackage[a4paper,margin=3.2cm]{geometry}
\usepackage[T1]{fontenc}
\usepackage{amssymb,amsmath,amsthm,bbm,setspace,xcolor,tikz, float, graphicx, mathdots,amsfonts,faktor,tikz-cd,stmaryrd,mathrsfs,dynkin-diagrams}
\usepackage[shortlabels]{enumitem}

\usepackage{comment}
\SetSymbolFont{stmry}{bold}{U}{stmry}{m}{n}

\makeatletter
\providecommand*{\xmapstofill@}{%
  \arrowfill@{\mapstochar\relbar}\relbar\rightarrow
}
\providecommand*{\xmapsto}[2][]{%
  \ext@arrow 0395\xmapstofill@{#1}{#2}%
}
\makeatother

\usepackage{subcaption}
\usetikzlibrary{positioning,patterns,nfold}
\usepackage{multirow}

\usepackage{tkz-euclide}

\usepackage{hyperref}
\hypersetup{
	colorlinks=true,
	linkcolor=blue,
	hypertexnames=false,
	filecolor=black,      
	urlcolor=blue,
	citecolor=blue
}
\usepackage[nameinlink]{cleveref}

\usepackage{drawingPCA}

\theoremstyle{plain}
\newtheorem{theorem}{Theorem}[section]
\newtheorem*{theorem*}{Theorem}
\newtheorem{proposition}[theorem]{Proposition}
\newtheorem{lemma}[theorem]{Lemma}
\newtheorem{corollary}[theorem]{Corollary}
\newtheorem{construction}[theorem]{Construction}

\theoremstyle{definition}
\newtheorem{definition}[theorem]{Definition}
\newtheorem*{convention}{Convention}
\theoremstyle{remark}
\newtheorem{remark}[theorem]{Remark}
\newtheorem{example}[theorem]{Example}
\newtheorem*{example*}{Example}
\newtheorem*{question*}{Question}

\newcommand{\R}{\mathbb{R}}

\newcommand{\Z}{\mathbb{Z}}
\newcommand{\C}{\mathbb{C}}
\newcommand{\Quat}{\mathbb{H}}
\newcommand{\Octo}{\mathbb{O}}
\newcommand{\Hyp}{\mathbb{H}}
\newcommand{\A}{\mathcal{A}}

\newcommand{\divisionAlg}{\mathbb{D}}
\newcommand{\SO}{\mathrm{SO}}

\newcommand{\Spin}{\mathrm{Spin}}
\newcommand{\GL}{\mathrm{GL}}
\newcommand{\SL}{\mathrm{SL}}

\newcommand{\SP}{\mathrm{Sp}}
\newcommand{\SU}{\mathrm{SU}}
\renewcommand{\sl}{\mathfrak{sl}}
\newcommand{\so}{\mathfrak{so}}
\renewcommand{\sp}{\mathfrak{sp}}
\newcommand{\su}{\mathfrak{su}}

\newcommand{\SLJ}{\SL_2(\jordan{J})}
\newcommand{\SLJi}{\SL_2(\jordan{J}_i)}
\newcommand{\slj}{\sl_2(\jordan{J})}

\newcommand{\CL}{\mathrm{Cl}}
\newcommand{\Mat}{\mathrm{Mat}}
\newcommand{\Id}{\mathbbm{1}}

\newcommand{\doverline}[1]{\overline{\overline{#1}}}
\newcommand{\opp}{\mathrm{opp}}
\newcommand{\Conf}{\mathrm{Conf}}
\newcommand{\Confx}{\mathrm{Conf}^\times}

\newcommand{\Hom}{\mathrm{Hom}}
\newcommand{\flag}{\mathscr{B}}
\newcommand{\decFlag}{\mathscr{A}}
\newcommand{\pot}{\delta}
\newcommand{\oppInv}[1]{{#1}^{*}}
\newcommand{\flagTheta}{\mathscr{P}_\Theta}
\newcommand{\decFlagTheta}{\decFlag_\Theta}
\newcommand{\WeylLength}{\ell}
\newcommand{\poscone}[1]{\overset{\circ}{c}\,(\mathfrak{u}_{\beta_{#1}})}

\newcommand{\tw}{\mathrm{tw}}
\newcommand{\puncTB}{T'S}
\newcommand{\ASpace}[2]{\mathcal{A}_{#1,#2}}
\newcommand{\Loc}[2]{\mathscr{A}_{#1,#2}}
\newcommand{\LocT}[3]{\mathscr{A}_{#1,#2}^{#3}}

\newcommand{\tr}{\mathrm{tr}}
\newcommand{\End}{\mathrm{End}}
\newcommand{\ed}{\mathrm{ed}}
\newcommand{\farey}{\mathcal{F}}
\newcommand{\pruned}{\mathrm{pr}}
\newcommand{\sau}{\sigma\hspace{-.8mm}\tau}

\newcommand{\Sym}{\mathrm{Sym}}
\newcommand{\splitt}{\mathrm{spl}}
\newcommand{\ad}{\mathrm{ad}}
\newcommand{\algclosure}[1]{#1^{\mathrm{alg}}}

\newcommand{\reduced}{\mathrm{red}}

\newcommand{\jordan}[1]{\mathbf{#1}}
\newcommand{\Inn}{\mathrm{Inn}}
\newcommand{\Aut}{\mathrm{Aut}}
\newcommand{\K}{\mathbb{K}}
\newcommand{\TitsCons}{\mathcal{T}}
\newcommand{\ThetaCons}[1]{\TitsCons(#1,\{\jordan{J}_i\})}
\newcommand{\ThetaConsJs}[2]{\TitsCons(#1,\{#2\})}

\newcommand{\In}{\mathrm{In}}
\newcommand{\Out}{\mathrm{Out}}

\newcommand{\keyword}[1]{\textbf{#1}}

\definecolor{Blue}{RGB}{0, 122, 255}
\definecolor{Green}{RGB}{12, 163, 2}
\definecolor{Orange}{RGB}{245, 154, 35}
\definecolor{DarkBlue}{RGB}{5, 70, 143}
\definecolor{Pink}{RGB}{255, 0, 162}
\definecolor{Red}{RGB}{209, 0, 0}
\definecolor{Purple}{RGB}{130, 33, 139}
\definecolor{BrightPurple}{RGB}{191, 85, 255}
\definecolor{Teal}{RGB}{51, 186, 150}
\definecolor{Yellow}{RGB}{255, 201, 0}

\newcommand{\st}{\vert}
\renewcommand{\tilde}{\widetilde}
\renewcommand{\hat}{\widehat}

\renewcommand{\vec}[1]{\mathbf{#1}}

\newcommand{\quiverAlgebra}{\mathcal{Q}}

\newcommand{\normalizer}[1]{\mathcal{N}(#1)}
\newcommand{\centralizer}[1]{\mathcal{C}(#1)}
\newcommand{\elemSeed}{E}
\newcommand{\elemNetwork}{E}

\makeatletter
\renewcommand\part{%
  \par
  \vspace*{4ex \@plus 1ex \@minus .2ex}
  \@startsection{part}{0}%
    \z@
    {0pt}
    {8ex}
    {\normalfont\huge\centering}
}
\makeatother

\title{Noncommutative Cluster Varieties and Moduli Spaces of Local Systems}

\author{Zachary Greenberg}
\address{Max Planck Institute for Mathematics in the Sciences\\
Inselstr. 22\\04103 Leipzig, 
Germany }
\email{greenberg@mis.mpg.de \url{https://zngzag42.github.io}}

\author{Dani Kaufman}
\address{Max Planck Institute for Mathematics in the Sciences\\
Inselstr. 22\\04103 Leipzig, 
Germany }
\email{kaufman@mis.mpg.de \url{https://sites.google.com/view/danikaufman/jordan}}

\author{Merik Niemeyer}
\address{Max Planck Institute for Mathematics in the Sciences\\
Inselstr. 22\\04103 Leipzig, 
Germany }
\email{niemeyer@mis.mpg.de}

\author{Anna Wienhard}
\address{Max Planck Institute for Mathematics in the Sciences\\
Inselstr. 22\\04103 Leipzig, 
Germany }
\email{wienhard@mis.mpg.de}

\date{}
\thanks{D.K. was supported by the Danish National Research Foundation (CPH-GEOTOP-DNRF151) and the Alexander von Humboldt Foundation. Z.G., M.N. and A.W. were supported by the European Research Council under ERC-Advanced Grant 101018839. A.W. thanks the Hector Fellow Academy for support. Part of this work was done at the Institute for Advanced Study in Princeton, and at the Erwin Schrödinger International Institute for Mathematics and Physics. We would also like to thank Max Riestenberg, Clarence Kineider, Kaitao Xie, Linhui Shen, and Alexander Goncharov for many helpful discussions. }

\begin{document}

    \begin{abstract}
In this article, we construct noncommutative cluster varieties, $\mathcal{A}_{R,S}$, for each reduced root system $R$ and marked surface $S$ simultaneously generalizing the cluster varieties of Fock-Goncharov, Li, Goncharov-Shen, Berenstein-Retakh, Goncharov-Kontsevich, and our previously introduced polygonal cluster algebras.
Additionally, we define a large class of algebraic groups, we call Jordan split groups. Given a reduced root system $R$ and a family of Jordan algebras, the Lie algebra for $G$ is constructed by unifying the Tits-Kantor-Koecher construction for a single Jordan algebra with the construction of a split Lie algebra. 

The notion of Jordan split groups is closely related to a grading of its Lie algebra  by the root system $R$. We show that these gradings are usually induced by a choice of standard parabolic subalgebra $\mathfrak{p}_\Theta$ and we classify $R$-graded pairs $(G,\Theta)$ via a condition depending only on the subset $\Theta\subset \Delta$ of the set of simple roots, providing many groups to which our results apply. 
Jordan split groups can be pinned using Jordan algebras, similarly to the pinnings of classically split groups. This allows us to develop a corresponding root system calculus. Using this, we compare the Weyl groups of $R$ and $G$, study properties of the $P_\Theta$-Bruhat decomposition, give a Lusztig style parameterization of the unipotent subgroup $U_\Theta$ based on words in the Weyl group $W(R)$, study decomposition relations between group elements, and define the noncommutative cluster varieties $\mathcal{A}_{R,S}$. 

We define Jordan algebra points of $\mathcal{A}_{R,S}$  which parameterize $G$-local systems on $S$ with boundary decoration related to cosets $G/U_\Theta$ when $G$ is Jordan split of type $R$. The cluster modular group $\Gamma_{R,S}$ considered by Goncharov-Shen, which in particular contains the mapping class group of $S$, acts on this noncommutative cluster variety by noncommutative quasi-cluster automorphisms. When $S$ is a disk, points of $\mathcal{A}_{R,S}$ parameterize configurations of decorated flags. We use this  to give noncommutative cluster structures on the double $R$-Bruhat cells of $G$, generalizing the cluster algebras of Berenstein-Fomin-Zelevinsky.

When each Jordan algebra is formally real, we say that $G$ has a positive structure with respect to $\Theta$. We define positive semigroups  $U_\Theta^{>0} \subset U_\Theta$ and $G^{>0}_\Theta \subset G$. For real algebraic groups, the pairs $(G,\Theta)$ which have positive structures are exactly those which admit a positive structure as defined by Guichard-Wienhard. 
As a consequence, our cluster structure gives explicit positivity tests, and can be used to give algebraic proofs of many of the properties of positive configurations of flags and of positive representations. 
    \end{abstract}
    
	\maketitle

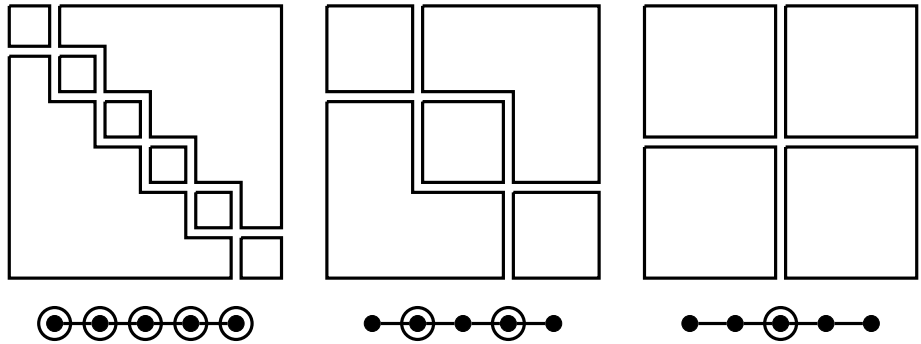
    \setcounter{tocdepth}{2}
    \makeatletter
    \def\l@subsection{\@tocline{2}{0pt}{2.5pc}{5pc}{}}
    \makeatother
\input{Figures/figureCoolIntro}
 \newpage 
	\tableofcontents
\newpage

\section{Introduction}
Split real semisimple Lie groups are the simplest among all real semisimple Lie groups. One of their key features is that they can be elegantly assembled from special three dimensional subgroups that arise from the $\sl(2)$-subalgebras associated to a set of positive simple roots. This leads to a characterization of split real Lie groups by their Dynkin diagrams, which in turn are classified by type $A_n, B_n,C_n, D_n, E_6, E_7, E_8, F_4,$ or $ G_2$. The property that split real Lie groups are built up from ${\SL}_2(\R)$'s plays a key role in their structure theory and underlies many interesting developments in recent years. Some examples are
\begin{enumerate}
\item Decomposition theorems generalizing the Gauss decomposition of matrices. 
\item The structure theory of totally positive semigroups by Lusztig \cite{lusztig1994total}, their cluster structure and explicit positivity criteria developed in \cite{berenstein2005cluster}.
\item Parameterizations of configuration spaces of flags and applications 
to cluster coordinates on higher Teichm\"uller spaces \cite{fock2006moduli, goncharov2019quantum}. 
\end{enumerate}
In this paper, we show that when we allow more general building blocks, many real semisimple Lie groups, and more generally semisimple algebraic groups, can be built up in a similar way. This allows us to draw interesting conclusions on their structure theory and give vast generalizations of the above mentioned applications to a broad family of real semisimple Lie groups. In particular, this includes all real semisimple Lie groups which admit a positive structure relative to a subset of simple roots $\Theta$ as introduced in \cite{guichard2022generalizing}. 

To motivate the basic building blocks we use, let us consider the simplest split real simple Lie group, $\SL_2(\R)$, which is the building block for split real semisimple groups. $\SL_2(\R)$ is itself built from more elementary pieces; two abelian additive groups, or ``root groups'', glued together along a multiplicative group, or a ``torus''. This manifest as the decomposition of a generic element $g$ in $\SL_2(\R)$ as:
$$g=\begin{bmatrix} a&b\\c&d \end{bmatrix}= \begin{bmatrix}
    1&0\\ca^{-1}&1
\end{bmatrix}\begin{bmatrix} a&0\\0&a^{-1} \end{bmatrix}\begin{bmatrix} 1&a^{-1}b\\0&1 \end{bmatrix}=[g]_-[g]_0[g]_+.$$ 
This is called the Gauss decomposition.
The matrices on the left and right sides belong to the root groups, while the matrix in the middle is in the torus. The root groups are isomorphic to the additive group $(\R,+)$, while the torus is the multiplicative group $(\R^\times, \cdot)$. The torus acts on the root groups by conjugation: 
$$\begin{bmatrix}
    a&0\\0&a^{-1}
\end{bmatrix}\begin{bmatrix} 1&b\\0&1 \end{bmatrix}\begin{bmatrix} a^{-1}&0\\0&a \end{bmatrix}= \begin{bmatrix} 1&a^2b\\0&1 \end{bmatrix} .$$ 
There is one essential symmetry of $\SL_2(\R)$ which swaps the two root groups; it is realized via conjugation by a special group element $s =\begin{bmatrix} 0&-1\\1&0 \end{bmatrix}$.  
If the entry $a$ is zero, then we cannot take the Gauss decomposition of $g$. However, in this case the element $gs$ admits a Gauss decomposition, a fact related to the Bruhat decomposition. 

The symmetry $s$ allows to define a map from an open subset of the root group to the torus: 
$$\iota\left(\begin{bmatrix} 1&b\\0&1 \end{bmatrix}\right) = \begin{bmatrix} b&0\\0&b^{-1} \end{bmatrix},$$
which is inspired by the calculation:
$$\begin{bmatrix} 1&b\\0&1 \end{bmatrix}\begin{bmatrix} 0&-1\\1&0 \end{bmatrix}= \begin{bmatrix} b&-1\\1&0 \end{bmatrix} = \begin{bmatrix}
    1&0\\b^{-1}&1
\end{bmatrix}\begin{bmatrix} b&0\\0&b^{-1} \end{bmatrix}\begin{bmatrix} 1&-b^{-1}\\0&1 \end{bmatrix}.$$ 
The map $\iota$ in fact provides a quadratic action of the root group on itself given by $\iota(a)(b)=a^{2}b$. 

This structure of $\SL_2(\R)$ can be generalized in a much broader context, since we used very few properties of $\R$. In particular we did not even use that $\R$ is commutative. The key structure is the map $\iota$ and the quadratic action it defines  of $\R$ on itself. 

There is a well known generalization of this structure called a \keyword{quadratic Jordan algebra}. 
A quadratic Jordan algebra is given by a triple $\jordan{J}=(V,\iota,\Id)$ of a vector space $V$, a map $\iota:V \to \End(V)$ and a special element $\Id\in V$ which is sent to the identity. The image of $\iota$ generates a subgroup of $\GL(V)$ called the \keyword{structure group} of $\jordan{J}$ denoted by $\Gamma(\jordan{J})$.
An example is the vector space of $k\times k$ matrices over a field $\K$ with $\Id$ being the identity matrix and $\iota$ given by $\iota(A)(B)= ABA$. In this case, the structure group is a subgroup of $\GL_k\times\GL_k$.
Given a quadratic Jordan algebra $\jordan{J}$, the Tits-Kantor-Koecher construction creates a Lie algebra $\slj$, from which we obtain an associated group $\SLJ$. This group can be thought of as two copies of the vector space $V$, replacing the root groups of $\SL_2$, glued along the (noncommutative) structure group, replacing the torus of $\SL_2$. 

These $\SLJ$'s are the basic building blocks we use. We call the groups which are built from these building blocks in the same way as split groups are built from $\SL_2$'s \keyword{Jordan split  groups}. 

We show that many semisimple algebraic groups can be realized as Jordan split groups, e.g. any semisimple algebraic group with a reduced root system is a Jordan split group. Viewing them this way provides interesting new structure on these groups and allows us to construct generalizations of the above mentioned features. 

In particular we obtain: 
\begin{enumerate}
\item Decomposition theorems generalizing the Gauss decomposition.
\item Non-commutative cluster structures on particular double Bruhat cells of $G$ generalizing those of \cite{berenstein2005cluster}.
\item A cluster structure and explicit positivity criteria for groups admitting a positive structure with respect to $\Theta$ as defined in \cite{guichard2022generalizing}. 
\item Parameterizations of configuration spaces of flags via non-commutative cluster structures, generalizing at the same time \cite{fock2006moduli,le2019cluster, goncharov2019quantum} as well as \cite{goncharov2021spectral}.
\item Non-commutative cluster parameterizations of moduli spaces of $G$ local systems. In particular, this gives a cluster parametrization of higher Teichm\"uller spaces for all groups admitting a positive structure with respect to $\Theta$ and a characterization of them as Jordan positive points. 
\end{enumerate}
In the following we describe the results of this paper in more detail, starting with the structure theory of Jordan split  groups $G$, describing the cluster structure and its applications for configurations spaces of flags, for the group $G$ and its special unipotent subgroup $U_\Theta$, and finally on the moduli spaces of decorated local systems on a surface $S$.

\subsection{Jordan split Lie groups and their structure theory} 
In the first part of the paper we develop the structure theory of semisimple algebraic groups defined over Jordan algebras. 

Given $G$, a semisimple algebraic group over a field $\K$, and $\Theta \subset \Delta$, a subset of the set of positive simple roots $\Delta$, we consider the associated parabolic subgroup $P_\Theta$, its Levi subgroup $L_\Theta$ and unipotent subgroup $U_\Theta$ along with their associated Lie algebras $\mathfrak{p}_\Theta,\mathfrak{l}_\Theta$ and $\mathfrak{u}_\Theta$. The action of the center, $\mathfrak{z}_\Theta$, of $\mathfrak{l}_\Theta$ on $\mathfrak{u}_\Theta$ decomposes $\mathfrak{u}_\Theta$ into weight spaces, giving rise to the parabolic root system associated to $\Theta$. 

We say that a pair $(G,\Theta)$ induces a \keyword{root system grading} when this parabolic root system is a reduced root system $R$, whose Dynkin diagram we denote by $\Delta_R$. The parabolic Weyl group $  W(R) = \normalizer{Z_\Theta}/\centralizer{Z_\Theta} $, where $Z_\Theta$ is the center of $L_\Theta$ is the Weyl group of this reduced root system. The Lie algebras with an $R$-grading have been classified in \cite{bm-lieAlgebrasGradedSimplyLaced, benkhart_lieAlgebrasGradedNonSimplyLaced}. 

If $(G,\Theta)$ induces a root system grading with $|\Theta|=1$ then the Lie algebra $\mathfrak{g}$ has an $A_1$ grading which implies that $U_\Theta$ is abelian. Lie algebras with an $A_1$-grading are exactly the Lie algebras $\slj$ which are built by applying the Tits-Kantor-Koecher construction to a Jordan algebra. In this situation, the vector space $\mathfrak{u}_\Theta$ is the vector space underlying the Jordan algebra $\jordan{J}$. The group $\SLJ$ is the simply connected algebraic group over $\K$ with Lie algebra $\slj$. 
The list of Dynkin diagrams of reduced restricted root systems and choices of $\Theta$ which give such an $A_1$ grading is given in \Cref{fig:possibleThetadynks}. 

Our first result is a classification of those $(G,\Theta)$ which induce root system gradings via a condition on the subset $\Theta$. To state the classification, we introduce the following notion. A subset $\Theta \subset \Delta$ is  \keyword{Jordan compatible} if for each root $\beta$ in $\Theta$ the connected subdiagram of $\Delta \setminus\Theta \cup \{\beta\}$ containing $\beta$ is one of the diagrams  associated to an $A_1$-grading  (see \Cref{fig:possibleThetadynks}).  This in fact associates a Jordan algebra to each element of $\Theta$. 

\begin{theorem}
A pair $(G, \Theta)$ induces a root system grading if and only if $\Theta$ is Jordan compatible. 
\end{theorem}
Note that in particular $(G ,\Delta)$ always induces a root system grading when $G$ has a reduced root system. In this case each root space is associated with a Jordan division algebra, i.e. a Jordan algebra where every nonzero element of $V$ is mapped by $\iota$ to an automorphism. If $\Delta_R$ is not equal to $\Delta$, then $\Delta_R$ has to be of type $A_p,B_p,C_p,F_4$ or $G_2$. In \Cref{thm:Identification}, we identify the possible classes of Jordan algebras that can appear in each type by studying the Jordan compatible subsets. For example, in type $A_p$ with $p>2$ each root corresponds to the same associative Jordan algebra. In type $B_p$, the short roots are assigned commutative Jordan algebras and the unique long root has a Jordan algebra of Clifford type. This is a reinterpretation of the results of \cite{bm-lieAlgebrasGradedSimplyLaced, benkhart_lieAlgebrasGradedNonSimplyLaced} in the classification of $R$ graded Lie algebras.

Next we reverse perspective, and initiate the study of groups $(G,\Theta)$, for $\Theta$ Jordan compatible, as ``split groups'' defined over a family of Jordan algebras.
For this we consider the root systems $R$ as the scaffold that gives a prescribed way to ``plug in'' a collection of Jordan algebras $\{\jordan{J}_i\}$ which recreate a Lie algebra in a similar fashion to the Tits-Kantor-Koecher construction. Each $R$ gives very strict rules for what collections $\{\jordan{J}_i\}$ can be used in this construction, as outlined in \Cref{thm:Identification}.

We call these groups \keyword{Jordan split} of type $R$ and we define a Jordan pinning as a collection of homomorphisms $\SLJi \to G$ such that the image of the vector space of $\jordan{J}_i$ lands in $U_{\beta_i}$. 
We investigate the algebraic structure of Jordan pinnings and establish the following theorem. 
\begin{theorem}
    The key algebraic relations among the images of Jordan pinnings depend only on the type of  $R$ and not the particular set of Jordan algebras $\{\jordan{J}_i\}$. 
\end{theorem}

We call these algebraic relations the \keyword{root system calculus} for any Jordan split group $G$ of type $R$.  In particular, we establish commutation and redecomposition relations among root systems and actions of the Weyl group both on the Levi and the root spaces.
This is key for later applications. 

The pinning also identifies the parabolic subgroup $P_\Theta$ with unipotent radical $U_\Theta$ as the subgroup generated by the images of the positive roots. Furthermore, we deduce properties of the Bruhat decomposition of $G$ with respect to $P_\Theta$, which mirror properties of the standard Bruhat decomposition with respect to the minimal parabolic subgroup. Here, the Weyl group $W(R)$ gives the organizing principle. For example we prove:  
\begin{theorem}
    For each word $w \in W(R)$, there is a unique cell in the $P_\Theta$-Bruhat decomposition of $G$, $P_\Theta w P_\Theta$. Moreover if $w_1,w_2$ are two words in $W(R)$ whose lengths add,
    \begin{equation*}
        P_\Theta w_1 P_\Theta w_2 P_\Theta = P_\Theta w_1w_2P_\Theta \,.
    \end{equation*}
\end{theorem}

We further explore the action of $W(R)$ on the full Levi subgroup $L_\Theta$. For example we prove:
\begin{theorem} Fix a Jordan pinning of $G$ by $\{\jordan{J}_i\}$. This induces coroot maps $\check{\beta}_i: \Gamma(\jordan{J}_i) \rightarrow L_\Theta \subset G$ generating $L_\Theta$ and lifts $\overline{\varsigma}_1 \cdots \overline{\varsigma}_p\in G$ of the simple reflections generating $W(R)$. When $\jordan{J}_i$ is commutative and $i$ is connected to $j$ in $\Delta_R$ the action of $W(R)$ on $L_\Theta$ is given by
    \begin{align*}
        \varsigma_i\left(\check{\beta}_j(A)\right) &\coloneq \overline{\varsigma}_i^{-1} \check{\beta}_j(\iota(A)) \overline{\varsigma}_i = \check{\beta}_i(N(A))\check{\beta}_j(A) \\
        \varsigma_j\left(\check{\beta}_i(a)\right) &\coloneq \overline{\varsigma}_j^{-1} \check{\beta}_i(\iota(a)) \overline{\varsigma}_j = \check{\beta}_i(a)\check{\beta}_j(a\Id_j) 
    \end{align*}
\end{theorem}
Here $N : \jordan{J}_j \rightarrow \K$ is a norm map defined from the $R$-grading which has degree $1,2$ or $3$ matching the connection between $i$ and $j$. This realizes the asymmetry of a non simply laced uniformly across the choice of $\jordan{J}_j$.

In the split case, the longest word in the Weyl group can be used to give a parametrization of an open and dense set in the unipotent radical of the Borel subgroup. This map was used prominently by Lusztig to give a parametrization of the totally positive unipotent semigroup. 

Here we give a Lusztig style parameterization of the unipotent subgroup $U_\Theta$ using the longest word in $W(R)$.
For $\vec{i}$ a reduced expression of the longest word in $W(R)$, the Lusztig map 
\begin{equation*}
    F_\vec{i} : \jordan{J}_{i_1}\times \cdots \jordan{J}_{i_N} \rightarrow U_\Theta
\end{equation*}
sends a list of invertible elements of the Jordan algebras onto an open dense subset of $U_\Theta$.

This map is inspired by the parametrization of the positive unipotent semigroup with respect to a subset $\Theta \subset \Delta$ in the generalization of Lusztig's total positivity introduced in \cite{guichard2022generalizing}.

Although the root system calculus and Jordan pinning only depend on the type of the root system, our later constructions will depend on the choice of Jordan algebras. We introduce the 
\keyword{noncommutative rank} of Jordan split groups $G$ of type $R$, which is the number of $\jordan{J}_i$ which are noncommutative Jordan algebras. 
The noncommutative rank is $0$ precisely the split or quasi-split groups which were previously studied. It turns out there is an important distinction between higher noncommutative rank (greater than or equal to 2) and noncommutative rank 1. Noncommutative rank 1 will share many structural properties with rank 0 and we refer to the combined cases as lower noncommutative rank.
We prove that non-quasi-split $G$ are noncommutative rank 1 if and only if $\Delta_R\simeq A_1, B_p ,\text{ or } G_2$ and higher non-commutative rank if and only if $\Delta_R\simeq A_p, C_p, \text{ or } F_4$.

The framework we develop gives a strong connection between groups of lower noncommutative rank and positive structures in Lie groups. 
When each $\jordan{J}_i$ is a formally real Jordan algebra we say that $G$ has a \keyword{positive structure with respect to $R$}. Given a Jordan pinning by formally real Jordan algebras, we define positive semigroups  
$U_\Theta^{>0} \subset U_\Theta$  as the image of the positive cones in each $\jordan{J}_i$ under a Lusztig map for any reduced expression of the longest word in $W(R)$.
When $\K=\R$ the pairs $(G,\Theta)$ which induce a positive structure on $G$ with respect to some root system $R$ are exactly those pairs $(G,\Theta)$ which admit a \keyword{positive structure with respect to $\Theta$}  as introduced in \cite{guichard2022generalizing}. The positive semigroups are the very same positive semigroups as the ones defined in \cite{guichard2022generalizing} and studied in \cite{guichard2026geometric}. Our approach here gives several new insights, which we discuss in more detail below. 

\subsection{Non-commutative cluster varieties and applications}
In the second part of the paper we use the structure theory of Jordan split groups of type $R$ and a marked surface $S$ to construct a noncommutative cluster variety $\ASpace{R}{S}$.  This noncommutative space is obtained by gluing together a collection of cluster charts. Each cluster chart consists of a noncommutative \emph{group} of functions called the \keyword{seed group} along with some extra combinatorial data. The seed group is a free group modulo some simple relations among a set of distinguished monomials called \keyword{angles}. The seed group takes the place the ring of functions on a classical torus used to define a classic cluster variety.

These cluster charts are glued together via mutations which are defined to have the property that the pullback of each angle is a (formal) sum of angles. These mutations are best understood via the \keyword{Jordan points} of the cluster variety $\ASpace{R}{S}$. The restriction on the collection $\{\jordan{J}_i\}$ imposed by $R$ identifies a single Jordan algebra $\jordan{J}$ which contains each $\jordan{J}_i$. Then, a Jordan point of a cluster chart is a map from the seed group to the structure group $\Gamma(\jordan{J})$ of this Jordan algebra so that every angle $A_i$ is mapped to $\iota(v_i)$  for some $v_i\in V$. Then the mutation of Jordan points is defined so that $\mu^*(A') = \iota(\sum v_k)$. In this way, a Jordan point consist of elements of the group $\Gamma(\jordan{J})$, but the addition appearing in the mutation happens in the vector space $V$.

Abelianinzing our construction recovers the $\mathcal{A}$-cluster varieties of Fock and Goncharov \cite{fock2006moduli} and of Goncharov and Shen \cite{goncharov2019quantum}. In the $A_p$ case this is the noncommutative cluster variety of Goncharov and Kontsevich \cite{goncharov2021spectral} with some minor modification. In the $B_p$ case this gives the noncommutative cluster variety structure conjectured in our previous work \cite{greenberg2024noncommutative}, while the $G_2$ case is a small variation of this. The $C_p$ case is entirely new, and we describe it by realizing it as a symmetric folding of the $A_{2p-1}$ case. The $F_4$ case and some exceptional cases are entirely new as well, and we do not explore all of their properties here. 

For any allowable choice of Jordan algebras $\{\jordan{J}_i\}$, $\ASpace{R}{S}$ has Jordan algebra points which parametrize $G$-local systems on $\Sigma$ which preserve decorated $\Theta$-flags, i.e. cosets in $G/U_\Theta$ at the marked points. 
Considering the non-commutative cluster varieties for particularly simple surfaces, give us noncommutative cluster structures on $G$ and the unipotent group $U_\Theta$. 
When $(G,\Theta)$ induces a positive structure, this cluster structure allows us to draw important conclusions. It can be used to give algebraic proofs of many of the properties of positive semigroups relative to $\Theta$, and of positive configurations of flags. Moreover it gives us explicit positivity test, generalizing the well known positivity test from split real Lie groups \cite{fock2006moduli,berenstein2005cluster}. 
Applied to the space of decorated $G$-local systems on $\Sigma$, the non-commutative cluster structures we define provide a unifying perspective  of the $\mathcal{A}$-theory of Fock and Goncharov for all groups admitting a positive structure, and thus for all (known) higher Teichm\"uller spaces associated to positive representations.

\subsubsection{Flag varieties, configuration spaces, and cluster coordinates} 
The noncommutative cluster varieties  $\ASpace{R}{S}$ associated to a Jordan-split type $R$ and a marked surface $S$ rely on the study of flag varieties and their configuration spaces as demonstrated in \cite{fock2006moduli} and \cite{goncharov2019quantum}.
Given a Jordan split group $G$ over a root system $R$ there is an associated unipotent subgroup, $U_\Theta$, corresponding to the positive roots of $R$. We consider the decorated flag variety $\decFlag:=G/U_\Theta$ and the configuration space
$\Confx_n(\decFlag)$ of pairwise transverse $n$-tuples of flags. The key point is to develop appropriate coordinate systems on $\Confx_n(\decFlag)$. 

 Many special cases of Jordan split groups have already been studied. For $G = \SL_n(\R)$ and $R$ of type $A_{n-1}$, Fock and Goncharov \cite{fock2006moduli} developed a set of coordinates for $\Confx_n(\decFlag)$. General split real Lie groups $G$ were studied in \cite{goncharov2019quantum}, where the authors describe a family of coordinate systems depending on a series of choices involving data from the root system. They then prove that these choices are a result of the cluster structure of $\Confx_n(\decFlag)$.

    The case of $G= \SP_{2n}(\R)$ with an $A_1$-grading given by the parabolic associated $\Theta = \{\alpha_n\}$, the unique long root, was studied in \cite{alessandrini2019noncommutative}. The case of $G = \SL_{kn}$ with an $A_{n-1}$-grading induced by $\Theta = \{k,2k,\cdots,(n-1)k\}$ was studied in \cite{goncharov2021spectral}. In both these cases, a generalized noncommutative cluster structure was developed to describe the coordinates. However, the generalizations appear very different, \cite{alessandrini2019noncommutative} are a realization of the noncommutative surfaces introduced by Berenstein and Retakh \cite{berenstein2018noncommutative}, while Goncharov and Kontsevich  \cite{goncharov2021spectral} used noncommutative plabic graphs. 

When the noncommutative rank is less than or equal to 1, the coordinates have the structure of a polygonal cluster algebra developed in \cite{greenberg2024noncommutative}. 
Both the usual cluster algebra and noncommutative surface cluster algebras are special cases of polygonal cluster algebras. We review and slightly extend the notion of polygonal cluster algebra developed in \cite{greenberg2024noncommutative} in \Cref{sec:PolygonalClusterAlgebras}. Polygonal cluster algebras come with a Weyl type, which can be $A_1$, $B_p$, or $G_2$, and a  parameter $r$ which can be any natural number for $A_1$,  has to be $r=2$ for $B_p$ and $r=3$ for $G_2$. 

 In higher noncommutative rank, a different generalization is needed which we call \keyword{grounded wiring networks}. These networks are introduced in \Cref{sec:NoncomNetworks} and they include the networks used in \cite{goncharov2021spectral}. 
We summarize our results in the following theorem. 

\begin{theorem}\label{thm:cluster_intro}
Let $\K$ be a field which is not of characteristic $2$ or $3$. 
Let $(G,\Theta)$ be a Jordan split semisimple algebraic group over $\K$ with $|\Delta_R|=:p$. Let $\Confx(\decFlagTheta)$ be the configuration space of decorated flags. Then there are cluster coordinates on $\Confx(\decFlagTheta)$, whose structure depend on the noncommutative rank of $(G,\Theta)$. 
\begin{enumerate}
  \item When $(G,\Theta)$ is of non-commutative rank 0, the cluster coordinates are obtained by evaluating the commutative cluster varieties of \cite{goncharov2019quantum} in field extensions of $\K$.
   \item When $(G,\Theta)$ is of non-commutative rank 1, the cluster coordinates are given by evaluations of non-commutative polygonal cluster algebras with general $r$ for $A_1$, $r=2$ for $B_p$ and $r=3$ for $G_2$ in Jordan algebras. 
   \item When $(G,\Theta)$ is of higher non-commutative rank and of type $A_p$ or $C_p$ the cluster coordinates are given by evaluations of cluster algebras of grounded wiring networks of type $A_p$ and $C_p$ in Jordan algebras. 
     \end{enumerate}
\end{theorem}

\subsubsection{Cluster structure on G, Gauss decompositions and positivity tests}
When $G$ is a split group and $\Theta=\Delta$, Berenstein, Fomin and Zelevinsky described a cluster structure on double Bruhat cells in \cite{berenstein2005cluster}. 
For two elements $w_1,w_2\in W$ of the Weyl group, the \keyword{double Bruhat cell} is given by
\begin{equation*}
    C^{w_1,w_2}:=P_\Delta w_1 P_\Delta\cap P_\Delta^\opp w_2 P_\Delta^\opp\,.
\end{equation*}
The group $G$ decomposes as a disjoint union of these cells. Moreover, if $w_0\in W$ is the longest element, the cell
\begin{equation*}
    C^{e,w_0}= P_\Delta\cap P_\Delta^\opp w_0 P_\Delta^\opp
\end{equation*}
is open and dense in $P_\Delta$. 
Thus, the cluster coordinates on $C^{e,w_0}$ describe an open dense subset of $P_\Delta$. As explained in \cite{le2019cluster,gilles2021fock} the cluster structure on $\Confx_3(\decFlag)$ is an extension of the cluster structure on $P_\Delta$, and we have
\begin{equation*}
    \Confx_3(\decFlag)\cong C^{e,w_0} \times H ,  
\end{equation*}
where $H=L_\Delta$ is a Cartan group.
Of course,  the same construction can be carried out for the cell $C^{w_0,e}$, which is open dense in $P_\Delta^\opp$. As $P_\Delta$, $P_\Delta^{opp}$ are semidirect products of $H$ and $U_\Delta$, resp. $U_\Delta^{opp}$, these correspondences can glued as in \Cref{fig:BorelTrianglesIntro}. 
\input{Figures/figureBorelTrianglesIntro}

Using cluster amalgamation, we obtain charts on an open dense subset $G_0$ of $G$ as indicated in \Cref{fig:BorelTrianglesIntro}: Namely, there is a map
\begin{equation*}
    G_0 = U_\Delta^{\opp}HU_\Delta \supset C^{w_0,e}C^{e,w_0} \to\Conf_4(\decFlag)
\end{equation*}

Thus the cluster coordinates on the configuration space $ \Confx_3(\decFlag)$ and on $\Confx_4(\decFlag)$ provide cluster coordinates on an open dense set of the unipotent group $U_\Delta$, the parabolic group $P_\Delta$ and the group $G$. 
The decomposition of the quadrilateral on the right hand side in \Cref{fig:BorelTrianglesIntro} corresponds to the Gauss decomposition. A  flip of the diagonal in this quadrilateral gives the opposite Gauss decomposition. Writing the flip as a sequence of cluster mutations thus provides an explicit way to compute the opposite Gauss decomposition. 

When $(G,\Theta)$ is Jordan split of type $R$, the cluster structure we develop allows us to give a generalization of this picture, where the relevant decompositions are now the Bruhat decompositions with respect to $P_\Theta$ instead of $P_\Delta$. 
\begin{theorem}
Let $(G,\Theta)$ be Jordan split of type $R$. Let $\decFlag = G/U_\Theta$ be the associated space of decorated flags. 
Then 
\begin{enumerate}
\item 
There are embeddings of open dense subsets $U_{\Theta,0}$ of 
$U_\Theta$ and $P_{\Theta,0}$ of $P_\Theta$ into $ \Confx_3(\decFlag)$ which induce noncommutative polygonal cluster structure on  $U_{\Theta,0}$ respectively on $P_{\Theta,0}$. An analogous statement holds for $U_\Theta^{\opp}$ and $P_\Theta^{\opp}$.
\item 
There is an embedding of an open dense subset $G_0 = U_\Theta^{opp}L_\Theta U_\Theta$ into 
$\Confx_4(\decFlag)$, which induces a noncommutative cluster structure on $G_0$, which can be depicted in a similar way as in \Cref{fig:BorelTrianglesIntro}. An analogous statement holds for $G^{\opp}_0 = U_\Theta L_\Theta U_\Theta^{\opp}$
\item 
The flip in the quadrilateral is a sequence of cluster mutations. This gives an explicit way to compute the change of coordinates then writing an element $g \in G_0 \cap G_0^{\opp}$ in the opposite Gauss decomposition. 
\item More generally for $(u,v) \in W(R)\times W(R)$ there is noncommutative cluster structure on the double $R$-Bruhat cell $P_\Theta u P_\Theta\cap P_\Theta^\opp v P_\Theta^\opp$.
\end{enumerate}
\end{theorem}

Note that also the multiplication of two elements in $G_0$ can be interpreted in terms of the cluster structure, see \Cref{fig:GroupProductIntro}. \input{Figures/figureGroupProductIntro}

When $(G, \Theta)$ admits a positive structure with respect to $R$, one can observe that these structures are compatible with the positive structure, which gives us the following theorems as corollaries. 

\begin{theorem}
The positive semigroup $U_\Theta^{>0}$ is contained in $U_{\Theta,0}$. The semigroup $U_\Theta^{>0}$ can be characterized as the set of  positive $\jordan{J}$-points of the noncommutative polygonal cluster structure on 
$U_{\Theta,0}$. Moreover 
\begin{enumerate}
\item Each cluster seed gives a precise set of positivity conditions, which are expressed as generalized minors of the group $G$ of type $R$ over the Jordan algebra $\jordan{J}$. 
\item There is a sequence of cluster mutations, which give an explicit way to compute the change of coordinates between different Lusztig parameterizations of $U_\Theta^{>0}$. 
\end{enumerate}
\end{theorem}
\begin{remark}
Note that getting explicit expression for the change of parameterizations without the cluster structure is quite hard in general. In \cite{guichard2026algebraic}, explicit polynomial type equations are derived following the approach of Berenstein-Zelevinsky \cite{Berenstein1997} using the universal enveloping algebra, and formulas for change of coordinates were obtained with a computer assisted proof. 
\end{remark}

\begin{theorem}
The positive semigroup $G_\Theta^{>0}$ can be characterized as the set of  positive $\jordan{J}$-points of the noncommutative cluster structure on 
$G_{0}$. Moreover 
\begin{enumerate}
\item Each cluster seed gives a precise set of positivity conditions, which are expressed as generalized minors of the group $G$ of type $R$ over the Jordan algebra $\jordan{J}$. 
\item There is a sequence of cluster mutations for the flip of the quadrilateral, this implies in particular that $G_\Theta^{>0} = U_\Theta^{>0} L^\circ_\Theta U_\Theta^{opp, >0}  = U_\Theta^{opp, >0} L^\circ_\Theta U_\Theta^{ >0} $, and gives an explicit way to compute the change of coordinates from one Gauss decomposition to the opposite one. 
\item This also gives an independent proof that the set $U_\Theta^{>0} L^\circ_\Theta U_\Theta^{opp, >0} $ is a semigroup. 
\end{enumerate}
\end{theorem}

 For the detailed formulas we refer the reader to Section~\ref{sec:NoncomRank1Examples}. The generalized minors that appear in the positivity conditions are precisely the same as for the split case, interpreted in the right way when looking at G as group of type $R$ over $\{\jordan{J}_i\}$. 

Whereas there is a lot of work on the nonnegative semigroups $U^{\geq 0}$ and $G^{\geq0}$ in the split case, much less is known for the nonnegative semigroups $U_\Theta^{\geq 0}$ and $G_\Theta^{\geq0}$ (see \cite{Wienhard2025} for some discussion). We expect the cluster parametrization to help understand the nonnegative semigroups.

\subsubsection{Cluster structure on moduli spaces of flat bundles and higher Teichmüller spaces} 

The cluster structure on the configuration spaces of $n$-tuples of flags, allow us to deduce a cluster variety structure on certain moduli spaces of local systems on the surface $S$ of negative Euler characteristic with marked points and punctures. This goes back to Penner's work on decorated Teichm\"uller space \cite{Penner}, but has been pioneered by Fock and Goncharov \cite{fock2006moduli} in the case of split real Lie groups. 
They introduced a pair of positive varieties, the $\mathcal{A}$-variety and the $\mathcal{X}$-variety, and showed that their positive points correspond to higher Teichm\"uller spaces. 
We adapt their strategy to the situation of Jordan split groups $G$, focusing on the $\mathcal{A}$-variety side of the story. We further explain that when $G$ carries a positive structure these local systems are related to positive representations \cite{guichard2021, beyrer2024}.

Given a surface $S$ with $k$-punctures and Jordan split group $G$ of type $R$ with associated parabolic subgroup $P_\Theta$, we consider the space of decorated (twisted) local systems on $S$, which depends on a representation $\rho: \pi_1(S) \rightarrow G$, with a decoration  by an affine flag $A_k$, i.e. an element of $G/U_\Theta$ for any puncture, which is fixed by the holonomy around this puncture. 

Building on the cluster structure on the space of configuration spaces, we build a cluster variety structure $\ASpace{R}{S}$ whose points parameterize the space of decorated local systems. 

\begin{theorem}
   The Jordan points of the cluster variety $\ASpace{R}{S}$ parametrize an open and dense set of the space of decorated local systems. 
   When $G$ of type $R$ carries a positive structure, the set of positive Jordan points, $\ASpace{R}{S}^{>0}$, of $\ASpace{R}{S}$ is independent of all choices and corresponds to the space of decorated positive representations of $\pi_1(S)$. 
\end{theorem}

The cluster structure on the space of decorated local systems has strong consequences. 
\begin{enumerate}
\item Given a curve $\gamma$ on $S$, the holonomy along $\gamma$ can be written explicitly in terms of the cluster coordinates. 
\item When $G$ of type $R$ is of non-commutative rank 1 and carries a positive structure, for any decorated local system in $\ASpace{R}{S}^{>0}$, the holonomy around any non-peripheral curve $\gamma$ is conjugate to an element in $G_\Theta^{>0}$.
\item Since the mapping class group action on triangulations of $S$ is generated by sequences of flips, the generalized mutations of the cluster structure gives explicit realizations of the action of the mapping class group. 
\end{enumerate}

Expanding on this last point, we further consider the subgroup of the cluster modular group defined  in \cite{goncharov2019quantum}.   This group, which we denote by $\Gamma_{R,S}$, depends on the marked surface and the Jordan split type $R$ of $G$; it is isomorphic to the group associated to the moduli space of decorated local systems for a split group of type $R$. It contains
\begin{enumerate}
    \item The mapping class group of $S$. 
    \item The outer automorphisms of $G$ which arise from automorphisms of the Dynkin diagram of $R$.
    \item For each puncture on $S$, a copy of the Weyl group $W(R)$.
    \item For each boundary component, a copy of a sub-quotient of the Artin-Tits Braid group of type $R$.
\end{enumerate}

In the commutative setting, the first three actions act by usual cluster automorphisms, i.e. maps which send cluster variables to cluster variables, while the last acts  by what are called cluster quasi-automorphisms, i.e. a cluster automorphism followed by a monomial transformation by frozen cluster variables. 

We define a noncommutative quasi-cluster automorphism to be a map which sends cluster variables to cluster variables followed by a noncommutative monomial transformation by frozen cluster variables. We show in \Cref{sec:ClusterModularGroupAction} the following:
\begin{theorem}
        The group $\Gamma_{R,S}$ acts on $\ASpace{R}{S}$ by noncommutative (quasi) cluster automorphisms.
\end{theorem}

Fixing a noncommutative ring $\mathcal{R}$, Goncharov and Kontsevich \cite{goncharov2021spectral} give a \keyword{spectral description} of the moduli space of fully twisted decorated rank $m$ local systems on $S$. They show that this space is birationally equivalent to the space of twisted decorated local rank 1 local systems on a more complicated surface called the \keyword{spectral surface} which is associated to a particular plabic graph embedded on $S$. 

The noncommutative cluster charts we construct for Jordan split groups of type $A_p$ are a natural generalization of those considered by Goncharov and Kontsevich, although there is a small difference in signs. This sign difference amounts to an untwisting of the local systems on the spectral surface. To this end we have:

\begin{theorem}
The space $\ASpace{R}{S}$ is birationally equivalent to the moduli space of untwisted twisted rank 1 local systems on the spectral surface.
\end{theorem}

When $R$ is $C_p$, there is a spectral description of the the space $\ASpace{R}{S}$ in terms of local systems which are invariant under a natural symmetry of the spectral surface.

\subsection{Structure of the paper}
In \Cref{part:RootGradedLieAlgebras}, we discuss the structure of Lie algebras and Lie groups with a root system grading. This begins with a review of the theory of quadratic Jordan algebras in \Cref{sec:jordan_algebras}. In \Cref{sec:RootGradedLieAlgebras} we classify which choices of $\Theta$ induce root system gradings and compare the Weyl group of the grading to the original Weyl group of the restricted root system. Finally in \Cref{sec:LieGroupsWithRootGradings}, we introduce the notion of a Jordan pinning which we use to study the $\Theta$-Bruhat decomposition and develop a ``root system calculus'' for computing in any group with a Jordan pinning.

In \Cref{part:NoncomClusterVarieties}, we provide the noncommutative generalizations of cluster varieties needed to parameterizes configurations of flags and moduli spaces of local systems. We begin in \Cref{sec:ClassicClusterAlgebra} with a review of classic cluster algebra theory. In \Cref{sec:PolygonalClusterAlgebras} we review the construction of polygonal cluster algebras we previously developed in \cite{greenberg2024noncommutative}. Finally in \Cref{sec:NoncomNetworks}, we introduce the notion of a grounded wiring network which provides the noncommutative cluster varieties in higher noncommutative rank. 

In \Cref{part:ClusterlikeCoordinates}, we describe the construction of noncommutative cluster varieties which parameterize configurations of flags. \Cref{sec:ClusterCoordinatesSetup} outlines the strategy from \cite{goncharov2019quantum} and describes the basic invariant of a pair of decorated flags, the $L$-distance. Following the general strategy, \Cref{sec:triplesToElementaryConfigs} describes how to break a configuration of a triple of flags into elementary configurations. Then in \Cref{sec:CoordsOnElementaryCofigurations}, we construct noncommutative cluster varieties for any elementary configuration. The construction depends on the noncommutative rank of the Jordan-split group, noncommutative rank less than 1 is discussed in \Cref{sec:elementarySeedsLowerNCrank} and higher noncommutative rank in \Cref{sec:elementaryNetworks}. In \Cref{sec:flagTriplesTheta}, we build the cluster varieties for triples of flags and in \Cref{sec:TuplesOfThetaFlags} we extend the construction to arbitrary tuples. In \Cref{sec:NoncomRank1Examples}, we carefully apply the theory in each noncommutative rank 1 example ($A_1, B_p, G_2$-gradings). \Cref{sec:HigherRankExamples} is the analogous section for higher noncommutative rank ($A_p$, $C_p$, $F_4$-gradings).

In \Cref{part:Applications}, we explore the natural applications of theory. First, we produce cluster-like coordinates on the algebraic groups themselves and their double Bruhat cells in \Cref{sec:coordinatesGroup}. Next, in \Cref{sec:ThetaPositivity}, we explore positivity and the natural semigroup structures that arise when all the Jordan algebras are Euclidean (formally real). In \Cref{sec:localSystems}, we describe how to extend the cluster varieties for tuples of flags to decorated twisted local systems. The positive points of these systems are closely related to the higher Teichm\"uller spaces of surfaces. 

\part{Root-graded Lie Algebras and Algebraic Groups}\label{part:RootGradedLieAlgebras}
In this part we develop the theory of Jordan split groups. First, we provide the necessary background on quadratic Jordan algebras.

\section{Quadratic Jordan algebras}\label{sec:jordan_algebras}
Let $\K$ be a field, not of characteristic 2 (usually we take $\R$ or $\C$).
\begin{definition}
    A \keyword{quadratic Jordan algebra}, $\jordan{J}=(V,\iota,\Id)$, consists of a vector space $V$ over $\K$, a quadratic map $\iota:V\to\End(V)$, and a special element $\Id\in V$ such that 
    \begin{enumerate}
        \item $\iota(\Id)$ is the identity transformation,
        \item\label{def:QuadraticJordanAlgebra_Composition} $\iota\big(\iota(a)(b)\big) = \iota(a)\iota(b)\iota(a)$,
    \end{enumerate}
    and defining $\iota(a,b)(c) \coloneq \big(\iota(a+c)-\iota(a)-\iota(c)\big)(b)$,
    \begin{enumerate}[resume]
        \item $\iota(x) \iota(y,x) = \iota(x,y) \iota(x)$.
    \end{enumerate}
\end{definition}

That $\iota$ is quadratic means by definition that the map $b_\iota:V\times V\to \End(V)$, defined by
\begin{equation*}
    b_\iota(v,w):=\iota(v+w)-\iota(v)-\iota(w)
\end{equation*}
is bilinear. This implies that $\iota(\cdot,\cdot)$ is a bilinear map (which actually maps to $\End(V)$), and thus induces a linear map $V\otimes V\to \End(V)$.

\begin{lemma}
    The map $\iota(\cdot,\cdot)$ determines the map $\iota$.
\end{lemma}
\begin{proof}
Using that $\K$ is not characteristic 2 we find
\begin{equation*}
    \iota(a)(c) = \frac{1}{2}\,\iota(a,c)(a)\,.
\end{equation*}
\end{proof}

\begin{remark}\label{rem:classical_jordan}
    A classical Jordan algebra over $\K$ is a vector space over $\K$ with a commutative but non associative product $V\circ V\to V$ satisfying some axioms.  
    The classical Jordan product is recovered from the quadratic definition as $a \circ b := \frac{1}{2}\iota(a,b)(\Id)$. On the other hand, the classical product determines the map $\iota$ by $$\iota(a)(b) = 2 a\circ (a\circ b)-(a\circ a)\circ b.$$ Thus the two theories are the same when the characteristic of the underlying field $\K$ is not 2.
\end{remark}

\begin{example}\label{ex:JordanHermitionType}
    Let $\mathcal{R}$ be a (possibly noncommutative) $\K$-algebra. We denote by $\jordan{J}(\mathcal{R})$ the Jordan algebra with vector space given by $\mathcal{R}$,  $\iota(a)(b)= aba$ and $\Id = 1_\mathcal{R}$. We call such a Jordan algebra \keyword{associative type}.
    
    If $\mathcal{R}$ is equipped with an anti-involution $\dagger:\mathcal{R}\to\mathcal{R}$ then $\mathcal{R}^\dagger$ denotes the subset of elements fixed by the anti-involution. There is a Jordan subalgebra $\jordan{H}(\mathcal{R},\dagger) \subset \jordan{J}(\mathcal{R}) $ with vector space $\mathcal{R}^\dagger \subset \mathcal{R}$, and $\iota(a)(b)=aba^\dagger$. We call these Jordan algebras \keyword{Hermitian type}.
\end{example}

\begin{example}\label{ex:JordanCliffordType}
    Let $\mathcal{R}$ be a \textit{commutative} $\K$-algebra and $M$ be an $\mathcal{R}$-module with a bilinear form $b:M\otimes M \to  \mathcal{R}$. Define $r_w(v) = \frac{1}{2}b(w,w)v- b(v,w)w$ and fix an element $e\in M$ with $b(e,e)=1$. We define a Jordan algebra $\jordan{J}(M,b,e)$ by $\iota(m)(n) = r_m(r_e(n))$. One easily checks that $e$ acts as the identity element and $\iota $ satisfies the axioms for a Jordan algebra. We call these Jordan algebras \keyword{Clifford type} as they are subalgebras of the associative Jordan algebra made from the Clifford algebra associated to $(M,b)$.
    \end{example}

\begin{example}
    Let $\mathbb{O}$ be an octonion algebra over $\K$, e.g. made from the Cayley-Dickson construction \cite{dickson-OnQuaternionsAndTheirGeneralizations}. Let $\jordan{H}_3(\mathbb{O})$ be the Jordan algebra with $V$ given by Hermitian $3\times  3$-matrices over $\mathbb{O}$, $\iota$ defined most easily by using the classical Jordan product $a\circ b =\frac{1}{2}(ab+ba)$ and \Cref{rem:classical_jordan}, and $\Id$ given by the identity matrix. Such a Jordan algebra is called an \keyword{Albert algebra} and is 27-dimensional over $\K$. 
\end{example}

\begin{definition}
    A Jordan algebra is called \keyword{special} if it is isomorphic to a subalgebra of $\jordan{J}(\mathcal{R})$ for some ring $\mathcal{R}$. We call a Jordan algebra \keyword{commutative} if it is special with $\mathcal{R}$ commutative and \keyword{noncommutative} otherwise. When $\jordan{J}$ is special there is a minimal specializing ring which contains it called the \keyword{universal specializing envelope}. 
\end{definition}

\noindent A proof of the following theorem is given in \cite[Appendix B]{mccrimmon-TasteJordanAlgebras}.
\begin{theorem}[Shirshov–Cohn]\label{thm:Shirshov-Cohn2GeneratorsSpecial}
    Any Jordan algebra with 2 generators is special. 
\end{theorem}

All Jordan algebras are special, other than the 27-dimensional exceptional Albert algebras. Over an algebraically closed field there is one exceptional Albert algebra. 

\subsection{The inner structure group}

\begin{definition}
    An element $a\in V$ is called \keyword{invertible} if its image $\iota(a)$ is invertible and there exists $b\in V$ such that $\iota(a)b = a$. We write $V^\times$ for the subset of invertible elements.  
    When all nonzero elements are invertible, the Jordan algebra is called a \keyword{Jordan division algebra}.
\end{definition}
The element $a^{-1}:=b$ is unique when it exists. 

\begin{definition}
  The \keyword{inner structure group} $\Inn(\jordan{J})$ is the subgroup of $\GL(V)$ generated by the image of $\iota(V^\times)$. The \keyword{inner automorphism group} $\Aut(\jordan{J})$ is the subgroup of $\Inn(\jordan{J})$ which fixes $\Id\in V$.  
\end{definition}

Generally, the map $\iota$ is clearly not injective since $\iota(a)=\iota(-a)$. However, we can construct a group $\Gamma(\jordan{J})$  which is a  covering $\Gamma(\jordan{J})\to\Inn(\jordan{J})$, so that there exists a map $\tilde{\iota}:V^\times \to \Gamma(\jordan{J})$ that extends $\iota$ and is injective, we will see this in \Cref{sec:A1gradedGroups}.

When $\jordan{J}$ is special, there is an action of the group $\mathcal{R}^\times\times\mathcal{R}^\times$ on $\mathcal{R}$ by $(a,b)\cdot r= arb^{-1}$. The elements $\{(z,z)\st z\in Z(\mathcal{R})\}$ act trivially. The inner structure group is the subgroup of $\mathcal{R}^\times\times\mathcal{R}^\times / Z(\mathcal{R})$ which is generated by the image of $\iota(v) = (v,v^{-1})$ for $v\in V^\times$.

\begin{definition}\label{def:structureInvolution}
    The \keyword{structure involution}, $\sigma:\Inn(\jordan{J}) \to \Inn(\jordan{J}) $, is the anti-involution defined on generators by $\sigma(\iota(v)) = \iota(v)$ so that if $g =\iota(v_1)\dots\tilde{\iota}(v_r) $ we have  $\sigma(g) = \tilde{\iota}(v_r)\dots\tilde{\iota}(v_1)$.
\end{definition}
One checks that this anti-involution is well-defined.

Notice that the structure involution satisfies that for $A\in \Inn(\jordan{J})$ we have $\iota(A(v))=A\iota(v)\sigma(A)$. For special algebras, the structure involution is induced by $\sigma(a,b)= (b^{-1},a^{-1})$.

\subsection{Real Jordan algebras}

The Jordan algebras over $\K=\R$ are classified using the classification of finite dimensional algebras over $\R$. We recall the classification here.

We denote by $\divisionAlg$ any of the rings $\R, \C$, the quaternions $\Quat$, the split complex numbers, $\C^\splitt$, or the split quaternions, $\Quat^\splitt$.
\begin{theorem}[\cite{mccrimmon-TasteJordanAlgebras}]
    A simple real Jordan algebra belongs (up to isotopy) to (at least) one of the following families:
    \begin{itemize}[topsep=0pt]
        \item (\keyword{Associative type}) $\jordan{J}$ is isomorphic to $\jordan{M}_n(\divisionAlg)$, the Jordan algebra of $n\times n$-matrices.
        \item (\keyword{Hermitian type}) $\jordan{J}$ is isomorphic to $\jordan{H}_n(\divisionAlg)$, the Jordan algebra of Hermitian matrices.\\
        Note that when $\divisionAlg = \Quat$ or $\Quat^\splitt$, there are two involutions which can be used to identify the set of Hermitian matrices, standard conjugation and $(a+b\textbf{i}+c\textbf{j}+d\textbf{k})^* = (a+b\textbf{i}+c\textbf{j}-d\textbf{k})$. We write $\jordan{H}_n^*(\Quat)$ for the Jordan algebra with the second form.
        \item (\keyword{Albert type}) $\jordan{J}$ can be isomorphic to $\jordan{H}_3(\Octo)$ or $\jordan{H}_3(\Octo^\splitt)$, Hermitian matrices over the octonions or split-octonions.
        \item (\keyword{Clifford type}) For $s>0$ and $r\geq 0$, $\jordan{J}$ can be isomorphic to $\jordan{J}(s,r)$ the Jordan algebra of Clifford type whose vector space has a signature $(s,r)$ quadratic form.
    \end{itemize}
\end{theorem}

These real Jordan algebras serve as key examples in the remainder of the paper, so it is critical to understand their structure. 

First, we observe the connection between the structure groups of Hermitian and associative type Jordan algebras. When there is no confusion, we will refer to $\jordan{M}_1(\divisionAlg)$ simply by $\divisionAlg$. For $\jordan{H}_n(\divisionAlg)$, the vector space is the set of $n\times n$ matrices over $\divisionAlg$ which are fixed by an involution $\dagger$. The structure group can be identified with $\GL_n(\divisionAlg)$ which acts on $V$ by $M \mapsto A M A^{\dagger}$. The structure involution is simply $\dagger$ in this case. 

The corresponding associative type algebra $\jordan{M}_n(\divisionAlg)$ has structure group $\mathrm{S}(\GL_n(\divisionAlg)\times \GL_n(\divisionAlg))$ $\subset \GL_n(\divisionAlg)\times \GL_n(\divisionAlg)$. To compare the structure group, we embed $\GL_n(\divisionAlg)$ in $\GL_n(\divisionAlg)\times \GL_n(\divisionAlg)$ via $A \mapsto (A,A^{-\dagger})$. This preserves the action of $\GL_n(\divisionAlg)$ on $\jordan{H}_n(\divisionAlg)$ when viewed as a sub-Jordan algebra of $\jordan{M}_n(\divisionAlg)$.

Next we observe that when $n \leq 2$, the $\jordan{H}_n(\Octo)$ and $\jordan{H}_n(\Octo^\splitt)$ are actually isomorphic to algebras in the other families:
\begin{align*}
    \jordan{H}_1(\Octo) \simeq \R, \hspace{2pc} \jordan{H}_1(\Octo^\splitt) \simeq \R, \hspace{2pc} \jordan{H}_2(\Octo)\simeq \jordan{J}(1,9), \hspace{2pc} \jordan{H}_2(\Octo^\splitt)\simeq \jordan{J}(5,5). 
\end{align*}
Therefore, these Jordan algebras are special. However, for $n=3$ the Albert type algebras are not special and their inner structure group is different from the other cases. Instead the inner structure groups are related to central extensions of the real Lie groups $E_6^{(-26)}$ or $E_6^{(6)}$ respectively. 

There are several other exceptional isomorphisms between small members of each family.
\begin{example}
    There are several examples involving the split division algebras, $\C^\splitt$ and $\Quat^\splitt$:
    \begin{align*}
        \jordan{H}_n(\C^\splitt) \simeq \jordan{M}_n(\R), \hspace{2pc} \jordan{H}_n^*(\Quat^\splitt) \simeq \jordan{H}_{2n}(\R), \hspace{2pc} \jordan{M}_{n}(\Quat^\splitt) \simeq \jordan{M}_{2n}(\R).
    \end{align*}
    The first isomorphism is given by sending a matrix $M \in \Mat_n(\R)$ to the Hermitian matrix $N$ with entries
    \begin{align*}
        N_{k\ell} = \begin{cases}
            \frac{1}{2}\left(M_{k\ell}(1-i) + M_{\ell k}(1+i)\right) & k < \ell\\
            M_{kk} & k=\ell\\
            \frac{1}{2}\left(M_{k\ell}(1+i) + M_{\ell k}(1-i)\right) & l > \ell.
        \end{cases}
    \end{align*}
    The last two isomorphisms follow from the fact that $\Quat^\splitt \simeq \Mat_2(\R)$. 
\end{example}

Finally we focus on the real Jordan algebras of Clifford type, $\jordan{J}(s,r)$. We consider a real vector space $V$ equipped with a quadratic form $q$ of signature $(s,r)$ with $s>0$. We write $q(v,w) = \frac{1}{2}(q(v+w)-q(v)-q(w))$ for the associated inner product. Let $e \in V$ such that $q(e) = 1$. Then $\jordan{J}(s,r)$ is defined to be $\jordan{J}(V,q,e)$ using the construction in \Cref{ex:JordanCliffordType}.

It is easier to understand $\jordan{J}(s,r)$ as a subalgebra of the even part of the Clifford algebra $\CL^0(s,r)$. Let $\sigma$ be the involution given by the Clifford algebra transpose followed by conjugation by $e$. Then the vector space $V$ sits inside the fixed points of $\sigma$ in $\CL^0(s,r)$ via $v \mapsto ve$. Using this construction, we see that the inner structure group is a quotient of the even Clifford group, $\Gamma^0(s,r)$.

There are many exceptional isomorphisms of Clifford type algebras with other families: 
    \begin{align*}
        \begin{aligned}
            \jordan{J}(2,0)\simeq~& \C\\
            \jordan{J}(1,2)\simeq~& \jordan{H}_2(\R)
        \end{aligned} \qquad
        \begin{aligned}
            \jordan{J}(1,1)\simeq~& \C^\splitt\\
            \jordan{J}(1,3)\simeq~& \jordan{H}_2(\C)
        \end{aligned} \qquad
        \begin{aligned}
            \jordan{J}(4,0)\simeq~& \Quat \\
            \jordan{J}(1,5)\simeq~& \jordan{H}_2(\Quat)
        \end{aligned}
        \qquad 
        \begin{aligned}
            \jordan{J}(2,2)\simeq~& \Quat^\splitt \\
            \jordan{J}(3,3)\simeq~& \jordan{H}_2(\Quat^\splitt)\,.
        \end{aligned}
    \end{align*}

\subsection{Isotopy and norm maps}

We will generally be interested in a slightly weaker notion of isomorphism of Jordan algebras called \keyword{isotopy}.

\begin{definition}
    Let $u$ be an invertible element of $\jordan{J}=(V,\iota,\Id)$. The \keyword{$u$-isotope} of $\jordan{J}$ is the Jordan algebra $\jordan{J}^{(u)}= (V,\iota^{(u)},u^{-1})$ where $\iota^{(u)}(x) =\iota(x)\iota(u)$. We write $\jordan{J}\sim \jordan{J'}$ if there is an isotopy $\jordan{J}$ to $\jordan{J'}$.
\end{definition}

Isotopy is different from isomorphism, which requires that the identity element is sent to the identity element. 
If $u$ has a square root, $u^{\frac{1}{2}}$, then the the element of the inner structure group $\iota(u^{-\frac{1}{2}})$ provides an isomorphism $\jordan{J} \to \jordan{J}^{(u)}$. However, if $u$ does not have a square root then there is no isomorphism between $\jordan{J}$ and its $u$-isotope. 

\begin{definition}
    The \keyword{twisted forms} of $\jordan{J}$ are the isotopy classes of $\jordan{J}$ modulo isomorphism.
\end{definition}

For a vector space $V$, let $\Sym^r(V^*)$ denote the $r^{\mathrm{th}}$ symmetric power of the dual vector space.
\begin{definition}\label{def:JordanNorm}
    We say a quadratic Jordan algebra has a \keyword{norm map of degree $r$} if it is equipped with a map $N \colon V \rightarrow \K$ which satisfies
    \begin{enumerate}
        \item $N(v)=\tilde{N}(v,\dots,v)$ for some $\tilde{N}\in\Sym^r(V^*)$,
        \item $N(1)=1$,
        \item $v \in V^\times$ if and only if $N(v) \neq 0$,
        \item $N(\iota(v)(w))=N(v)^2N(w)$\,.
    \end{enumerate}
\end{definition}
In \Cref{sec:A1gradedGroups}, we will see the norm map lifts to a character of a cover of $\Inn(\jordan{J})$ defined on generators by $N(\iota(v))=N(v)$. 

\begin{proposition}[\cite{mcrimmon-norms}]
    Every (finite dimensional) Jordan algebra is normed of degree $r$ for some $r$.
\end{proposition}
\begin{proof}
    In \cite{mcrimmon-norms}, they construct a \textit{generic norm} for any finite dimensional Jordan using a generic minimal polynomial. 
\end{proof}

Note that in the $u$-isotope, we have $N^{(u)}(x) = N(x)N(u)$.

\begin{definition}\label{def:adjugateJordanAlgebra}
    Every degree $r$ normed Jordan algebra has an \keyword{adjugate map}, $\tau : V^\times \rightarrow V^\times$ defined by  $\tau(v)= v^{-1}N(v)$.
\end{definition}

\noindent The adjugate gives an anti-involution on (a cover of) $\Inn(\jordan{J})$, defined on generators by  
\begin{equation*}
    \tau(\iota(v)) = N(v) \iota(v)^{-1}\,.
\end{equation*}

\begin{example}\leavevmode
    \begin{itemize}
        \item For $\jordan{J}=M_n(\R)$ a norm of degree $n$ is given by the determinant.
        \item If $\jordan{J}=\jordan{J}(s,r)$ is a real Jordan algebra of Clifford type, a norm of degree 2 is given by the quadratic form $q$. On the inner structure group $\Inn(\jordan{J})=\Gamma^0(s,r)/{\pm 1}$ the norm can be expressed as $v_1\cdots v_k\mapsto q(v_1)\cdots q(v_k)$, using that any element in the Clifford group can be written as a product of vectors. Using the defining relation of the Clifford algebra, $vv=q(v)$, the adjugate map is given by the Clifford transpose, $\tau(v_1\cdots v_k)=v_k\cdots v_1$.
    \end{itemize}
\end{example}

\begin{remark}
    The adjugate map can be extended to a map  on all of $V$ satisfying $\iota(v)(\tau v) = N(v) v$. We will not need to take the adjugate of non-invertible elements in which case the definition via the norm is sufficient.
\end{remark}

When $\K$ is formally real we can define a notion of a positivity on a Jordan algebra.
\begin{definition}\label{def:JordanAlgebraPositive}
    A \keyword{positive structure} on a quadratic Jordan algebra over a formally real field $\K$ is a subset $V^{>0}\subset V^\times$ called the \keyword{positive cone} which is compatible with the Jordan algebra structure as follows:
    \begin{enumerate}
        \item $V^{>0}$ is closed under addition and scalar multiplication by positive elements of $\K$.
        \item $V^{>0}$ contains 1 and is preserved by the action of $\Inn(\jordan{J})$.
        \item $V^{>0}\subset V^\times$.
        \item $V^{>0}$ is closed under inverse. 
    \end{enumerate}    
\end{definition}

\begin{definition}
    An \keyword{ideal} of a Jordan algebra $\jordan{J}$ is a subset $I \subset V$ which satisfies $\iota(V)(I) \subset I$.
    A Jordan algebra is called \keyword{simple} if it has no nontrivial ideals. 
\end{definition}
The simple finite dimensional real Jordan algebras with positive structure are also known as \keyword{Euclidean Jordan algebras}.
\begin{proposition}[Jordan-von Neumann-Wigner Classification]\label{thm:posJordanAlgsClassification}
     The finite dimensional simple real Jordan algebras with positive structure are given by the following list:  
    \begin{center}
        \begin{tabular}{c | c c c c c }
             Algebra& $H_n(\R)$ & $H_n(\C)$ & $H_n(\Quat)$ & $\jordan{J}(1,n)$ & $H_3(\Octo)$\\
             \hline
             Degree & $n$ & $n$ & $n$ & $2$ & $3$  
        \end{tabular}    
    \end{center}
    
\end{proposition}
\begin{proof}
    This is Fundamental Theorem 2 in \cite{jordan_EuclideanJordanAlgebra}. 
\end{proof}

The classification of simple real and complex Jordan algebras is well known. After extending scalars to $\C$, each real simple Jordan algebra becomes a simple complex Jordan algebra. We remark that each simple complex Jordan algebra has a unique real form (up to isotopy) which is Euclidean.

\subsection{Tits-Kantor-Koecher Lie algebra}
Let $\mathfrak{g}_\jordan{J}$ denote the Lie algebra of the  inner structure group $\Inn(\jordan{J})$. It is called the \keyword{inner structure Lie algebra} and is spanned by the elements $\iota(a,b)$ for $a,b\in V$. Since $\mathfrak{g}_\jordan{J}\subset\mathfrak{gl}(V)$, we can combine two copies of the vector space with the inner structure algebra $\mathfrak{g}_\jordan{J}$ to obtain a new $(-1,0,1)$ graded Lie algebra.
\begin{definition}[Tits-Kantor-Koecher Construction]\label{def:TKKconstruction}
    The \keyword{Tits-Kantor-Koecher-Lie algebra} $\slj$ associated to the Jordan algebra $\jordan{J}$ is the graded Lie algebra $V_{-1}\oplus \mathfrak{g}_\jordan{J}\oplus V_1$ with Lie brackets 
    $$ [A_{(0)},v_{(1)}] = A(v)_{(1)} \quad [A_{(0)},v_{(-1)}] = -\sigma(A)(v)_{(-1)} \quad [v_{(1)},w_{(-1)}] = \iota(v,w)_{(0)}.$$
\end{definition}

\begin{lemma}\label{thm:Simple_slj}
    If $\jordan{J}$ is special, then $\slj$ is a subalgebra of $\End_2(\mathcal{R})$.
\end{lemma} 
\begin{proof}
    We embed the generators of $\slj$ in  $\End_2(\mathcal{R})$ as follows:
    \begin{align*}
        v_{(1)} = \begin{bmatrix}
        0 & v \\ 0& 0
    \end{bmatrix}, \quad v_{(-1)} = \begin{bmatrix}
        0 & 0 \\ v& 0
    \end{bmatrix} \text{ and } \iota(a,b)_{(0)} = \begin{bmatrix}
        ab & 0 \\ 0 & -ba
    \end{bmatrix}\,.
    \end{align*}
    It remains to verify the bracket relations in $\End_2(\mathcal{R})$, for example
    \begin{align*}
        [\iota(a,b)_{(0)},c_{(1)}] = \begin{bmatrix}ab&0\\0 & -ba\end{bmatrix}\begin{bmatrix} 0&c \\0 & 0\end{bmatrix} - \begin{bmatrix}0&c\\ 0&0 \end{bmatrix}\begin{bmatrix}ab&0\\ 0& -ba\end{bmatrix} = \begin{bmatrix}0&abc+cba\\0 &0 \end{bmatrix} = \big(\iota(a,b)(c)\big)_{(1)}\,.
    \end{align*}
\end{proof}

\begin{lemma}[\cite{koecher-JordanAlgebrasInLieAlgebras} Theorem 3]\leavevmode
\begin{enumerate}
    \item  If $\jordan{J}$ is isotopic to $\jordan{J'}$ then $\slj\simeq \mathfrak{sl}_2(\jordan{J'})$.
    \item  If $\jordan{J}$ is simple, then $\slj$ is a simple Lie algebra.
\end{enumerate}
\end{lemma}

Our philosophy is that the Lie algebra $\slj$ encodes all of the properties of $\jordan{J}$ which are invariant under isotopy. In \Cref{sec:A1gradedLieAlgebras}, we will give a combinatorial classification for which Lie algebras can be realized as $\slj$ for some Jordan algebra.

\section{Lie algebras graded by root systems}\label{sec:RootGradedLieAlgebras}

We now describe an additional structure on a Lie algebra, called a grading by a root system, which allows us to realize the algebra as a ``split algebra'' over a smaller root system. The simplest case, an $A_1$ grading, is equivalent to realizing the algebra as $\slj$ for some Jordan algebra. 

\subsection{Preliminaries}
We first provide some background on abstract root systems, see \cite{knapp1996lie}. Let $E$ be a Euclidean space with inner product denoted $(\cdot,\cdot)$. We define for $v,w\in E$
\begin{equation*}
    \langle v,w\rangle := 2\frac{(v,w)}{(w,w)}\,.
\end{equation*}
Recall the definition of a  root system:
\begin{definition}\label{def:RootSystem}
    A finite subset $R\subset E\setminus\{0\}$ is called a \keyword{root system} if
    \begin{enumerate}
        \item $R$ spans $E$.
        \item For $\alpha\in R$, the \keyword{reflection} $s_\alpha$ defined by $s_\alpha(v)=v-\langle v,\alpha\rangle \alpha$, leaves $R$ invariant.
        \item For $\alpha,\beta\in R$, $\langle \beta,\alpha\rangle\in\Z$.
    \end{enumerate}
    We call $\dim(E)=:\mathrm{rk}(R)$ the \keyword{rank} of the root system $R$. The subgroup of $\GL(E)$ generated by the reflections $s_\alpha$ for $\alpha\in R$ is called the \keyword{Weyl group} of $R$. The root system $R$ is called \keyword{reduced} if the only scalar multiples of a root $\alpha$ which are contained in $R$ are $\alpha$ and $-\alpha$.
\end{definition}
A subset $\Delta = \{\alpha_1,\cdots, \alpha_p\} \subset R$ is called a set of \keyword{simple roots} if it is a basis of $E$ such that any element $\alpha\in R$ can be written as a linear combination with integer coefficients all of which are $\geq 0$ or all of which are $\leq 0$. It is a standard fact that a root system admits a set of simple roots. This defines sets of \keyword{positive} and \keyword{negative} roots, depending on the sign of the coefficients. We use the notation $[\alpha : \alpha_i]$ for the coefficient of $\alpha_i$ in the decomposition of $\alpha$ so that $\alpha = \sum [\alpha : \alpha_i] \alpha_i$.

\begin{definition}\label{def:highestRoot}
    Given a set of simple roots, every other root $\alpha$ has a well defined \keyword{height} given by $\sum [\alpha : \alpha_i]$. There is a unique root with maximum height called the \keyword{highest root}, $\alpha_M$. 
\end{definition}
The highest root is also the maximum element under the partial order given by $\alpha> \beta$ if and only if $\alpha-\beta$ is a nonnegative sum of simple roots. This further implies that for every simple root $\alpha_i$ and root $\alpha$, $[\alpha_M: \alpha_i]\geq [\alpha: \alpha_i]$. The simple roots are the minimal positive roots under this order.

\begin{lemma}\label{thm:RootSystemChains}
    For each simple root $i$, there is an (non-unique) ascending chain of roots $\alpha_i = \beta^1 ,\cdots, \beta^k = \alpha_M$ such that $\beta^{j}-\beta^{j-1}$ is a simple root.  Moreover if $\beta^0 = 0$, each simple root $\alpha$ appears $[\alpha_M : \alpha]$ many times as $\beta^{j}-\beta^{j-1}$. 
\end{lemma}
\begin{proof}
    The chain can be built inductively. By \cite[Proposition 4.9]{kac_infinite_1984}, the highest root is the unique root such that $\alpha_M + \alpha_i$ is not a root for any $\alpha_i$. So for each $\beta^{j}$ there is some simple root $\alpha$ so $\beta^{j+1} = \beta^j + \alpha$ is a root.
    
    We observe that at step $j$, the coefficient, $[\beta^{j} : \alpha_i]$ is the number of times $\alpha_i$ was used. Applying this to the final step, $\beta^k =\alpha_M$ proves the final statement. 
\end{proof}

The reflections in the simple roots  generate the Weyl group. The \keyword{length} of an element $w$ in the Weyl group is the minimal number of simple reflections required to write $w$. The \keyword{longest word} in the Weyl group, $w_0$ is the unique element with maximal length. It sends all positive roots to negative roots. The action of $-w_0$, fixes the set of simple roots and is called the \keyword{opposition involution}. 

Finally, if we fix a set of simple roots $\Delta=\{\alpha_1,\dots,\alpha_p\}$, we obtain a $p\times p$-matrix $C$ with entries
\begin{equation*}
    c_{ij}=\langle\alpha_i,\alpha_j\rangle 
\end{equation*}
which is called the \keyword{Cartan matrix} of $R$. This defines the \keyword{coroots} $\check{\alpha}_i\in E^*$ in the dual of $E$ via the condition
\begin{equation*}
    \check{\alpha}_i(\alpha_j)=c_{ij}\,.
\end{equation*}
The Weyl group acts on $E^*$ via the dual maps of the reflections $s_\alpha$, and we thus obtain the dual root system $\check{R}\subset E^*$. In particular, we obtain a coroot associated to any root, not just the simple ones.

Importantly, root systems are used to study Lie algebras: Let $G$ be a reductive group over $\K$ and $\mathfrak{g}$ its Lie algebra. Let $T$ be a maximal split $\K$-torus in $G$, $\mathfrak{t}$ its Lie algebra. The adjoint action of $\mathfrak{t}$, or equivalently $T$, decomposes $\mathfrak{g}$ into root spaces as
\begin{equation*}
    \mathfrak{g}=\mathfrak{g}_0\oplus\bigoplus_{\alpha\in R} \mathfrak{g}_\alpha
\end{equation*}
for some subset $R\subset \mathfrak{t}^*$, called the \keyword{restricted root system} of $(G,T)$. For any $t\in \mathfrak{t}$ and $x\in\mathfrak{g}_\alpha$, we have $[t,x]=\alpha(t)x$. The subset $R$ is an abstract root system in $\mathfrak{t}^*$ where the inner product comes from endowing $\mathfrak{t}$ with the Killing form.

The Weyl group of $(G,T)$ is defined as $W(G,T):=\mathcal{N}_G(T)/\mathcal{C}_G(T)$ where $\mathcal{N}_G(T)$ is the normalizer of $T$ in $G$, $\mathcal{C}_G(T)$ is the centralizer. This acts by conjugation on $T$, and via the induced action on $\mathfrak{t}$ by reflections on $R$. In fact $W(G,T)$ is isomorphic to the Weyl group of the abstract root system.

For the subsequent discussion it is useful to interpret the roots as characters of $T$, i.e. as elements of $\Hom(T,\K^\times)$, rather than as linear maps $\mathfrak{t}\to\K$. The latter arise simply as differentials of the former at the identity in $T$.

For each root $\alpha$ we write $T_\alpha\subset T$ for the kernel of the root $\alpha$, $G_\alpha:= \mathcal{C}_G(T_\alpha)$ for the centralizer of $T_\alpha$, and $S_\alpha := [G_\alpha,G_\alpha]$ is its derived subgroup. $S_\alpha\cap T$ is the image of the coroot $\check{\alpha}:\K^\times \to G$. Just like the roots, the coroots can be viewed in different ways: If we interpret the roots as elements of $\mathfrak{t}^*$, the coroots are elements of $\mathfrak{t}$ or equivalently of $\Hom(\K,\mathfrak{t})$. Viewing the roots as characters of $T$ instead, the coroots should be thought of as cocharacters, i.e. elements of $\Hom(\K^\times, T)$, which is the perspective we will usually take. If we pick a set of simple roots $\{a_1,\dots,\alpha_r\}$, the simple coroots are defined by
\begin{equation*}
    \alpha_i\circ\check{\alpha}_j=\left(k\mapsto k^{c_{ji}}\right)\in\Aut(\K^\times)\,.
\end{equation*}

Fix a choice of simple roots $\Delta$ and write $R^+$ for the set of positive roots. To each choice $\Theta\subset\Delta$ we associate a \keyword{standard parabolic subgroup} $P_\Theta$. It is the subgroup whose Lie algebra is given by
\begin{equation*}
    \mathfrak{p}_\Theta := \mathfrak{g}_0\oplus\bigoplus_{\alpha\in R^+} \mathfrak{g}_\alpha\bigoplus_{\alpha\in \mathrm{Span}(\Delta\setminus\Theta)\cap R^+}\mathfrak{g_{-\alpha}}\,.
\end{equation*}
The \keyword{opposite standard parabolic} $P_\Theta^\opp$ is obtained by switching the roles of the positive and negative roots. The \keyword{Levi subgroup}, $L_\Theta$, is defined to be $P_\Theta \cap P_\Theta^\opp$. The center of the Levi subgroup is denoted $Z_\Theta$. The \keyword{unipotent subgroups} $U_\Theta$ and $U_\Theta^\opp$ are the unipotent radicals of $P_\Theta$ and $P_\Theta^\opp$ respectively. We have $P_\Theta\simeq U_\Theta\rtimes L_\Theta$.

 The subgroups defined above are called standard; any other parabolic subgroup is conjugate in $G$ to some $P_\Theta$.

\begin{example}\label{ex:LusztigSLn}
    Consider the special linear group $G = \SL_n(\R)$. Its Lie algebra $\mathfrak{g}$ is the set of traceless $n\times n$-matrices. A maximal split torus is the set $T$ of diagonal matrices of determinant 1. The corresponding root system is of type $A_{n-1}$ with the Weyl group acting by permuting the entries of the elements of $T$. Let us fix the set $\Delta=\{\alpha_1,\dots,\alpha_{n-1}\}$ of simple roots such that $\mathfrak{g}_{\alpha_i}= \R E_{i,i+1}$ is spanned by the matrix with a $1$ in the $i^{\text{th}}$ entry of the upper off-diagonal. That means that the roots, written as characters,  are given by
    \begin{equation*}
        \alpha_i\big(\mathrm{diag}(t_1,\dots,t_n)\big)=\frac{t_i}{t_{i+1}}
    \end{equation*}
    with corresponding coroots
    \begin{equation*}
        \check{\alpha}_i(t):=\mathrm{diag}\big(1,\dots,1,\underset{i}{t},\underset{i+1}{t^{-1}},1,\dots,1\big)\,.
    \end{equation*}
    The standard parabolic subgroups which are specified by $\Theta\subset\Delta$ are block upper triangular with the size of the blocks determined by $\Theta$. For instance, when $n=4$, we have the following parabolic subgroups:
    \begin{align*}
    \begin{aligned}
         P_{\{\alpha_1\}}&=\left\{
            \begin{bsmallmatrix}
                * & * & * & *\\
                0 & * & * & *\\
                0 & * & * & *\\
                0 & * & * & *
            \end{bsmallmatrix}
        \right\}\,,\\
         P_{\{\alpha_1,\alpha_2\}}&=\left\{
            \begin{bsmallmatrix}
                * & * & * & *\\
                0 & * & * & *\\
                0 & 0 & * & *\\
                0 & 0 & * & *
            \end{bsmallmatrix}
        \right\}\,,
    \end{aligned}
       \qquad
       \begin{aligned}
            P_{\{\alpha_2\}}&=\left\{
            \begin{bsmallmatrix}
                * & * & * & *\\
                * & * & * & *\\
                0 & 0 & * & *\\
                0 & 0 & * & *
            \end{bsmallmatrix}
        \right\}\,,\\
        P_{\{\alpha_1,\alpha_3\}}&=\left\{
            \begin{bsmallmatrix}
                * & * & * & *\\
                0 & * & * & *\\
                0 & * & * & *\\
                0 & 0 & 0 & *
            \end{bsmallmatrix}
        \right\}\,,
       \end{aligned}
       \qquad
       \begin{aligned}
           P_{\{\alpha_3\}}&=\left\{
            \begin{bsmallmatrix}
                * & * & * & *\\
                * & * & * & *\\
                * & * & * & *\\
                0 & 0 & 0 & *
            \end{bsmallmatrix}
        \right\}\,,\\
         P_{\{\alpha_2,\alpha_3\}}&=\left\{
            \begin{bsmallmatrix}
                * & * & * & *\\
                * & * & * & *\\
                0 & 0 & * & *\\
                0 & 0 & 0 & *
            \end{bsmallmatrix}
        \right\}
       \end{aligned}       
    \end{align*}
    and also $P_\emptyset=G$ and the minimal parabolic $P_\Delta$, the set of upper triangular matrices.
\end{example}

Write $p=|\Theta|$. For each root $\alpha_i$ write $A_i$ for the set of positive roots for which have $[\alpha:\alpha_i]>0$  but $[\alpha:\alpha_j]=0$ for any other $\alpha_j\in \Theta$. The unipotent subgroup $U_\Theta$ is generated by abelian subgroups $U_{\beta_1}\dots,U_{\beta_p}$ Where $U_{\beta_i} := \exp(\mathfrak{u}_{\beta_i}) = \exp(\bigoplus_{\alpha\in A_i} \mathfrak{g}_{\alpha})$. These subgroups are abelian and the quotient $U_\Theta/[U_\Theta,U_\Theta]$ is isomorphic to $\prod_iU_{\beta_i}$.

In the next section we will understand when the groups $\mathfrak{u}_{\beta_i}$ behave like the simple root spaces of a root space decomposition of different root system than the one of $(G,T)$.

\subsection{Root Grading}
We have seen that the adjoint action of a maximal split torus on a Lie algebra decomposes it into root spaces. We now consider a similar decomposition, where we start with an abstract root system.
\begin{definition}
    Let $R$ be a reduced root system. A Lie algebra $\mathfrak{g}$ is \keyword{graded by $R$} if it has a decomposition
    \begin{equation*}
        \mathfrak{g} = \bigoplus_{\beta \in R \cup \{0\}} \mathfrak{g}_\beta
    \end{equation*}
    such that 
    \begin{enumerate}
        \item\label{prop:RootGradeGenerated} $\mathfrak{g}$ is generated as a Lie algebra by $\mathfrak{g}_\beta$ for $\beta \in R$. In particular every element of $\mathfrak{g}_0$ is given by a sum of iterated brackets of elements in the other graded summands.
        \item\label{prop:RootGradeSplit} The split Lie algebra $\mathfrak{g}^\splitt$ with root system $R$ is a subalgebra of $\mathfrak{g}$ such that there is a choice of a split Cartan subalgebra $\mathfrak{h} \subset \mathfrak{g}^\splitt$ such that $\mathfrak{g}^\splitt_\beta \subset \mathfrak{g}_\beta$ for all $\beta \in R\cup \{0\}$.
        \item\label{prop:RootGradeFullAction} The adjoint of action of $\mathfrak{h}$ on $\mathfrak{g}_\beta$ is the diagonal action with eigenvalue given by $\beta$. Symbolically, 
        \begin{equation*}
            \forall h\in \mathfrak{h}\subset \mathfrak{g}_0, \forall g \in \mathfrak{g}_\beta,~  [h,g] = \beta(h)\cdot g \,.
        \end{equation*}
    
    \end{enumerate}
\end{definition}

$R$-graded Lie algebras were introduced in \cite{bm-lieAlgebrasGradedSimplyLaced}, who classified all Lie algebras which can be graded by simply laced root systems of rank at least 2. Later, \cite{benkhart_lieAlgebrasGradedNonSimplyLaced} extended this classification to all reduced finite root systems and beyond.

Given a choice of Chevalley generators for $\mathfrak{g}^\splitt$, we obtain an element $\Id_\beta \in \mathfrak{g}_\beta$ with the following property: If $z\in \mathfrak{g}_\gamma$ then
\begin{equation*}
    \big[[\Id_{-\beta}, \Id_{\beta}], z\big] = 2 \frac{\langle \gamma, \beta\rangle}{\langle \beta,\beta\rangle} z\,,
\end{equation*}
where $\langle \cdot, \cdot \rangle$ is the inner product for the root system $R$. Conversely a choice of identity elements for each $\beta \in R$ along with $[\Id_{-\beta},\Id_{\beta}]$ generate a split subalgebra with root system $R$ \cite{neher_LieAlgebrasGraded}.

\begin{proposition}
    For any Jordan algebra $\jordan{J}$, the Tits-Kantor-Koecher Lie algebra, $\slj$ is $A_1$-graded.     
\end{proposition}
\begin{proof}
    By construction $\slj = V_{-1} \oplus \mathfrak{g}_{\jordan{J}} \oplus V_{1}$ already splits into three pieces indexed by the roots of an $A_1$-root system with $\mathfrak{g}_0 = \mathfrak{g}_\jordan{J}$. Property (\ref{prop:RootGradeGenerated}) holds since the inner structure algebra is spanned by elements of the form $\iota(v,w)_0 = [v_{(1)},w_{(-1)}]$. In order to verify Property (\ref{prop:RootGradeSplit}), we observe the elements $(\Id_{(1)},[\Id_{(1)},\Id_{(-1)}],\Id_{(-1)})$ form an $\sl_2$-triple. This gives an embedding of the split $A_1$-Lie algebra, $\sl_2(\K)$, into $\slj$. This identifies the rank 1 split Cartan subalgebra $\mathfrak{h}$ as $\K \cdot \iota(\Id,\Id)_0$. To verify Property (\ref{prop:RootGradeFullAction}), we check the action of $\mathfrak{h}$ on arbitrary elements $v_{(1)} \in \mathfrak{g}_{(1)}$ and $w_{(-1)} \in \mathfrak{g}_{(-1)}$ is by the $+1$ or $-1$ root. This calculation takes place in the two generator Jordan algebra generated by $v$ and $w$ which is special by \Cref{thm:Shirshov-Cohn2GeneratorsSpecial}. We compute
    \begin{align*}
        [\iota(\Id,\Id),v_{(1)}] &= \iota(\Id,\Id)(v)_{(1)} = \big( (v+\Id)^2 -\Id - v^2\big)_{(1)} =  2 v_{(1)}\,,\\
        [\iota(\Id,\Id),w_{(-1)}] &= - \sigma(\iota(\Id,\Id))(w)_{(-1)} = -2 w_{(-1)}\,.
    \end{align*}
\end{proof}

\noindent One source of a root system grading comes from a choice of a standard parabolic subalgebra $\mathfrak{p}_\Theta$.
\begin{definition}\label{def:ParabolicGrading}
    Let $\mathfrak{p}_\Theta$ be a parabolic subalgebra of $\mathfrak{g}$. Let $\mathfrak{z}_\Theta$ be the center of the associated Levi subalgebra $\mathfrak{l}_\Theta$. We say $\mathfrak{p}_\Theta$ \keyword{induces a root system grading} if the simultaneous eigen-decomposition of the action of $\mathfrak{t}_\Theta=\mathfrak{z}_\Theta\cap\mathfrak{t}$ is a root system grading.
\end{definition}

\begin{lemma}\label{thm:UnderstandingRootGrading}
    Let $\Theta \subset \Delta$ be a subset of simple roots. Then \begin{equation*}
        \mathfrak{t}_\Theta=\bigcap_{\alpha\in\Delta\setminus\Theta}\ker(\alpha)\,.
    \end{equation*}
    Moreover, the simultaneous eigenspaces of $\mathfrak{t}_\Theta$ are the union of root spaces whose roots agree when restricted to the span of $\Theta$. Equivalently the eigenspaces are indexed by the possible lists of coefficients of $([\beta: \alpha ])_{\alpha\in \Theta}$.
\end{lemma}
\begin{proof}
    The expression for $\mathfrak{t}_\Theta$ is shown in \cite[Section 2.9]{guichard2022generalizing}. This implies, for $t \in \mathfrak{t}_{\{\alpha\}}$ and $\beta$ any root, 
    \begin{equation*}
        [t,\mathfrak{g}_\beta] = \beta(t) \mathfrak{g}_\beta = \left(\sum_{\alpha \in \Delta}[\beta:\alpha] \alpha(t) \right) \mathfrak{g}_\beta = \left(\sum_{\alpha \in \Theta}[\beta:\alpha] \alpha(t) \right) \mathfrak{g}_\beta\,.
    \end{equation*} 
    Therefore any roots which agree on the span of $\Theta$ belong to the same eigenspace as claimed. 

\end{proof}

\begin{example}
    Let $\mathfrak{g} = \sl_4(\R)$ with restricted root system of type $A_3$. We choose the simple roots so that $P_\Delta$ is the set of upper triangular matrices. Let $\Theta = \{\alpha_2\}$ the middle simple root. We will see $\mathfrak{p}_\Theta$ induces an $A_1$ grading on $G$. We have
    \begin{align*}
        \mathfrak{l}_\Theta = \left\{\begin{bmatrix}
                    * & *& &\\
                    * & * & &\\
                    & & * & *\\
                    & & * & *
                    \end{bmatrix}\right\}
        \hspace{2pc}
        \mathfrak{t}_\Theta = \mathfrak{z}_\Theta = \left\{\begin{bmatrix}
            \lambda_1 & & & \\
            & \lambda_1 & & \\
            & & \lambda_2 & \\
            & & & \lambda_2
        \end{bmatrix}\right\}
    \end{align*}
    with $\lambda_1+\lambda_2 = 0$. Then action of $\mathfrak{z}_\Theta$ on $\mathfrak{g}$ has two nonzero weights $\beta= \lambda_1 - \lambda_2$ and $-\beta=\lambda_2 - \lambda_1$. The corresponding weight spaces are
    \begin{align*}
        \mathfrak{g}_\beta = \left\{\begin{bmatrix}
                     & & *& *\\
                     & & *& *\\
                     & &  & \\
                     & &  & 
                    \end{bmatrix}\right\}
        \hspace{2pc}
        \mathfrak{g}_{0} = \left\{\begin{bmatrix}
                    * & *& &\\
                    * & * & &\\
                    & & * & *\\
                    & & * & *
                    \end{bmatrix}\right\}
        \hspace{2pc}
        \mathfrak{g}_{-\beta} = \left\{\begin{bmatrix}
                     & & & \\
                     & & & \\
                     *& *&  & \\
                     *& *&  & 
                    \end{bmatrix}\right\}
    \end{align*}
    We check Property (\ref{prop:RootGradeSplit}) of an $A_1$-grading, that $\mathfrak{g}^\splitt \simeq \sl_2(\R)$ is contained as a subalgebra with
    \begin{align*}
        \begin{bmatrix}0 & 1\\0 & 0\end{bmatrix} \mapsto \begin{bmatrix} & & 0 & 1\\ & & 1 & 0\\ & & & \\& & & \end{bmatrix} \qquad
        \begin{bmatrix}1 & 0\\0 & -1\end{bmatrix} \mapsto \begin{bmatrix} 1& 0&  & \\ 0& 1& &\\ & &-1&0\\ & & 0& -1 \end{bmatrix}\qquad
        \begin{bmatrix}0 & 0\\1 & 0\end{bmatrix} \mapsto \begin{bmatrix} & &  & \\ & & &\\ 0& 1& &\\ 1& 0& &  \end{bmatrix}.
    \end{align*}
    The matrices above are a choice of Chevalley generators and so define the elements $\Id_\beta$ and $\Id_{-\beta}$. This also identifies $\mathfrak{h}$ with $\mathfrak{z}_\Theta$, verifying property (\ref{prop:RootGradeFullAction}) that $\mathfrak{h}$ acts correctly on all of $\mathfrak{g}_\beta$ not just $\mathfrak{g}^\splitt_\beta$.
    Property (\ref{prop:RootGradeGenerated}) that $\mathfrak{g}_0$ is generated by iterated brackets in $\mathfrak{g}_\beta$ and $\mathfrak{g}_{-\beta}$ is a linear algebra exercise. 
\end{example}

Even when a root system grading is not induced by a parabolic subalgebra, the root system grading identifies a parabolic subalgebra.
\begin{definition}\label{def:ParabolicAssociatedToGrading}
    Let $\mathfrak{g}$ be a simple Lie algebra graded by a root system $R$. We can chose the maximal split torus $\mathfrak{t}$ defining the restricted root system to be in $\mathfrak{g}_0$ so the restricted root system of $\mathfrak{g}$ refines the root system grading. Choosing a set of simple roots for $\mathfrak{g}$ induces a decomposition of $R$ into positive root $R_{>0}$ and $R_{<0}$.
    
    The \keyword{associated parabolic subalgebra}, $\mathfrak{p}_\Theta$, and \keyword{associated opposite parabolic subalgebra}, $\mathfrak{p}^\opp_\Theta$ are given by
    \begin{equation*}
        \mathfrak{p}_\Theta = \mathfrak{g}_0 \oplus \bigoplus_{\beta \in R_{>0}}\mathfrak{g}_\beta \hspace{2pc} \mathfrak{p}^\opp_\Theta = \mathfrak{g}_0 \oplus \bigoplus_{\beta \in R_{<0}}\mathfrak{g}_\beta\,.
    \end{equation*}
    These parabolic subalgebras are standard and so the index $\Theta$ can be interpreted as the associated subset of simple restricted roots of $\mathfrak{g}$. 

    The \keyword{associated Levi subalgebra}, $\mathfrak{l}_\Theta = \mathfrak{p}_\Theta \cap \mathfrak{p}_\Theta^\opp = \mathfrak{g}_0$. 
\end{definition}
When the root system grading on $\mathfrak{g}$ is induced by a parabolic subalgebra, $\mathfrak{p}_\Theta$, the associated parabolic subalgebra recovers $\mathfrak{p}_\Theta$. Thus there is no conflict in notation. 

\subsection{Recognition theorems}\label{sec:RecognitionTheorems}

In \cite{bm-lieAlgebrasGradedSimplyLaced,benkhart_lieAlgebrasGradedNonSimplyLaced} the authors proved a \keyword{recognition theorem} for each finite simple Weyl type. These theorems show that each Lie algebra with an $R$-grading is ``built'' out of copies of a fixed algebra $A$ organized by the root system $R$. As the Dynkin diagram for $R$ becomes more complicated, there are more restrictions on the algebras used to build it. 

\begin{theorem}\cite{benkhart_lieAlgebrasGradedNonSimplyLaced}\label{thm:RecognitionA1}
    Any Lie algebra with an $A_1$-grading is (a central extension of) $\slj$ for some quadratic Jordan algebra $\jordan{J}$.
\end{theorem}

This implies that every $A_1$-graded Lie algebra is essentially given by the Tits-Kantor-Koecher construction (\Cref{def:TKKconstruction}). In \Cref{tab:A1graded} we see all simple real Lie algebras with an $A_1$-grading, the corresponding Jordan algebra $\jordan{J}$ and the type of the restricted root system $\Delta$. 

\begin{table}[H]
    \centering
    \begin{tabular}{c||c|c|c|c|c|c|c}
      $\mathfrak{g}$   & $\sl_{2n}(\R)$ & $\sl_{2n}(\C)$ & $\sl_{2n}(\Quat)$ & $\sp_{2n}(\R)$ & $\su(n,n)$ & $\so^*(4n)$ & $\sp(2n,2n)$ \\ \hline
      $\jordan{J}$   &  $\jordan{M}_n(\R)$ & $\jordan{M}_n(\C)$ & $\jordan{M}_n(\Quat)$ &  $\jordan{H}_n(\R)$ & $\jordan{H}_n(\C)$& $\jordan{H}_n(\Quat)$  & $\jordan{H}_n^*(\Quat)$ \\ \hline 
      $\Delta$ & $A_{2n-1}$ & $A_{2n-1}$ & $A_{2n-1}$ & $C_n$ & $C_{n}$ & $C_{n}$ & $C_{n}$  \\ \hline\hline
      
      $\mathfrak{g}$ & \multicolumn{2}{|c|}{$\so(p+1,p+1)$} & \multicolumn{2}{|c|}{$\so(p+1,q+1)$} &  $\sl_{4n}(\R)$ & $\sl_{2n}(\R)$ & $\so(2n,2n)$ \\ \hline
      $\jordan{J}$   & \multicolumn{2}{|c|}{$\jordan{J}(p,p)$} & \multicolumn{2}{|c|}{$\jordan{J}(p,q)$} & $\jordan{M}_n(\Quat^\splitt)$ & $\jordan{H} _n(\C^\splitt)$ &  $\jordan{H}_n(\Quat^\splitt)$  \\ \hline
      $\Delta$       & \multicolumn{2}{|c|}{$D_{p+1}$} & \multicolumn{2}{|c|}{$B_{p+1}$}  & $A_{4n-1}$ & $A_{2n-1}$ & $D_{2n}$\\ \hline\hline
      
      $\mathfrak{g}$ &  \multicolumn{2}{|c|}{$\so(2p+1)^{(\C)}$} &  $\sp_{2n}(\C)$ & $\so(2p)^{(\C)}$ & $\mathfrak{e}_7^{(-25)}$ & $\mathfrak{e}_7^{(7)}$ & $\mathfrak{e}_7^{(\C)}$  \\ \hline
      $\jordan{J}$   & \multicolumn{2}{|c|}{$\jordan{J}(2p-1)$} & $\jordan{H}_n(\R)^{(\C)}$ & $\jordan{J}(2p-2)$ & $\jordan{H}_3(\Octo)$ & $\jordan{H}_3(\Octo^\splitt)$ & $\jordan{H}_3(\Octo^{(\C)})$  \\ \hline
      $\Delta$       & \multicolumn{2}{|c|}{$B_{p}$}  & $C_{p}$ & $D_{p}$ & $C_3$ & $E_7$ & $E_7$
    \end{tabular}
    \caption{Simple real Lie algebras with an $A_1$-grading.}
    \label{tab:A1graded}
\end{table}

\begin{theorem}[\cite{bm-lieAlgebrasGradedSimplyLaced}]
    Any Lie algebra with an $A_2$-grading is a (central extension of) $\sl_3(\divisionAlg)$ for some alternative unital algebra  $\divisionAlg$.
\end{theorem}

\noindent As an example of this theorem, we present \Cref{tab:A2graded} of all real Lie algebras with an $A_2$-grading.
\begin{table}[H]
    \centering
    \begin{tabular}{c|| c|c|c|c|c|c|c | c}
        $\mathfrak{g}$ & $\mathfrak{sl}_{3n}(\R)$ & $\mathfrak{sl}_{3n}(\C)$ & $\mathfrak{sl}_{3n}(\Quat)$& $\mathfrak{sl}_{3n}(\C^\splitt)$ & $\mathfrak{sl}_{3n}(\Quat^\splitt)$&$\mathfrak{e}_6^{(-26)}$&$\mathfrak{e}_6^{(6)}$ & $\mathfrak{e}_6^{(\C)}$\\ \hline
        $\divisionAlg$ & $\jordan{M}_n(\R)$ & $\jordan{M}_n(\C)$ & $\jordan{M}_n(\Quat)$ & $\jordan{M}_n(\C^\splitt)$ & $\jordan{M}_n(\Quat^\splitt)$ &  $\Octo$ & $\Octo^\splitt$ & $\Octo^{(\C)}$\\ \hline
        $\Delta$ & $A_{3k-1}$ & $A_{3k-1}$& $A_{3k-1}$& $A_{3k-1}\times A_{3k-1}$& $A_{6k-1}$& $A_{2}$ & $E_6$ & $E_6$ 
    \end{tabular}
    \caption{Simple real Lie algebras with an $A_2$-grading.}
    \label{tab:A2graded}
\end{table}

\begin{theorem}[\cite{bm-lieAlgebrasGradedSimplyLaced}]\label{thm:RecognitionAp}
    Any Lie algebra with an $A_p$-grading for $p > 2$ is (a central extension of) $\sl_{p+1}(\divisionAlg)$ for some associative unital algebra  $\divisionAlg$.
\end{theorem}

\begin{theorem}[\cite{bm-lieAlgebrasGradedSimplyLaced}]
    Any Lie algebra with a $D_p$, $E_6$, $E_7$, or $E_8$-grading (for $p > 4$) is (a central extension of) $\mathcal{G} \otimes A$, with $\mathcal{G}$ the split algebra of corresponding type over $\mathbb{K}$ and $A$ a commutative associate unital $\mathbb{K}$-algebra. 
\end{theorem}

The recognition theorems for gradings by root systems with Dynkin diagrams which are not simply laced are slightly more complicated, as they impose more interesting conditions on the underlying algebra. 

\begin{theorem}[\cite{benkhart_lieAlgebrasGradedNonSimplyLaced}]\label{thm:RecognitionCp}
    Any Lie algebra with a $C_p$-grading for $p \geq 4$ is (a central extension of) the algebra $\sp_{2p}(\mathcal{R},\sigma)$ where $\mathcal{R}$ is a unital associative algebra with an anti-involution $\sigma : \mathcal{R} \rightarrow \mathcal{R}$.

    When $p=3$, the assumption that $\mathcal{R}$ is associative can be reduced to the assumption that only the fixed points of $\sigma$ associate. 

    When $p=2$, any $C_2$ graded Lie algebra is (a central extension of) to $\sp_2(\jordan{J})$ for $\jordan{J}$ a quadratic Jordan algebra whose identity element is contained in a sub-Jordan algebra isomorphic to $\jordan{H}_2(\R)$. 
\end{theorem}
\noindent Once again, we observe \Cref{tab:C3graded} of all simple real Lie algebras with a $C_3$-grading.
\begin{table}[H]
    \centering
    \begin{tabular}{c||c | c |c |c |c |c |c }
         $\mathfrak{g}$ & $\sp_{6k}(\R)$ & $\su(3k,3k)$ & $\so^*(12k)$ & $\so(6k,6k)$ & $\sp(6k,6k)$ & $E_7^{(-25)}$ & $E_7^{(7)}$ \\ \hline
         $\mathcal{R}$  & $\jordan{M}_k(\R)$ & $\jordan{M}_k(\C)$ & $\jordan{M}_k(\Quat)$ &  $\jordan{M}_k(\Quat^\splitt)$ &  $\jordan{M}^{*}_k(\Quat)$ & $\Octo$ & $\Octo^\splitt$\\ \hline
         $\mathcal{R}^\sigma$ & $\jordan{H}_k(\R)$ & $\jordan{H}_k(\C)$ & $\jordan{H}_k(\Quat)$ & $\jordan{H}_k(\Quat^\splitt)$ & $\jordan{H}^{*}_k(\Quat)$ & $\R$ & $ \R$\\  \hline
         $\Delta$       & $C_{3k}$ & $C_{3k}$ & $C_{3k}$ & $D_{6k}$ & $C_{6k}$ & $C_3$ & $E_7$
    \end{tabular}
    \caption{Simple real Lie algebras with a $C_3$-grading.}
    \label{tab:C3graded}
\end{table}
For $p\geq 4$, the real simple Lie algebras with a $C_p$-grading are given in \Cref{tab:Cpgraded}.
\begin{table}[H]
    \centering
    \begin{tabular}{c||c | c |c |c |c  }
         $\mathfrak{g}$ & $\sp_{2pk}(\R)$ & $\su(pk,pk)$ & $\so^*(4pk)$ & $\so(2pk,2pk)$ & $\sp(2pk,2pk)$ \\ \hline
         $\mathcal{R}$  & $\jordan{M}_k(\R)$ & $\jordan{M}_k(\C)$ & $\jordan{M}_k(\Quat)$ &  $\jordan{M}_k(\Quat^\splitt)$ &  $\jordan{M}^{*}_k(\Quat)$ \\ \hline
         $\mathcal{R}^\sigma$ & $\jordan{H}_k(\R)$ & $\jordan{H}_k(\C)$ & $\jordan{H}_k(\Quat)$ & $\jordan{H}_k(\Quat^\splitt)$ & $\jordan{H}^{*}_k(\Quat)$ \\  \hline
         $\Delta$       & $C_{pk}$ & $C_{pk}$ & $C_{pk}$ & $D_{2pk}$ & $C_{pk}$
    \end{tabular}
    \caption{Simple real Lie algebras with a $C_p$-grading, $p\geq 4$.}
    \label{tab:Cpgraded}
\end{table}

In order to recognize the Lie algebras for types $B_p$, $F_4$, or $G_2$, Benkart and Zelmanov constructed a Lie algebra using what they call the \keyword{generalized Tits construction} \cite{benkhart_lieAlgebrasGradedNonSimplyLaced}. This construction produces a Lie algebra, $\TitsCons(\mathcal{R}, \jordan{J})$, from an alternative algebra $\mathcal{R}$ and a Jordan algebra $\jordan{J}$. The construction is inspired by a construction of the exceptional Lie algebra $e_8$ by Tits as $e_8^{(8)} = \TitsCons\big(\Octo^\splitt, \jordan{H}_3(\Octo^\splitt)\big)$.   

\begin{theorem}[\cite{benkhart_lieAlgebrasGradedNonSimplyLaced}]
    Any Lie algebra with a $G_2$-grading is (a central extension of) the algebra $\TitsCons(\Octo^\splitt, \jordan{J})$ where $\jordan{J}$ is a degree 3 Jordan algebra.
\end{theorem}
\noindent The simple real Lie algebras with a $G_2$-grading are given in \Cref{tab:G2graded}.

\begin{table}[H]
    \centering
    \begin{tabular}{c||c|c|c|c|c}
        $\mathfrak{g}$ & $f_4^{(4)}$ & $e_6^{(2)}$ & $e_6^{(6)}$ & $e_7^{(-5)}$ & $e_7^{(7)}$ \\ \hline
        $\jordan{J}$ & $\jordan{H}_3(\R)$ & $\jordan{H}_3(\C)$ & $\jordan{H}_3(\C^\splitt)$ & $\jordan{H}_3(\Quat)$ & $\jordan{H}_3(\Quat^\splitt)$\\ \hline
        $\Delta$ & $F_4$ & $F_4$ & $E_6$ & $F_4$ & $E_7$ \\ \hline\hline
        $\mathfrak{g}$ & $e_8^{(-24)}$ & $e_8^{(8)}$ & $\so(p+3,p+4)$ & $\so(p+3,p+3)$ & $g_2^{(2)}$\\ \hline
        $\jordan{J}$  & $\jordan{H}_3(\Octo)$ & $\jordan{H}_3(\Octo^\splitt)$ & $\R \oplus \jordan{J}(p,p+1)$ &$\R \oplus \jordan{J}(p,p)$ &$\R  $ \\ \hline
        $\Delta$  & $F_4$ & $E_8$ & $B_{p+3}$ & $D_{p+3}$ &$G_2$
    \end{tabular}
    \caption{Simple real Lie algebras with a $G_2$-grading.}
    \label{tab:G2graded}
\end{table}

\begin{theorem}[\cite{benkhart_lieAlgebrasGradedNonSimplyLaced}]
    Any Lie algebra with an $F_4$-grading is (a central extension of) the algebra $\TitsCons\big(\mathcal{R}, \jordan{H}_3(\Octo^\splitt)\big)$ where $\mathcal{R}$ is an alternative algebra with a trace map satisfying the degree 2 Cayley Hamilton equation.
\end{theorem}
Over $\R$, the common examples of alternative algebras are composition algebras. These algebras have a norm compatible with multiplication. The simple real Lie algebras with an $F_4$-grading are as follows.
\begin{table}[H]
    \centering
    \begin{tabular}{c||c|c|c|c|c|c|c}
        $\mathfrak{g}$ & $f_4^{(4)}$ & $e_6^{(2)}$ & $e_6^{(6)}$ & $e_7^{(-5)}$ & $e_7^{(7)}$ & $e_8^{(-24)}$ & $e_8^{(8)}$\\ \hline
        $\mathcal{R}$ & $\R$ & $\C$ & $\C^\splitt$ & $\Quat$ & $\Quat^\splitt$ & $\Octo$ & $\Octo^\splitt$ \\ \hline
        $\Delta$ & $F_4$ & $F_4$ & $E_6$ & $F_4$ & $E_7$ & $F_4$ & $E_8$
    \end{tabular}
    \caption{Simple real Lie algebras with a $F_4$-grading.}
    \label{tab:F4graded}
\end{table}

\begin{theorem}[\cite{benkhart_lieAlgebrasGradedNonSimplyLaced}]
     Any Lie algebra with a $B_p$-grading is (a central extension of) the algebra $\TitsCons\big(\jordan{J}(p,p+1), \jordan{J}\big)$ where $\jordan{J}$ is a degree 2 Jordan algebra over a commutative associative $\K$-algebra $A$.
\end{theorem}

Over $\R$ there are two choices of commutative associative algebras, $A = \R$ or $A = \C$. The following table shows all possible simple real Lie algebras with a $B_p$-grading.
\begin{table}[H]
    \centering
    \begin{tabular}{c||c|c|c|c}
         $\mathfrak{g}$& $\so(p+r,p+1+q)$ & $\so(p+r,p+r)$ &$\so\big(2(p+r)+1\big)^\C$ & $\so\big(2(p+r)\big)^\C$\\ \hline
         $\jordan{J}$ & $\jordan{J}(r,q)$ & $\jordan{J}(r,r-1)$ & $\jordan{J}(2r)$ & $\jordan{J}(2r-1)$ \\ \hline
         $\Delta$ & $B_{p+r}$ & $D_{p+r}$ & $B_{p+r}$ & $D_{p+r}$
    \end{tabular}
    \caption{Simple real Lie algebras with a $B_p$-grading.}
    \label{tab:Bpgraded}
\end{table}
In the next section, we will see that most of these gradings are induced by a choice of parabolic subalgebra.

\subsection{\texorpdfstring{$A_1$}{A1}-grading}\label{sec:A1gradedLieAlgebras}
\begin{figure}[htb]
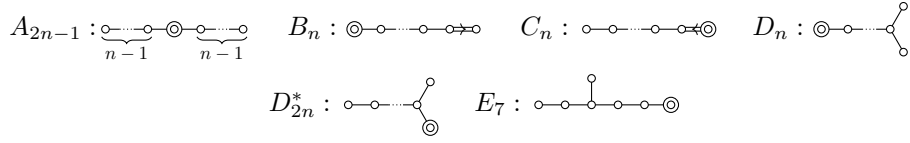

    \centering
    $A_{2n-1}:$\begin{dynkinDiagram}A{o.oOo.o}
        \dynkinBrace[n-1]12
        \dynkinBrace[n-1]45
    \end{dynkinDiagram}
    \quad
    $B_n:$ \dynkin B{Oo.ooo}
    \quad
    $C_n:$ \dynkin C{oo.ooO}
    \quad
    $D_n:$ \dynkin D{Oo.ooo}\\
    \quad
    $D^*_{2n}:$ \begin{dynkinDiagram}D{oo.ooO}
    \end{dynkinDiagram}
    \quad 
    $E_7:$ \dynkin E{ooooooO}
    \caption{Choice of $\Theta=\{\alpha\}$ giving an $A_1$-grading.}
    \label{fig:possibleThetadynks}
\end{figure}

There is a simple classification of $A_1$-gradings. We will show in \Cref{thm:A1GradingsInducedByTheta} that all $A_1$ gradings are induced by a parabolic subgroup $P_{\{\alpha\}}$ for some simple root $\alpha$ as in \Cref{def:ParabolicGrading}. Moreover, there are explicit necessary and sufficient conditions on $\alpha$ for $P_{\{\alpha\}}$ to induce an $A_1$-grading. All possible choices of $\alpha$ that result in an $A_1$-grading are shown in \Cref{fig:possibleThetadynks} by the circled node. 

First, we classify 3-graded simple Lie algebras, which are almost $A_1$-graded.
\begin{definition}
    A Lie algebra $\mathfrak{g}$ is \keyword{3-graded} if $\mathfrak{g}$ decomposes as $\mathfrak{g}_{-1}\oplus \mathfrak{g}_0\oplus\mathfrak{g}_1$ such that 
    \begin{equation*}
        [\mathfrak{g}_0, \mathfrak{g}_i] \subset \mathfrak{g}_i, \hspace{2pc} [\mathfrak{g}_1,\mathfrak{g}_1] = 0, \hspace{2pc} [\mathfrak{g}_{-1},\mathfrak{g}_{-1}]=0, \hspace{2pc} [\mathfrak{g}_1,\mathfrak{g}_{-1}] \subset \mathfrak{g}_0\,.
    \end{equation*}
\end{definition}
A 3-grading gives the appropriate decomposition for an $A_1$-grading, but does not guarantee the existence of a ``good'' split $A_1$-subalgebra satisfying Property \ref{prop:RootGradeFullAction} of a root grading. 

\begin{lemma}\label{thm:3gradingCharacterization}
    Let $\mathfrak{g}$ be a simple Lie algebra. Then $\mathfrak{g}$ is 3-graded if and only if there is a simple root $\alpha$ such that the parabolic subgroup $P_{\{\alpha\}}$ induces the grading. Moreover the coefficient of $\alpha$ in the longest root, $[\alpha_M : \alpha] =1$. 
\end{lemma}
\begin{proof}
    We first show that when $[\alpha_M : \alpha] = 1$, $P_{\{\alpha\}}$ induces a 3-grading. By \Cref{thm:UnderstandingRootGrading}, the eigenspaces of the decomposition are indexed by the possible coefficients of $\alpha$ for any root. Since $[\alpha_M : \alpha] = 1$, there are three options $\{-1,0,1\}$, splitting $\mathfrak{g} = \mathfrak{g}_{-1} \oplus \mathfrak{g}_0 \oplus \mathfrak{g}_1$.
    \Cref{thm:UnderstandingRootGrading} also states that each $\mathfrak{g}_i$ further splits into root spaces, so it suffices to understand the brackets between root spaces to check the bracket relations for 3-grading. For example, to check $[\mathfrak{g}_1,\mathfrak{g}_1] = 0$  it suffices to consider the bracket between $\mathfrak{g}_\beta \subset \mathfrak{g}_1$ and $\mathfrak{g}_\gamma \subset \mathfrak{g}_1$.  The coefficient $[\beta+\gamma : \alpha] = [\beta: \alpha] + [\gamma : \alpha] = 2$ and so $\beta + \gamma$ is not a root and $[\mathfrak{g}_\beta,\mathfrak{g}_\gamma] = 0$ as needed. The remaining bracket relations follow from analogous calculations of coefficients of $\alpha$ among root spaces living in the corresponding graded pieces.
    
    For the converse, suppose that $\mathfrak{g} = \mathfrak{g}_{-1} \oplus \mathfrak{g}_0 \oplus \mathfrak{g}_1$ is a 3-grading. We can choose the maximal split torus, $\mathfrak{t}$, to be contained in $\mathfrak{g}_0$ so that the root space decomposition refines the 3-grading. Therefore we can choose a set of simple roots so that the subalgebra $\mathfrak{p} = \mathfrak{g}_0\oplus \mathfrak{g}_1$ is a standard parabolic subalgebra.

    Let $\Theta$ be the subset of simple roots which do not vanish uniformly on $\mathfrak{t}\cap \mathfrak{z}$. Since $\mathfrak{g}_1$ is nonempty, $\Theta$ is nonempty as well. If $|\Theta| > 1$ then let $\alpha_1,\alpha_2$ be distinct simple roots in $\Theta$. Consider a chain from $\alpha_1$ to the highest root $\alpha_M$. Since $[\alpha_M, \alpha_2] > 0$ there is a positive root on the chain $\alpha'$ with $[\alpha': \alpha_2] = 0$ and $\alpha' + \alpha_2$ a root. The root space $\mathfrak{g}_{\alpha'}$ is in $\mathfrak{g}_1$ and has nontrivial bracket with $\mathfrak{g}_{\alpha_2} \subset \mathfrak{g}_1$ contradicting the 3-grading. Similarly if $[\alpha_M: \alpha] > 1$, there is a positive root on the chain from $\alpha_1$ to $\alpha_M$ with $[\alpha':\alpha_1] = 1$ and $\alpha' + \alpha_1$ a root. The same contradiction applies in this case. 
    Therefore the 3-grading is induced by $\mathfrak{p}_{\{\alpha\}}$ with $[\alpha_M: \alpha] = 1$ as claimed.
\end{proof}

In order to extend the classification of 3-gradings to $A_1$-gradings, we need to build a good $\sl_2$ triple in $\mathfrak{g}$ compatible with the 3-grading. This requires some explicit computations using the highest root in a finite root system. 

The coefficients of the highest root in a root system $R$ are related to a special labeling for the nodes of the Dynkin diagram of the associated \keyword{affine root system}, $\tilde{R}$. This Dynkin diagram for $\tilde{R}$ can be obtained by adding one node, called the \keyword{extending node}, and edges to the Dynkin diagram for $R$. We refer to \cite{kac_infinite_1984} for details about affine root systems and this labeling. 

\begin{lemma}[\cite{kac_infinite_1984}, Proposition 4.7]\label{thm:affineLabel}
    The Dynkin diagrams of an affine root system are characterized by the existence of a labeling of the nodes by positive integers such that the smallest label is 1 and twice each label is the sum of the adjacent labels ``with multiplicity.'' Equivalently, there is a positive integer vector $(d_1,\cdots,d_p)$ in the cokernel of the Cartan matrix of $\tilde{R}$ whose minimal entry is 1, $\sum c_{ij}d_i = 0$.
\end{lemma}

\begin{figure}[htb]
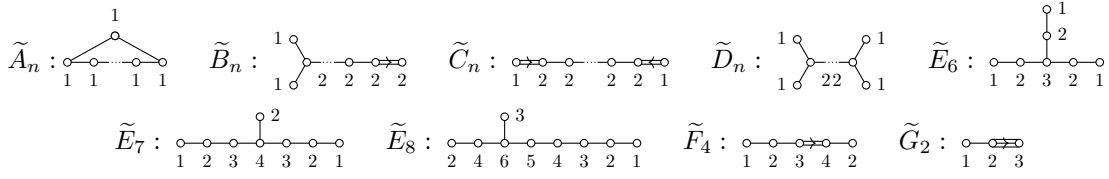

    \centering
    $\tilde{A}_{n}:$\dynkin[labels={1,1,1,1,1}] A[1]{oo.oo}
    \quad
    $\tilde{B}_n:$ \dynkin[labels={1,1,2,2,2,2}] B[1]{oo.ooo}
    \quad
    $\tilde{C}_n:$ \dynkin[labels={1,2,2,2,2,1}] C[1]{oo.ooo}
    \quad
    $\tilde{D}_n:$ \dynkin[labels={1,1,2,2,1,1}] D[1]{oo.ooo}
    \quad
    $\tilde{E}_6:$ \dynkin[labels={1,1,2,2,3,2,1}] E[1]{oooooo}\\
    \quad
    $\tilde{E}_7:$ \dynkin[labels={1,2,2,3,4,3,2,1}] E[1]{ooooooo}
    \quad 
    $\tilde{E}_8:$ \dynkin[labels={1,2,3,4,6,5,4,3,2}] E[1]{oooooooo}
    \quad
    $\tilde{F}_4:$ \dynkin[labels={1,2,3,4,2}] F[1]{oooo}
    \quad
    $\tilde{G}_2:$ \dynkin[labels={1,2,3}] G[1]{oo}
    \caption{Affine Dynkin diagrams and labeling as in \Cref{thm:affineLabel}.}
    \label{fig:affineDykin}
\end{figure}

The Dynkin diagrams and labelings are shown in \Cref{fig:affineDykin}. Removing a node labeled $1$ from the Dynkin diagram of $\tilde{R}$ results in a labeled Dynkin diagram for $R$. Let $d_i$ be the label of the node corresponding to $\alpha_i$. This is independent of the choice of node labeled 1 to remove. 

\begin{lemma}[\cite{kac_infinite_1984}, Remark 4.9]
    The coefficient $[\alpha_M:\alpha_i]$ of the highest root is equal to the label $d_i$.
\end{lemma}

We now decompose the longest word in the Weyl group, $w_0$, into commuting reflections about positive roots. This decomposition will help build the split $A_1$-subalgebra inside $\mathfrak{g}$. 
Recall from \Cref{def:RootSystem}, that each root, $\alpha$, defines a reflection ,$s_\alpha$, by $s_\alpha(v) = v- \langle v,\alpha\rangle \alpha$.
\begin{lemma}\label{lem:longest_element_orthogonal}
    Let $X$ be a reduced root system and $r$ be the number of orbits of simple roots under the opposition involution. There exists a family of $r$ strongly orthogonal roots $\delta_1,\cdots, \delta_r$ such that $w_0 = \prod s_{\delta_i} $ and $\delta_i$ is the highest root of a connected subdiagram of $X_i$.
\end{lemma}
\begin{proof}
    We begin with $\delta_1 = \sum d_i \alpha_i$ the highest positive root. The other $\delta_i$ will be constructed inductively as the highest roots of subdiagrams fixed by reflection in all previous roots. Now we consider the action of $\delta_1$ on $\alpha_j$ an arbitrary simple root,  
    \begin{equation*}
        s_{\delta_1} \alpha_j = \alpha_j - \langle \alpha_j,\delta_1\rangle \delta_1 = \alpha_j - \langle \delta_1, \alpha_j\rangle \frac{(\alpha_j,\alpha_j)}{(\delta_1,\delta_1)} \delta_1 = \alpha_j - \left(\sum_i  c_{ij}d_i\right) \frac{(\alpha_j,\alpha_j)}{(\delta_1,\delta_1)} \delta_1\,.
    \end{equation*}
    By \Cref{thm:affineLabel}, when $j$ is not adjacent to the extending node $\sum c_{ij}d_i = 0$ and $\alpha_j$ is fixed by $\delta_1$. Otherwise $\sum c_{ij}d_i > 0$ and we obtain a negative root.  

    We now consider the strict subdiagram of nodes fixed by $s_{\delta_1}$. By case by case observation, we see the nonfixed nodes form a single orbit of the opposition involution. This leaves $r-1$ orbits in the subdiagram, whose longest word inductively splits into a product of $c-1$ strongly orthogonal reflections about $\delta_2,\cdots, \delta_r$. Furthermore each $\delta_i$ for $i > 1$ is a sum of nodes fixed by $\delta_1$ and thus is strongly orthogonal to $\delta_1$ as well. Note that when the fixed diagram of $\delta_1$ is disconnected the highest roots in each subdiagram are strongly orthogonal and so should both be chosen to obtain $r$ roots.
\end{proof}

\begin{example}\label{ex:e7orthogonal}
    The longest word in $W(E_7)$ is given by reflections in the following roots
    \begin{align*}
        \delta_1=\substack{2\hphantom{0}\\234321} \quad \delta_2=\substack{1\hphantom{0}\\012221}\quad \delta_3 = \substack{0\hphantom{0}\\000001} \quad \delta_4=\substack{1\hphantom{0}\\012100} \quad \delta_5=\substack{1\hphantom{0}\\000000} \quad \delta_6=\substack{0\hphantom{0}\\010000} \quad \delta_7=\substack{0\hphantom{0}\\000100}\,.
    \end{align*}
    Here we graphically labeled the nodes of the Dynkin diagram by the coefficient of the simple root. So $\delta_1 = 2\alpha_1 + 2\alpha_2 + 3\alpha_3 + 4\alpha_4 + 3\alpha_5 + 2\alpha_6+\alpha_7$ using the indexing for $E_7$ in \cite{knapp1996lie}.
\end{example}

\begin{theorem}\label{thm:A1GradingsInducedByTheta}
    Let $\mathfrak{g}$ be a simple Lie algebra. Then $\mathfrak{g}$ is $A_1$-graded if and only if there is a root $\alpha$ such that $P_{\{\alpha\}}$ induces the grading. Moreover, $\alpha$ is fixed by the opposition involution and has $[\alpha_M : \alpha] =1$.
\end{theorem}
\begin{proof}
    We begin by showing that whenever $\alpha$ is fixed by the opposition involution and $[\alpha_M:\alpha] = 1$, the parabolic $P_{\{\alpha\}}$ induces an $A_1$ grading. By \Cref{thm:3gradingCharacterization}, we know $P_{\{\alpha\}}$ induces a 3-grading, $\mathfrak{g} = \mathfrak{g}_{-1} \oplus \mathfrak{g}_{0} \oplus \mathfrak{g}_1$. In order to upgrade to an $A_1$-grading we need to check the remaining three properties.
    
    First we check $\mathfrak{g}_0$ is generated by iterated brackets in $\mathfrak{g}_1$ and $\mathfrak{g}_{-1}$. It suffices to hit every simple root space in $\mathfrak{g}_0$ as these generate simple Lie algebras. Let $\beta\neq \alpha$ be a simple root and consider a chain from $\alpha$ to the highest root $\alpha_M$. At some point there is a root $\alpha'$ with $[\alpha':\beta] = 0$ and $\alpha' + \beta$ a root. Moreover $[\alpha':\alpha] = 1$ so $\mathfrak{g}_{\alpha'+\beta} \in \mathfrak{g}_1$ and $\mathfrak{g}_{-\alpha'}\in \mathfrak{g}_{-1}$. Then $[\mathfrak{g}_{\alpha'+\beta},\mathfrak{g}_{-\alpha'}] = \mathfrak{g}_\beta$ as needed.
    
    Next we check there is a good split $A_1$ subalgebra in $\mathfrak{g}$ satisfying Properties \ref{prop:RootGradeSplit} and \ref{prop:RootGradeFullAction}. Concretely, we find an $\sl_2$ triple $(E,F,H)$ with $E \in \mathfrak{g}_1$, $F \in \mathfrak{g}_{-1}$, and $H \in \mathfrak{t}\cap \mathfrak{z}_{\{\alpha\}} \subset \mathfrak{g}_0$. We build $E$ using the root spaces corresponding to the set of strongly orthogonal roots, $\{\delta_i\}$ constructed in \Cref{lem:longest_element_orthogonal}. First relabel the roots such that the roots containing $\alpha$ are the first $c$ roots and that for $i<j<c$, $\delta_i > \delta_j$ in the root partial order. This is possible since each $\delta_i$ is the highest root for some nested subdiagram. For $i = 1,\cdots,c$, chose a nonzero $ E_i \in \mathfrak{g}_{\delta_i}$ and complete it to an $\sl_2$-triple $(E_i,F_i,H_i=\check{\delta}_i)$.
    Then $E=\sum E_i$, $F_i = \sum F_i$ and $H=\sum H_i$. Since the $\delta_i$ are strongly orthogonal, for $i\neq j$, $[E_i,F_j] = [H_i,E_j] = [H_i,F_j] =0$ and  we have
    \begin{equation*}
        [E,F] = \left[\sum_{i=1}^c E_i, \sum_{j=1^c}F_j\right] = \sum_{i=1}^c\sum_{j=1}^c[E_i,F_j] = \sum_{i=1}^c H_i =H\,.
    \end{equation*}
    Thus $(E,F,H)$ is an $\sl_2$ triple which gives an embedding of a split $A_1$-Lie algebra into $\mathfrak{g}$ respecting the 3-grading. It remains to check that $H$ acts correctly on all of $\mathfrak{g}_1$ (and $\mathfrak{g}_{-1}$), in other words, $H \in \mathfrak{t}\cap \mathfrak{z}_{\{\alpha\}}$. By  \Cref{thm:UnderstandingRootGrading} it suffices to check that all simple roots other than $\alpha$ vanish on $H$. 

    We now verify this fact in each case. Let $\alpha_j$ be a simple root and compute
    \begin{equation*}
        \alpha_j(H) = \alpha_j\left(\sum_{k=1}^c \check{\delta}_k\right) = \left\langle \alpha_j, \sum_{k=1}^c \delta_k \right\rangle = \sum \frac{(\alpha_j,\delta_k)}{(\delta_k,\delta_k)} \,.
    \end{equation*}
By construction $\delta_k$ is the highest root for a subdiagram $X_k$. Thus there are two possible cases for $(\delta_k,\alpha_j)$ to be nonzero: $\alpha_j$ is the extending node of $X_k$ or $\alpha_j$ is adjacent to $X_k$. These two cases are related in the construction of the set of highest roots; whenever $\alpha_j$ is the extending node for some $X_{k'}$ it is then adjacent to $X_{k'+1}$. Moreover, each node can only be the extending node of a subdiagram $X_k$ at most once. 

Next, we analyze each of the finitely many cases of $\alpha$ that are fixed by the opposition involution and have $[\alpha_M:\alpha] = 1$ that are shown in \Cref{fig:possibleThetadynks}. For each $\alpha_j \neq \alpha$ that is the extending node for $X_{k'}$, we observe that $X_{k'+1}$ is the only diagram adjacent to $\alpha_j$. So $k'$ and $k'+1$ are the only times $(\delta_k,\alpha_j) \neq 0$. Moreover we compute in each case that     
\begin{equation*}
        \frac{(\alpha_j,\delta_k)}{(\delta_k,\delta_k)} = - \frac{(\alpha_j,\delta_{k+1})}{(\delta_{k+1},\delta_{k+1})} \,.
    \end{equation*} and that $(\alpha_j,\delta_k') = -(\alpha_j,\delta_{k'+1})$.
This implies that $\alpha_j(H) = 0$ for all $\alpha_j \neq \alpha$ that are the extending node for some $X_{k'}$. If $\alpha_j \neq \alpha$ is never an extending node, then $(\alpha_j, \delta_k)=0$ for all $k$ and $\alpha_j(H)$ vanishes as needed.

Finally we consider the case $\alpha_j = \alpha$. In this case, $\delta_c = \alpha$ and the only nonzero pairing is $(\alpha,\delta_c) = 2$ as needed.

    For illustration, we give the explicit calculation for the two possible cases in $D_6$ and $D_6^*$ as in \Cref{fig:possibleThetadynks}. The set of orthogonal roots from \Cref{lem:longest_element_orthogonal} is
    \begin{equation*}
        1222\substack{1\\1} \qquad 10000\substack{0\\0} \qquad 0012\substack{1\\1} \qquad 0000\substack{0\\1} \qquad 0000\substack{1\\0} \qquad 0001\substack{1\\1} \,.
    \end{equation*}
    When $\alpha$ is $1000\substack{0\\0}$ as in case $D_6$, the $\delta$ used to build $E$ are
    \begin{equation*}
        \delta_1 = 1222\substack{1\\1} \qquad \delta_2 = 10000\substack{0\\0}\,.
    \end{equation*}
    The only interesting simple root is $\alpha_2 = 0100\substack{0\\0}$. Since $(\delta_1,\delta_1) = (\delta_2,\delta_2)$ it suffices to compute $(\delta_1,\alpha_2) = 1$ and $(\delta_2,\alpha_2) = -1$.

    When $\alpha = 0000\substack{0\\1}$ as in case $D_6^*$, the $\delta$ used to build $E$ are
    \begin{equation*}
        \delta_1 = 1222\substack{1\\1} \qquad \delta_2 = 0012\substack{1\\1} \qquad \delta_3 = 0000\substack{0\\1}\,.
    \end{equation*}
    In this case there are two interesting simple roots $\alpha_2 = 0100\substack{0\\0}$ and $\alpha_4 = 0100\substack{0\\0}$. To check the result we compute $(\delta_1,\alpha_2) = 1$, $(\delta_2,\alpha_2) = -1$ and $(\delta_2,\alpha_4) = 1$, $(\delta_3,\alpha_4) = -1$ as needed.
  
    \medskip
    For the opposite direction we assume that $\mathfrak{g}$ is $A_1$-graded. Then in particular, $\mathfrak{g}$ is 3-graded and so there is a simple root $\alpha$ with $[\alpha_M: \alpha]=1$ and $P_{\{\alpha\}}$ induces the 3-grading. We can take the split torus $\mathfrak{t}$ which induces the root system grading to include $\mathfrak{h}$, the Cartan subalgebra of the split $A_1$ subalgebra of the root grading.
    
    Next let $G$ be a algebraic group with Lie algebra $\mathfrak{g}$ and chose a lift $\overline{s}$ of the unique generator of the Weyl group of the split $A_1$ subalgebra to $G$. By definition, the adjoint action of $\overline{s}$ conjugates $P_{\{\alpha\}}^\opp$ to $P_{\{\alpha\}}$. Similarly, a lift $\overline{w_0}$ of the longest word in the full Weyl group conjugates $P_{\{\alpha\}}^\opp$ to $P_{\{-w_0\alpha\}}$. This implies that $P_{\{\alpha\}}$ and $P_{\{-w_0\alpha\}}$ are conjugate and thus are the same standard parabolic. So $\alpha=-w_0\alpha$ and $\alpha$ is fixed by the opposition involution as needed. 
\end{proof}

\begin{remark}
This classification of $A_1$-gradings is very similar to the classification of Hermitian symmetric domains of tube type proved in \cite{deligne-HermitianTubeType}.     
\end{remark}

\subsection{Larger gradings}\label{sec:largerGrading}
Let $\mathfrak{g}$ be a simple Lie algebra of Weyl type $X$ with simple roots $\Delta= \{\alpha_1, \cdots, \alpha_n\}$ and let $\Theta = \{ \alpha_{i_1},\cdots \alpha_{i_p}\}\subset \Delta$. We label the simple $\Theta$-root corresponding to the $\alpha_{i_j}$-space by $\beta_j$. For each $\alpha_{i_j} \in \Theta$, let $F_j$ be the connected component containing $\alpha_{i_j}$ of the Dynkin diagram given by removing $\alpha_{i_k}$ for $k \neq j$.

\begin{definition}
    We say that $\Theta$ is \keyword{Jordan compatible} if for each $\alpha_{i_j} \in \Theta$, the root $\alpha_{i_j}$ satisfies the condition to induce an $A_1$-grading in $F_j$, i.e is labeled 1 and fixed by the opposition involution for $F_j$. 
\end{definition}

\begin{figure}[htb]
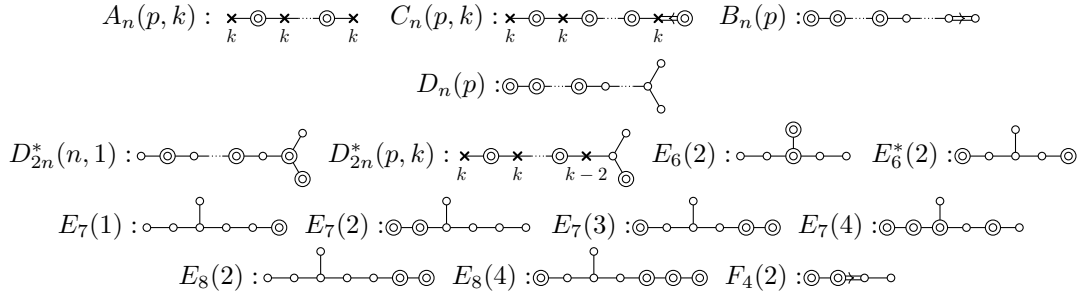

    \centering
    $A_n(p,k):$ \dynkin[labels={k,,k,,k}] A{xOx.Ox} \hspace{1mm}
    $C_n(p,k):$\dynkin[labels={k,,k,,,k,}] C{xOxO.OxO} \hspace{1mm}
    $B_n(p):$\dynkin B{OO.Oo.oo}\hspace{1mm} 
    $D_n(p):$\dynkin D{OO.Oo.ooo} \hspace{1mm}\\
    $D^*_{2n}(n,1):$\dynkin D{oOo.OoOoO}\hspace{1mm}
    $D^*_{2n}(p,k):$\dynkin [labels={k~,,k,,k-2,,}]D{xOx.OxooO}\hspace{1mm}
    $E_6(2):$\dynkin E{oOoOoo}\hspace{1mm}
    $E^*_6(2):$\dynkin E{OooooO}\hspace{1mm}\\
    $E_7(1):$\dynkin E{ooooooO}\hspace{1mm}
    $E_7(2):$\dynkin E{OoOoooo}\hspace{1mm}
    $E_7(3):$\dynkin E{OooooOO}\hspace{1mm}
    $E_7(4):$\dynkin E{OoOOoOo}\hspace{1mm}\\
    $E_8(2):$\dynkin E{ooooooOO}\hspace{1mm}
    $E_8(4):$\dynkin E{OooooOOO}\hspace{1mm}
    $F_4(2):$\dynkin F{OOoo}\hspace{1mm}
    \caption{Jordan compatible choices of $\Theta$.}
    \label{fig:possibleWeylTypedynks_all}
\end{figure}

\Cref{fig:possibleWeylTypedynks_all} shows all possible Jordan compatible choices of $\Theta$, other than $\Theta=\Delta$. The roots in $\Theta$ are circled, and the $\times$ labeled $k$ denotes an $A_k$-subdiagram with no circled nodes. This list is obtained by combining the definition with the classification of $A_1$-gradings as shown in \Cref{fig:possibleThetadynks}.

\begin{theorem}\label{thm:ClassifyParabolicInduceRootGrading}
      Let $\Theta \subset \Delta$ be a subset of simple roots and $\mathfrak{p}_\Theta$ the associated parabolic subalgebra. Then the grading induced by $\mathfrak{p}_\Theta$ is a grading by a root system of some type $Y$ if and only if $\Theta$ is Jordan compatible. 
\end{theorem}

\begin{proof}
    First we assume that $\mathfrak{p}_\Theta$ induces a root grading of type $Y$. Consider the subalgebra generated by the roots in the subdiagram $F_i$ corresponding to $\beta_i$. The restriction of $P_\Theta$ to this subalgebra induces an $A_1$-grading and thus is given by a parabolic associated to a root $\gamma_i$ with $[\alpha_{M_i}: \gamma_i] = 1$ for $\alpha_{M_i}$ the highest root in $F_i$ and $\gamma_i$ fixed by the opposition involution in $F_i$. Moreover this parabolic is the restriction of $P_\Theta$ and so $\gamma_i = \beta_i$ as needed for $\Theta$ to be Jordan compatible.\smallskip
    
    We now show the converse. Assume $\Theta$ is Jordan compatible and thus that each diagram $F_j$ has an induced $A_1$-grading.
    This gives elements $\Id_{\pm\beta_j}$ for each simple root in $\Theta$ as the element $E$ constructed from a sum of strongly orthogonal root spaces as in \Cref{thm:A1GradingsInducedByTheta}. The parabolic $\mathfrak{p}_\Theta$ divides $\mathfrak{g}$ into the union of root spaces which agree on $\mathfrak{h} = \mathfrak{z}_\Theta \cap \mathfrak{t}$. Thus each graded piece is labeled by a $Y$-root, the projection the restricted root space onto the subspace spanned by $\Theta$. It remains to check this decomposition satisfies the three properties of a root grading.
    
    Property \ref{prop:RootGradeGenerated} of being $Y$-graded follows directly from the $A_1$-grading on each subdiagram. The sums of iterated brackets of $\Id_{\pm\beta_j}$ generate the part of the Levi subgroup corresponding to each diagram. Therefore all of the Levi ($\mathfrak{g}_0$) is generated by brackets in the proposed $Y$-roots. 
    
    Property \ref{prop:RootGradeSplit} is the bulk of the proof and requires a case by case check that the subalgebra generated by the identity elements $\Id_{\pm \beta_j}$ is a split Lie algebra of type $Y$. This subalgebra will have Cartan subspace $\mathfrak{h}$.
    
    Property \ref{prop:RootGradeFullAction} is then an easy consequence of the parabolic grading. The action of the split torus, $\mathfrak{h}$, acts by the $Y$-roots on each restricted root space, not only on $\K \Id_{\beta}$ as needed.   

    We now focus on Property \ref{prop:RootGradeSplit} and check the identity elements generate a split subalgebra. In order to determine the type of the split subalgebra we focus on the subdiagram $F_{ij} = F_i \cup F_j$ for two simple roots $\beta_i$ and $\beta_j$. There are only a few possibilities for $F_{ij}$:
    \begin{enumerate}
    \setcounter{enumi}{-1}
        \item $F_{ij}$ is disconnected.
        \item $F_{ij}$ is $A_n(2,k)$ or $E^*_6(2)$.
        \item $F_{ij}$ is $B_n(2),C_n(2,k),D_n(2),$ or $D^*_{2n}(2,k)$.
        \item $F_{ij}$ is $F_4(2),E_6(2),E_7(2),$ or $E_8(2)$.
    \end{enumerate}
    The Dynkin diagram for $Y$ is constructed by examining $F_{ij}$ for each pair of roots in $\Theta$ and assigning $k$ edges between node $i$ and $j$ when $F_{ij}$ is in case $(k)$. It suffices to show the Serre relations hold for each subalgebra generated by $\Id_{\pm \beta_i}$ and $\Id_{\pm \beta_j}$ for the groups of type $A_1\times A_1, A_2, B_2,G_2$ is each case. 
    
    The case $(0)$ for $A_1\times A_1$ is trivial as the disconnected diagrams result in commuting $A_1$ subalgebras.

    For the remaining cases we start by verifying that $[h_i, \Id_{\beta_j}] = \langle \beta_i, \beta_j \rangle \Id_{\beta_j}$ with $h_i = [\Id_{\beta_i},\Id_{-\beta_i}]$.  By construction, the element $h_j = [\Id_{\beta_j},\Id_{-\beta_j}] = \sum_{k=1}^{c_j} (\check{\delta}_k)_j$ for the set $\{(\delta_k)_j\}$ of strongly orthogonal roots with $[(\delta_k)_i:\alpha_{i_j}] > 0$ constructed in \Cref{thm:A1GradingsInducedByTheta}. Generically we have
    \begin{equation*}
        [h_1, \Id_{\beta_2}] = \left[\sum_{k=1}^{c_1} (\check{\delta}_k)_1, \sum_{l=1}^{c_2} \Id_{(\delta_l)_2}\right] = \sum_{k=1}^{c_1}\sum_{l=1}^{c_2} \langle (\delta_l)_2, (\delta_k)_1 \rangle (\delta_l)_2\,.
    \end{equation*}
    We will compute the final sum, using the explicit family of roots $\delta_k$ in each case.

    We begin by examining the choices of $F_{ij}$ which we claim result in an $A_2$-grading. In \Cref{tab:A2rootDecomposition}, we present the strongly orthogonal roots which decompose $\Id_{\beta_1}$ and $\Id_{\beta_2}$. Each root is encoded by writing the coefficient of each positive root over the respective node of the Dynkin diagram of the original restricted root system. The underlined entries correspond to the choice of $\Theta$. The roots in \Cref{tab:A2rootDecomposition} are organized so the only nontrivial pairing is among roots in the same row.
    \begin{table}[H]
    \centering
        \begin{equation*}
            \begin{array}{c||c |c}
                 & \beta_1 & \beta_2\\ \hline
                 \multirow{5}{*}{$A_n(2,k)$}&  11\cdots 11\underline{1}11\cdots11\underline{0}00\cdots00& 000\cdots 00\underline{0}00\cdots00\underline{1}00\cdots00  \\
                     &  01\cdots 11\underline{1}11\cdots10\underline{0}00\cdots00& 000\cdots 00\underline{0}00\cdots01\underline{1}10\cdots00  \\
                     & \vdots & \vdots \\
                     &  00\cdots 01\underline{1}10\cdots00\underline{0}00\cdots00& 000\cdots 00\underline{0}01\cdots11\underline{1}11\cdots10  \\
                     &  00\cdots 00\underline{1}00\cdots00\underline{0}00\cdots00& 000\cdots 00\underline{0}11\cdots11\underline{1}11\cdots11  \\ \hline
                 \multirow{2}{*}{$E_6^*(2)$} & \substack{1\\\underline{1}221\underline{0}} &\substack{0\\ \underline{0}000\underline{1}} \\
                     & \substack{0\\ \underline{1}000\underline{0}} & \substack{1\\ \underline{0}221\underline{1}} 
            \end{array}
        \end{equation*}
        \caption{Decompositions of simple roots for $A_2$-gradings}\label{tab:A2rootDecomposition}
    \end{table}
    For example, we consider $A_n(2,k)$.  The set of orthogonal roots for $\beta_1$ are $\alpha_{[1,2k+1]}$, $\alpha_{[2,2k]}$,$\cdots$, $\alpha_{[k,k+2]}$,$\alpha_{k+1}$ where $\alpha_{[i,j]} = \sum_{\ell=i}^j \alpha_\ell$. Similarly the roots for $\beta_2$ are $\alpha_{2k+2}$, $\alpha_{[2k+1,2k+3]}$,$\cdots$ $\alpha_{[k+3,3k+1]}$, $\alpha_{[k+2,3k+2]}$. Since the pairing $\langle \alpha,\beta\rangle$ is nonzero if and only if $\alpha\pm \beta$ is a root, we check the only nonzero pairing is across rows of \Cref{tab:A2rootDecomposition}.\medskip
    
    For the choices of $F_{ij}$ with a $B_2$-grading, \Cref{tab:B2rootDecomposition} shows the analogous sets of strongly orthogonal roots. In this case, two roots decomposing $\beta_1$ can have nontrivial pairing with a root of $\beta_2$. When this occurs, the roots decomposing $\beta_2$ are placed between the rows for nontrivial paired roots decomposing $\beta_1$ in \Cref{tab:B2rootDecomposition}. All pairings outside these ``rows'' are trivial.
    \begin{table}[hb]
        \begin{equation*}
            \begin{array}{c||c |c}
                 & \beta_1 & \beta_2\\ \hline
                 \multirow{2}{*}{$B_n(2)$}&  \underline{01}2\cdots2\Rightarrow2 & \multirow{2}{*}{$\underline{10}\cdots0\Rightarrow0$}  \\
                     & \underline{01}0\cdots0\Rightarrow0\\ \hline                
                 \multirow{5}{*}{$C_n(2,k)$} & 11\cdots11\underline{1}11\cdots11\Leftarrow\underline{0} & 00\cdots00\underline{0}00\cdots00\Leftarrow\underline{1}\\
                     & 01\cdots11\underline{1}11\cdots10\Leftarrow\underline{0}& 00\cdots00\underline{0}00\cdots02\Leftarrow\underline{1}\\
                     & \vdots & \vdots\\
                     & 00\cdots01\underline{1}10\cdots00\Leftarrow\underline{0}& 00\cdots00\underline{0}02\cdots22\Leftarrow\underline{1}\\
                     & 00\cdots00\underline{1}00\cdots00\Leftarrow\underline{0}& 00\cdots00\underline{0}22\cdots22\Leftarrow\underline{1}\\ \hline
                 \multirow{2}{*}{$D_n(2)$} & \underline{0}\underline{1}0\cdots0\substack{0\\0} & \multirow{2}{*}{$\underline{10}0\cdots0\substack{0\\0}$}\\
                     & \underline{0}\underline{1}2\cdots2\substack{1\\1} \\ \hline
                 \multirow{5}{*}{$D_{2n}^*(2,2k)$} & 11\cdots 1\underline{1}11\cdots1\substack{1\\ \underline{0}}& \multirow{2}{*}{$00\cdots00\underline{0}00\cdots0\substack{0\\ \underline{1}}$}\\
                     &01\cdots 1\underline{1}11\cdots1\substack{0\\ \underline{0}} &  \\
                     &\vdots & \vdots  \\
                     &00\cdots01\underline{1}10\cdots0\substack{0\\ \underline{0}} & \multirow{2}{*}{$00\cdots00\underline{0}12\cdots2\substack{1\\ \underline{1}}$}  \\
                     &00\cdots00\underline{1}00\cdots0\substack{0\\ \underline{0}} & 
            \end{array}
        \end{equation*}
        \caption{Decompositions of simple roots for $B_2$-gradings}\label{tab:B2rootDecomposition}
    \end{table}
    
    In this case the pairing between $\beta_1$ and $\beta_2$ is asymmetric, $[h_1,\Id_{\beta_2}] = -2\Id_{\beta_2}$ while $[h_2,\Id_{\beta_1}] = -\Id_{\beta_1}$. In type $C_n(2,k)$ the asymmetry in the pairing comes from the asymmetry of the underlying pairing $\langle \delta_i,\delta_j\rangle$. In every other case the asymmetry comes from the fact that the root $\beta_1$ decomposes into twice as many root spaces as $\beta_2$. Each root of $\beta_2$ pairs nontrivially with two roots of $\beta_1$ with value $-1$. The asymmetry in sizes results in $[h_1,\Id_{\beta_2}] = -2\Id_{\beta_2}$ while $[h_2,\Id_{\beta_1}] = -\Id_{\beta_1}$ as needed. 

    The final case is the subdiagrams which generate a $G_2$ subsystem. As before, \Cref{tab:G2rootDecomposition} shows the strongly orthogonal roots used to construct $\Id_{\beta_1}$ and $\Id_{\beta_2}$. In each case, $\Id_{\beta_2}$ decomposes into a single root space which pairs nontrivially with each space in $\Id_{\beta_1}$. As in $B_2$, this asymmetry of dimensions corresponds to the required asymmetry of pairing.
    \begin{table}[hb]
        \begin{equation*}
            \begin{array}{c||c|c}
                  & \beta_1 & \beta_2\\ \hline
                 \multirow{3}{*}{$F_4(2)$} & \underline{01}\Rightarrow22 &  \multirow{3}{*}{$ \underline{10}\Rightarrow00$} \\
                      & \underline{01}\Rightarrow20&   \\
                      & \underline{01}\Rightarrow00 &  \\ \hline
                  \multirow{3}{*}{$E_6(2)$} &\substack{\underline{0}\\11\underline{1}11}& \multirow{3}{*}{$\substack{\underline{1}\\00\underline{0}00}$} \\
                      &\substack{\underline{0}\\01\underline{1}10}&  \\
                      &\substack{\underline{0}\\00\underline{1}00}&  
            \end{array}
            \hspace{2pc}
            \begin{array}{c||c|c}
                  & \beta_1 & \beta_2\\ \hline
                  \multirow{3}{*}{$E_7(2)$} & \substack{1\hphantom{0}\\ \underline{01}2221} & \multirow{3}{*}{$ \substack{0\hphantom{0}\\ \underline{10}0000}$} \\
                      &\substack{1\hphantom{0}\\ \underline{01}2100} & \\
                      &\substack{0\hphantom{0}\\ \underline{01}0000} & \\ \hline
                  \multirow{3}{*}{$E_8(2)$} & \substack{2\hphantom{00}\\23432\underline{10}} & \multirow{3}{*}{$\substack{0\hphantom{00}\\00000\underline{01}}$} \\
                      &\substack{1\hphantom{00}\\01222\underline{10}} &  \\
                      &\substack{0\hphantom{00}\\00000\underline{10}} & 
            \end{array}
        \end{equation*}
        \caption{Decompositions of simple roots for $G_2$-gradings}\label{tab:G2rootDecomposition}
    \end{table}

    The finiteness of the Serre relation, $\ad (\Id_{\beta_i})^{-\langle \beta_j,\beta_i\rangle+1}(\Id_{\beta_j})=0 $, is clear by observing the coefficients of the nodes in $\Theta$ in the highest root. In each connected case, the coefficients of $(\beta_1,\beta_2)$ are $(1,1)$, $(2,1)$, or $(3,2)$ when the subdiagram is type $A_2,B_2$ or $G_2$ respectively. Since the adjoint action maps root spaces to root spaces, the maximum numbers of times $\ad (\Id_{\beta_i})$ can be nonzero is the coefficient, $[\alpha_M : \beta_i]$. Moreover, $[\alpha_M: \beta_i]$ is the coefficient of the highest root in a root system of type $Y$. For example, in $F_4(2)$, the highest root is $\alpha_M = 23\Rightarrow43 = 2\alpha_1+3\alpha_2+4\alpha_3+2\alpha_4$. Then $\beta_1 = 10\Rightarrow00$ and $[\alpha_M:\beta_1] = 2$ and $\beta_2 = 01\Rightarrow00$ has $[\alpha_M:\beta_2]=3$ matching $G_2$ as needed. Moreover these coefficients match the coefficients of the highest root of the relevant root system. This implies that these final Serre relations match the proposed root system.
    
\end{proof}

\subsection{Grading by folding}
There are several gradings which arise from an alternate construction, a folding of the Dynkin diagram. Recall that each non-simply laced diagram can be obtained by \keyword{folding} a simply laced diagram via a symmetry of the diagram as follows:
\begin{equation*}
    A_{2p-1}\rightarrowtail C_{p}, \qquad D_{p+1} \rightarrowtail B_p, \qquad E_6 \rightarrowtail F_4, \qquad D_4 \rightarrowtail G_2\,.
\end{equation*}
There is one special folding of $B_3$ onto $G_2$, which does not arise from an automorphism of the diagram. Instead, the symmetry identifies the root spaces corresponding to the two ends. This makes sense as it ``finishes'' the folding of $D_4$ onto $G_2$ that began by folding $D_4$ onto $B_3$.

\begin{proposition}
    Let $R$ be a root system which folds onto a root system $R'$. Then any Lie algebra which is $R$-graded is also $R'$-graded.
\end{proposition}
\begin{proof}
    The root spaces of the folded algebra are given by the orbits of the symmetry used to fold it. Thus, the $R'$-grading is given by identifying the pieces of the $R$-grading under the folding symmetry. 
\end{proof}

\begin{proposition}\label{thm:ClassifyingRGradingConstructions}
    Every $R$-grading of a simple Lie algebra comes from a grading by a parabolic subgroup possibly followed by a folding. 
\end{proposition}
\begin{proof}
    Suppose $\mathfrak{g}$ is simple and $R$-graded. 
    By the recognition theorems (\Cref{sec:RecognitionTheorems}) we can assume that $\mathfrak{g}$ has a reduced restricted root system (i.e. is not of type $BC$), as otherwise it would not be graded by a reduced root system. 

    Denote by $R^+$ the set of positive roots in $R$. Consider the associated parabolic subalgebra $\mathfrak{p}=\bigoplus_{\beta \in R^+ \cup \{0\}} \mathfrak{g}_\beta$. We wish to prove that $\mathfrak{p}$ induces a grading by a root system $R'$ which comes with a folding $R'\to R$. Write $\mathfrak{z}$ for the center of $\mathfrak{g}_0$. We have $\mathfrak{h}\subset\mathfrak{z}$ where $\mathfrak{h}$ is the Cartan subgroup of the split group coming from the $R$-grading. Consider a simple $R$-root space $\mathfrak{g}_\beta$. The action of $\mathfrak{z}$ on it decomposes it into $\mathfrak{g}_\beta = \bigoplus_i\mathfrak{g}_{\gamma^i}$, since $\mathfrak{g}_\beta$ is abelian these commute with each other. As $\mathfrak{p}$ is a parabolic subalgebra, each $\gamma^i$ corresponds to a root $\alpha^i$ of $\mathfrak{g}$. 
    
    Elements of the normalizer of $\mathfrak{h}$ must map $\mathfrak{g}_0$ into itself, and hence must map $\mathfrak{z}$ into itself as well. Therefore the $W(R)$ action respects the decomposition of $\mathfrak{g}_\beta$ by $\mathfrak{z}$, i.e. any root space for a root in the $W(R)$ orbit of $\beta$ also decomposes into the same number of $\mathfrak{z}$-root spaces. In particular $\mathfrak{g}_{-\beta}$ decomposes into $\bigoplus_i\mathfrak{g}_{-\gamma^i}$. Since these commute with each other, we have that $[\mathfrak{g}_{\gamma^i},\mathfrak{g}_{-\gamma^j}] = 0 $ for $i\neq j$ as well.
    
    Now we observe that the projection of the identity elements $\Id_{\pm\beta}$ into the $\mathfrak{g}_{\pm\gamma^i}$ give $\sl_2$-triples in $\mathfrak{g}_{\gamma^i}\oplus\mathfrak{g}_{-\gamma^j}\oplus [\mathfrak{g}_{\gamma^i},\mathfrak{g}_{-\gamma^j}]$. This implies that these subalgebras are $A_1$-graded and hence the collection of roots $\alpha^i$ must be Jordan compatible so that this parabolic induces a root grading $R'$. 

    The folding from $R'$ to $R$ is given by grouping all the simple $\gamma^j$ whose root spaces are contained in the same $\mathfrak{g}_{\beta}$. 
\end{proof}

\subsection{Identifying Jordan factors}\label{sec:IdentificationOfJordanFactors}

Our classification of Jordan compatible subsets of roots of a Dynkin diagram (\Cref{thm:ClassifyParabolicInduceRootGrading}) recovers most of the features of the recognition theorems given in \Cref{sec:RecognitionTheorems}. In particular, if $\mathfrak{g}$ is a simple Lie algebra over $\K$ and the parabolic subalgebra associated to a subset $\Theta\subset \Delta$ induces an $R$-grading we obtain restrictions on the possible Jordan algebras associated to the roots in $R$.

First, we reinterpret what we have established about $A_1$-graded algebras to give a classification analogous to the classification of Jordan algebras over an algebraically closed field given in \cite[Renaissance Structure Theorem]{mccrimmon-TasteJordanAlgebras}.

\begin{proposition}\label{prop:classification_of_jords} Let $\mathfrak{g}$ be an $A_1$-graded simple Lie algebra over an algebraically closed field $\algclosure{\K}$. Then $\mathfrak{g} = \slj$ and belongs to one of the following cases:
    \begin{itemize}
        \item \emph{\textbf{Associative type:}} $\slj$ is  type $A_{2n-1}(1)$ and  $\jordan{J}\sim\jordan{M}_n(\algclosure{\K})$.
        \item \emph{\textbf{Clifford type:}} $\slj$ is  type $B_n(1)/D_{n}(1)$ and $\jordan{J}\sim\jordan{J}(V,b)$, \\
        for $V$ a $2n-1$ or $2n-2$ dimensional vector space and $b$ a nondegenerate form on $V$.
        \item \emph{\textbf{Symmetric Hermitian type:}} $\slj$ is  type $C_{n}(1)$ and $\jordan{J}\sim\jordan{H}_n(\algclosure{\K})$.
        \item \emph{\textbf{Quaternionic Hermitian type:}} $\slj$ is  type $D_{2n}^*(1)$ and $\jordan{J}\sim\jordan{H}_n(\mathbb{D})$, \\
        for $\mathbb{D}$  the split quaternions over $\algclosure{\K}$.
        \item \emph{\textbf{Exceptional Type:}} $\slj$ is  type $E_{7}(1)$ and  $\jordan{J}\sim\jordan{H}_3(\mathbb{O})$,\\
        for $\mathbb{O}$  the split octonions over $\algclosure{\K}$.
    \end{itemize}

\end{proposition}
\begin{proof}
    By \Cref{thm:ClassifyingRGradingConstructions}, every $A_1$-grading is induced by a parabolic subgroup. This parabolic subgroup corresponds to an isotopy class of Jordan algebras over $\algclosure{\K}$. Since $\algclosure{\K}$ is algebraically closed, there is one simple Lie algebra for each Dynkin diagram in our classification of parabolic subgroups which induce $A_1$-gradings (\Cref{thm:A1GradingsInducedByTheta}).
\end{proof}
In order to extend the result to non algebraically closed $\K$ we need to recall the classification of simple Lie algebras by the $\keyword{Satake-Tits Index}$: 

\begin{definition}
    Suppose $\mathfrak{g}$ is a Lie algebra over $\K$ with restricted root system $\Delta$. We write $\Delta^{\mathrm{alg}}$ for the Dynkin diagram of $\algclosure{\mathfrak{g}}:=\mathfrak{g}\otimes_\K\algclosure{\K}$ over the algebraic closure $\algclosure{\K}$.
    There is a subset $\Delta_\bullet \subset \Delta^{\mathrm{alg}}$ of roots whose restriction to the maximal torus of $\mathfrak{g}$ is trivial. There is also an induced action, $\diamond$, of the Galois group $\mathrm{Gal}(\algclosure{\K}/\K)$ on the non trivial roots, $\Delta_\circ:= \Delta^{\mathrm{alg}}\setminus\Delta_\bullet $. \\
    The triple $(\Delta^{\mathrm{alg}},\Delta_\circ ,\diamond) $ is called the \keyword{Satake-Tits index} of $\mathfrak{g}$.
\end{definition}

\begin{lemma}
    Suppose $\mathfrak{g}$ has restricted root system $A_1$. Then $\mathfrak{g}\simeq\slj$ for $\jordan{J}$ a division Jordan algebra over $\K$. Also, the Satake-Tits index is such that $\Delta^{\mathrm{alg}}$ is a product of identical copies of identical $A_1$ graded diagrams, $\Delta_\circ$ is the set of $A_1$ grading nodes, and $\diamond$ acts transitively on $\Delta_\circ$. 
\end{lemma}
\begin{proof}
    The fact that $\jordan{J}$ is a division Jordan algebra is proved later by analyzing the $\Theta$-Bruhat decomposition of $\mathfrak{g}$.  Since $\mathfrak{g}$ has restricted root system $A_1$, we will see that the corresponding Jordan algebra has capacity 1 (see \Cref{def:capacity}) and so every nonzero element is in the top Bruhat cell by \Cref{thm:SP2JBruhatDecomposition}. Then \Cref{thm:BruhatSigmaIntersection} will show that all of these elements are invertible.
    
    For the second part, first assume that $\diamond$ acts trivially. In this case $\Delta_\circ$ must be one element, its easy to see from the definition of the Tits index that this induces an $A_1$-grading on $\algclosure{\mathfrak{g}}$. If the $\diamond$ action is nontrivial, we can first pass to the smallest extension of $\K$ which makes the $\diamond$ action trivial, denote it by $\K'$. The number of embeddings of this field into $\algclosure{\K}$ gives the number of copies of this $A_1$ graded diagram which are permuted by the $\diamond$ action.
    \end{proof}

We call $\jordan{J}$ associative, Clifford, Hermitian, or exceptional type depending on the corresponding type of the Jordan algebra $\algclosure{\jordan{J}}$ obtained by extending scalars.

\begin{lemma}\label{thm:LiftingGradingToAlgClosure}
    If $\Delta$ is a reduced root system, then $\mathfrak{g}^{\mathrm{alg}}$ is also $\Delta$-graded. Moreover $\Theta = \Delta_\circ$ is a Jordan compatible subset of $\Delta^{\mathrm{alg}}$. The $\diamond$ action descends to $\algclosure{\Delta}_\Theta$, folding by this action gives the $\Delta$-grading. 
\end{lemma}
\begin{proof} 
    Since $\cdot \otimes_\K \algclosure{\K}$ commutes with direct sum, the $\Delta$-grading on $\algclosure{\mathfrak{g}}$ is given by applying $\cdot \otimes_\K \algclosure{\K}$ to each graded piece of $\mathfrak{g}$. The definition of the Satake-Tits index implies that the center of the Levi of the parabolic subgroup determined by $\Theta$ contains the maximal split torus of $\mathfrak{g}$. We can see that the maximal split torus is exactly the fixed points of the Galois group. Thus if the star action is trivial, then each node in $\Delta_\circ$ corresponds directly to a node in $\Delta$, each of which correspond to $A_1$ graded subalgebras, which implies that $\Delta_\circ$ is Jordan compatible. In the case that the $\diamond$ action is nontrivial, each orbit must correspond to a number of copies of the same $A_1$ graded subalgebras. This gives a folding of $\algclosure{\Delta}_\Theta$. 
\end{proof}

The previous lemma allows us to prove an ``identification theorem'' for any $R$-graded algebra which identifies which Jordan algebras can appear in each type.

\begin{theorem}[Identification]\label{thm:Identification} Let $\mathfrak{g}$ be a semisimple $R$-graded Lie algebra over $\K$. The grading identifies a Jordan algebra for each simple root $\alpha_j$ of $R$, with the following restrictions given by the Dynkin type of the simple factor containing $\alpha_j$: 
    \begin{description}
        \item [$\mathbf{A_p}$] We identify $p$ copies of an associative Jordan algebra, one for each simple root.\\
        \hspace*{-1.4pc} When $p=2$ we could instead identify two copies of an octonion algebra over $\K'$, where $\K'$ is a finite field extension of $\K$.\\
        \hspace*{-1.4pc} When $p=1$, we can additionally identify any semisimple Jordan algebra over $\K$.
        \item [$\mathbf{B_p}$] We identify $p-1$ copies of $\K'$ and one copy of a Clifford type algebra over $\K'$.
        \item [$\mathbf{C_p}$] We identify $p-1$ copies of $\jordan{J}(\mathcal{R)}$ associated to each short root and one copy of $\jordan{H}(\mathcal{R},\star)$ to the unique long root, for $\mathcal{R}$ a finite-dimensional algebra over $\K'$.\\
        \hspace*{-1.4pc} When $p=3$, we could instead identify two copies of an octonion algebra over $\K'$ for each short root and a copy of $\K'$ for the long root. 
        \item [$\mathbf{F_4}$] We identify 2 copies of $\K'$ associated to the long roots and 2 copies of $\mathcal{C}$, a Cayley-Dickson algebra over $\K'$ associated to the short roots.
        \item [$\mathbf{G_2}$] We identify one copy of $\K'$ for the long root. The short root either identifies a copy of $\jordan{H}_3(\mathcal{C})$ or a semisimple Jordan algebra consisting of $\K'$ plus a Clifford type algebra.  
        \item [$\mathbf{D_p}$] We identify $p$ copies of $\K'$, one for each simple root.
        \item[$\mathbf{E_6,E_7,E_8}$] We identify a copy of $\K'$ for each simple root.
    \end{description}
\end{theorem}
\begin{proof}
    The key idea is that the Dynkin diagrams which occur in the Satake-Tits index are in the same category of diagram as the original restricted root system. So the classification of the Jordan algebra associated to each node can be read from either the restricted root system or the Satake-Tits index. 
    For example consider the case $A_p$. By \Cref{thm:LiftingGradingToAlgClosure}, the subset $\Delta_\circ \subset \algclosure{\Delta}$ is Jordan compatible with an associated $A_p$ grading. There are only two possible choices for this: $A_n(p,k)$ and $E_6^*(2)$. In the first case, each Jordan algebra has the same Dynkin diagram and is associative type. The other case can only occur when $p=2$, and $\mathfrak{g}$ is a form of the algebra $\mathfrak{e}_6$ over $\algclosure{\K}$. To see that the Jordan algebras are octonion algebras over a field extension of $\K$, we observe the Jordan algebras must be 8-dimensional with a product that respects a quadratic norm.

    Another interesting case is type $C_p$. As before the Satake-Tits index must be one of: $C_n(p,k)$, $D_{2n}^*(p,k)$, or $E_7(3)$ with trivial $\diamond$ action, or $A_n(2p-1,k)$ via folding. In both, $C_n(p,k)$ and $D_{2n}^*(p,k)$, we see $p-1$ roots are the same associative type Jordan algebra and the last is Hermitian type given by the fixed points of an involution on the associative type Jordan algebra. Similarly the folding on $A_n(2p-1,k)$ identifies the paired nodes with an associative type Jordan algebra and the unpaired node with the Hermitian type given by the fixed points of the folding. When $p>3$, these are the only cases, but when $p=3$ there is also the case $E_7(3)$. Two roots correspond to the exceptional $A_2$ and identify octonion algebras over $\K'$ and the last root is identified by the fixed point under conjugation, the base algebra $\K'$.

    The final case we consider explicitly is type $G_2$. The Satake-Tits index must be one of $F_4(2)$, $E_6(2)$, $E_7(2)$, or $E_8(2)$ with no folding or a folding of a $B_3$ or $D_4$-graded diagram. The first cases straightforwardly identify $\K'$ for the long root and $\jordan{H}_3(\K')$, $\jordan{H}_3(\C_{\K'}^\splitt)\simeq \jordan{M}_3(\K')$, $\jordan{H}_3(\Quat^\splitt_{\K'})$, or $\jordan{H}_3(\Octo_{\K'}^\splitt)$ respectively for the short root. These are all $\jordan{H}_3(\mathcal{C})$ for some Cayley-Dickson algebra as claimed. 
    In the folded cases, we observe the singleton orbit, (which becomes the long root) is always identified to be $\K'$. In type $B_3$, the other orbit consists of a node identified to be $\K'$ and Clifford type algebra which are combined to be the Jordan algebra of the short root. In type $D_4$, three copies of $\K'$ are combined identifying the short root as $\K' \oplus \K' \oplus \K' \simeq \K' \oplus \jordan{J}(1,1)$ as claimed.

\end{proof}

\begin{remark}
    Regardless of Dynkin-type, the Jordan algebra for $\K'$, a finite field extension of $\K$ can be identified for every simple root. For example, in type $A_p$, $\K'$ is an associative Jordan algebra or in type $B_p$, $\K'$ is a one dimensional Clifford algebra over $\K'$.
\end{remark}

This identification theorem is equivalent to the recognition theorems for simple graded Lie algebras. It is interesting to compare the identification theorem with the recognition theorem in type $C_2/B_2$. In the identification theorem, we separated the Lie algebras into two families
\begin{align*}
    B_2: \K' \Leftarrow \jordan{J}(V_{\K'},b) \qquad \text{and} \qquad C_2: \jordan{J}(\mathcal{R}) \Rightarrow \jordan{H}(\mathcal{R},\star).
\end{align*}
When $\mathcal{R}^\star \simeq \K'$ is commutative, these families agree as $x \mapsto x x^\star$ defines a quadratic form on $\mathcal{R}$ realizing $\jordan{J}(\mathcal{R})$ as a Clifford type Jordan algebra. However when $\mathcal{R}^\star$ is noncommutative these families are distinct. This helps to explain the complications in the $C_2$-recognition theorem (\Cref{thm:RecognitionCp}).

\subsection{Noncommutative Lie algebras over Jordan algebras}
In the previous sections, we have started with a simple Lie algebra, $\mathfrak{g}$, and classified the possible root systems $R$ which can grade $\mathfrak{g}$. This classification identified an isotopy class of Jordan algebras for each node of the Dynkin diagram for $R$ and gave restrictions on the types of Jordan algebras that appear. 

We now reverse the process and begin with a simple root system $R$ and a compatible family of Jordan algebras $\{\jordan{J}_i\}$. Each Jordan algebra determines a Satake-Tits index of type $A_1$ which can be glued to obtain the Satake-Tits index of an $R$-graded Lie algebra, $\mathfrak{g}$, which identifies the chosen $\jordan{J}_i$. 

\begin{definition}
    We write $\ThetaCons{R}$ for the (adjoint form) of Lie algebra with $R$-grading and the collection  $\{\jordan{J}_i\}$ of Jordan algebras assigned to its nodes. If $\mathfrak{g}$ is a finite central extension of $\ThetaCons{R}$, we call it a \keyword{Jordan split algebra} of type $R$ over the Jordan algebras $\{\jordan{J}_i\}$. 
\end{definition}

$\ThetaCons{R}$ is a unification of the Tits construction and the Tits-Kantor-Koecher construction which spells out the dependence on the Jordan algebras associated to each root space.

\begin{remark}
    Gradings by folding correspond to using non-simple Jordan algebras in $\ThetaCons{R}$. For example, the $F_4$ grading of $e_6^{(6)}$ is given by $\ThetaCons{F_4}$ with $\jordan{J}_1  = \jordan{J}_2= \R$ and $\jordan{J}_{3}=\jordan{J}_4 = \C^\splitt$.

    Similarly there is a $G_2$-grading of $\so(p+3,p+4+q)$ given by $\ThetaCons{G_2}$ with long root space $\jordan{J}_1 = \R$, and short root space $\R \oplus \jordan{J}(p,p+1+q)$. This folds the $B_3$-grading given by $\ThetaConsJs{B_3}{\R,\R,\jordan{J}(p,p+1+q)}$. 
\end{remark}

\begin{definition}
    Let $\mathfrak{g}$ be Jordan split of type $R$, so that $\mathfrak{g}\simeq \ThetaCons{R}$. The number of noncommutative Jordan algebras in $\{\jordan{J}_i\}$  is the \keyword{noncommutative rank} or \keyword{NC Rank} of $\mathfrak{g}$ with respect to $R$. 
\end{definition}

By the identification theorem \Cref{thm:Identification}, there are at most two different isotopy classes of Jordan algebras in the collection $\{\jordan{J}_i\}$. Each of these isotopy classes is assigned to a short or long root. 

\begin{proposition}
    We have the following classification of simple Jordan split $\mathfrak{g}$ of type $R$ by noncommutative rank:
        \begin{enumerate}
        \item NC Rank $=0$: $\mathfrak{g}$ is split or quasi-split  and $R$ is the root system of $\mathfrak{g}$.
        \item NC Rank $=1$: $R$ is type $A_1$, $B_p$, or $G_2$ and not split or quasi-split.
        \item NC Rank $>1$: $R$ is type $A_p$, $C_p$, with $p>1$, or $F_4$ and not split or quasi-split. 
        \end{enumerate}
\end{proposition}
\begin{proof}
    First, if $\mathfrak{g}$ is split or quasi-split, then by definition there is a Borel subalgebra in $\mathfrak{g}$, this means that the structure groups of all of the Jordan algebras appearing when we choose $\Theta=\Delta$ are abelian, which means they are commutative Jordan algebras. 

    For the second statement, we look at the conditions for $B_p$ and $G_2$ and see that in either case a commutative Jordan algebra must be assigned to the long roots of the root system. Since there is only one short root remaining in each type, the non split or quasi-split type $A_1$, $B_p$, and $G_2$ have noncommutative rank 1. 

    There is a slight ambiguity since $B_2\simeq C_2$, but we will consider NC rank 1 groups of this type as $B_2$ and higher NC rank groups as $C_2$.

    For the last part we only have to consider the $C_p$ case and see that the Hermitian type algebra is the fixed points of an involution on the associative type algebra, meaning that if the associative algebra is commutative, then all are and we are back in the quasi-split case. Since $p>2$ this implies that the NC rank is bigger than 1.

    Every possible simple root system is covered by these cases, so the classification is complete.
\end{proof}
 
\begin{example}
    Over the real numbers, the noncommutative rank 0 Jordan split Lie algebras are the split real Lie algebras, the complex Lie algebras viewed as real Lie algebras and the real forms $\mathfrak{su}(n,n)$, $\so(n,n+2)$ and  $\mathfrak{e}_6^{(2)}$. 
\end{example}

\subsection{Weyl groups}\label{sec:ThetaWeylGroup_real}
When $\mathfrak{g}$ is graded by a root system $R$ there are two Weyl groups that act on $G$, the Weyl group of $R$, $W(R)$, and the Weyl group of the restricted root system of $\mathfrak{g}$. Let $\Delta$ be a choice of simple restricted roots for $\mathfrak{g}$ and let $W(\Delta)$ be the Weyl group of the restricted root system.
In this section we define two injections $W(R) \rightarrow W(\Delta)$ allowing us to understand $W(R)$ using $W(\Delta)$. 

\begin{convention}
    We use $s_1,\dots,s_n$ for the generators of $W(\Delta)$, and $\varsigma_1\dots,\varsigma_p$ for those of $W(R)$.
\end{convention}

\begin{definition}
    For $w\in W(\Delta)$ we denote by $\ell(w)$ the \keyword{length of $w$ in the generators $s_i$}.
    For $F\subset\Delta$, any subset of simple roots, we denote by $w_F$ the \keyword{longest element of the Weyl group $W(F)<W(\Delta)$} generated by $F$. We will denote by $w_0$ the longest element of $W(R)$.
\end{definition}

We begin by assuming the $R$-grading on $\mathfrak{g}$ is induced by the associated parabolic subalgebra $\mathfrak{p}_\Theta$ for some $\Theta \subset \Delta$. For $1\leq i\leq p$ we denote by $F_i$ the subset of $\Delta$ consisting of the connected component of $\Delta\setminus\Theta \cup \beta_i$ containing $\beta_i$ as in \Cref{sec:largerGrading}.\medskip

Now suppose that $\mathfrak{g}$ is $A_1$-graded and so $\Theta=\{\beta\}$ has one element. In this case we decompose $w_\Delta= s_{\delta_1} \dots s_{\delta_c} \dots s_{\delta_r}  $ where the $\delta_i$ are strongly orthogonal roots and $\delta_1,\dots,\delta_c$ are the roots which contain $\beta$ as in \Cref{lem:longest_element_orthogonal}.  We write $\varsigma^{(j)} =\prod_{i=1}^j s_{\delta_i}$.

\begin{definition}\label{def:capacity}
     The number of $\delta_i$ with $[\delta_i : \beta] > 0$ denoted by $c$ above is called the \keyword{capacity} of the Jordan algebra $\jordan{J}$.
\end{definition}

\begin{remark}
    This capacity is actually the same as the capacity defined as the maximal number of orthogonal idempotent elements you can write the identity element as a sum as, see \cite{mccrimmon-TasteJordanAlgebras}.
\end{remark}

\begin{proposition}\label{prop:double_cosets_A1}
    Let $\Theta = \{\beta\}$ be a single element. The double cosets $$[w] \in W(\Delta \setminus \Theta) \backslash W(\Delta) / W(\Delta \setminus \Theta)$$ are indexed by the elements $\varsigma^{(j)}$ for $0\leq j \leq c$. 
\end{proposition}

\begin{proof}
    We calculate the double cosets of $W(\Delta)$ by the subgroup generated by the reflections $s_{\alpha_j}$ for $\alpha_j\neq \beta$ in each possible finite Weyl type in \Cref{fig:possibleThetadynks}.
    
    In type $A_{2n-1}$ and $C_n$, $W(\Delta)$ is a subgroup of $S_{2n}$, with the $C_n$ Weyl group acting by symmetric permutations. For $f\in S_{2n}$, let $|f| =|\{1 \leq i \leq n\,\st\, f(i)>n\}|$, be the number of the left half of indices which are mapped to the right half by $f$.  Then $f,g $ are in the same double coset if and only if $|f|=|g|$. The maximum value of $|f|$ is $n$ which is the capacity, $c$. Since $|\varsigma^{(i)}| = i$, each $\varsigma^{(i)}$ is in a distinct double cosets and they cover all possible cosets. \medskip

    For the $D^*_{2n}$ case, we realize the Weyl group as a subgroup of $(\Z/2\Z)^{2n}\rtimes S_{2n}$ where $s_\beta$ is the element $\big((1,\dots,1,-1,-1),(2n-1\hspace{1mm}2n)\big)$ and the other $s_j$ act only as permutations. Therefore, the invariant of a double coset is the number of ones in the $(\Z/2\Z)^{2n}$ factor. The elements $\varsigma^{(i)}$ have $2i$ minus ones and thus are in distinct double cosets. Moreover every coset contains some $\varsigma^{(i)}$ because the number of minus ones is an even integer, $2k$, between $0$ and $2n$ and the capacity $c$ is equal to $n$ in this case.\medskip

    Next we consider the cases $D_{n+1}$ and $B_{n}$ where $\beta = s_1$. In both cases the capacity is 2 and the set of strongly orthogonal roots are $\delta_1 = \alpha_M$ and $\delta_2=\alpha_1$. The Weyl group for $B_n$ is realized by taking the realization of Weyl group of $D_{n+1}$ and ``folding'' the generators for the tails $s_n$ and $s_{n+1}$ into one generator. Concretely, the Weyl group for $B_n$ is the subgroup of $(\Z/2\Z)^{n+1}\rtimes S_{n+1}$ generated by the simple transpositions $s_1,\cdots s_{n-1}$ and $s_n = \big((1,\cdots,1,-1,-1),())$. There is an injective map $W(B_n) \rightarrow W(D_{n+1}$ given by sending $s_i \mapsto s_i$ for $i < n$ and $s_n \mapsto s_n s_{n+1}$. This map remains injective after taking the double quotient so we can study the cosets in $D_n$ and pullback the result to $B_n$. 
    
    The three cosets are identified by the sign/image of $1$: the identity coset has image $+1$, the coset of $s_1$ has image $\pm k$ for $k\neq 1$ and the coset for $w_0$ has image $-1$. We then compute $\varsigma^{(0)}$ is the identity and is in the coset for $+1$. Next, $\varsigma^{(1)} = \big((-1,-1,\cdots,1,1),(12)\big)$ and so send $1$ to $-2$ is in the second coset. Finally, $\varsigma^{(2)} = \big((-1,-1,\cdots,1,1),()\big)$ which sends $1$ to $-1$ and is in the coset of $w_0$ as claimed.\medskip

    The final case is for $E_7$ with $\beta= s_7$. It has capacity 3 with $\delta_1,\delta_2,\delta_3$ the first 3 orthogonal roots from \Cref{ex:e7orthogonal}. The double cosets can be computed in two steps. First we identify $W(\Delta)/W(\delta\setminus\Theta)$ with the orbit of the fundamental weight $\varphi_\beta$, dual to $\beta$. This orbit contains 56 elements. The full double cosets can then be identified by computing the left action of $W(\Delta\setminus \Theta)$ on these points, resulting in 4 orbits. The  ``angle'' $\langle \varphi_\beta , w \varphi_\beta\rangle$ is well defined on each double coset and takes a different value on each of the four orbits. So it suffices to check that $\langle \varphi_\beta, \varsigma^{(i)}\rangle$ takes four distinct values to prove the theorem.  
\end{proof}

Note that $\varsigma^{(c)}$ is in the same double coset as $w_\Delta$. There is a second distinguished element $\varsigma^\reduced$ in this double coset defined by the equation
\begin{equation*}
    w_\Delta = \varsigma^\reduced w_{\Delta\setminus\Theta}\,.
\end{equation*}
This element satisfies $\WeylLength(w_\Delta) = \WeylLength(\varsigma^\reduced) + \WeylLength(w_{\Delta\setminus\Theta})$ and is the minimal length element in the coset of $w_\Delta$ in $W(\Delta)/W(\Delta\setminus\Theta)$ \cite{borel-reductiveGroups,guichard2022generalizing}.

\medskip

Now suppose that $\Theta$ is  an arbitrary Jordan compatible subset. We have elements $\varsigma_i^{(j)} \in W(\Delta)$ for each $\beta_i$ obtained by applying the above construction to $F_i\subset\Delta$. We denote by $c_i$ the capacity of $\beta_i$.

\begin{proposition}
    We have an injective group homomorphism defined by
    \begin{align*}
        \phi:W(R) &\to W(\Delta)\\
        \varsigma_i &\mapsto \varsigma_i^{(c_i)}\,.
    \end{align*}
\end{proposition}
\begin{proof}
    This follows from the fact that $\mathfrak{g}$ is graded by a root system with Dynkin diagram $\Delta_\Theta$. Our explicit embedding of the split Lie algebra $\mathfrak{g}_\Theta$ satisfies that the simple root spaces of $\mathfrak{g}_\Theta$ are mapped to the sum of orthogonal root spaces. Simultaneous reflection about these root spaces in the Weyl group $W(\Delta)$ gives the element $\phi (\varsigma_i)$  for each simple root space $\beta_i$. 
\end{proof}

\begin{lemma}\label{lem:weyl_group_normalizes}
    For any $w \in W(R)$, $\phi(w)$  normalizes $W(\Delta\setminus \Theta)$. Moreover, $\phi(w_0)$ is in the same coset as $w_\Delta$ in $W(\Delta)/ W(\Delta\setminus \Theta)$.
\end{lemma}
\begin{proof}
    We check that this is true for every generator $\varsigma_i$. Conjugation by $\phi(\varsigma_i)$ acts by the opposition involution on $F_i$ and trivially on the other connected components of $\Delta\setminus \Theta$, which normalizes $W(\Delta\setminus \Theta)$ since $\beta_i$ is fixed by the opposition involution on $F_i$. The second statement follows since both $\phi(w_0)$ and $ w_\Delta$ act the same way on the split torus of $\mathfrak{g}_\Theta$. 
\end{proof}

\begin{corollary}\label{lem:product_double_cosets}
    Let $w\in W(R)$.  The product of the double cosets of $\phi(w)$ and any other word $s \in W(\Delta)$ is the double coset of $\phi(w)s$. 
\end{corollary}
\begin{proof}
    Chose representatives $u_1 \phi(w) u_2$ and $u_3 \dot{s} u_4$ for the double cosets of $\phi(w)$ and $s$. Then since $\phi(w)$ normalizes $W(\Delta\setminus \Theta)$, $\phi(w) u_2u_3 = u' \phi(w)$ for some $u' \in W(\Delta \setminus\Theta)$ and  we calculate
    \begin{align*}
        u_1\phi(w)u_2 \cdot u_3 s u_4 = u_1 u'~ \phi(w) s ~u_4 \in [\phi(w)s]\,. 
    \end{align*}
\end{proof}

 Next we will construct a second injection $W(R) \rightarrow W(\Delta)$. The image will land in the same double coset as the previous map, but will be better behaved with respect to the length of words. For each subdiagram $F_i$, we define a reduced word $\varsigma_i^\reduced$ via the following equation
\begin{equation*}
    w_{F_i} = \varsigma_i^\reduced w_{F_i\setminus\{\beta_i\}}\,.
\end{equation*}

\begin{lemma}\label{thm:varsigmaRedProperties}
    The element $\varsigma_i^\reduced$ satisfies $\WeylLength(\varsigma_i^\reduced)+\WeylLength(w_{F_i\setminus\{\beta_i\}}) =\WeylLength(w_{F_i})$ and is the minimal length element in the coset of $w_{F_i}$ in $W(\Delta)/W(\Delta\setminus F)$. Furthermore $\varsigma_i^\reduced$ is order 2 in $W(\Delta)$ and normalizes $W(\Delta\setminus\Theta)$.
\end{lemma}
\begin{proof}
    The statements about the length are proved in \cite{borel-reductiveGroups,guichard2022generalizing}. From this we know that $\phi^\reduced(\varsigma_i)$ is the shortest representative of $w_{F_i}$ in $W(F_i)/W(F_i\setminus\{\beta_i\})$. But since all nodes in $\Delta\setminus\Theta$ which are not in $F_i\setminus\{\beta_i\}$ are not connected to $F_i$, $\phi^\reduced(\varsigma_i)$ has to be the shortest representative of $w_{F_i}$ in $W(\Delta)/W(\Delta\setminus\Theta)$ as claimed. 
    
    The remaining statements are also proved by restricting to the $A_1$ case by considering $F_i\subset\Delta$: First observe that the elements $\varsigma_i^\reduced $ and $w_{F_i\setminus\{\beta_i\}}$ commute. Namely, we can move $w_{F_i\setminus\{\beta_i\}}$ past $\varsigma_i^\reduced$ letter by letter at the cost of applying the opposition involution of $F_i$. But since the opposition involution fixes $\beta_i$ and the longest element of $F_i\setminus\{\beta_i\}$ is unique, the elements commute as claimed. Since the longest elements $w_{F_i}$ and $w_{F_i\setminus\{\beta_i\}}$ are unique, they must be of order 2, and thus $\varsigma_i^\reduced$ is order 2.

    To see that $\varsigma_i^\reduced$ normalizes $W(\Delta\setminus\Theta)$ it suffices to show that it normalizes $W(F_i\setminus\{\beta_i\})$. But this is clear:
    \begin{equation*}
        \varsigma_i^\reduced W(F_i\setminus\{\beta_i\}) = w_{F_i} W(F_i\setminus\{\beta_i\}) = W(F_i\setminus\{\beta_i\}) w_{F_i} = W(F_i\setminus\{\beta_i\}) \varsigma_i^\reduced\,.
    \end{equation*}
    Again, we have used that the opposition involution of $F_i$ fixes $\beta_i$.
\end{proof}

\begin{proposition}
    We have an injective group homomorphism defined by
    \begin{align*}
        \phi^{\reduced}:W(R) &\to W(\Delta)\\
        \varsigma_i &\mapsto \varsigma_i^\reduced\,.
    \end{align*}
    This homomorphism satisfies that $\ell(\phi^{\reduced}(w_1w_2)) = \ell(\phi^{\reduced}(w_1)) +\ell(\phi^{\reduced}(w_2))$ for $w_i\in W(R)$ for which $w_1w_2$ is a reduced word in $W(R)$.
\end{proposition}
\begin{proof}
    We first prove that the image of the longest element $w_0\in W(R)$ is reduced. By definition $\phi(\varsigma_i)$, $\phi^\reduced(\varsigma_i)$ and $w_{F_i}$ are in the same coset in $W(\Delta)/W(\Delta\setminus\Theta)$. Using \Cref{thm:varsigmaRedProperties} we find that $\phi^{\reduced}(w_0)$ is in the same coset as $\phi(w_0)$ which is in the same coset as $w_\Delta$ in $W(\Delta)/ W(\Delta\setminus\Theta)$ by \Cref{lem:weyl_group_normalizes}. Moreover the minimal length element in this coset has length $\ell(w_\Delta) -\ell(w_{\Delta\setminus \Theta})$.
    
    Pick a reduced expression $w_0 = \varsigma_{i_1}\dots \varsigma_{i_N} $. We will compute the sum $\sum_j \ell\big(\phi^{\reduced}(\varsigma_{i_j})\big)$ and check that this is equal to the minimal possible length $\ell(w_\Delta) -\ell(w_{\Delta\setminus \Theta})$. Then
    \begin{equation*}
        \WeylLength\big(\phi^\reduced(w_0)\big)\leq\sum_j \ell\big(\phi^{\reduced}(\varsigma_{i_j})\big) = \ell(w_\Delta) -\ell(w_{\Delta\setminus \Theta})\leq\WeylLength\big(\phi^\reduced(w_0)\big)
    \end{equation*}
    which implies that $\phi^\reduced(w_0)$ is reduced. We will leave its explicit calculation to the end of the proof. To show that $\phi^{\reduced}$ is a homomorphism, we just need to show that it satisfies the braid relations. This follows by realizing that the braid relations are exactly different expressions of $w_0$ when $p=2$. But clearly these are satisfied since the image of $w_0$ does not depend of the choice of reduced expression we took.

    The calculation that shows $\sum_j \ell(\phi^{\reduced}(\varsigma_{i_j})) = \ell(w_\Delta) -\ell(w_{\Delta\setminus \Theta})$ can be done easily in each case. First we recall the length of the longest word in each root system:
    \begin{center}
        \begin{tabular}{l|| c|c|c|c|c|c|c|c|c}
            Type & $A_n$ & $B_n$ & $C_n$ & $D_n$ & $E_6$ & $E_7$ &$E_8$ & $F_4$ & $G_2$  \\
            \hline
            $\ell(w_0)$ & $n(n+1)/2$ & $n^2$ & $n^2$ & $n(n-1)$ & $36$ & $63$ & $120$ & $24$ & $6$
        \end{tabular}
    \end{center}
    Next we calculate all possible values of $\ell(\phi^{\reduced}(\varsigma_{i_j}))$ depending on the type of $F_{i_j}$. These are given by the following table.
    \begin{center}
        \begin{tabular}{l|| c|c|c|c|c|c}
            Type & $A_{2k+1}(1,k)$ & $B_n(1)$ & $C_{k+1}(1,k)$ & $D_n(1)$ & $D_{2k}^*(1,2k-1)$ & $E_7(1)$   \\
            \hline
            $\ell(\phi^\reduced(\varsigma_{i_j}))$ & $(k+1)^2$ & $2n-1$ & $(k+1)(k+2)/2$ & $2n-2$ & $k(2k-1)$ & $27$ 
        \end{tabular}
    \end{center}
    Finally we proceed case by case to compute the sum:
    \begin{description}
        \item [Case $A_n(p,k)$] Recall that in this case we have $(p+1)(k+1)=n+1$ since $n-k = p(k+1)$. In addition, for each $i$, $\ell(\phi^{\reduced}(\varsigma_i)) = (k+1)^2$. So\\
        \begin{minipage}{\linewidth}
            \begin{align*}
                \ell(w_\Delta) -\ell(w_{\Delta\setminus\Theta}) &= \frac{n(n+1)}{2}-\frac{(p+1)k(k+1)}{2}\\
                &= \frac{p(p+1)(k+1)^2}{2}= \sum_{j=1}^{p(p+1)/2}\ell(\phi^{\reduced}(\varsigma_{i_j}))\,.
            \end{align*}
        \end{minipage}
        \item [Case $B_n(p)$] In this case, $\Delta_\Theta$ is type $B_p$ and any reduced expression of the longest word consists of $p(p-1)$ reflections about long simple roots with $\ell(\phi^\reduced(\varsigma_i)) = 1$ and $p$ reflections about the short simple root with $\ell(\phi^\reduced(\varsigma_p)) = 2(n-p)+1$. Thus\\
        \begin{minipage}{\linewidth}
            \begin{align*}
                \ell(w_\Delta)-\ell(w_{\Delta\setminus\Theta}) &= n^2-(n-p)^2 = 2np-p^2\\
                &= p\big((p-1)+(2(n-p)+1)\big)= \sum_{j}\ell\big(\phi^{\reduced}(\varsigma_{i_j})\big)\,.
            \end{align*}
        \end{minipage}
        \item [Case $C_n(p,k)$] In this case, $n=(k+1)p$ and $\Delta_\Theta$ is a $C_p$ diagram, with long root labeled $p$. Then,\\
        \begin{minipage}{\linewidth}
            \begin{equation*}
                n^2-p\frac{k(k+1)}{2} = n(n-k/2) = p(p-1)(k+1)^2+p(k+1)(k+2)/2\,.
            \end{equation*}
        \end{minipage}\medskip
        \item [Case $D_n(p)$] In this case $\Delta_\Theta$ is type $B_p$ and the calculation is analagous to the $B_n(p)$ case.\\
        \begin{minipage}{\linewidth}
            \begin{equation*}
                n(n-1)-(n-p)(n-p-1) = 2np-p^2-p = p\Big((p-1)+\big(2(n-p+1)-2\big)\Big)\,.
            \end{equation*}
        \end{minipage}\medskip
        \item [Case $D^*_{2n}(p,2k-1)$]  In this case,  $2kp= 2n$ and $\Delta_\Theta$ is a $C_p$ diagram, with the long root labeled $p$. Then\\
        \begin{minipage}{\linewidth}
            \begin{align*}
                \begin{aligned}
                    \ell(w_\Delta) &= 2n(2n-1)\\
                    \ell(w_{\Delta\setminus \Theta}) &= pk(2k-1)=n(2k-1)
                \end{aligned}
                \qquad    
                \ell\big(\phi^{\reduced}(\varsigma_i)\big) = \begin{cases}
                    4k^2 & i\neq p\\
                    k(2k-1) & i = p \,.
                \end{cases}
            \end{align*}
        \end{minipage}
        In a reduced expression of $w_0$, there are $p^2$ terms total and $p$ terms which are $\varsigma_p$. This means we have\\
        \begin{minipage}{\linewidth}
            \begin{equation*}
                (p^2-p)(4k^2) + pk(2k-1) = 4n^2 - 4nk+2nk-n =2n(2n-1) - n(2k-1)
            \end{equation*}
        \end{minipage}\medskip\\
        as desired. 
    \end{description}
\end{proof}

\begin{lemma}\label{lem:reduced_words_inthetaweylgroup}
    For any reduced word $w\in W(R)$ and any reduced word $s\in W(\Delta\setminus\Theta)$, $\phi^{\reduced}(w)s$ is a reduced word in $W(\Delta)$. In particular if $w\varsigma_i$ is reduced in $W(R)$, then $\phi^{\reduced}(w)\varsigma_i^{(j)}$ is reduced for any $1\leq j\leq c_i$.
\end{lemma}
\begin{proof}
    We can complete $w$ to $w_0$ by $w'w=w_0$ and $s$ to $w_{\Delta\setminus\Theta}$ by $ss' =w_{\Delta\setminus\Theta} $. Since $\phi^{\reduced}(w_0)w_{\Delta\setminus\Theta}$ is reduced this proves the result. 
    For the second statement we notice that $\phi(\varsigma_i) = \varsigma_i^{(c_i)} = \phi^{\reduced}(\varsigma_i)s $ for some $s\in W(\Delta\setminus\Theta)$. Now we apply the first part to $w\varsigma_i$ to see that  $\phi^{\reduced}(w\varsigma_i)s =\phi^{\reduced}(w)\varsigma_i^{(c_i)} $ is reduced, and finally see that $\phi^{\reduced}(w)\varsigma_i^{(j)}$ is reduced since it is a subword of $\phi^{\reduced}(w)\varsigma_i^{(c_i)}$.
\end{proof}

Finally we return to the case of arbitrary $R$-gradings. By \Cref{thm:ClassifyingRGradingConstructions}, at worst the associated parabolic $\mathfrak{p}_\Theta$ induces an $R'$-grading on $\mathfrak{g}$ with $R'$ folding onto $R$. The simple roots of $R'$ that are grouped by the folding have commuting reflections and so we have a natural injective map $W(R) \rightarrow W(R')$ sending $\varsigma^R_\beta$ to $\prod \varsigma^R \alpha_i$. Post composing with the injections for $R'$ define the required injections for any grading. 

\section{Jordan split groups}\label{sec:LieGroupsWithRootGradings}
We now focus on the structure of algebraic groups with Lie algebra graded by a root system $R$. We identify a ``good'' set of generators to understand $G$ via the root system $R$ in analogy to a choice of Chevalley basis of a split group. These generators have simple relations which can be used to lift the Weyl group of $R$ to $G$ and compute the Gauss decomposition of products of elements. Finally we generalize the fundamental weights of a root system to functions which take values in the Jordan algebras. These Jordan weights will provide the coordinate functions needed to define cluster structures in \Cref{part:ClusterlikeCoordinates}.

\begin{definition}\label{def:JordanSplitGroup}
    An algebraic group $G$ is \keyword{Jordan split} of type $R$ if its Lie algebra is (a central extension of) $\ThetaCons{R}$ for some collection of Jordan algebras $\{\jordan{J}_i\}$.
\end{definition}
Since Lie algebras can be graded by root systems of several types, the same group $G$ can be Jordan split in several ways. 
\begin{example}\leavevmode
\begin{enumerate}
    \item The group $\SL_6(\R)$ is Jordan split in three ways. Its Lie algebra $\sl_6(\R)$ is isomorphic to $\ThetaConsJs{A_5}{\R,\R,\R,\R,\R}$, $\ThetaConsJs{A_3}{\jordan{M}_2(\R),\jordan{M}_2(\R),\jordan{M}_2(\R)}$, and $\ThetaConsJs{A_2}{\jordan{M}_3(\R),\jordan{M}_3(\R)}$.
    \item Also $\SP_6(\R)$ is Jordan split in two ways, $\sp_6(\R) \simeq \ThetaConsJs{C_3}{\R,\R,\R} \simeq \ThetaConsJs{A_1}{\jordan{H}_3(\R)}$.
    \item For a final example consider $\SO(3,8)$. It is Jordan split in three ways as well, $\so(3,8) = \ThetaConsJs{B_3}{\R,\R,\jordan{J}(0,5)} \simeq \ThetaConsJs{B_2}{\R,\jordan{J}(1,6)} \simeq \ThetaConsJs{A_1}{\jordan{J}(2,7)}$.
\end{enumerate}
\end{example}

Given a Jordan split group $G$ of type $R$, our key construction is what we call the Jordan pinning of the group.  As motivation, we first look at groups with an $A_1$-graded Lie algebra.

\subsection{\texorpdfstring{$A_1$-graded groups}{A1 graded groups}}\label{sec:A1gradedGroups}
The recognition theorem for $A_1$ allows us to reduce the study of $A_1$-graded Lie algebras to the Lie algebras $\slj$ for some Jordan algebra $\jordan{J}= (V,\iota,\Id)$. 

\begin{definition}
     We write $\SLJ$ for the simply connected algebraic group over $\K$ with Lie algebra $\slj$.
\end{definition}

In \Cref{tab:SimpleJordanAlgebras}, we list the complex and real algebraic groups which are of the form $\SLJ$ for some $\jordan{J}$.
\begin{table}[htb]
    \centering
    \begin{tabular}{c|c|c||c|c|c}
        Degree & $\SL_2(\jordan{J}^\C)$  & $\jordan{J}^\C$ & $\SL_2(\jordan{J}^\R)$ & $\jordan{J}^\R$   & Restricted Type\\
        \hline \hline
         1 & $\SL_2$ & $\C$ & $\SL_2(\R)$ & $\R$   & $A_1$ \\
         \hline
         2 & $\Spin(s+q+2)$ & $\jordan{J}(s+q)$  & $\Spin(s+1,q+1)$ & $\jordan{J}(s,q)$  & $B_n/D_n$ \\
         \hline
         \multirow{2}{*}{3} & \multirow{2}{*}{$E_7$} & \multirow{2}{*}{$H_3(\Octo^\C)$} & $E_7^{(-25)}$ &  $H_3(\Octo)$  & \multirow{2}{*}{$E_7$} \\   
         &  & & $E_7^{(7)}$ & $H_3(\Octo^\splitt)$& \\
         \hline
         n  & $\SP_{2n}$ & $H_n(\C)$ & $\SP_{2n}(\R)$ & $H_n(\R)$   & $C_n$\\
           \hline
          \multirow{2}{*}{n}  &  \multirow{2}{*}{$\SL_{2n}$} &  \multirow{2}{*}{$M_n(\C)$} & $\mathrm{SU}(n,n)$ & $H_n(\C)$  &  \multirow{2}{*}{$A_{2n-1}$} \\ 
           & &  & $\SL_{2n}(\R)$ & $M_n(\R)$  & \\
           \hline
          \multirow{2}{*}{n} & \multirow{2}{*}{$\Spin(4n)$} & \multirow{2}{*}{$H_n(\Quat^\C)$} & $\Spin^*(4n)$ & $H_n(\Quat)$   & \multirow{2}{*}{$D_{2n}$} \\
          & & & $\Spin(2n,2n)$ &  $H_n(\Quat^\splitt)$  & \\
          \hline
         2n& $\SP_{4n}$& $H_{2n}(\C)$ & $\SP(2n,2n)$ & $H_n^*(\Quat)$   & $C_{2n}$ \\
         \hline
         2n& $\SL_{2n}\times\SL_{2n}$& $M_{n}(\C)\times M_n(\C)$ & $\SL_{2n}(\C)$ & $M_n(\mathbb{C})$   & $A_{2n-1}\times A_{2n-1}$\\
         \hline
         2n& $\SL_{4n}$& $M_{2n}(\C)$ & $\SL_{2n}(\Quat)$& $M_n(\Quat)$   & $A_{4n-1}$\\
    \end{tabular}
    \caption{$\SLJ$ for real and complex Jordan algebras.}
    \label{tab:SimpleJordanAlgebras}
\end{table}

\begin{definition}\label{def:structureGroup}
    The \keyword{structure group}, denoted by $\Gamma(\jordan{J})$, is the Levi factor of the parabolic subgroup associated to the $A_1$ grading of $\SLJ$. 
\end{definition}

The structure group covers the inner structure group and will replace the inner structure group moving forward. The map $\iota$ from $\jordan{J}$ will extend to an injective map $\tilde{\iota} : V^\times \rightarrow \Gamma(\jordan{J})$. 

When $\jordan{J}$ is special, we can relate the group $\SLJ$ to (a covering of) a subgroup of a general linear group over the ring $\mathcal{R}$ containing $\jordan{J}$. When $\jordan{J}$ is of Hermitian type, $\SLJ$  covers $ \SP_2(\mathcal{R},\dagger)$ the generalized symplectic group studied in \cite{rogozinnikov2020symplectic,alessandrini2022-SP2Asigma,GreenbergEtAl2024MathrmSL_2}. For $\jordan{J}$ of Clifford type, $\SLJ$ actually is a subgroup of the generalized symplectic group \cite{rogozinnikov2025}.

\begin{lemma}
    Suppose that $\jordan{J} \subset\jordan{J}(\mathcal{R})$ is a special Jordan algebra for some finite dimensional algebra $\mathcal{R}$ over $\K$. Then $\SLJ$ covers the subgroup of $\GL_2(\mathcal{R})$ which is generated by the elements 
    $ \begin{bmatrix}
        1 & v\\ 0 &1
    \end{bmatrix}$ and $ \begin{bmatrix}
        1 & 0\\ v &1
    \end{bmatrix}$ for $v\in V$.
\end{lemma}
\begin{proof}
    The Lie algebra $\slj$ is a subalgebra of $\End_2(\mathcal{R})$ when $\jordan{J}$ is special by \Cref{thm:Simple_slj}. It is generated by the elements $v_{(1)}, w_{(-1)}$ since 
    \begin{equation*}
        \left[v_{(1)},w_{(-1)} \right] = \begin{bmatrix}&v\\ & \end{bmatrix} \begin{bmatrix}&\\w & \end{bmatrix} -\begin{bmatrix}&\\w& \end{bmatrix} \begin{bmatrix}&v\\ & \end{bmatrix} = \begin{bmatrix}vw&\\ &-wv \end{bmatrix} = \iota(v,w)_{(0)}\,.
    \end{equation*}
    Thus the subgroup generated by upper and lower triangular unipotent matrices has Lie algebra $\slj$ and is covered by $\SLJ$.
\end{proof}

Next observe that for each $v \in V^\times$, the diagonal matrix $\begin{bsmallmatrix}v&\\&v^{-1}\end{bsmallmatrix}$ is in $\SLJ$ by the calculation
\begin{equation*}
        \begin{bmatrix}v & \\ & v^{-1} \end{bmatrix} = \begin{bmatrix}1 & \\ -v^{-1}& 1 \end{bmatrix} 
        \begin{bmatrix}1 & v\\  & 1 \end{bmatrix} \begin{bmatrix}1 & -1 \\ & 1 \end{bmatrix} \begin{bmatrix}1 & \\1& 1 \end{bmatrix} \begin{bmatrix}1 & -1 \\ & 1 \end{bmatrix} \begin{bmatrix}1 & v^{-1}\\ & 1 \end{bmatrix}\,.
    \end{equation*}

In order to describe non special Jordan algebras, we consider two natural injections, $x: V \to \SLJ$ and $y: V\to \SLJ$ given by the $A_1$-grading of $\slj$:
\begin{equation*}
 x(v) = \exp\big((v)_1\big)  \qquad y(v) = \exp\big((v)_{-1}\big)\,.
\end{equation*}
Then the previous calculation inspires the following definitions:

\begin{definition}\label{def:SLJ_JordanThings}\leavevmode
\begin{itemize}
    \item The \keyword{lifted reflection}, $\bar\varsigma_\beta$, is the element $x(-\Id)y(\Id)x(-\Id)$. 
    \item The \keyword{structure map}  $\tilde{\iota}:V^\times \to \Gamma(\jordan{J})$ is defined by $\tilde{\iota}(v)= y(-v^{-1})x(v)\bar\varsigma_\beta x(v^{-1})$.
    \item The \keyword{structure involution} is the anti-involution $\sigma:\Gamma(\jordan{J})\to\Gamma(\jordan{J})$ given by 
    \begin{equation*}
        \sigma(A)= \bar\varsigma_\beta {A}^{-1} \bar\varsigma^{-1}_\beta\,.
    \end{equation*}
\end{itemize}
    
\end{definition}

The structure involution here extends the involution on the inner structure group given in \Cref{def:structureInvolution}. The norm map defined in \Cref{def:JordanNorm} extends to the structure group as a character given on elements of the form $\tilde{\iota}(v)$ by $N(v)$.

\begin{definition}
    Suppose that $\jordan{J}$ is special and let $\mathcal{R}$ be its universal specializing envelope.
\begin{itemize}
    \item The \keyword{general structure group} $\mathrm{G}\Gamma(\jordan{J)}$ is the subgroup of $\mathcal{R}^\times \times \mathcal{R}^\times$ which preserves the vector space $V$ by the action $(a,b)\cdot v = avb^{-1}$. 
    In the case of an exceptional Jordan algebra, there is a form of the group $E_8$ which contains the group $\SLJ$ as a subgroup. In this case, we define the general structure group to be the Levi subgroup of this group which is a central extension of $\Gamma(\jordan{J})$.
    \item  The \keyword{adjoint structure group} $\mathrm{P}\Gamma(\jordan{J})$ is the general structure group modulo the subgroup of trivially acting elements. 
\end{itemize}
\end{definition}

When $\jordan{J}$ is special the map $\iota$ on $V^\times$ can be recovered from the composition of $\tilde{\iota}$ with passing to the quotient $\mathcal{R}^\times\times \mathcal{R}^\times\to\mathcal{R}^\times\times \mathcal{R}^\times/\{(-1,-1)\}$.

\begin{remark}
    What we are calling the structure group is different than what is referred to as the structure group in \cite{koecher-JordanAlgebras}. 
\end{remark}

\begin{definition}
    There is a cocharacter $\check{\beta}:\K^\times \to \Gamma(\jordan{J})$ defined by $\check{\beta}(k):=\tilde{\iota}(k\Id)$.
\end{definition}

The cocharacter extends to an injective map from $\Gamma(\jordan{J})$ into $\SLJ$ whose image is the associated Levi subgroup. Then together with $x: V \to \SLJ$ and $y: V\to\SLJ$, these maps ``pin'' $\SLJ$. This is analogous to the pinning of $\SL_2(\K)$ and when $\jordan{J}$ is special, they have the same matrix form
\begin{align*}
    x(v) = \begin{bmatrix}1 & v\\& 1\end{bmatrix} \qquad \check{\beta}(a) = \begin{bmatrix}a & \\& a^{-1}\end{bmatrix} \qquad y(v) = \begin{bmatrix}1 & \\v & 1\end{bmatrix}\,.
\end{align*}
In fact the triple $(x(\Id), \check{\beta}(\Id), y(\Id))$ pin a distinguished $\SL_2(\K)$ subgroup.

\begin{definition}
    Let $G$ be Jordan split of type $A_1$. So there is an isomorphism $\phi: \slj \to \mathfrak{g}$ for some $\jordan{J} = (V,\iota,\Id)$. Since $\SLJ$ is simply connected $\phi$ induces a map $\phi_*: \SLJ \rightarrow G$. The \keyword{$A_1$-Jordan pinning} of $G$ is given by the collection of maps $\phi_* \circ x: V\to G$, $\phi_* \circ \check{\beta} : \Gamma(\jordan{J}) \to G$, $\phi_*\circ y: V \to G$.
\end{definition}

Since $\slj$ is independent of isotopy we can freely choose the identity element of $\jordan{J}$ when choosing an isomorphism $\mathfrak{g} \simeq \slj$ or equivalently choosing the pinning. We define a  \keyword{normalized pinning} using the decomposition of $w_0$ into strongly orthogonal roots $\{\delta_i\}$ constructed \Cref{lem:longest_element_orthogonal}. A normalized pinning is one where $x(\Id) = \sum_{i=1}^c e_i$ for $e_i$  the root space $\mathfrak{g}_{\delta_i}$. For the remainder of the paper we only consider normalized pinnings.\medskip

Next we lift the Weyl group of the $A_1$-grading to $G$. This is an element $\overline{\varsigma}_\beta$ with $\overline{\varsigma}_\beta^2$ in the center of $G$. The element $\overline{\varsigma}$ from \Cref{def:SLJ_JordanThings} is such a lift,
\begin{align*}
    \overline{\varsigma}_\beta = x\left(-\Id\right)y\left(\Id\right)x\left(-\Id\right) \,.
\end{align*}
Since $\big(x(\Id), \check{\beta}(\Id), y(\Id)\big)$ form an $\SL_2$-triple by construction, $\overline{\varsigma}_\beta$ is the classic lift of the Weyl group for $\SL_2$ \cite[Section 3.1]{goncharov2019quantum}.\medskip

\begin{lemma}
    The lifted reflection, $\overline{\varsigma}_\beta$, conjugates $x(v)$ into $y(-v)$.
\end{lemma}
\begin{proof}
    The sub Jordan algebra generated by $v$ is special by \Cref{thm:Shirshov-Cohn2GeneratorsSpecial}, so the calculation can take place in $\GL_2(\mathcal{R})$ for some specializing ring $\mathcal{R}$. Then
    \begin{align*}
        \overline{\varsigma}_\beta = \begin{bmatrix}
            0 & -1\\
            1 & 0
        \end{bmatrix}\qquad \text{and} \qquad \overline{\varsigma}_\beta^{-1}x(v)\overline{\varsigma}_\beta = \begin{bmatrix}0 & 1\\ -1 & 0\end{bmatrix} \begin{bmatrix} 1 & v\\ 0 & 1\end{bmatrix} \begin{bmatrix}0 & -1\\ 1 & 0\end{bmatrix} = \begin{bmatrix} 1 & 0\\ -v & 1 \end{bmatrix} = y(-v)\,.
    \end{align*}
\end{proof}

 Our next goal is to understand when a generic product of images of pinning maps and $\overline{\varsigma}_\beta$ can be decomposed as $y(w)\check{\beta}(a)x(v)$. We call the subset of elements which can be decomposed this way \keyword{Gauss decomposable}.

\begin{proposition}\label{thm:decompose_A1}
Let $v \in \jordan{J}$ be invertible.
\begin{enumerate}
    \item $x(v)\overline{\varsigma_\beta} = \overline{\varsigma_\beta}^{-1}y(-v) = y(v^{-1})\check{\beta}(\iota(v))x(-v^{-1}) $
\end{enumerate}
    If $w\in \jordan{J}$ is also invertible such that $(w+v)$ is invertible as well, then 
    \begin{enumerate}[resume]
        \item $(w^{-1}+v^{-1})$ is also invertible. 
        \item $\iota(w+v)\iota(w)^{-1}=\iota(v)\iota(w^{-1}+v^{-1})$
        \item $ x(v) y(w^{-1}) = 
y\big((v+w)^{-1}\big)~\check{\beta}\big(\iota(w+v)\iota(w)^{-1}\big) x\left((w^{-1}+v^{-1})^{-1}\right)$.
    \end{enumerate}
\end{proposition}

\begin{proof}
Denote by $\jordan{J}_{v,w}$ the sub Jordan algebra generated by $v,w$ which is special by \Cref{thm:Shirshov-Cohn2GeneratorsSpecial}. 
    These relations live entirely in the subgroup $\SL_2(\jordan{J}_{v,w})$ of $\SLJ$. Thus each statement can be checked in $\GL_2(\mathcal{R})$ for $\mathcal{R}$ the ring containing the subgroup generated by $v,w$.
    For example in $\mathcal{R}$, we have $(w^{-1}+v^{-1}) = v^{-1}(w+v)w^{-1}$ and so $w^{-1}+v^{-1}$ is invertible when $v,w$ and $v+w$ are. The final computation reduces to confirming 
    \begin{align*}
        \begin{bmatrix}1 & v\\ & 1\end{bmatrix}\begin{bmatrix}1 & \\ w^{-1}& 1\end{bmatrix} =
        \begin{bmatrix}1 & \\(v+w)^{-1} & 1\end{bmatrix}\begin{bmatrix}(w+v)w^{-1} & \\ & (w+v)^{-1}w\end{bmatrix}\begin{bmatrix}1 & (w^{-1}+v^{-1})^{-1}\\ & 1\end{bmatrix}\,. 
    \end{align*}
\end{proof}

This proposition gives an interpretation of the $\iota$ map in terms of the group $\SLJ$; $\iota(v)$ is the element of the structure group which appears when decomposing $x(v)\overline{\varsigma_\beta}$.

\medskip
Next we want to understand which elements of $G$ are Gauss decomposable.  This will be closely related to the double cosets of $P_{\{\beta\}}^\opp$ in $G$, which is called the $\Theta$-Bruhat decomposition of $G$ for $\Theta = \{\beta\}$. We will see the coset of $\overline{\varsigma}_\beta$ is dense in $G$. Moreover for any $g \in P_\Theta^\opp \overline{\varsigma}_\beta P_\Theta^\opp$, the element $g \overline{\varsigma}_\beta $ is Gauss decomposable. 

We construct elements $\overline{\varsigma}_\beta^{(j)}$ for $j \in [0,c]$ which will represent each coset. Recall that $\Id = \sum_{i=1}^c e_i$ for a choice of $e_i$ in each root space $\mathfrak{g}_{\delta_i}$.

\begin{lemma}
    For any $1\leq i \leq c$, $\{(e_i)_1, (e_i)_{-1}, \iota(e_i,e_i)\}$ is an $\sl_2$-triple. The $\sl_2$-subalgebras generated by these commute with each other.
\end{lemma}

\begin{proof}
    By abuse of notation, we write $e_i:=(e_i)_1$, and $f_i:=(e_i)_{-1},\, h_i:=[e_i,f_i]$. By construction, the subspace spanned by any such triple is closed under the Lie bracket and therefore a copy of $\sl_2$.
    
    These commute since the roots $\delta_i$ are strongly orthogonal. Namely, for $i\neq j$, we have $0 = [e_i,e_j] = [f_i,f_j] =[e_i,f_j]$ and thus
    \begin{equation*}
        [h_i,e_j] = \big[[e_i,f_i],e_j\big] = - \big[[f_i,e_j],e_i\big] - \big[[e_j,e_i],f_i] = 0\,.
    \end{equation*}
    Again since the roots $\delta_i$ are strongly orthogonal, we have
    \begin{equation*}
        \left[\sum_i e_i,\sum_j f_j\right] = \sum_i [e_i,f_i] = \sum_i h_i\,.
    \end{equation*}
    By assumption $\left(\sum_i e_i, \sum_j f_j, [\sum_i e_i,\sum_j f_j]\right)$ is an $\sl_2$-triple. Thus
    \begin{equation*}
        2\sum_i e_i = \left[\sum_j h_j, \sum_i e_i\right] = \sum_i [h_i,e_i]\,.
    \end{equation*}
    The equality has to hold for the individual summands since the triples span disjoint subalgebras. Analogously, we get $[h_i,f_i]=-2f_i$.
\end{proof}

We now set
\begin{equation*}
    e^{(j)} = \sum_{i=1}^j e_i
\end{equation*}
so that $\Id = e^{(c)}$. Again, the elements $(e^{(j)})_1$, $ (e^{(j)})_{-1}$, $ \iota(e^{(j)},e^{(j)})$ form an $\sl_2$-triple. For any algebraic group $G$ with Lie algebra $\sl_2(\jordan{J})$ we define
\begin{equation*}
    \overline{\varsigma_{\beta}}^{(j)}=x\left(-e^{(j)}\right)y\left(e^{(j)}\right)x\left(-e^{(j)}\right)\,.
\end{equation*}
Note that $\overline{\varsigma_{\beta}} =\overline{\varsigma_{\beta}}^{(c)}$. 

\noindent We can use the elements $\overline{\varsigma_\beta}^{(j)}$ to decompose any algebraic group which is locally isomorphic to $\SLJ$:
\begin{proposition}[$\Theta$-Bruhat decomposition]\label{thm:SP2JBruhatDecomposition}
    Let $G$ be a algebraic group whose Lie algebra is $\slj$. Denote by $\delta_1,\dots,\delta_c$ the strongly orthogonal roots of the decomposition in \Cref{lem:longest_element_orthogonal} which contain $\beta$. Then the Bruhat decomposition of $G$ with respect to  $P_\Theta^\opp$ is given by
    \begin{equation*}
        G= \coprod_{j=0}^c P_\Theta^\opp \overline{\varsigma_\beta}^{(j)} P_\Theta^\opp
    \end{equation*}
\end{proposition}
\begin{proof}
    It is a standard fact \cite{warner2012harmonic} that
    \begin{equation}\label{eq:thetaBruhatDecomp}
        G = \coprod_{w\in W(\Delta\setminus\Theta)\backslash W(\Delta)/W(\Delta\setminus\Theta)} P_\Theta^\opp \dot{w} P_\Theta^\opp\,,
    \end{equation}
    where $\dot{w}$ is any lift of a representative in $W(\Delta)$ of the double coset to $G$. This is independent of the lift.
    
    From \Cref{prop:double_cosets_A1} we know that the elements $\varsigma_\beta^{(j)}:=\prod_{i=1}^j s_{\delta_i}$ for $0\leq i\leq c$ represent the double cosets. Thus, it suffices to prove that $\overline{\varsigma_\beta}^{(j)}$ is indeed a lift. Because the roots $\delta_i$ are strongly orthogonal we can rearrange the product to obtain
    \begin{equation*}
        \overline{\varsigma_\beta}^{(j)} = \prod_{i=1}^j x(-e_i) \cdot \prod_{i=1}^j y(e_i) \cdot \prod_{i=1}^j x(-e_i) = \prod_{i=1}^j x(-e_i)y(e_i)x(-e_i)\,.
    \end{equation*}
    Each factor in the final expression is a lift of the reflection $s_{\delta_i}$ so $\overline{\varsigma}_\beta^{(i)}$ lifts the correct word.
\end{proof}

\begin{definition}
    The \keyword{rank} of an element $v\in V$ is the $j$ such that $x(v) \in P_\Theta^\opp\varsigma_\beta^{(j)}P_\Theta^\opp$. We write $V^{(j)}$ for the subset of rank $j$.
\end{definition}
In \Cref{thm:BruhatSigmaIntersection} we will see that $V^\times= V^{(c)}$.

\begin{remark}
    Recall the maps $\phi$ and $\phi^\reduced: W(R)\to W(\Delta)$ defined in \Cref{sec:ThetaWeylGroup_real}. For $G$ as before, we have $W(R) = \langle\varsigma_\beta\rangle \cong \Z/2\Z$. The top cell of the Bruhat decomposition is
    \begin{equation*}
        P_\Theta^\opp\varsigma_\beta^{(c)}P_\Theta^\opp = P_\Theta^\opp\phi(\varsigma_\beta)P_\Theta^\opp = P_\Theta^\opp\phi^\reduced(\varsigma_\beta)P_\Theta^\opp\,.
    \end{equation*}
\end{remark}

In the following sections, we will generalize everything above to a Jordan split group over an arbitrary simple root system. 

\subsection{Jordan pinning}
\label{sec:ThetaPinning}
In this section, we define a Jordan pinning for a algebraic group, which plays a key role in the remainder of the paper. We will use it to  understand generic elements of the group and eventually to define coordinates on the associated flag varieties and local systems. 
\begin{definition}
    Let $G$ be a Jordan split group over a root system $R$ indexed by simple roots $\Theta = \{\beta_1,\cdots \beta_p\}$. Let $\jordan{J}_i$ be a representative in the isotopy class of Jordan algebras such that $\mathfrak{g} = \ThetaCons{R}$.  
    A \keyword{Jordan pinning} of $G$ is a collection of homomorphisms 
    $\check{\beta}_i:\Gamma(\jordan{J}_i) \to G$ and maps $x_i,y_i : V_i \rightarrow G$ indexed by $i \in \Theta$
    such that 
        for $i \in \Theta$ the maps $x_i,y_i,\check{\beta}_i$ form a homomorphism $\varphi_i : \SLJi \rightarrow G$ by
        \begin{equation*}
                \varphi_i\left(A\right) = \check{\beta}(A), \quad \varphi_i\left(
                \exp(v_{(1)})\right) = x_i(v),\quad \varphi_i\left(\exp(v_{(-1)}\right)) = y_i(v)\,.
        \end{equation*}
    We call the Jordan pinning \keyword{normalized} if it furthermore satisfies that for each $\beta_i \in \Theta$ we have
    \begin{equation*}
        x_i(\Id_i) \in \exp(\mathfrak{g}_{\delta_1} \oplus \cdots \oplus \mathfrak{g}_{\delta_{c_i}}) \quad\text{and}\quad y_i(\Id_i) \in \exp(\mathfrak{g}_{-\delta_1} \oplus \cdots \oplus \mathfrak{g}_{-\delta_{c_i}})\,.
    \end{equation*}
\end{definition} 

If a Jordan pinning is not normalized we can always change $\jordan{J}_i$ by an isotopy to make it so. This does not affect the homomorphisms from $\SLJi$, so we will always take normalized Jordan pinnings.

\begin{definition}
     The \keyword{Levi factor}, $L_\Theta$, of a Jordan pinning of $G$ of type $R$ is the subgroup generated by $\{\check{\beta}_i(A) \st A \in \Gamma(\jordan{J}_i), i \in \Theta\}$. This is also the Levi subgroup of the associated parabolic subgroup, $P_\Theta$, of the grading (\Cref{def:ParabolicAssociatedToGrading}). 
\end{definition}

\begin{definition}\label{def:BetaMaps}
    The Jordan pinning induces \keyword{extended root maps} $\beta_i: L_\Theta\to \mathrm{P}\Gamma(\jordan{J}_i)$ from the Levi factor to the adjoint structure group for each $i$ by $lx_i(v)l^{-1} = x_i\left(\beta_i(l)(v)\right) $. These extend the usual root maps from $Z( L_\Theta) \to \K^\times$ given by the inclusion of the split group of type $R$ into $G$.
\end{definition}

\begin{remark}
    $L_\Theta$ acts on the set of Jordan pinnings by conjugation. We can view this action as changing each $\jordan{J}_i$ by an isomorphism. Therefore, the action of $L_\Theta$ on the space of pinnings identifies a twisted form for each $\jordan{J}_i$. We expect that the set of pinnings with respect to a given set of twisted forms is an $L_\Theta$-torsor as in the classic pinnings of a split group.
\end{remark}

Suppose that $\beta_i$ and $\beta_j$ are connected in $R$ and $\jordan{J}_i$ is one dimensional (and hence isomorphic to $\K$).
Now fix a Jordan pinning of $G$. The composition $\beta_i\circ \check{\beta}_j$ is a character of $\Gamma(\jordan{J}_j)$. We can use this to define a map $N:V_j^\times \to \K$ by $N(v) = \beta_i(\check{\beta}_j(\iota(v)^{-1}))$. In other words $N(v)$ satisfies the following equation for all $x$ in $\mathfrak{g}_{\beta_i}=V_i$:
\begin{equation*}
    \forall x\in\mathfrak{g}_{\beta_i}:\quad\check{\beta}_j\big(\iota(v)\big)^{-1} \cdot x \cdot \check{\beta}_j\big(\iota(v)\big) = N(v)x\,. 
\end{equation*}
We can extend $N$ to all of $V_j$ by taking $N(v) = 0$ when $v \notin V^\times$.

\begin{lemma}\label{def:norm_on_rootspace}
    The map $N$ above is a norm map of the Jordan algebra $\jordan{J}_i$ of degree 1,2,3 if the subdiagram $F_{ij}$ is $A_2$, $B_2$ or $G_2$. 
\end{lemma}
\begin{proof}
    To check that $N$ is a norm map, we need to verify the properties from \Cref{def:JordanNorm}. By definition and the fact that $\iota(\Id) = 1 \in \Gamma(\jordan{J}_j)$,
    \begin{align*}
        N(\Id) = \beta_i(\check{\beta}_j(\iota(\Id))) = \beta_i(\check{\beta}_j(1)) = \beta_i(1) = 1\,.
    \end{align*}
    On $V_j^\times$, the map is nonzero because the roots $\beta_i$ and $\beta_j$ are connected in $R$. Furthermore $N(v)$ is defined to be zero for noninvertible elements, so $N(v) \neq 0$ if and only if $v\in V^\times$ as needed.
    Finally since $\beta_i \circ \check{\beta}_j$ is a character, $N$ is multiplicative. So
    \begin{align*}
        N(\iota(v)(w)) = (\beta_i\circ \check{\beta}_j)\big(\iota(\iota(v)(w))\big) = (\beta_i\circ \check{\beta}_j)(\iota(v)\iota(w)\iota(v)) = N(v)\cdot N(w)\cdot N(v)\,,
    \end{align*}
    and $N(\iota(v)(w)) = N(v)^2N(w)$ as needed.

    The final property to check is that $N$ is given by $\tilde{N}(v,\cdots,v)$ for some $\tilde{N} \in \Sym^r(V^*)$ for $r=2,3$. To see this, define $\tilde{N}$ using the Lie algebra of $G$. By abuse of notation denote by $v  = \log(x_j(v)), 1_i = \log(x_i(\Id_i))$. The brackets $\ad^r((\Id_j)_j)(1)_i = I$ gives an element in the 1 dimensional root space $\beta_i +r\beta_j$. Define $\tilde{N}$ by $\tilde{N}(v_1,\dots,v_r) \cdot I =[v_1,[v_2,\dots,[v_r,1]\dots]] $. The Jacobi identity implies that this is an element in $\Sym^r(V^*)$. To see that $\tilde{N}(v,\dots,v) = N(v)$ we conjugate everything by $\check{\beta}_j(\iota(w))$. This acts on the root space $\beta_i$ by multiplication by $N(w)^{-1}$ and on the root space $\beta_i+r\beta_j$ by multiplication by $N(w)$, and on the root space $\beta_j$ by $v\to \iota(w)(v)$. Therefore we find $\tilde{N}(\iota(w)(v),\dots,\iota(w)(v))\cdot N(w)^{-1} =\tilde{N}(v,\dots,v)\cdot N(w)$ for all $w,v$. In particular for $v=\Id_j$ we have 
    $\tilde{N}(w^2,\dots,w^2) =  N(w)^2 = N(w^2)$. Since over an algebraically closed field every invertible $w$ has a square root, we can conclude that $\tilde{N}$ gives $N$ in this case, and taking the $\K$-points gives the result generally.
\end{proof}

\begin{definition}
    We call $G$ \keyword{simply connected with respect to $\{\jordan{J}_i\}$} if it has a Jordan pinning where every map is an injection. 
\end{definition}

\begin{definition}\label{def:GaussDecompTheta}
    We say an element $g \in G$ is \keyword{Gauss-decomposable with respect to $\Theta$} if $g \in U_\Theta^\opp L_\Theta U_\Theta$. In other words $g$ factors uniquely as 
    \begin{equation}
    g=[g]_{-}[g]_0[g]_{+} \,,
    \end{equation}
    where $[g]_{-} \in U_\Theta^\opp$, $[g]_0\in L_\Theta$ and $[g]_{+}\in U_\Theta$.
\end{definition}
The set of Gauss decomposable elements $G_0=U_\Theta^\opp L_\Theta U_\Theta$ is an open dense subset of $G$. \medskip

\begin{example}\label{ex:Sl6Pinning}
    The group $\SL_6(\R)$ is Jordan split over $A_2$ with $\jordan{J}_i = \jordan{M}_2(\R)$. The structure group $\Gamma(\jordan{J}_i)$ is generated by pairs of matrices $(A,A^{-1})$ sitting inside $\GL_2(\R)\times \GL_2(\R)$. One choice of Jordan pinning is given by 
    \begin{equation*}
        \begin{aligned}
            x_1(A) &= \begin{bsmallmatrix}1 & A & \\& 1 & \\ & & 1\end{bsmallmatrix} \\
            x_2(A) &= \begin{bsmallmatrix}1 &  & \\& 1 & A\\ & & 1\end{bsmallmatrix}
        \end{aligned}\qquad
        \begin{aligned}
            \check{\beta}_1(A,B) &=  \begin{bsmallmatrix}A &  & \\& B & \\ & & 1\end{bsmallmatrix}\\
            \check{\beta}_2(B,C) &=  \begin{bsmallmatrix}1 &  & \\& B & \\ & & C\end{bsmallmatrix}
        \end{aligned} \qquad  
        \begin{aligned}
            y_1(A) &= \begin{bsmallmatrix}1 &  & \\ A & 1 & \\ &  & 1\end{bsmallmatrix}\\
             y_2(A) &= \begin{bsmallmatrix}1 &  & \\& 1 & \\ & A & 1\end{bsmallmatrix}\,.
        \end{aligned}
    \end{equation*}
\end{example}

\begin{example}\label{ex:SP4pinning}
    Before we give a nontrivial pinning of $\Spin(3,4)$, we recall the $C_2$-pinning of $\SP_4(\R)$. This pinning is split so each Jordan algebra is $\R$ with structure group $\R^\times$.
\begin{align*}
        \begin{array}{ccc}
           x_1(t) = \begin{bsmallmatrix}
                1 & 0 & 0 & 0\\
                0  & 1 & 0 & t\\
                0  & 0  & 1 & 0 \\ 
                0  & 0 & 0 &1
            \end{bsmallmatrix}
            & y_1(t) = \begin{bsmallmatrix}
                1 & 0 & 0 & 0\\
                0  & 1 & 0 & 0\\
                0  & 0  & 1 & 0 \\ 
                0  & t & 0 &1
            \end{bsmallmatrix}
            & \check{\beta}_1(a) = \begin{bsmallmatrix}
                1 & 0 &  0 & 0\\
                0  & a & 0 & 0\\
                0  & 0  & 1 & 0 \\ 
                0  & 0 & 0 &a^{-1}
            \end{bsmallmatrix}\\
           x_2(t) = \begin{bsmallmatrix}
                1 & t & 0 & 0\\
                0  & 1 & 0 & 0\\
                0  & 0  & 1 & 0 \\ 
                0  & 0 & -t &1
            \end{bsmallmatrix}
            & y_2(t) = \begin{bsmallmatrix}
                1 & 0 & 0 & 0\\
                t  & 1 & 0 & 0\\
                0  & 0  & 1 & -t \\ 
                0  & 0 & 0 &1
            \end{bsmallmatrix}
             & \check{\beta}_2(a) = \begin{bsmallmatrix}
                a & 0 & 0 & 0\\
                0 & a^{-1} & 0 & 0\\
                0  & 0  & a^{-1} & 0 \\ 
                0  & 0 & 0 &a
            \end{bsmallmatrix}\,.
        \end{array}
    \end{align*}
\end{example}

\begin{example}\label{ex:Spin34pinning}
    A more complicated example is pinning $G=\Spin(3,4)$ as a Jordan split group of type $B_2$. The short root space is labeled $\beta_2$ and corresponds to a copy of $\Spin(3,2)$ inside $G$. We use an exotic isomorphism to identify $\Spin(3,2)\simeq\SP_4(\R)\simeq\SLJ$ for the Jordan algebra $\Sym_2(\R)$. Using the fact that $B_2\simeq C_2$, we can upgrade the pinning for $\SP_4$ of type $C_2$ from \Cref{ex:SP4pinning} to provide a Jordan pinning for $G$ of type $B_2$. Let $M\in \Sym_2(\R)$ and $t\in \R$. Recall that the adjugate $\tau$ is a linear map on $\Sym_2(\R)$ since $\jordan{J}$ is of degree 2. We use the  8-dimensional faithful representation of $\Spin(3,4)$ to write these pinning maps as block $4\times4$ matrices over $M_2(\R)$. Scalar entries are interpreted as scalar multiples of $\Id_2$. We have:
    \begin{align*}
        \begin{array}{ccc}
           x_1(t) = \begin{bsmallmatrix}
                1 & 0 & 0 & 0\\
                0  & 1 & 0 & t\\
                0  & 0  & 1 & 0 \\ 
                0  & 0 & 0 &1
            \end{bsmallmatrix}
            & y_1(t) = \begin{bsmallmatrix}
                1 & 0 & 0 & 0\\
                0  & 1 & 0 & 0\\
                0  & 0  & 1 & 0 \\ 
                0  & t & 0 &1
            \end{bsmallmatrix}
            & \check{\beta}_1(t) = \begin{bsmallmatrix}
                1 & 0 &  0 & 0\\
                0  & t & 0 & 0\\
                0  & 0  & 1 & 0 \\ 
                0  & 0 & 0 &t^{-1}
            \end{bsmallmatrix}\\
           x_2(M) = \begin{bsmallmatrix}
                1 & M & 0 & 0\\
                0  & 1 & 0 & 0\\
                0  & 0  & 1 & 0 \\ 
                0  & 0 & -\tau(M) &1
            \end{bsmallmatrix}
            & y_2(M) = \begin{bsmallmatrix}
                1 & 0 & 0 & 0\\
                M  & 1 & 0 & 0\\
                0  & 0  & 1 & -\tau(M) \\ 
                0  & 0 & 0 &1
            \end{bsmallmatrix}
             & \check{\beta}_2(A) = \begin{bsmallmatrix}
                A & 0 & 0 & 0\\
                0 & \sigma(A)^{-1} & 0 & 0\\
                0  & 0  & \tau(A)^{-1} & 0 \\ 
                0  & 0 & 0 &\sau(A)
            \end{bsmallmatrix}\,.
        \end{array}
    \end{align*}

    We can verify some of the statements about this Jordan pinning, for example \Cref{def:norm_on_rootspace}:
    $$ \check{\beta}_2(A)x_1(t)\check{\beta}_2(A)^{-1} = \begin{bsmallmatrix}
            1 & 0 & 0 & 0\\
            0  & 1 & 0 & t\sigma(A)^{-1}\sau(A)^{-1}\\
            0  & 0  & 1 & 0 \\ 
            0  & 0 & 0 &1
        \end{bsmallmatrix} =
        \begin{bsmallmatrix}
            1 & 0 & 0 & 0\\
            0  & 1 & 0 & tN(A)^{-1}\\
            0  & 0  & 1 & 0 \\ 
            0  & 0 & 0 &1
        \end{bsmallmatrix} \,.$$
    
\end{example}

\subsection{\texorpdfstring{Jordan weights}{Jordan Weights}}
Now we assume that $G$ is such that $G$ is simply connected with respect to a root system $R$. We wish to define a collection of functions dual to the functions $\check{\beta}_i$ to generalize the fundamental weights of split algebraic group. In  noncommutative rank less than or equal to one, what we also call lower noncommutative rank,  there is a natural system generalizing the fundamental weights. However in noncommutative rank bigger or equal to two, what we also call higher noncommutative rank, we must weaken the duality between coroot and fundamental weights. With this weakened duality there will be many possible systems of weights and we will prove our constructions are independent of the choice.

\begin{definition}\label{def:ThetaWeightSystem}
    A system of \keyword{Jordan weights} is a set of functions $\Lambda_i:L_\Theta\to \Gamma(\jordan{J}_i) $ such that:
    \begin{enumerate}
        \item $\Lambda_i(\check{\beta}_j(\Lambda_j(l))) = \begin{cases}
            1 & i\neq j\\
            \Lambda_j(l) & i = j. 
        \end{cases}$\hfill (weak duality)
        \item For all $i \neq j$, $\check{\beta}_i(\Lambda_i(l))\check{\beta}_j(\Lambda_j(l))=\check{\beta}_j(\Lambda_j(l))\check{\beta}_i(\Lambda_i(l))$.\hfill (commutativity)
        \item\label{ThetaWeightReconstruct} For all $l \in L_\Theta$, $\prod_i\check{\beta}_i(\Lambda_i(l))=l$.\hfill (reconstructability)
    \end{enumerate}
    We say a system is \keyword{faithful on $\beta_i$} if $\Lambda_i(\check{\beta}_i(a)) = a$. 
\end{definition}

We note that the Jordan algebras, $\{\jordan{J}_i\}$ are only well-defined up to isotopy. In every case it is possible to choose the isotopy class so that for any $i,j$ one of $\jordan{J}_i$ or $\jordan{J}_j$ is a Jordan subalgebra of the other. This implies that we can choose a pining which identifies a single identity element for all the Jordan algebras and we can consider the target of $\Lambda_i$ to be the same for all $i$.
\begin{example}
    When pinning $\SL_6(\R)$ with $\Theta = \{\alpha_2,\alpha_4\}$, both $\jordan{J}_1$ and $\jordan{J}_2$ are isotopic to $\jordan{M}_2(\R)$. We will take the identity matrix as the identity element in both cases. 

    Similarly when pinning $\SP_{12}(\R)$ with $\Theta = \{\alpha_2,\alpha_4,\alpha_6\}$ we have $\jordan{J}_1$ and $\jordan{J}_2$ isotopic to $\jordan{M}_2(\R)$ and $\jordan{J}_3$ isotopic to $\jordan{H}_2(\R)$. Choosing the identity matrix as the identity in all three cases realizes the required subalgebra conditions.  
\end{example}

\begin{remark}
    The Jordan weights naturally extend to Gauss decomposable elements of $G$, by applying the weight to the Levi factor $[g]_0$.
\end{remark}

\begin{proposition}
    In noncommutative rank 0 or 1, there is a unique system of Jordan weights that is faithful on all simple roots. 
\end{proposition}
\begin{proof}
    In rank 0, the usual fundamental weights provide such a system. In rank 1, property (\ref{ThetaWeightReconstruct}) determines $\Lambda_p$ by forcing $\Lambda_p\big(\check{\beta}_p(A)\big) = A$ and all other $\Lambda_i$ to be fundamental weights.
\end{proof}

\begin{example}
    A system of Jordan weights for the pinning of $\Spin(3,4)$ in \Cref{ex:Spin34pinning} is 
    \begin{align*}
        \begin{aligned}
        &\Lambda_1 : L_\Theta \rightarrow \R \\
        &\Lambda_1 \left(\begin{bsmallmatrix}
            A & & & \\ & B & & \\ & & \tau(A)^{-1} & \\ & & & \tau(B)^{-1}
        \end{bsmallmatrix}\right) = \sqrt{\det(A)\det(B)}
        \end{aligned}
        \hspace{2pc}
        \begin{aligned}
            &\Lambda_2 : L_\Theta \rightarrow \Sym_2(\R) \\
            &\Lambda_2 \left(\begin{bsmallmatrix}
            A & & & \\ & B & & \\ & & \tau(A)^{-1} & \\ & & & \tau(B)^{-1}
            \end{bsmallmatrix}\right) = A\,.
        \end{aligned}
    \end{align*}
\end{example}

We now focus on higher noncommutative rank. As a motivating example, consider the group $\SL_6(\C)$ the $A_2$-pinning analogous to \Cref{ex:Sl6Pinning}. The Levi subgroup is identified with the triple of matrices $(A,B,C)$ with $\det(ABC) = 1$. So any element of the Levi can be decomposed as $\check{\beta}_1(A,A^{-1})\check{\beta}_2(AB,C)$. The natural guess of Jordan weights is  $\Lambda_1(A,B,C) = (A,A^{-1})$ and $\Lambda_2(A,B,C) = (AB,C)$ which has weak duality and reconstructability but not commutativity.  Instead there are two choices of weights which are faithful either on $\check{\beta}_1$ or $\check{\beta}_2$ but not both:
\begin{equation*}
    \Lambda_1(A,B,C) = (A,\sqrt{\det(AB)}^{-1}B) \qquad \Lambda_2(A,B,C) = (\sqrt{\det(AB)}\Id_2,C) \text{ or}
\end{equation*}
\begin{equation*}
    \Lambda_1(A,B,C) = (A,\sqrt{\det(A)}^{-1}\Id_2) \qquad \Lambda_2(A,B,C) = (\sqrt{\det(A)}B,C)\,.
\end{equation*}
Over non algebraically closed fields, the square roots we used might not be well defined. Instead we use a \keyword{weak system of Jordan weights} which lands in a central extension of the structure group by $\K$. The coroot maps can be similarly extended to map the extended structure group into a central extension of $G$. Under this construction the weak equivalent of the first set of weights above is 
\begin{equation*}
    \Lambda_1(A,B,C) = (A,B) \qquad \Lambda_2(A,B,C) = (\Id_2,C) \,.
\end{equation*}
The difference between weak and standard systems of Jordan weights is only important for some groups of higher noncommutative rank. As such, we will drop the word weak for the rest of the paper.

\begin{proposition}
    In higher noncommutative rank there exists many (weak) systems of Jordan weights. Moreover, given a root $\beta_i$ we can choose a system that is faithful on $\beta_i$.
\end{proposition}
\begin{proof}
    There are two families of root systems that admit higher noncommutative rank gradings, $A_p, C_p$ and three exceptional cases $A_2$, $C_3$, and $F_4$. 

    We begin with type $A_p$. Let $\mathcal{R}$ be a finite dimensional algebra over $\K$ and consider the group $G:= \GL_{p+1}(\mathcal{R})$. Let $P_\Theta$ be the parabolic subgroup consisting of upper triangular matrices with Levi subgroup $\GL_1(\mathcal{R})^{p+1}$. This parabolic subgroup induces an $A_p$-grading on $\SL_{p+1}(\mathcal{R})$ where each Jordan algebra is $\jordan{J}(\mathcal{R})$.

    There is a natural Jordan pinning of $G$ given by the $p$ standard embeddings of $\SL_2(\mathcal{R})\rightarrow G$ as $2\times2$ diagonal blocks. To define a weak system, we extend the coroot maps $\check{\beta}_i$ to have domain $\GL_1(\mathcal{R})\times \GL_1(\mathcal{R})$, so that these extended pinning maps are maps from $\GL_2(\mathcal{R})\to G$. This identifies $L_\Theta$ with $\GL_1(\mathcal{R})^{p+1}$ with the image of $\check{\beta}_i$ contained in the factors $i$ and $i+1$. 
    
    We can now construct a (weak) system of Jordan weights faithful on $\beta_i$ by defining the image of each $\SL_1(\mathcal{R})$ factor to either be the trivial map to $\Id$ or the identity map. Each $\Lambda_j$ will map factors other than $j$ or $j+1$ using the trivial map. To be faithful, we define $\Lambda_i$ to be the identity map on both factors $i$ and $i+1$. We then define for $j<i$, $\Lambda_j$ the trivial map on factor $j+1$ and the identity map on factor $j$. Similarly for $j>i$ we take $\Lambda_j$ trivial on factor $j$ and the identity on factor $j+1$. This construction maps each factor with identity map exactly once, so it satisfies commutativity and reconstructability.  \medskip

Next we focus on type $C_p$. Now suppose that $\mathcal{R}$ has an anti-involution $\dagger:\mathcal{R}\to \mathcal{R}$. Let $\Id_p'$ be the anti-diagonal matrix with $1$ along the antidiagonal and let $\Omega = \begin{bmatrix}
    0 & \Id_p' \\ -\Id_p' & 0
\end{bmatrix}$. The group $H:= \SP_{2p}(\mathcal{R},\dagger) = \{g\in \GL_{2p}(\mathcal{R}) | g \Omega (\dagger({g}))^T = \Omega\}$ has a $C_p$-grading, with $\Theta$ associated to the parabolic subgroup of upper triangular matrices in $H$. The Levi subgroup of this parabolic is $\GL_1(\mathcal{R})^p$. $H$ has a natural pinning obtained from the pinning of $G=\GL_{2p}(\mathcal{R})$ by setting $x_i^H(v)= x^G_i(v)x^G_{2p-i}$ for $i<p$ and $x_i^H(v)= x^G_i(v)$ for $i=p$. A faithful choice of Jordan weights is given by a symmetric choice of faithful weights on $\GL_{2p}(\mathcal{R})$.\medskip

Finally, consider the exceptional groups of higher noncommutative rank over $A_2,C_3,$ or $F_4$. In every case the absolute root system is type $F_4$, $E_6$, $E_7$ or $E_8$. Over an algebraically close field, their Levi subgroups are central extensions of $\Spin(8)$ or $(\SL_2)^{3}$, while the structure groups are central extensions of $\Spin(8)$ and $(\SL_2)^{2}$ respectively. The weight functions are then constructed analogously by choosing isomorphisms on the factors corresponding to $\beta_i$ for $\Lambda_i$ and the identity map on those factors for overlapping weights. 

In all cases we see that there is a choice of Jordan weights faithful on the given noncommutative root $\beta_i$. 
\end{proof}

\subsection{Lifts of the grading Weyl group}\label{sec:LiftsOfThetaWeylGroup}
The Jordan pinning over the root system $R$ will allow us to lift the Weyl group $W(R)$ of $R$ to $G$. This lift will act on the decomposition of $G$ into $R$-root spaces as the ordinary Weyl group acts on the split algebraic group associated to $R$. Let $\Delta_R$ be the Dynkin diagram for the root system $R$.
\begin{definition}\label{def:ThetaWeylLifts}
 We define two lifts of each generator $\varsigma_i$ of $W(R)$ 
    \begin{align*}
        \overline{\varsigma_i} &:= y_i(\Id_i)x_i(-\Id_i)y_i(\Id_i) \hspace{2pc} \text{and} \hspace{2pc }\doverline{\varsigma_i} := x_i(\Id_i)y_i(-\Id_i)x_i(\Id_i) \,. 
    \end{align*}
\end{definition}
As in the classical split case, $\overline{\varsigma_i}^{-1} = \doverline{\varsigma_i}$ and $\overline{\varsigma_i}^2\in Z_\Theta$ is in the center of the Levi. This implies that the lifts can be equivalently defined as 
\begin{align*}
    \overline{\varsigma_i} &:= x_i(-\Id_i)y_i(\Id_i)x_i(-\Id_i) \hspace{2pc} \text{and} \hspace{2pc }\doverline{\varsigma_i} := y_i(-\Id_i)x_i(\Id_i)y_i(-\Id_i) \,. 
\end{align*}
Since the lifts fulfill the braid relations for the Dynkin diagram of type $R$, we obtain unique lifts $\overline{w}$ and $\doverline{w}$ of any element $W(R)$ by lifting a reduced expression. These lifts define an ``action'' of $W(R)$ on $G$ (action of the positive braid group) by $\varsigma_i(g) = \doverline{\varsigma_i}g \overline{\varsigma_i}$.
\medskip

Recall from \Cref{def:structureInvolution} that there is an anti-involution on the structure group $\Gamma(\jordan{J}_i)$ we will denote by $\sigma_i$. Using this involution we compute:
\begin{lemma}\label{thm:ActionOfSigmaPonSpecialRoot}
    For any $i$, $$\varsigma_i \left(\check{\beta}_i(A)\right) = \check{\beta}_i\big(\sigma_i(A)^{-1}\big).$$
\end{lemma}
\begin{proof}
    By construction this computation lives entirely in a single $\SLJi$. This group is the exponential of a Lie algebra constructed via the Tits-Kantor-Koecher construction (\Cref{def:TKKconstruction}). Therefore we can compute $\varsigma_i(\check{\beta}_i(A))$ by its adjoint action on $x_i(\Id_i)$.
    \begin{align*}
        \varsigma_i \left(\check{\beta}_i(A)\right). x_i(\Id_i) &= \doverline{\varsigma_i}\check{\beta}_i(A)\overline{\varsigma}_i x_i(\Id_i) \doverline{\varsigma_i}\check{\beta}_i(A)^{-1} \overline{\varsigma_i}\\
        &=\doverline{\varsigma_i}\check{\beta}_i(A)~x_i(-\Id_i)y_i(\Id_i)x_i(-\Id_i)~x_i(\Id_i)~y_i(-\Id_i)x_i(\Id_i)y_i(-\Id_i)~\check{\beta}_i(A)^{-1} \overline{\varsigma_i}\\
        &= \doverline{\varsigma_i}\check{\beta}_i(A)~y_i(-\Id_i)~\check{\beta}_i(A)^{-1} \overline{\varsigma_i}\\
        &= \doverline{\varsigma_i}y_i\big(\sigma_i(A)^{-1}(-\Id_i)\big) \overline{\varsigma_i}\\
        &= x_i\big(\sigma_i(A)^{-1}(\Id_i)\big)\,.
    \end{align*}
    Thus $\varsigma_i\left(\check{\beta}_i(A)\right)   =\check{\beta}_i\big(\sigma_i(A)^{-1}\big)$ as claimed.
\end{proof}

\begin{lemma}\label{thm:ActionOfSigmaPonOtherRoots}
Assume that $\jordan{J}_i$ is a commutative Jordan algebra. Then 
\begin{align*}
    \varsigma_i\left(\check{\beta}_j(A)\right) &= \begin{cases}
        \check{\beta}_i(N(A))\check{\beta}_j(A) & \text{$i$ connected to $j$ in $\Delta_R$}\\
        \check{\beta}_j(A) & \text{otherwise}
    \end{cases}
\\
    \varsigma_j\left(\check{\beta}_i(a)\right) &= \begin{cases}
        \check{\beta}_i(a)\check{\beta}_j(a\Id_j) & \text{$i$ connected to $j$ in $\Delta_R$}\\
        \check{\beta}_j(A) & \text{otherwise}
    \end{cases}
\end{align*}
\end{lemma}
\begin{proof}
For $j$ connected to $i$, the norm on $\jordan{J}_j$ was defined in \Cref{def:norm_on_rootspace} so that
\begin{align*}
    \check{\beta}_j(A).x_i(\Id_i) = x_i\big(N(A)^{-1}\big) \hspace{2pc} \check{\beta}_j(A).y_i(\Id_i) = y_i\big(N(A)\big)\,.
\end{align*}
It is then a simple computation to see
\begin{align*}
    \varsigma_j\left(\check{\beta}_j(A)\right) &= x_i(\Id_i)y_i(-\Id_i)x_i(\Id_i)\check{\beta}_j(A) y_i(\Id_i)x_i(-\Id_i)y_i(\Id_i) \\
    &= x_i(\Id_i)y_i(-\Id_i)x_i(\Id_i)y_i\big(N(A)\big)x_i\big(-N(A)^{-1}\big)y_i\big(N(A)\big)\check{\beta}_j(A) \\
    &= \check{\beta}_i(N(A))\check{\beta}_j(A) \,. 
\end{align*}
Similarly, since the roots are adjacent in the original Dynkin diagram $\Delta$ we have 
\begin{align*}
    \check{\beta}_i(a).x_j(\Id_j) = x_j(a^{-1}\Id_j) \hspace{2pc} \check{\beta}_i(a).y_j(\Id_j) = y_j(a\Id_j) \,.
\end{align*}
The computation is then analogous as it restricts to a calculation in the split subgroup. Finally, when $\beta_i$ and $\beta_j$ are not adjacent in $\Delta_R$, the corresponding roots are not adjacent in the original Dynkin diagram. Therefore the actions of $\varsigma_i$ and $\varsigma_j$ are trivial.
\end{proof}

With the previous two lemmas (\Cref{thm:ActionOfSigmaPonSpecialRoot,thm:ActionOfSigmaPonOtherRoots}) we have computed the full action of the Weyl group in lower noncommutative rank. We observe that when $\Delta_R$ is not simply laced, the degree of the norm map induces the asymmetry given by the high weight edge of $\Delta_R$.

\subsection{\texorpdfstring{$R$}{R}-Bruhat decomposition}\label{sec:ThetaBruhatDecomp}
The Weyl group of the grading root system can be used to define a family of disjoint subsets of $G$. Although these ``cells'' do not cover all of $G$, their combinatorics mirror the combinatorics of a true Bruhat decomposition. First, recall from \Cref{eq:thetaBruhatDecomp} that the Bruhat decomposition with respect to $P_\Theta^\opp$ is given by
\begin{equation*}
    G =\coprod_{w \in W(\Delta\setminus\Theta)\setminus W(\Delta)/W(\Delta\setminus\Theta)} P_\Theta^\opp \dot{w} P_\Theta^\opp\,.
\end{equation*}

We now restrict to the cells which are indexed by elements of $W(R)$. These will be double cosets of $P_\Theta^\opp$, which is the opposite parabolic subgroup associated to the $R$-grading of $G$. Recall the map $\phi^\reduced:W(R)\to W(\Delta)$ defined in \Cref{sec:ThetaWeylGroup_real}.
\begin{definition}
    The \keyword{$R$-Bruhat cell} associated to the element $w \in W(R)$ is $P_\Theta^\opp \phi^\reduced(w) P_\Theta^\opp$.
\end{definition}

We can represent this cell in several other ways, using either the map $\phi: W(R)\to W(\Delta)$ from \Cref{sec:ThetaWeylGroup_real}, or the lifts $\overline{w}$ and $\doverline{w}$ from \Cref{def:ThetaWeylLifts}.
\begin{lemma}
    For any $w\in W(R)$ we have
    \begin{equation*}
        P_\Theta^\opp\phi^\reduced(w)P_\Theta^\opp = P_\Theta^\opp\phi(w)P_\Theta^\opp = P_\Theta^\opp\overline{w}P_\Theta^\opp = P_\Theta^\opp\doverline{w}P_\Theta^\opp\,.
    \end{equation*}
\end{lemma}

We will often simply write $P_\Theta^\opp w P_\Theta^\opp$ for $w\in W(R)$ for this cell, and switch between the different expressions as needed.

\begin{proof}
    By construction $\phi^\reduced(w)$ and $\phi(w)$ lie in the same coset in $W(\Delta)/W(\Delta\setminus\Theta)$. For a single reflection we proved that $P_\Theta^\opp\phi(\varsigma_i)P_\Theta^\opp = P_\Theta^\opp\overline{\varsigma}_iP_\Theta^\opp$ in the proof of \Cref{thm:SP2JBruhatDecomposition}. The lemma then follows from \Cref{thm:BruhatCellsMultiply}.
\end{proof}

\begin{proposition}\label{thm:ThetaLargestCellDense}
    The largest cell $P_\Theta^\opp \phi^\reduced(w_0) P_\Theta^\opp$ is dense in $G$.
\end{proposition}
\begin{proof}
    This cell can equivalently be represented by $w_\Delta$, the longest element in $W(\Delta)$. Thus, left multiplication by $w_\Delta$ identifies this cell with $P_\Theta P_\Theta^\opp$, the set of Gauss decomposable elements. This set is dense $G$.
\end{proof}

\begin{proposition}\label{thm:BruhatSigmaIntersection}
    The Bruhat cells associated to a single generator $\varsigma_i$ intersect $U_\Theta$ exactly in the invertible elements of the pinned root space $x_i$. In symbols:
     \begin{equation}
         U_\Theta \cap P_\Theta^{\opp}\varsigma_i P_\Theta^{\opp} = x_i(V_i^{\times})\,.
     \end{equation}
\end{proposition}
\begin{proof}
    This was proved in \cite{guichard2022generalizing} for groups with a $\Theta$-positive structure. We give an alternate proof that only uses the properties of the Jordan pinning and not the original Weyl group. We instead prove the equivalent statement $U_\Theta\overline{\varsigma}_i \cap P_\Theta^\opp \varsigma_i P_\Theta^\opp \overline{\varsigma}_i = x_i(V_i^\times)\overline{\varsigma}_i$. First we note that $P_\Theta^\opp \varsigma_i P_\Theta^\opp \overline{\varsigma}_i = P_\Theta^\opp x_i(V_i)$ since acting by $\varsigma_i$ sends all negative roots to negative roots except for $-\beta_i$. Let $u\overline{\varsigma}_i \in U_\Theta\overline{\varsigma}_i$ and consider the projection of $u$ onto any positive root space other than $\beta_i$. When multiplied by $\overline{\varsigma}_i$, this projection only intersects $P_\Theta^\opp x_i(V)$ at the identity element. So the only nontrivial intersection occurs when $u = x_i(v)$ for some $v \in V$. Then by \Cref{thm:decompose_A1} elements of the form $x_i(v)\overline{\varsigma}_i$ only decompose as $y_i(v^{-1})\check{\beta}_i(\iota(v))x_i(v)$ when $v$ is invertible. Since all such elements are clearly in $P_\Theta^\opp x_i(V)$, we have the intersection $U_\Theta\overline{\varsigma}_i \cap P_\Theta^\opp \varsigma_i P_\Theta^\opp \overline{\varsigma}_i$ is exactly $x_i(V^\times)\overline{\varsigma}_i$ as claimed.
\end{proof}

The final property of Bruhat decompositions that extends to the $R$-graded decomposition is the result of multiplying two cells whose index word $v$ and $w$ add lengths when multiplying. Without this assumption the result is much more complicated to describe. 
\begin{proposition}\label{thm:BruhatCellsMultiply}
    Let $w_1,w_2 \in W(R)$ such that $\WeylLength(w_1w_2) = \WeylLength(w_1) + \WeylLength(w_2)$. Then
    \begin{equation*}
        P_\Theta^\opp w_1 P_\Theta^\opp w_2 P_\Theta^\opp = P_\Theta^\opp w_1w_2 P_\Theta^\opp\,.
    \end{equation*}
\end{proposition}
\begin{proof}
    This statement follows by induction on the length of $w_2$. If $l(w_2)> 1$ then $w_2 = \varsigma_i w_2'$ for some $w_2'$ of shorter length. Then by induction
    \begin{align*}
        P_\Theta^\opp w_1 P_\Theta^\opp w_2 P_\Theta^\opp &= P_\Theta^\opp w_1 \varsigma_i \varsigma_i P_\Theta^\opp \varsigma_i w_2' P_\Theta^\opp = P_\Theta^\opp w_1 \varsigma_i \varsigma_i P_\Theta^\opp \varsigma_i P_\Theta^\opp w_2' P_\Theta^\opp \,.\\
        \intertext{Moreover $\varsigma_i P_\Theta^\opp \varsigma_i P_\Theta^\opp = x_i(V_i) P_\Theta^\opp$ as in the previous proof. So the expression becomes}
        & P_\Theta^\opp w_1 \varsigma_i x_i(V_i) P_\Theta^\opp w_2' P_\Theta^\opp\,.
    \end{align*}
   When $w_1 \varsigma_i$ acts on $x_i (V_i)$ it transforms into some element in $P_\Theta^\opp$ because $\WeylLength(w_1 \varsigma_i) = \WeylLength(w_1) + 1$. So $P_\Theta^\opp w_1\varsigma_i x_i(V_i) = P_\Theta^\opp w_1\varsigma_i$ and we obtain $P_\Theta^\opp w_1 \varsigma_i P_\Theta^\opp  w_2' P_\Theta^\opp$. By induction we can move the rest of $w_2'$ across proving the theorem. 
\end{proof}

Finally, let us consider some `intermediate' cells which are not indexed by elements of $W(R)$.
\begin{corollary}\label{thm:bruhat_cells_mulitply_partial}
    Suppose that $w$ is a word in $W(R)$ and $w\varsigma_i$ is a reduced word in $W(R)$. Then
    \begin{equation*}
        P_\Theta^{\opp}\phi^\reduced(w)\varsigma_i^{(j)} P_\Theta^{\opp} = P_\Theta^{\opp}\phi^\reduced(w) P_\Theta^{\opp}\cdot P_\Theta^{\opp}\varsigma_i^{(j)} P_\Theta^{\opp}\,,
    \end{equation*}
    for all $1\leq j\leq c_i$. Moreover, for $j_1\neq j_2$, $P_\Theta^{\opp}\phi^\reduced(w)\varsigma_i^{(j_1)} P_\Theta^{\opp}$ and $P_\Theta^{\opp}\phi^\reduced(w)\varsigma_i^{(j_2)} P_\Theta^{\opp}$ are disjoint.
\end{corollary}

\begin{proof}
    We first generalize a result from \cite{guichard2022generalizing}, that $P_\Theta^\opp \phi^\reduced(w)P_\Theta^\opp = P_\Theta^\opp\phi^\reduced P_\Delta^\opp$. 
    By \cite{warner2012harmonic}, for any $w\in W(R)$ we have
    \begin{align*}
        P_\Theta^\opp\phi^\reduced(w)P_\Theta^\opp &= \bigcup_{u\in W(\Delta\setminus\Theta)} P_\Theta^\opp\phi^\reduced(w)P_\Delta^\opp uP_\Delta^\opp\,.
    \intertext{Since $\phi^\reduced(w)u$ is reduced by \Cref{lem:reduced_words_inthetaweylgroup} we use standard Bruhat theory to obtain}    
        &= \bigcup_{u\in W(\Delta\setminus\Theta)} P_\Theta^\opp\phi^\reduced(w)u P_\Delta^\opp = P_\Theta^\opp\phi^\reduced(w)P_\Delta^\opp\,.
    \end{align*}
    The final equality holds because $\phi^\reduced(w)$ normalizes $W(\Delta\setminus\Theta)$ by \Cref{thm:varsigmaRedProperties}. \\
    By assumption, $w\varsigma_i$ is reduced, so \Cref{lem:reduced_words_inthetaweylgroup} implies that $\phi^\reduced(w)\varsigma_i^{(j)}$ is reduced as well. This allows us to conclude
    \begin{equation*}
        P_\Theta^{\opp}\phi^\reduced(w) P_\Theta^{\opp}\cdot P_\Theta^{\opp}\varsigma_i^{(j)} P_\Theta^{\opp} = P_\Theta^{\opp}\phi^\reduced(w) P_\Delta^{\opp}\varsigma_i^{(j)} P_\Theta^{\opp} = P_\Theta^{\opp}\phi^\reduced(w) \varsigma_i^{(j)} P_\Theta^{\opp}.
    \end{equation*}
    By \Cref{lem:product_double_cosets} we see that the double cosets of $w\varsigma_i^{(j)}$ are all distinct which implies that the corresponding cells are disjoint.
\end{proof}

We can parameterize unipotent elements in each cell using a generalization of the Lusztig map, which parameterizes the positive unipotent semigroup in split real Lie groups.
\begin{definition}\label{def:LuszitigMap}
    Let $D = \varsigma_{i_1}\cdots,\varsigma_{i_N}$ be a reduced expression for the longest word in the Weyl group for $R$. A \keyword{prefix $R$-Lusztig map}, $F_D^k : \jordan{J}_{i_1}\times \cdots \times \jordan{J}_{i_k}\rightarrow U_\Theta$ is given by 
    \begin{equation*}
        F_D^k(v_1,\cdots,v_k) = x_{i_1}(v_1)\cdots x_{i_k}(x_k)\,.
    \end{equation*}
    The \keyword{$R$-Lusztig map} is given by $F_D^N$ for $N$ the length of the longest word.
\end{definition}
\begin{corollary}
    For any choice of invertible elements $v_j\in \jordan{J}_{i_j}$, the element $F_D^k(v_1,\cdots,v_k)$ is in the $R$-Bruhat cell $P_\Theta^\opp w P_\Theta^\opp$ where $w = \varsigma_{i_1}\cdots\varsigma_{i_k}$.
\end{corollary}
\begin{proof}
    This follows by induction on $k$ from \Cref{thm:BruhatCellsMultiply}. Let $w_1 = \varsigma_{i_1}\cdots \varsigma_{i_{k-1}}$. We know that $F_D^k(v_1) = F_D^{k-1}(v_1,\cdots,v_{k-1}) x_{i_k}(v_k)$ and by induction that $F_D^{k-1}(v_1,\cdots,v_{k-1}) \in P_\Theta^\opp w_1 P_\Theta^\opp$. We also know $x_{i_k}\in P_\Theta^\opp \varsigma_{i_k} P_\Theta^\opp$ and that $\ell(w_1\varsigma_{i_k}) = \ell(w_1) + \ell(\varsigma_{i_k})$ since $\varsigma_{i_1}\cdots\varsigma_{i_k}$ is reduced. So \Cref{thm:BruhatCellsMultiply} implies that 
    \begin{align*}
        F_D^k(v_1) = F_D^{k-1}(v_1,\cdots,v_{k-1})\cdot x_{i_k}(v_k) \in P_\Theta^\opp w_1 P_\Theta^\opp \varsigma_{i_k}P_\Theta^\opp = P_\Theta^\opp w P_\Theta^\opp \,.
    \end{align*}
\end{proof}

\subsection{Root system calculus}\label{sec:RootSystemCalculus}
The above results enable us to make calculations in Jordan split groups of type $R$ that are independent of the choice of Jordan algebras. The calculations depend only on the root system, and we call this the \keyword{root system calculus}.

Let $G$ be a Jordan split group of type $R$ and fix a normalized Jordan pinning of $G$. This defines elements $x_i(v),y_i(v), \check{\beta}_i(a)$ for $v\in \jordan{J}_i$ and $a\in \Gamma(\jordan{J}_i)$. In \Cref{sec:LiftsOfThetaWeylGroup}, we also defined elements $\overline{\varsigma}_i$ and $\doverline{\varsigma}_i$ which lift the Weyl group of type $R$ to $G$.  The following lemmas provide all relations among these generators needed for the remainder of the paper. Many of the lemmas in this section are restatements of earlier lemmas gathered here for easy reference.

We begin with relations that are entirely within a single pinned $\SLJi$.
\begin{lemma}[\Cref{thm:decompose_A1}]\label{thm:DecomposeXSigma}\label{thm:DecomposeXY}
    Let $v\in V_i^\times$. Then 
    \begin{equation*}
    x_i(v) \overline{\varsigma}_i = \doverline{\varsigma}_i y_i(-v) = y_i(v^{-1})\check{\beta}_i(\iota(v))x_i(-v^{-1}) \,.
\end{equation*}
Let $v,w,v+w \in V_i^\times$. Then
    \begin{equation*}
    x_i(v) y_i(w^{-1}) = y_i\left((v+w)^{-1}\right)\check{\beta}_i\left(\iota(w+v)\iota(w^{-1})\right)x_i\left((w^{-1}+v^{-1})^{-1}\right) \,.
\end{equation*}
\end{lemma}

Next we consider the relations among pinning maps from different roots.
\begin{lemma}
    Let $i\neq j$ index two simple roots in $R$ and consider their connection in the Dynkin diagram. 
    \begin{itemize}
        \item Regardless of connection, $x_i(v) y_j(w) = y_j(w) x_i(v)$.
        \item If $i$ is disconnected from $j$ then $x_i(v)x_j(w) = x_j(w)x_i(v)$ and $y_i(v)y_j(w) = y_j(w)y_i(v)$.
    \end{itemize}
\end{lemma}
\begin{proof}
    The first and second statement follow from computing the brackets of the corresponding root spaces in $\mathfrak{g}$. The brackets vanish as $\beta_i -\beta_j$ is never a root and $\beta_i + \beta_j$ isn't a root when $i$ is not connected to $j$ in a simple Dynkin diagram. So by the Baker-Campbell-Hausdorff theorem the corresponding group elements commute. 
\end{proof}

Now we recall the ``braid action'' of the Weyl group on $G$ given by $\varsigma_i(g) = \doverline{\varsigma}_ig\overline{\varsigma}_i$. The following lemma shows that the action on the Levi subgroup mirrors the usual action on the Cartan subgroup with the  addition of structure involution $\sigma_i$ on $\Gamma(\jordan{J}_i)$ (\Cref{def:structureInvolution}) and the norm map (\Cref{def:norm_on_rootspace}).
\begin{lemma}[\Cref{thm:ActionOfSigmaPonOtherRoots,thm:ActionOfSigmaPonSpecialRoot}]\label{thm:RootCalculusWeyOnLevi}\leavevmode
    \begin{itemize}
        \item For any $i$, we have $\varsigma_i \left(\check{\beta}_i(A)\right) = \check{\beta}_i\big(\sigma_i(A)^{-1}\big)$.
        \item For $i\neq j$ with $i$ disconnected from $j$ in $R$ we have $\varsigma_i(\check{\beta}_j(A)) = \check{\beta}_j(A)$.
        \item For $i\neq j$, with $i$ connected to $j$ in $R$ and $\jordan{J}_i$ commutative we have
            \begin{align*}
                \varsigma_i\left(\check{\beta}_j(A)\right) &= \check{\beta}_i(N(A))\check{\beta}_j(A) \\
                \varsigma_j\left(\check{\beta}_i(a)\right) &= \check{\beta}_i(a)\check{\beta}_j(a\Id_j) 
            \end{align*}
    \end{itemize}
\end{lemma}

Moreover the Weyl group action sends root spaces to root spaces. There is not a unique expression for elements of nonsimple root spaces as products of pinning maps. However it will suffice to know the following lemma.
\begin{lemma}\label{thm:RootCalculusWeylOnRoots}
    For $i \neq j$, $\varsigma_j(x_i(v)) \in P_\Theta$ and $\varsigma_j(y_i(v)) \in P_\Theta^\opp$.
\end{lemma}
\begin{proof}
    By construction $\varsigma_j$ acts by reflection on the root grading. Since the only reflection which sends a simple positive root $\beta_i$ to a negative root is $\varsigma_i$ the statement holds.
\end{proof}

The next lemma describes the action of the Levi factor on the root spaces and on itself. This requires the map $\beta_i : L_\Theta \rightarrow P\Gamma(\jordan{J}_i)$ defined in \Cref{def:BetaMaps} via $l x_i(v) l^{-1} = x_i\big(\beta_i(l)(v)\big)$.
\begin{lemma}\label{thm:RootCalculusLeviOnRoots}
    For any $i$, we have
    \begin{itemize}
        \item $\beta_i(\check{\beta}_i(A))(v) =A(v)$
        \item $ly_i(v)l^{-1} = y_i(\beta_i(\varsigma_i(l))(v)) = y_i(\sigma_i(\beta_i(l))^{-1}(v))$
        \item $l\check{\beta_i}(\iota(v))\varsigma_i(l)^{-1} = \check{\beta_i}\big(\iota\big(\beta_i(l)(v)\big)\big)$
    \end{itemize}
\end{lemma}
\begin{proof}
    The first statement combines definition of $\beta$ with the action of the coroot $\check{\beta}(A)$ on $x(v)$. 
    For the second statement, we conjugate everything by $\overline{\varsigma}_i$ to change $y_i(v)$ into $x_i(-v)$ and recover the definition of $\beta_i$. 
    For the final statement write $w=\beta_i(l)(v)$. On the one hand 
    \begin{align*}
        lx_i(v)l^{-1}\overline{\varsigma}_i = x_i(w) \overline{\varsigma}_i = y_i(w^{-1})\cdot \check{\beta}_i(\iota(w))\cdot x_i(-w^{-1})\,.
    \end{align*}
    On the other hand, we observe that $w^{-1} = \sigma_i(\beta_i(l))^{-1}(v^{-1})$ and compute
    \begin{align*}
        lx_i(v)l^{-1}\overline{\varsigma}_i &= l x_i(v) \overline{\varsigma}_i ~\varsigma_i(l^{-1}) \\
        &=  l \big(y_i(v^{-1}) \check{\beta}_i(\iota(v))x_i(-v^{-1})\big) \varsigma_i(l)^{-1}\\
        &= l y_i(v^{-1})l^{-1} \cdot l \check{\beta}_i(\iota(v)) \varsigma_i(l)^{-1} \cdot \varsigma_i(l)x_i(-v^{-1}) \varsigma_i(l)^{-1}\\
        &=y_i(w^{-1}) \cdot l \check{\beta}_i(\iota(v)) \varsigma_i(l)^{-1} \cdot x_i(-w^{-1})\,.
    \end{align*}
    Equating the $\check{\beta}_i$ terms proves the final statement. 
\end{proof}

\part{Noncommutative Cluster Varieties}\label{part:NoncomClusterVarieties}
Our next goal is to develop a theory of noncommutative cluster varieties rich enough to provide cluster structures for the natural geometric objects related to Jordan split groups. In noncommutative rank 0, the algebraic groups are quasi-split and the classical cluster structure will suffice as described in \cite{goncharov2019quantum}. In our previous work, we defined a noncommutative cluster framework that suffices for noncommutative rank 1 \cite{greenberg2024noncommutative}. In higher noncommutative rank, we instead require a further generalization of the noncommutative wiring diagrams defined in \cite{goncharov2021spectral}. In this part, we review the classical cluster theory and describe the two noncommutative generalizations needed for the remainder of the paper.

\section{Classic cluster algebras}\label{sec:ClassicClusterAlgebra}
We first recall the construction of a cluster algebra. It is \emph{grown} from an initial \keyword{seed} which consists of a quiver and a set of variables associated to each node of the quiver.

\subsection{Quivers and seed mutation}
\begin{definition}
    A \keyword{quiver} $Q$ is a directed graph on $n$ nodes with no self loops or two cycles. The nodes of $Q$ are partitioned into two sets, the \keyword{mutable/unfrozen} nodes and \keyword{frozen} nodes. 
    
    A \keyword{quiver with weighted nodes} is a quiver where each node is assigned a positive integer weight $d_k$. The \keyword{adjacency matrix} of a weighted quiver is a skew-symmetrizable matrix $\epsilon_{ij}$ with entries $\epsilon_{ij}=e_{ij}d_j/\gcd(d_j,d_i)$, where $e_{ij}$ is the number of arrows from node $i$ to node $j$ counted with signs.
\end{definition}
Note that the diagonal matrix $D$ with $D_{ii} = d_i$ is the skew symmetrizer of $\epsilon$, in other words $D\epsilon$ is skew symmetric. 

\begin{example}
     We will only need quivers with nodes of two different weights. These will either be in $\{1,2\}$ or in $\{1,3\}$. In either case we draw nodes of higher weight larger than nodes of weight 1. For example in \Cref{fig:WeightedQuiverExample} the first two nodes are weight 1 while the rightmost node is weight 2. The skew-symmetrizable adjacency matrix is also given.
\begin{figure}[htb]
    \centering
    \begin{tikzpicture}
        \node[mutable]	(1)	at (0,0)		[label = above: $1$]	{};
        \node[mutable]	(2)	at (\distL,0)		[label = above: $2$]	{};
        \node[mutableBig]	(3)	at (2*\distL,0)		[label = above: $3$]	{};
        \draw[arrow] (1) to (2);
        \draw[arrow] (2) to (3);
    \end{tikzpicture}
    \hspace{2pc}
    $\begin{bmatrix}
        0 & 1 & 0\\
        -1 & 0 & 2\\
        0 & -1 & 0
    \end{bmatrix}$
    \caption{Quiver with weighted nodes.}
    \label{fig:WeightedQuiverExample}
\end{figure}
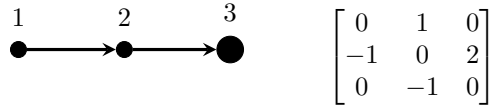

We would draw the same quiver if node 3 had weight 3 instead. In this case the adjacency matrix would be
\begin{equation*}
    \begin{bmatrix}
        0 & 1 & 0 \\ -1 & 0 & 3 \\ 0 & -1 & 0
    \end{bmatrix}\,.
\end{equation*}
It will always be clear from the context which case we are dealing with.
\end{example}

\begin{definition}
    A \keyword{seed}, $(Q,\vec{a})$ consists of a weighted quiver $Q$ and a list of algebraically independent variables $\vec{a}=(a_1,\dots,a_n)$ in $\K(x_1,\dots,x_n)$ associated to the $n$ nodes of $Q$. The variables $a_i$ are called \keyword{cluster variables}.
\end{definition}

Given a seed $(Q,\vec{a})$ and a mutable node $k$, we can produce a new seed $\mu_k(Q,\vec{a}) = (Q',\vec{a'})$  by a process called \keyword{mutation} which consists of mutating the quiver $Q$ and the cluster variables. 
\begin{definition}
    The \keyword{mutation of $Q$ at node $k$} is the quiver $\mu_k(Q)$ formed by
    \begin{enumerate}
        \item For each pair of arrows $i \rightarrow k \rightarrow j $ add $d_k \gcd(d_i,d_j) / (\gcd (d_i,d_k) \gcd(d_k,d_j))$ arrows from $i$ to $j$.
        \item Remove any 2-cycles created by the previous step.
        \item Reverse all arrows incident to $k$.
    \end{enumerate}
    This procedure is equivalent to the following operation on the adjacency matrix:
    \begin{equation}
		b'_{ij}=\begin{cases}
			-b_{ij} & k\in\{i,j\}\\
			b_{ij}+b_{ik}^+b_{kj}^+-b_{ik}^-b_{kj}^- & k\notin\{i,j\}
		\end{cases}
    \end{equation}
    where $x^\pm = \pm\max(0,\pm x)$.
\end{definition}

\begin{example}\label{ex:BasicClusterAlgebra}
    Consider the quiver $Q$ in \Cref{fig:exMutationStart}. It has one mutable node drawn with a circle and 4 frozen nodes drawn as squares. Mutation at the mutable node produces the quiver \Cref{fig:exMutationEnd}.
    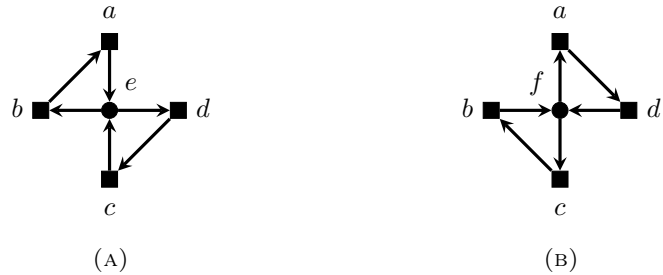
\begin{figure}[htb]
        \centering
        \begin{subfigure}{0.4\textwidth}
            \centering
            \begin{tikzpicture}[scale = 1.3]
                \node[frozen] (12) at (0.5*\distL, 0) [label = above: $a$] {};
                \node[frozen] (23) at (0*\distL, -0.5*\distL) [label = left: $b$] {};
                \node[frozen] (34) at (0.5*\distL, -1*\distL) [label = below: $c$] {};
                \node[frozen] (14) at (1*\distL, -0.5*\distL) [label = right: $d$] {};
                \node[mutable] (13) at (0.5*\distL, -0.5*\distL) [label = above right: $e$] {};
                \draw[arrow] (12) to (13);
                \draw[arrow] (13) to (23);
                \draw[arrow] (23) to (12);
                \draw[arrow] (34) to (13);
                \draw[arrow] (13) to (14);
                \draw[arrow] (14) to (34);
            \end{tikzpicture}
            \caption{}
            \label{fig:exMutationStart}
        \end{subfigure}
        \begin{subfigure}{0.4\textwidth}
            \centering
            \begin{tikzpicture}[scale = 1.3]
                \node[frozen] (12) at (0.5*\distL, 0) [label = above: $a$] {};
                \node[frozen] (23) at (0*\distL, -0.5*\distL) [label = left: $b$] {};
                \node[frozen] (34) at (0.5*\distL, -1*\distL) [label = below: $c$] {};
                \node[frozen] (14) at (1*\distL, -0.5*\distL) [label = right: $d$] {};
                \node[mutable] (24) at (0.5*\distL, -0.5*\distL) [label = above left: $f$] {};
                \draw[arrow] (24) to (12);
                \draw[arrow] (23) to (24);
                \draw[arrow] (24) to (34);
                \draw[arrow] (14) to (24);
                \draw[arrow] (34) to (23);
                \draw[arrow] (12) to (14);
            \end{tikzpicture}
            \caption{}
            \label{fig:exMutationEnd}
        \end{subfigure}
        \caption{Mutation at the mutable node in the middle relates the two quivers.}
        \label{fig:exampleMutation}
    \end{figure}
\end{example}

We are now ready to define seed mutation and the cluster algebra:
\begin{definition}
    The \keyword{mutation of a seed $(Q,\vec{a})$ at $k$} is the new seed $(\mu_k(Q),\vec{a'})$ whose quiver is given by mutation above. All the cluster variables of $\vec{a'}$ are the same as $\vec{a}$ except for $a_k'$ which satisfies the \keyword{exchange relation}
    \[ a_k a_k' = \prod_{b_{ik}>0} a_i^{b_{ik}} + \prod_{b_{ik}<0} a_i^{-b_{ik}}\]
    The \keyword{cluster algebra} generated by a seed $(Q,\vec{a})$ is the subalgebra of the field of fractions $\K(\vec{a})$ generated by all cluster variables obtained by arbitrary sequences of mutations.
\end{definition}

\begin{example}
    Consider again the quivers in \Cref{fig:exampleMutation}. We have labeled both with variables $a,\dots,f$. The exchange relation reads:
    \begin{equation}\label{eq:2in2outExchangeRelation}
        ef = ac + bd\,.
    \end{equation}
\end{example}

Cluster algebras are of particular interest because they endow the coordinate rings of some varieties with a combinatorial structure: We say a variety has a \keyword{cluster structure} if there is a seed which generates a cluster algebras isomorphic to the coordinate ring of the variety. Thus we will often implicitly identify the abstract cluster coordinates with elements of the given ring.
\begin{example}
    Consider the decorated flag variety $\decFlag$ for $\SL_2(\R)$. It consists of nonzero vectors in $\R^2$. The configuration space
    \begin{equation*}
        \Conf_4(\decFlag) = \SL_2(\R)\backslash\decFlag^4
    \end{equation*}
    consists of 4-tuples of decorated flags up to the diagonal action of $\SL_2(\R)$. The subspace $\Confx_4(\decFlag)$ where no two vectors lie on a common line has a cluster structure given by either quiver in \Cref{fig:exampleMutation}. Namely, if the configuration $x$ is represented by $(v_1,\dots,v_4)$, we can compute coordinate functions
    \begin{equation*}
        a_{ij}(x):=\det(v_i | v_j)
    \end{equation*}
    which are clearly invariant under the action of $\SL_2(\R)$ on $\decFlag^4$ and thus descend to $\Confx_4(\decFlag)$. We get $\binom{4}{2}=6$ coordinates since $a_{ij}=-a_{ji}$. One can check that these satisfy the equation
    \begin{equation*}
        a_{13}a_{42}=a_{14}a_{32}+a_{12}a_{43}
    \end{equation*}
    which has exactly the form of the exchange relation \Cref{eq:2in2outExchangeRelation}. Thus, we can assign the coordinate functions to the nodes of the quivers in \Cref{fig:exampleMutation} as $a = a_{14},\,b = a_{12},\,c = a_{32},\,d = a_{43},\,e = a_{13}$ and $f = a_{42}$.
    
    This is the most basic example of the theory we will discuss in \Cref{part:ClusterlikeCoordinates}.
\end{example}

\subsection{Cluster algebras from surfaces}
The quivers in \Cref{fig:exampleMutation} have another interesting property: They are obtained from the two triangulations of the square by drawing a node for every arc in the triangulation (where the boundary arcs correspond to frozen nodes) and inscribing a clockwise three-cycle of arrows into every triangle. Cluster algebras arising in this way were studied in \cite{fock2007dual, fomin2008clusterTri}.

For any quiver obtained from a triangulated surface in the way described above, any single mutation is equivalent to a flip in the triangulation. This is the operation where one diagonal of a square is removed and replaced with the other one.

The exchange relation for such a mutation can be read off in an interesting way: Consider a square as in \Cref{fig:zigZagCommutative}. To replace the diagonal with the other one, we consider the two zigzag paths which connect the endpoints of the new diagonal.

\input{Figures/figureZigZagCommutative}

By writing down the cluster variables we encounter along these two paths we get the expression
\begin{equation}\label{eq:zigZagExchangeCommutative}
    {\color{Red}aec}+{\color{Blue}bed} = fe^2
\end{equation}
where $f$ is the new cluster variable obtained by mutation and the equality simply follows by multiplying the exchange relation by $e$. In this commutative setting neither the orientation of the paths nor the order of the product matter. But we observe the surface can provide a natural order which is exactly the approach used by Berenstein-Retakh in their definition of noncommutative surfaces \cite{berenstein2018noncommutative}.

\subsection{Key properties}
We now recall a few key properties of cluster algebras \cite{fomin2002cluster}.
\begin{proposition}
    If two seeds $(Q,\vec{a})$ and $(Q',\vec{a'})$ are related by a sequence of mutations, they produce isomorphic cluster algebras.
\end{proposition}
\begin{proposition}[Laurent Phenomenon]
    Every cluster variable can be written as a Laurent polynomial in the cluster variables from a single seed $(Q,\vec{a})$. Moreover the coefficients of the Laurent polynomial are positive.
\end{proposition}
\begin{definition}
    Given a (semi)ring $R$, the \keyword{$R$-points} of a cluster algebra are the homomorphisms from the cluster algebra to $R$. When $R = \R$ this is equivalent to choosing a real number for each cluster variable subject to the exchange relations. An $\R$-point is called \keyword{positive} if all the cluster variables are evaluated at strictly positive real numbers. 
\end{definition}
\begin{proposition}[Positivity]
    The positive $\R$-points of a cluster algebra can be constructed by evaluating the cluster variables of any seed at positive real numbers and using the exchange relations to evaluate all other cluster variables.   
\end{proposition}
\begin{corollary}
    The positive points of a cluster algebra form a well-defined semi-algebraic set. Furthermore a positive point can be specified by choosing arbitrary positive numbers for all the cluster variables in a single seed $(Q,\vec{a})$.
\end{corollary}

\section{Polygonal cluster algebras}\label{sec:PolygonalClusterAlgebras}
We now recall some of the basic constructions from \cite{greenberg2024noncommutative} we will need to describe the cluster structure on our noncommutative cluster coordinates. We introduce a small modification to include the $G_2$ case. In this section we fix a parameter $r$ which will be $2$ in type $B_p$ and $3$ in type $G_2$.

\subsection{2-colored and ordered quivers}

All of the quivers we will consider are weighted quivers where the weights are either 1 or $r$. We denote by $Q_0^{(1)}$ and $Q_0^{(2)}$ the set of weight one and weight $r$ nodes respectively. We call the weight 1 nodes \keyword{commutative nodes} and the weight $r$ nodes \keyword{noncommutative nodes}.

\begin{definition}
    A \keyword{2-colored} quiver $Q$ is a weighted quiver where all of the arrows between noncommutative nodes are partitioned into two sets called \keyword{parallel} and \keyword{crossing arrows} respectively. The combination of one parallel and one crossing arrow with the same source and target is called a \keyword{paired arrow}.

    The \keyword{pruned quiver} denoted by $Q^{\pruned}$ is obtained by removing all of the commutative nodes and all paired arrows between weight $r$ nodes.     We call the arrows in $Q^{\pruned}$ \keyword{noncommutative arrows}.
\end{definition}

\begin{convention}
    We draw parallel arrows by black arrows with solid tips ($\solidrightarrowtriangle$) and crossing arrows by blue arrows with a hollow tip ($\textcolor{Blue}{\rightarrowtriangle}$). Paired arrows are green with a double arrow head ($\textcolor{Green}{\twoheadrightarrow}$).
\end{convention}

Crossing arrows are assigned weight $r-1$, so that during mutation they count as $r-1$ arrows. Consequently, paired arrows have weight $r$. This is reflected in the following mutation rule:

\begin{definition}
    Mutation of a 2-colored quiver is defined by the following rules:
    First suppose that $k \in Q_0^{(2)}$ is a weight $r$ node, then mutation at $k$ is determined by the following process:
     \begin{enumerate}
        \item For every sequence $i\rightarrow k\rightarrow j$ with $i,j$ weight $r$ nodes, introduce arrows $i\rightarrow j$, whose colors are determined by the \keyword{coloring rules}
	        \begin{align*}
	    		\textcolor{black}{\solidrightarrowtriangle} + \textcolor{black}{\solidrightarrowtriangle} =~& \textcolor{black}{\solidrightarrowtriangle}\\
	    		\textcolor{Blue}{\rightarrowtriangle} + \textcolor{Blue}{\rightarrowtriangle} =~& \textcolor{black}{\solidrightarrowtriangle} + (r-2)\textcolor{Green}{\twoheadrightarrow} \\
	    		\textcolor{Blue}{\rightarrowtriangle} + \textcolor{black}{\solidrightarrowtriangle} =~& \textcolor{Blue}{\rightarrowtriangle} \,.
			\end{align*}
            If instead at least one of $i,j$ is weight one, then the new arrows are always colored black.
        \item Cancel any resulting 2-cycles of the same color. Note that a paired arrow can cancel with either a parallel or a crossing arrow to result in an arrow of the opposite class:
            \begin{align*}
	    		\textcolor{black}{\solidrightarrowtriangle} + \textcolor{Green}{\twoheadleftarrow} =~& \textcolor{Blue}{\leftarrowtriangle}\\
	    		\textcolor{Blue}{\rightarrowtriangle} + \textcolor{Green}{\twoheadleftarrow} =~& \textcolor{black}{\text{\reflectbox{$\solidrightarrowtriangle$}}}
			\end{align*}
        \item Reverse any arrow connected to $k$.  
    \end{enumerate}
    If instead $k\in Q_0^{(1)}$ is weight one then mutation is determined as follows:
    \begin{enumerate}
        \item For every sequence $i\rightarrow k\rightarrow j$ with $i,j$ weight $r$ nodes, introduce a paired arrow $i \textcolor{Green}{\twoheadrightarrow} j$. If instead at least one of $i,j$ is weight one, then the new arrows are always colored black.
        \item Cancel any resulting 2-cycles of the same color.
        \item Reverse arrows connected to $k$.
    \end{enumerate}
\end{definition}

\begin{remark}
    The key motivation for these definitions comes from the type of coordinates we will define on $\Theta$-positive representations. In these cases, we assign elements of the structure group $\Gamma(\jordan{J})$ of a Jordan algebra to the non commutative nodes of $Q$ and real numbers to the commutative nodes. The number $r$  is the degree of the norm map of $\jordan{J}$. The commutative arrows will really `see' the norm of the elements while crossing arrows will `see' an involution $\tau$ applied to the corresponding elements.
\end{remark}

\begin{example}
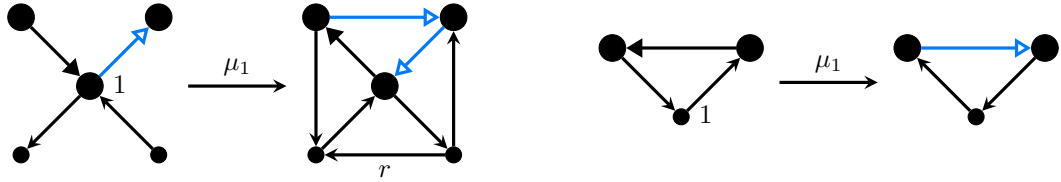
    Some examples of 2-colored quivers and their mutation are shown in \Cref{fig:2coloredMutationExample}.
    
    \input{Figures/figureTwoColoredMutation}
\end{example}

Adding color to the quivers is not sufficient to capture the combinatorics of the coordinates which we will discuss in \Cref{part:ClusterlikeCoordinates}. For this, we will also need to take the order into account since some of the functions will take values in a non-abelian group. This requires a further upgrade of the quivers:

\begin{definition}
    An \keyword{ordered (2-colored) quiver} consists of a 2-colored quiver $Q$ together with a splitting of the incoming and outgoing arrows in the pruned quiver $Q^{\pruned}$ at every node $k\in Q_0^{(2)}$ into two ordered tuples, denoted $\In_f(k)$ and $\Out_f(k)$ for $f\in\{0,1\}$.
\end{definition}

\begin{figure}[ht]
	\begin{center}
		\begin{tikzpicture}
			\node[vertex]	(A) at (-1.5,1)		[label = left: $\In_\bullet(k)$]	{};
			\node[vertex]	(B) at (-1.5,-0.8)	[]{};
			\node[vertex]	(C) at (-1.1,-1.1)	[]{};
			\node[vertex]	(AB) at ($0.5*(B)+0.5*(C)+(0,-0.2)$) [label = left: $\Out_\bullet(k)$]	{};
			\node[vertex]	(D) at (1.5,1)		[label = right: $\Out_\circ(k)$]{};
			\node[vertex]	(F) at (1.5,-1)		[label = right: $\In_\circ(k)$]	{};

			\pic[name=E] at (0,0) {NCnode};
			\node[vertex] (X) at (0,0.8)		[label = above: $k$]		{};
            \node[vertex]	(OutC) at (3.5,0.0)		[label = right: $\Out^c(k)$]	{};

			\draw[barrow] (A.-20) to [out=-20,in=140] (E-circle.140);
			\draw[xarrow] (E-circle.220) to [out=220,in=10] (B.10);
			\draw[barrow] (E-circle.240) to [out=250,in=60] (C.60);
			\draw[xarrow] (E-circle.45) to [out=45,in=-180] (D.-180);
			\draw[barrow] (F.180) to [out=180,in=-45] (E-circle.-45);
            \draw[carrow] (0,0) to[out=15,in=195] (OutC);
            
		\end{tikzpicture}
	\end{center}
 	\caption{Decoration of big nodes in an ordered quiver with partition of the arrows.}
	\label{fig:decNode}
\end{figure}
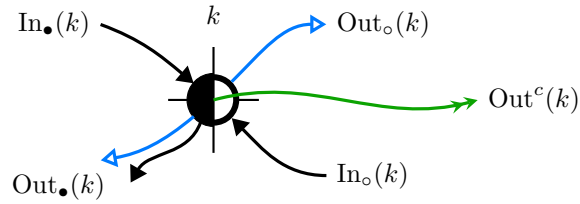 

\begin{convention}
    We draw the noncommutative nodes of an ordered quiver as shown in \Cref{fig:decNode}. The order of the tuples is given by going counterclockwise around the node. We call the two halves open and filled, respectively. To match the pictures we use $\circ$ for $0$ and $\bullet$ for $1$. The commutative arrows (arrows connected to weight 1 nodes and paired arrows) connected to a node $k\in Q_0^{(2)}$ are denoted by the sets $\In^c(k)$ and $\Out^c(k)$.   
\end{convention}

In order to mutate such quivers we restrict to a special subclass.

\subsection{Decorated polygons}
In \cite{greenberg2024noncommutative} we introduced the notion of source-target-compatible quivers using an algebraic compatibility condition, and proved that these are in correspondence with decorated polygonal tiling. To keep the exposition concise, we will simply consider decorated polygons and their corresponding quivers here.

\begin{definition}\label{def:decPoly}
    A \keyword{decorated polygon} consists of the following data:
    \begin{itemize}
        \item a polygon $P_n$,
        \item an orientation of its sides,
        \item a partition of its vertices into two disjoint sets, called \keyword{parallel} and \keyword{crossing angles} such that there is an \emph{odd number of parallel angles},
        \item a choice of vertex for every side such that the straight lines connecting the midpoint of each side to its corresponding vertex intersect pairwise (allowing intersection at the endpoints). We call these lines \keyword{angle indicators}.
    \end{itemize}
\end{definition}

Some examples of decorated polygons are visualized in \Cref{fig:decPolyEx}. We indicate a parallel angle by a black solid arc near the vertex, while a crossing angle is indicated by a blue dashed arc.

\input{Figures/figureDecPolyExamples}

\medskip
Given a decorated polygon we construct an ordered quiver as follows:
\begin{enumerate}
    \item For every side, we get a noncommutative node such that the filled half is on the left as given by the orientation.
    \item The halves of the nodes facing into the polygon are connected with arrows: For a given node $k$ all nodes corresponding to sides to the right of the angle indicator are sources of arrows with target $k$, while all nodes which correspond to sides to the left of the angle indicator are targets of arrows starting at $k$.
    \item The order of the arrows at every node is given by the clockwise order of the sides.
    \item To determine the color of the arrows, start by coloring the arrows between nodes corresponding to adjacent sides: Unless all angle indicators end at the vertex between the two sides, color the arrow the color of the angle. Otherwise take the opposite color. Color all remaining arrows such that every set of three nodes
    \begin{itemize}
        \item has an odd number of parallel arrows if the nodes form a 3-cycle,
        \item has an even number of parallel arrows otherwise.
    \end{itemize}
\end{enumerate}

\begin{remark}
    The coloring rules determine a unique coloring of the quiver \cite{greenberg2024noncommutative}.
\end{remark}

The ordered quivers associated to the polygons in \Cref{fig:decPolyEx} are shown in \Cref{fig:decPolyToQuiverEx}.

\input{Figures/figureDecPolyTranslation}

By gluing decorated polygons along their sides, we can use them to tile a surface.
\begin{definition}
    A \keyword{decorated tiling} of an oriented surface $S$ with marked points consists of the following data
    \begin{itemize}
        \item a set of non-intersecting oriented arcs between marked points on $S$, such that the complementary regions are polygons,
        \item an admissible decoration for every such polygon such that the orientation of the edges agree with the orientation of the arcs on the surface,
    \end{itemize}
\end{definition}

Using the same rules as above for every polygon in the tiling, we obtain an ordered quiver from every decorated tiling.
\begin{definition}
    An ordered quiver for which the pruned quiver can be represented as a decorated tiling is called a \keyword{polygonal quiver}.
\end{definition}

\begin{remark}
    This differs slightly from the  terminology in \cite{greenberg2024noncommutative}. There we gave an alternate classification of these quivers which we called  `source-target compatible'.
\end{remark}

We can now use the polygons to mutate polygonal quivers:
\begin{definition}
    Let $Q$ be an polygonal quiver corresponding to a decorated tiling $T$. The (admissible) \keyword{mutation} at a node $k$ produces a new polygonal quiver $\mu_k(Q)$ as follows:
    \begin{enumerate}
        \item Perform a mutation of the underlying 2-colored quiver at $k$.
        \item Order the arrows at every node such that the ordered quiver corresponds to a decorated tiling $T'$. The tiling $T'$ is determined by the following conditions:
        \begin{itemize}
            \item If $k\in Q_0^{(1)}$, then $T'$ has the same polygons but possibly different angle indicators.
            \item If $k\in Q_0^{(2)}$, then $T'$ is related to $T$ by a flip as shown in \Cref{fig:flipDecPoly}. All the other angle indicators and angle colors are determined by the mutated 2-colored quiver.
        \end{itemize}
    \end{enumerate}
    We denote by $|Q|$ the set of polygonal quivers which are obtained from $Q$ by sequences of admissible mutations.
\end{definition}

\input{Figures/figureFlipDecPoly}

\begin{remark}
    The second step is not always possible: There are polygonal quivers for which a mutation of the underlying 2-colored quiver does not lead to a quiver which can be obtained from the prescribed tiling. We do not consider mutations at such nodes.
\end{remark}

\begin{example}
     Let us illustrate the mutation procedure and possible inadmissibility of mutations with two examples: 
     
     \Cref{fig:mutationEx} illustrates an admissible mutation at a noncommutative node. The order of the arrows is deduced from the decorated tilings which are also shown

    \input{Figures/figureSTmutationExample}

    On the other hand, consider the polygonal quiver in \Cref{fig:nonAdmissibleMutationQuiver} with its pruned part represented by the decorated triangle as shown. Mutation at the commutative node is not admissible because the new 2-colored quiver cannot be represented by changing the angle indicators of the decorated polygon.

    \input{Figures/figureNonAdmissibleMutation}
\end{example}

Finally, we define two operations on polygonal quivers which will be used later:
\begin{definition}\label{def:weavingQuiver}
    Let $Q$ be an polygonal quiver with decorated tiling $T$, $k$ a noncommutative node of $Q$, respectively a side in $T$.  The \keyword{switch at $k$} gives an polygonal quiver $\sigma_k(Q)$ which corresponds to the tiling $\sigma_k(T)$ where the orientation of $k$ is reversed.

    When $r=2$, we moreover define: The \keyword{weaving at $k$} gives $w_k(Q)$ corresponding to $w_k(T)$ where the types of all angles adjacent to $k$ is switched, i.e. straight angles become crossing and vice versa.
\end{definition}
The first operation corresponds to reversing the fill of the node $k$, while the second operation switches the color of all arrows connected to $k$. Since we would like these operations to commute with mutation, the weaving only makes sense if $r=2$ because for $r=3$ the different kinds of arrows behave combinatorially quite differently.

\subsection{Seed group and angles}
Let us now introduce the algebraic framework needed to define mutations. Let $Q$ be an polygonal quiver with corresponding decorated tiling $P$. We assign a collection of variables to $Q$:
\begin{align*}
    & \left\{X_i^{(f)},N_i \Big |\,i\in Q_0^{(2)}, f\in \mathbb{Z}/2\mathbb{Z}\right\} \cup  \left\{N_i \Big |\,i\in Q_0^{(1)}\right\}\,.
\end{align*}
Let $F_Q$ be the free group generated by these. We define two anti-automorphisms of $F_Q$ called $\sigma$ and $\tau$ as follows: For $i\in Q_0^{(2)}, f\in\Z/2\Z$ and $j\in Q_0^{(1)}\cup Q_0^{(2)}$
\begin{align*}
    &\begin{matrix}
        & X_i^{(f)} & \xmapsto{\sigma} & X_i^{(f+1)}\\
        & N_j &\xmapsto{\sigma} &N_j
    \end{matrix} \hspace{2pc}
    &\begin{matrix}
        & X_i^{(f)} & \xmapsto{\tau} & \big(X_i^{(f)}\big)^{-1}N_i \\
        & N_j & \xmapsto{\tau} & N_j
    \end{matrix}\,.
\end{align*}
To lighten the notation, we write $X_i:=X_i^{(0)}$ and $\sau(X_i^{(f)}):=\sigma\big(\tau(X_i^{(f)})\big)$.\medskip

We now consider paths in the decorated tiling:
\begin{definition}
    An \keyword{angle path} $\gamma$ is a (continuous) path along the sides and angles of a decorated tiling, which avoids the vertices of the polygons, as shown in \Cref{fig:angleEx}. The \keyword{parity} $\ell$ on angle paths is the $\Z/2\Z$-valued additive function with
    \begin{equation*}
        \ell(e) = \begin{cases}
            1 & e \text{ is a side or a crossing angle}\\
            0 & e \text{ is a parallel angle.}
        \end{cases}
    \end{equation*}
\end{definition}

\input{Figures/figureAngleExample}

To every angle path we associate an expression in the variable from above:
\begin{definition}\label{def:anglePCA}
    Let $\gamma$ be an angle path in a tiling $P$. The \keyword{label} $L(\gamma)\in F_Q$ of $\gamma$ is computed recursively according to the following rules:
    \begin{itemize}
        \item The label of a path following an angle is 1.
        \item The label along a side $i$ is $\sau(X_i)$ if the path follows the orientation, $\tau(X_i)$ otherwise.
        \item The label of the path $\gamma * \delta$ obtained by concatenation is given by
        \begin{equation*}
            L(\gamma*\delta) = L(\gamma)\sau^{\ell(\gamma)}\big(L(\delta)\big)\,.
        \end{equation*}
    \end{itemize}
    An \keyword{angle} is the label associated to any loop around a polygon.
\end{definition}

\begin{example}\label{ex:surfaceAngle}
    Consider the polygon in \Cref{fig:angleEx} with the indicated angle path $\gamma$. Since this is a loop, the associated label is an angle, namely
    \begin{equation*}
        L(\gamma) = \tau(X_1)\sigma(X_2)\sau(X_3)\,.
    \end{equation*}
    
    The definition of the angles is motivated by the work of Berenstein-Retakh on noncommutative surfaces \cite{berenstein2018noncommutative}. Namely, label the vertices of the triangle as indicated in \Cref{fig:angleEx}, and rename the variables
    \begin{equation*}
        X_1 = x_{21}\qquad X_2 = x_{32}\qquad X_3 = x_{31}
    \end{equation*}
    and $x_{ji}:=\sigma(x_{ij})$ according to the orientation of the corresponding sides (this is compatible with our notion of switch). Then the angle reads
    \begin{equation*}
        L(\gamma) = \tau(x_{21})x_{23}\tau(x_{13}) = N_1N_3\cdot x_{21}^{-1}x_{23}x_{13}^{-1} =: N_1N_3\cdot T_1^{2,3}\,.
    \end{equation*}
    This is exactly the noncommutative angle $T_1^{2,3}$ from \cite{berenstein2018noncommutative} up to multiplication by the central variables.
\end{example}

\begin{definition}
	The \keyword{seed group} $\quiverAlgebra_Q$ is obtained by taking the quotient of the free group $F_Q$ by the relations:
	\begin{align*}
		& N_i  \quad  \text{is central} \\
        & \Delta = \sigma(\Delta) \text{ for any angle } \Delta\,.
	\end{align*}
    We write $\A_Q$ for the noncommutative space whose group of functions is given by the seed group $\quiverAlgebra_Q$. We call it the \keyword{seed torus}.
\end{definition}

Consider an angle $\Delta$ given by a loop $\gamma$ around a polygon. By reversing the loop $\gamma$, we obtain the label $L(\overline{\gamma})=\sigma\big(L(\gamma)\big)$, so the relation for the angles in the seed group can be understood as assigning the same label to a loop around a polygon and its reverse.

\begin{remark}
    In \cite{greenberg2024noncommutative}, we constructed an algebra associated to each quiver instead and extended $\sigma$ and $\tau$ linearly to the entire algebra. However, this approach only works for $r=2$, while we extend the theory here to include the case $r=3$.
\end{remark}

If two polygonal quivers are related by one of the operations (switching or weaving) from \Cref{def:weavingQuiver}, they have isomorphic seed groups:
\begin{proposition}\label{def:weavingGroup}
    Let $Q$ be an polygonal quiver with a noncommutative node $k$. Then we have an isomorphism defined by
    \begin{align*}
        &\begin{matrix}
            \sigma_k: & \quiverAlgebra_Q & \overset{\sim}{\to} & \quiverAlgebra_{\sigma_k(Q)}\\
            & X_i^{(f)} & \mapsto & \begin{cases}
                X_i^{(f)} & i\neq k\\
                X_k^{(f+1)} & i = k
            \end{cases}\\
            & N_i & \mapsto & N_i\,.
        \end{matrix}\\
        \intertext{When $r=2$ we have another isomorphism defined by}
        &\begin{matrix}
            w_k: & \quiverAlgebra_Q & \overset{\sim}{\to} & \quiverAlgebra_{w_k(Q)}\\
            & X_i^{(f)} & \mapsto & \begin{cases}
                X_i^{(f)} & i\neq k\\
                \big(X_k^{(f+1)}\big)^{-1}N_k & i = k
            \end{cases}\\
            & N_i & \mapsto & N_i\,.
        \end{matrix}
    \end{align*}
\end{proposition}

Thus, whenever we perform a switch at $k$, we will replace $X_k$ with $\sigma(X_k)$ in all formulas. Similarly, when we weave at $k$, we replace $X_k$ with $\sau(X_k)$. All other variables remain unchanged.

\begin{proof}
    It suffices to check that $\sigma_k,\, w_k$ map angles onto angles: For instance if an angle path goes through $k$, the angle in $\quiverAlgebra_{\sigma_k(Q)}$ is obtained from the angle for the same path in $\quiverAlgebra_Q$ by replacing $X_k$ with $\sigma(X_k)$, which is exactly the effect of the map $\sigma_k$ defined above.
\end{proof}

These maps define weaving and switch maps of the associated seed tori.

\subsection{Jordan algebra points}
So far we have given a description of the mutation of polygonal quivers, and we have assigned a group $\quiverAlgebra_Q$ and a torus $\A_Q$ to each one. When $r=2$ it is possible to define a \keyword{polygonal cluster algebra} given by gluing together quotients of a quiver \emph{algebra} \cite{greenberg2024noncommutative}. The defining relations were constructed so that the algebra could be evaluated in arbitrary Clifford algebras. 

However for $r=3$, such an algebra does not exist. Nevertheless, we can define ``evaluations'' of $\quiverAlgebra_Q$ into arbitrary Jordan algebras. These evaluations can be glued together to form a geometric object in analogy with the construction of a cluster variety for classic commutative cluster algebras.  

Now let $\jordan{J}=(V,\iota,\Id)$ be a normed Jordan algebra of degree $r$. 

\begin{definition}
    Let $Q$ be an polygonal quiver. A \keyword{$\jordan{J}$-point} of the seed tori $\A_Q$ is a group homomorphism $\phi:\quiverAlgebra_Q \to \Gamma(\jordan{J})$ such that:
    	\begin{enumerate}
		\item $\phi$ intertwines the $\sigma$ and $\tau$ defined for $\quiverAlgebra_Q$ with those defined for $\Gamma(\jordan{J})$.
		\item $\phi(\Delta) = \iota(v_\Delta)$ for some $v_\Delta\in V$ for any angle $\Delta$.
	\end{enumerate}
    We denote by $\A_Q(\jordan{J})$ the set of $\jordan{J}$-points of $\A_Q$. 
\end{definition}

\begin{remark}\label{rem:angleTransition}
    Because of the first requirement, we can define a $\jordan{J}$-point of $\A_Q$ by picking images $\phi(X_i)$ for $i\in Q_0^{(2)}$ and $\phi(N_j)$ for $j\in Q_0^{(1)}$. This determines the images of all other generating variables. Also, observe that $\phi(N_i)\in\K^\times$ for any $i\in Q_0$.\\
    For the second condition, it suffices to check one angle per polygon. This is because angles in the same polygon with different base points are related by applying $\sau$ and $\Delta\mapsto A\Delta\sigma(A)$ for $A\in\Gamma(\jordan{J})$. Both of these operations preserve $\iota(V)\subset\Gamma(\jordan{J})$.
\end{remark}

If $\jordan{J}$ moreover has a positive structure, we can define positive points:
\begin{definition}
    We say that a $\jordan{J}$-point $\phi:\quiverAlgebra_Q\to\jordan{J}$ is \keyword{positive} if 
    \begin{enumerate}
        \item $\phi(X_i) \in \Gamma(\jordan{J})^\circ$ for all $i\in Q_0^{(2)}$.
        \item $\phi(N_i) \in \K^{>0}$ for all $i \in Q_0^{(1)}$.
		\item $v_\Delta\in V^{>0}$ for any angle $\Delta$.
    \end{enumerate}
    The set of such points is denoted $\A_Q(\jordan{J}^{>0})$.
\end{definition}

\begin{remark}
    Again, it is sufficient to check condition (3) for one angle per polygon as the operations in \Cref{rem:angleTransition} preserve the cone if we restrict to the identity component of the structure group, which is condition (1).
\end{remark}

Our next goal is to define a mutation of Jordan points, i.e. a map $\A_Q(\jordan{J})\to\A_{Q'}(\jordan{J})$ that is defined whenever $Q'$ is obtained from $Q$ by (admissible) mutation. Note that this map will only be defined for an open subset of Jordan points. This gives rise to a birational map of seed tori $\A_Q \to \A_{Q'}$.

\begin{definition}\label{def:JordanPointMutation}
    Given an admissible mutation $\mu_k \colon Q \to Q'$ at a noncommutative node $k\in Q_0^{(2)}$ we define a mutation map $\mu_k:\A_Q(\jordan{J})\to\A_{Q'}(\jordan{J})$ by requiring that
    \begin{equation*}\label{eqn:FullMutationFormula}
        \big(\mu_k(\phi)\big)(X'_j) = \phi(X_j) \hspace{1cm}\text{for } j \neq k 
    \end{equation*}
    and that the angles as indicated in \Cref{fig:mutationAngles} follow the \keyword{angle summation formula}
    \begin{equation}
        \big(\mu_k(\phi)\big)\big(\sau^{\varepsilon_f}(\Delta_{f}'(k))\big) = \iota\left(\frac{\phi(O_f)}{\phi(N_k)}v_{\Delta_f^l(k)}+ \frac{\phi(I_f)}{\phi(N_k)}v_{\Delta_f^r(k)}\right)\,,
    \end{equation}
    for $f\in\{\circ,\bullet\}$ where the powers $\varepsilon_f$ are chosen such that $\sau^{\varepsilon_f}(\Delta_f'(k))$ contains $\sigma^f(X_k')$ (no $\tau$). Moreover, we have defined
    \begin{align*}
        O_f &= \prod_{i \in \Out^c(k)} N_i\cdot \prod_{\alpha \in \Out_{f+1}(k)}N_\alpha\,,\\
        I_f &= \prod_{i \in \In^c(k)} N_i\cdot \prod_{\alpha \in \In_{f}(k)}N_\alpha\,.
    \end{align*}
    This mutation is defined whenever the argument of $\iota$ in \Cref{eqn:FullMutationFormula} lies in $V^\times$, which is a generic condition on $\phi$.
\end{definition}

\input{Figures/figureMutationAngles}

\begin{remark}
    Note that in the mutation map, the only elements which are added are in the vector space $V$ (which gets mapped to $\Gamma(\jordan{J})$ via $\iota$), and the only elements which are multiplied belong to $\Gamma(\jordan{J})$. As such $\jordan{J}$ does not need to be contained in an algebra.
\end{remark}

If we consider a positive $\jordan{J}$-point of $\quiverAlgebra_Q$, the new angle always lies in $\iota(V^\times)$:
\begin{lemma}\label{thm:admissibleMutationPositivity}
    An admissible mutation $\mu_k : Q \rightarrow Q'$ sends positive $\jordan{J}$-points to positive $\jordan{J}$-points.
\end{lemma}
\begin{proof}
    Consider a positive $\jordan{J}$-point of $\A_Q$. Then, in particular, the elements $v_{\Delta_f^l(k)}$ and $v_{\Delta_f^r(k)}$ from the angle summation formulas are in $V^{>0}$. This positive cone is preserved by the action of elements in $\Gamma(\jordan{J})^\circ$ and $\R^{>0}$ and addition. Thus the new angles given by \Cref{eqn:FullMutationFormula} are again in $V^{>0}$. As all other angles are unchanged or equivalent to the new angles by actions of elements in $\Gamma(\jordan{J})^\circ$ the resulting point is also positive. 
\end{proof}

\begin{definition}
    We write $\A_{|Q|}$ for the noncommutative space which is obtained by gluing together the seed tori along all sequences of admissible mutation maps. Given a normed Jordan algebra of degree $r$ $\jordan{J}$ we write $\A_{|Q|}(\jordan{J})$ for the set of Jordan points of $\A_{|Q|}$.
\end{definition}

The previous lemma implies that the set of positive Jordan points $\A_{|Q|}(\jordan{J}^{>0})$ can be computed on any single seed torus, as the mutation maps are isomorphisms on the positive points. 

In \cite{greenberg2024noncommutative} for the case $r=2$, we construct an analog of the cluster algebra called the \keyword{noncommutative polygonal cluster algebra} which contains all of the cluster variables obtained by admissible mutation and prove a noncommutative Laurent phenomenon. This implies that the ring of global functions on $\A_{|Q|}$ contains the noncommutative polygonal cluster algebra associated to $Q$. 

\begin{example}
    To make the mutation formulas more explicit, let us give an example. Consider again, the mutation shown in \Cref{fig:mutationEx}. Suppose that $\phi$ is a $\jordan{J}$-point for the seed torus associated to the left quiver. Then $\phi':=\mu_6(\phi)$ is defined by the angle addition formulas
    \begin{align*}
        \phi'\big(\sau(\Delta_\circ'(6))\big) &= \phi'\big(\sigma(X_3)\tau(X_2)X_6'\sau(X_5)\big) = \iota\left(\phi\Bigg(\frac{N_2N_3}{N_6}\Bigg) v_{\Delta_\circ^l(6)} + \phi\Bigg(\frac{N_5}{N_6}\Bigg) v_{\Delta_\circ^r(6)}\right)\\
        \phi'\big(\Delta_\bullet'(6)\big) &= \phi'\big(\sau(X_4)\sigma(X_6')\tau(X_1)\big) = \iota\left(\phi\Bigg(\frac{N_4}{N_6}\Bigg)v_{\Delta_\bullet^l(6)} + \phi\Bigg(\frac{N_1}{N_6}\Bigg)v_{\Delta_\bullet^r(6)} \right)
    \end{align*}
    with the vectors determined by the equations
    \begin{align*}
    \begin{array}{lcl}
        \iota(v_{\Delta_\circ^l(6)}) = \phi\big(\sau(X_6)X_4\sau(X_5)\big) & & \iota(v_{\Delta_\circ^r(6)}) = \phi\big(\sigma(X_3)\tau(X_2)\sigma(X_1)\tau(X_6)\big) \\
        \iota(v_{\Delta_\bullet^l(6)}) = \phi\big(\tau(X_6)\tau(X_3)\sigma(X_2)\tau(X_1)\big) & & \iota(v_{\Delta_\bullet^r(6)}) = \phi\big(\sau(X_4)X_5\sau(X_6)\big)\,.
    \end{array}
    \end{align*}
\end{example}

\subsection{Mutation in special cases}
We will now motivate the mutation procedure for $\jordan{J}$-points and in particular the angle summation formula \Cref{eqn:FullMutationFormula}: Consider again the flip in \Cref{fig:flipDecPoly} and the zigzag mutation procedure for surface cluster algebras as visualized in \Cref{fig:zigZagCommutative}. In \Cref{fig:zigZagDecPoly1} we have combined both ideas and outlined two paths in a schematic picture of a decorated tiling.

\input{Figures/figureZigZagNoncomm}

As before, let the node at which we mutate be $k$. From the figure a first guess for a mutation formula would be
\begin{equation}\label{eq:naivePCAmutation}
    X_k'N_k = I_c\cdot\sau^{\varepsilon_1}\big({\color{Red}L(\gamma_1)}\big) + O_c\cdot\sau^{\varepsilon_2}\big({\color{Blue}L(\gamma_2)}\big)
\end{equation}
where $\gamma_1$ and $\gamma_2$ are the angle paths marked in the figure. We choose $\varepsilon_1,\,\varepsilon_2$ such that $\tau(X_k)$, respectively $\sau(X_k)$ appears in the summands. This mirrors \Cref{eq:zigZagExchangeCommutative} on the noncommutative part since $\tau(X_k)/N_k=X_k^{-1}$. We have also added the central variables coming from small nodes and paired arrows
\begin{equation*}
    I_c := \prod_{i \in \In^c(k)} N_i\,,\qquad O_c := \prod_{i \in \Out^c(k)}N_i\,.
\end{equation*}

This exchange relation is essentially the one used in \cite{greenberg2024noncommutative} to glue the quotients of the quiver algebras. As we already remarked above, this cannot be done when $r=3$, essentially since we cannot assign a new central element $N_k'$. This is why we only define mutation on $\jordan{J}$-points of the quiver algebra.\medskip

Thus, suppose that $\phi:\quiverAlgebra_Q\to\Gamma(\jordan{J})$ is a $\jordan{J}$-point and that we performed an admissible mutation at a noncommutative node $k$ to get $\mu_k(Q)=Q'$. To construct the $\jordan{J}$-point $\phi'$ of $\quiverAlgebra_{Q'}$, we keep the image of all variables $X_i,\,N_i$ for $i\neq k$ unchanged. A reasonable idea to define $\phi'(X_k')$ is to use \Cref{eq:naivePCAmutation} and just replace all variables by their image under $\phi$. However, the image of $\phi$ lies in the \emph{group} $\Gamma(\jordan{J})$, so it does not make sense to add the two terms unless $\Gamma(\jordan{J})$ is contained in some algebra (which will often be the case).\medskip

Instead, as we saw mutation should be performed using the angles: Consider \Cref{fig:zigZagDecPoly2}, the local picture of the decorated tiling after the flip. The label assigned to the indicated path is an angle in $\Delta'\in\quiverAlgebra_{Q'}$. If we write it down and plug in the naive mutation formula \Cref{eq:naivePCAmutation}, we find
\begin{equation}\label{eq:angleMutation}
    \sau^\varepsilon\big(\Delta'\big) N_k = I\cdot{\color{Red}\Delta_1} + O\cdot{\color{Blue}\Delta_2}
\end{equation}
where $\Delta_1$ is the angle corresponding to the red loop in the right polygon, and $\Delta_2$ corresponds to the blue loop in the left polygon in \Cref{fig:zigZagDecPoly3}. $\varepsilon$ is chosen such that $\sau^\varepsilon(\Delta')$ contains $X_k'$ (not $\sau(X_k')$). Moreover
\begin{equation*}
    I = I_c\cdot\prod_{i\in\In_\bullet(k)}N_i\,,\hspace{2cm} O = O_c\cdot\prod_{i\in\Out_\circ(k)}N_i
\end{equation*}
are central terms which correspond to going back and forth as indicated in \Cref{fig:zigZagDecPoly3}, as well as the commutative nodes and paired arrows as before. All loops used to compute angles have base points as above, see \Cref{fig:mutationAngles}. We could also deduce a similar equation for the other angle after the flip. If $\Gamma(\jordan{J})$ is contained in an algebra and $\iota$ provides a linear embedding of $V$, this approach starting with the zigzag paths, leads exactly to the mutation we introduced in \Cref{def:JordanPointMutation}. Importantly, this is the case for spin groups with Jordan split type $B_p$.

\begin{remark}
    Deducing \Cref{eq:angleMutation} from \Cref{eq:naivePCAmutation} requires a slightly more careful analysis of the powers of $\sau$. This uses the parity condition for decorated polygons and can be found in \cite{greenberg2024noncommutative}.
\end{remark}

\section{Noncommutative networks}\label{sec:NoncomNetworks}

In higher noncommutative rank, we use different combinatorial data to provide coordinates, a planar bipartite network. These networks were introduced by Postnikov to study the Grassmannian \cite{postnikov-TotalPositivityNetworks} and have many applications to cluster algebra theory \cite{fwz_IntroToClusters_PlabicGraphs}.  The noncommutative generalization we consider here is mostly the same as that of \cite{goncharov2021spectral}.

\begin{definition}
    A \keyword{plabic network} on a disk is a bipartite graph embedded in a disk with edge weights in an (associative) algebra $\mathcal{R}$ such that
    \begin{enumerate}
        \item A vertex is degree 1 if and only if it is on the boundary of the disk.
        \item The graph in the interior of the disk is planar.
        \item The counterclockwise product of edge weights around any vertex is $1$.
    \end{enumerate}
\end{definition}

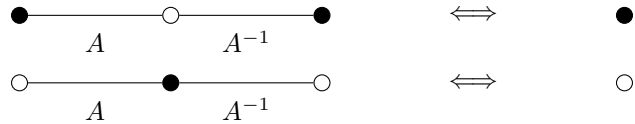
\begin{figure}[hbt]
    \centering
    \begin{tikzpicture}
        \node[mutable] (l1) at (-2, 0) [] {};
        \node[mutable,fill=white] (c1) at (0, 0) [] {};
        \node[mutable] (r1) at (2, 0) [] {};
        \draw (l1) to (c1) to (r1);
        \node[] (Al) at ($1/2*(l1)+1/2*(c1)$) [label=below:$A$] {};
        \node[] (Ar) at ($1/2*(r1)+1/2*(c1)$) [label=below:$A^{-1}$] {};
        \node[] (iff) at (4,0) [] {$\iff$};
        \node[mutable] (res) at (6,0) [] {};
    \end{tikzpicture}
    
    \begin{tikzpicture}
        \node[mutable,fill=white] (l1) at (-2, 0) [] {};
        \node[mutable] (c1) at (0, 0) [] {};
        \node[mutable,fill=white] (r1) at (2, 0) [] {};
        \draw (l1) to (c1) to (r1);
        \node[] (Al) at ($1/2*(l1)+1/2*(c1)$) [label=below:$A$] {};
        \node[] (Ar) at ($1/2*(r1)+1/2*(c1)$) [label=below:$A^{-1}$] {};
        \node[] (iff) at (4,0) [] {$\iff$};
        \node[mutable,fill=white] (res) at (6,0) [] {};
    \end{tikzpicture}
    \caption{Local equivalence moves of plabic networks.}
    \label{fig:NetworkMove}
\end{figure}

Plabic networks are considered up to natural equivalences given by isotopy of the disk fixing the boundary and contracting/adding vertices of degree 2 (\Cref{fig:NetworkMove}).

We will be interested in control over the boundary edges which cause the network to be nonplanar.
\begin{definition}
    A \keyword{triangular grounded wiring network} on $n$ strands is a bipartite network embedded in a triangle with a designated \keyword{ground side} such that
    \begin{enumerate}
        \item The left and right sides of the triangle have $n$ $\bullet$ vertices labeled $1$ to $n$.
        \item The edges incident to the ground are called \keyword{legs}. Removing all the legs results in a planar graph.
        \item The non leg edges can be partitioned into two groups, \keyword{strands} and \keyword{bridges} such that the side boundary nodes labeled $i$ are connected by a unique path of strands.
    \end{enumerate}
    A general \keyword{grounded wiring network} on a disk is a network embedded in disk with a triangulation of the disk such that the network restricted to each triangle is a triangular grounded wiring network. 
\end{definition}
A grounded wiring network can be obtained by amalgamating copies of the elementary networks in \Cref{fig:ElementaryNetworks} along the boundary vertices with matching strand numbers and potentially deleting legs.
\input{Figures/figureElementaryNetworksNoweights}

\begin{definition}\label{def:NetworkGroup}
    Each grounded wiring network $P$ has an associated \keyword{seed group} $\quiverAlgebra_P$ obtained by taking the quotient of the free group generated by the edges of the network by the relation that the counterclockwise product around each vertex of $P$ is $1$. The \keyword{seed torus}, $\A_P$ is the noncommutative torus with $\quiverAlgebra_P$ as its group of coordinate functions. 
 
\end{definition}

\begin{remark}
    The condition on the weights that the counterclockwise product around each vertex is equal to 1 is different from the convention in \cite{goncharov2021spectral}. We will expand on this distinction in \Cref{sec:spectral_discription}.
\end{remark}

\begin{definition}
    Let  $\jordan{J}$ be a special Jordan algebra with specializing ring $\mathcal{R}$. A \keyword{$\jordan{J}$-point} of a grounded wiring network $P$ is a group homomorphism $\quiverAlgebra_P \rightarrow \mathcal{R}^\times$. 
\end{definition}

We consider two grounded wiring networks equivalent if they are related by \keyword{leg slide} and identification of legs attached to the same strand. 
\begin{definition}\label{def:LegSlide}
    A \keyword{leg slide} is performed by moving a leg left or right horizontally until it travels across a bridge onto a neighboring strand (\Cref{fig:LegSlide}). \input{Figures/figureLegSlide}

    A leg slide is associated with the following monomial map on the weights:
    Let $w$ be the weight of the leg. Each edge that is traversed from $\circ$ to $\bullet$ has its weight left multiplied by $w$. Each edge traversed from  $\bullet$ to $\circ$  has its weight right multiplied by $w^{-1}$.
\end{definition}

\begin{definition}\label{def:SquareMove}
    There is a natural mutation of weighted plabic networks given by the \keyword{square move} with a noncommutative generalization given by \cite[Equation 68]{goncharov2021spectral}. \Cref{fig:SquareMove} shows the local transformation of the network and the birational map from $\A_P \rightarrow \A_{P'}$ and its inverse are given in \Cref{eqn:SquareMove}.
\end{definition}
\begin{align}\label{eqn:SquareMove}
    \begin{aligned}
        b_1&=(1+a_3a_4a_1a_2)a_2^{-1}a_1^{-1}a_4^{-1}\\
         &=a_2^{-1} a_1^{-1} a_4^{-1} (1+a_4 a_1 a_2 a_3)\\
        b_3&=(1+a_1a_2a_3a_4)a_4^{-1}a_3^{-1}a_2^{-1}\\
           &=a_4^{-1}a_3^{-1}a_2^{-1}(1+a_2 a_3 a_4 a_1)\\
        \vspace{1pc}\\
        a_1 &= b_3(1+b_4 b_1 b_2 b_3)^{-1}\\
            &= (1+b_3 b_4 b_1 b_2)^{-1} b_3\\
        a_3 &= b_1(1+b_2b_3b_4b_1)^{-1}\\
            &= (1+b_1 b_2 b_3 b_4)^{-1}b_1
    \end{aligned}&
    \hspace{2pc}
    \begin{aligned}
        b_2 &= a_4(1+a_1a_2a_3a_4)^{-1}\\
            &= (1+a_4a_1a_2a_3)^{-1}a_4\\
        b_4 &= a_2(1+a_3a_4a_1a_2)^{-1}\\
            &= (1+a_2a_3a_4a_1)^{-1}a_2\\
        \vspace{1pc}\\
        a_2 &= (1+b_4b_1b_2b_3)b_3^{-1}b_2^{-1}b_1^{-1}\\
            &= b_3^{-1}b_2^{-1}b_1^{-1}(1+b_1b_2b_3b_4)\\
        a_4 &= (1+b_2b_3b_4b_1)b_1^{-1}b_4^{-1}b_3^{-1}\\
            &= b_1^{-1}b_4^{-1}b_3^{-1}(1+b_3 b_4 b_1 b_2)
    \end{aligned}
\end{align}
\input{Figures/figureSquareMove}
Note that for $i \in \Z/4\Z$, $a_ia_{i+1} = (b_ib_{i+1})^{-1}$. Thus the products around each vertex in square are preserved, which ensures the result of a square move is a correctly weighted network. 

\begin{definition}\label{def:HoneycombNetwork}
    The \keyword{honeycomb networks}, $P_n$, are a family of networks embedded in a triangle following the pattern demonstrated in \Cref{fig:GlnNetwork}.
\end{definition}
The honeycomb networks are used to parameterize triples of $\GL_n(\R)$ flags in \cite{goncharov2021spectral}. We will recover the honeycomb networks when parameterizing groups with an $A_p$-grading in \Cref{sec:NoncomAp}. These networks are important as they are fully planar, so they can be regrounded.

\input{Figures/figureGLn_network}

\begin{definition}\label{def:ChangeOfGround}
    A full planar triangular grounded network can be \keyword{regrounded} by choosing a new side of the triangle to be the ground side. 
\end{definition}
In \Cref{fig:GlnNetwork}, we observe that any side of the triangle can be chosen to be the ground for a honeycomb network. The horizontal strands are identified by ``pulling the sides tight'' parallel to the ground. 

\begin{definition}
    We write $|P|$ for the set of grounded networks which are obtained from $P$ by equivalence moves, square moves, leg slides and regroundings. We write $\A_{|P|}$ for the noncommutative space obtained by gluing together the seed tori $\A_P$ along the birational maps defined by each of these operations.
\end{definition}

 From this definition we obtain noncommutative spaces from grounded wiring networks and from polygonal cluster algebras. We refer to both these constructions as \keyword{noncommutative cluster varieties}. Each space is built by gluing seed groups via noncommutative mutation relations. Although the construction of the seed group is different in each case, when we refer to the seed group the correct construction will be clear from context.

\part{Cluster-like Coordinates on Decorated Flag Varieties}\label{part:ClusterlikeCoordinates}
We will now use Jordan pinnings to define coordinates on configuration spaces of decorated flags. Thus, let $G$ be Jordan split over a root system $R$ over the $\K$ Jordan algebras $\{\jordan{J_i}\}$. We will write $p$ for the rank of $R$.  We recall that there are 3 cases based on the noncommutative rank of $G$ of type $R$.
\begin{itemize}
    \item When $G$ is non-commutative rank 0 of type $R$, we will see that our coordinates are obtained by evaluating the commutative cluster varieties of \cite{goncharov2019quantum} in field extensions of $\K$.
    \item When $G$ is non-commutative rank 1 of type $R$, our coordinates will give evaluations of non-commutative polygonal cluster algebras with general $r$ for $A_1$, $r=2$ for $B_p$ and $r=3$ for $G_2$.  In this case we will assume that the simple roots in $R$ are indexed such that the noncommutative root is the $p$-th one.
    \item When $G$ is higher noncommutative rank of type $R$, we will relate our coordinates to the description in \cite{goncharov2021spectral}.
\end{itemize}
Together this amounts to giving for each $R$ a noncommutative cluster variety $\A_{R,n}$ whose Jordan points parameterize configurations of $n$ decorated flags.

\section{Setup and strategy}\label{sec:ClusterCoordinatesSetup}
We consider the partial flag variety defined by the associated parabolic $P_\Theta$ to the $R$-grading. We denote the flag variety by $\flagTheta:=G/P_\Theta$ and the decorated flag variety by $\decFlagTheta:=G/U_\Theta$. Moreover, let $\pi:\decFlagTheta \to \flagTheta$ be the projection map induced by the inclusion $U_\Theta\subset P_\Theta$. In this case the fiber over a point in $\flagTheta$ is an $L_\Theta$-torsor. The \keyword{configuration space} of $n$ decorated flags for $n\in\mathbb{N}$ is given by
\begin{equation*}
    \Conf_n(\decFlagTheta):=G\backslash\decFlagTheta^n
\end{equation*}
using the diagonal action of $G$ on the product. Our ultimate goal is to endow this space with cluster-like coordinates.

\subsection{Pairs of flags}\label{sec:pairsFlagsTheta}
We start by considering pairs of flags. Let $w_0\in W(R)$ be the element of maximal length which we already encountered in \Cref{sec:LiftsOfThetaWeylGroup}.

\begin{definition}
    Two flags $B_1,B_2\in\flagTheta$ are called \keyword{transverse} if there exists $g\in G$ such that $(gB_1,gB_2)=(B,\overline{w_0}B)$. This is denoted $B_1\pitchfork B_2$. Two decorated flags $A_1,A_2\in\decFlagTheta$ are called transverse if $\pi(A_1)\pitchfork\pi(A_2)$. \\For an element $w\in W(R)$, we call a two-tuple of flags, resp. decorated flags, \keyword{$w$-transverse} if there exits $g\in G$ such that $(gB_1,gB_2)=(B,\overline{w}B)$ resp. $(g\pi(A_1),g\pi(A_2))=(B,\overline{w}B)$.
\end{definition}

\begin{remark}
    By \Cref{thm:ThetaLargestCellDense} transversality is a generic condition for pairs of (decorated) flags. But as we saw in \Cref{eq:thetaBruhatDecomp}, for $B_1,B_2\in\flagTheta$ there does not necessarily exist a $w$ in the grading Weyl group such that $(gB_1,gB_2)=(B,\overline{w}B)$. However, all flags that we will consider will be $w$-transverse for some $w$ in the grading Weyl group.
\end{remark}

\begin{remark}
    If $(A_1,A_2)$ is $w$-transverse, then $(A_2,A_1)$ is $w^{-1}$-transverse. To construct coordinates, we will only need the notion of $w_0$-transverse and $\varsigma_i$-transverse where this distinction vanishes. 
\end{remark}

We define the basic invariant for a pair of flags using that the fibers of $\pi$ are $L_\Theta$-torsors:
\begin{definition}
    For two $w$-transverse decorated flags, we define the \keyword{$L$-distance} $l(A_1,A_2)$ to be the unique element of $l \in L_\Theta$ such that there exists an element $g\in G$ for which 
    \[g(A_1,A_2)=(U_\Theta,\overline{w}U_\Theta l\,) \,.\] 
\end{definition}
This is constant on $G$-orbits and thus defines a map $\Conf_2(\decFlagTheta)\to L_\Theta$. Given a Jordan pinning and choice $\{\Lambda_i\}$ of Jordan weights, we can reconstruct this element of $L_\Theta$ as
\begin{equation*}
    l= \prod_i\check{\beta}_i(\Lambda_i(l))\,.
\end{equation*}
Therefore, our collection of coordinates on $\Confx_2(\decFlagTheta)$, the space of transverse configurations, is the set $\{\Lambda_i(l(A_1,A_2)\}\in \prod_i\Gamma(\jordan{J}_i)$. Before we take a closer look at the $L$-distance, let us describe the general strategy for the construction of cluster coordinates on $\Conf_n(\decFlagTheta)$.

\subsection{General strategy}\label{sec:coordinateStrategy}
Consider the space $\Confx_n(\decFlagTheta)\subset\Conf_n(\decFlagTheta)$ defined by the condition that any two flags in the configuration are transverse. For split groups $G$ coordinates on this space were previously constructed: Fock-Goncharov pioneered the approach presented below and described concrete coordinates for $G=\SL(n,\R)$ \cite{fock2006moduli}. This was extended to all classical split algebraic groups by Le \cite{le2019cluster} and further to include exceptional groups by Gilles \cite{gilles2021fock}. Finally, a very general construction was provided by Goncharov-Shen \cite{goncharov2019quantum}. Our goal will be to extend their construction to $(G,\Theta)$ reduced of Weyl type.

Thus, let $[(A_1,\dots,A_n)]\in\Confx_n(\decFlagTheta)$. Its coordinates are obtained as follows:\medskip

\noindent\textbf{\color{Orange}Step 1:} We \keyword{associate the flags to the vertices of an $n$-gon} in a counterclockwise manner.\medskip

\noindent\textbf{\color{Orange}Step 2:} Next, \keyword{fix a triangulation of the $n$-gon}. This naturally breaks the configuration of $n$ flags into $n-2$ configurations of 3 flags (associated to the triangles), which share sub-configurations of 2 flags (associated to the edges). This amalgamation approach was introduced by Fock and Goncharov and they proved that the configurations associated to the triangles determine the full configuration for split $G$ \cite[Thm. 8.1]{fock2006moduli}. We extend this in \Cref{thm:coordinatesDetectTuplePositivityTheta}.\medskip

Therefore, it suffices to \keyword{describe coordinates on $\Confx_3(\decFlagTheta)$}.\medskip

\noindent\textbf{\color{Orange}Step 3:} We further decompose a configuration of 3 flags into \keyword{elementary configurations} (see \Cref{sec:triplesToElementaryConfigs}). These are elements of $\Conf_3(\decFlagTheta)$ which are not fully transverse. Instead, one pair of flags is only $\varsigma_i$-transverse for $\varsigma_i$ a simple reflection in $W(R)$. See \Cref{def:ElementaryConfiguration}. These configurations are obtained from a configuration of 3 decorated flags by computing an \keyword{interpolating sequence} of decorated flags between one pair of flags in the configuration (or along an edge of the triangle), see \Cref{thm:flagSequenceTheta}. It depends on a choice of reduced expression for the longest word $w_0\in W(R)$. The set of elementary configurations determines the full configuration. This approach is due to Goncharov-Shen in the split case \cite{goncharov2019quantum}, and we extend it in \Cref{thm:coordinatesTriplesTheta}.\medskip

\noindent\textbf{\color{Orange}Step 4:} Invariants of elementary configurations can be obtained by computing pairwise $L$-distances between the flags. To \keyword{produce (noncommutative) coordinates} from this, one applies the Jordan weights to these $L$-distances as explained previously. These will be assigned to the nodes of a \keyword{elementary quiver} or edges of an \keyword{elementary network}, which are upgraded version of the ordinary quivers described in \cite{goncharov2019quantum}. Here it will be a polygonal quiver in noncommutative rank 1, see \Cref{sec:elementarySeedsLowerNCrank} or a grounded wiring network in higher noncommutative rank \Cref{sec:NoncomNetworks}. We will see in \Cref{thm:tripleFlagsPolyCAEvaluation} that the coordinate functions define a $\jordan{J}$-point of the corresponding noncommutative cluster variety.\medskip

\noindent\textbf{\color{Orange}Step 5:} All the elementary seeds are \keyword{amalgamated} to obtain a quiver, with corresponding coordinate functions, first for $\Confx_3(\decFlagTheta)$, and then for $\Confx_n(\decFlagTheta)$.\medskip

\input{Figures/figureGeneralStrategy}

Our main goal now is to understand how the choices we made effect these coordinates:
\begin{enumerate}
    \item the choice of triangulation,
    \item the choice of edge within each triangle along which to compute an interpolating sequence,
    \item the choice of reduced expression of $w_0$ for each edge, 
    \item the choice of Jordan weights used to compute coordinates.
\end{enumerate}
Specifically, we will show that making different choices correspond to noncommutative mutations. We only need to consider a few different elementary operations to make all of the changes above. Concretely, we will prove that we can use mutation to relate coordinate systems
\begin{itemize}
    \item on $\Confx_3(\decFlagTheta)$ with different reduced expressions (`2,3,4,6 moves', \Cref{sec:changingReducedExpressionTheta})
    \item on $\Confx_3(\decFlagTheta)$ which differ by the choice of interpolating edge (`rotation', \Cref{sec:rotationCharts})
    \item on $\Confx_4(\decFlagTheta)$ arising from the two triangulations of the square (`flip', \Cref{sec:FlipTriangulationGenericTheta}).
\end{itemize}
The choice of Jordan weights only appears in the higher noncommutative rank case and our coordinates will be naturally independent of this choice. 

\begin{remark}
    In \Cref{sec:ThetaPositivity}, we will see that, when the Jordan algebra $\jordan{J}_i$ has a positive structure, mutations sends positive points to positive points. This will recover the $\Theta$-positive structures on $\Confx_n(\decFlagTheta)$ defined in \cite{guichard2022generalizing}. Moreover, the above approach allows to define cluster-like coordinates on the moduli space of decorated twisted local systems and also yields a positive structure on this space as we will see in \Cref{sec:localSystems}.
\end{remark}

\subsection{Opposition involution}
We recall the \keyword{opposition involution} on a root system given by $\oppInv{\beta_i} = -w_0(\beta_i)$. For root systems of type $A_1,B_p$, $C_p$, $F_4$, or $G_2$, $\oppInv{\beta_i} = \beta_i$, but in $A_p$ we have $\oppInv{\beta_i}=\beta_{p-i+1}$ using the standard indexing of the simple roots. It is closely related to the following involution on $L_\Theta$:
\begin{definition}\label{def:OppositionInvolution}
    The \keyword{opposition involution} on $L_\Theta$ is the map
    \begin{equation*}
        (.)^*:L_\Theta\to L_\Theta\,,\quad l\mapsto w_0\big(l^{-1}\big)=: \oppInv{l}\,.
    \end{equation*}
    It induces opposition involutions $\check{\beta}_i\mapsto\check{\beta}_{\oppInv{i}}$ on the coroot system and $(.)^*:V_i\to V_{\oppInv{i}}$ on the vector spaces via
    \begin{equation*}
        \oppInv{\Big(\check{\beta}_i\big(\iota(v)\big)\Big)}=\check{\beta}_{\oppInv{i}}\big(\iota(\oppInv{v})\big)\,.
    \end{equation*}
\end{definition}

This map shows up naturally as we will see in \Cref{thm:HDistPropertiesTheta}(1). If $A_1\pitchfork A_2\in\decFlagTheta$, the $L$-distance satisfies
\begin{equation*}
    l(A_2,A_1)=s_G\oppInv{l(A_1,A_2)}\,.
\end{equation*}
Here $s_G= \overline{w}_0^2$ is a distinguished element of order 2 in the center of $G$. It depends only on the type of the root system $R$ and will be used in \Cref{sec:localSystems} to twist local systems.

\begin{lemma}\label{thm:w0conjugationPinning}
    We have the following identity for the Jordan pinning:
    \begin{equation*}
        \doverline{w_0}x_i(v)\overline{w_0}=y_{\oppInv{i}}(-\oppInv{v})\,.
    \end{equation*}
\end{lemma}
\begin{proof}
    Because of the action of $W(R)$ on the coroots, we know that $\overline{w_0}x_i(v)\doverline{w_0}=y_{\oppInv{i}}(v')$ for some $v'\in V_i$. Multiply both sides of this equation by $\overline{\varsigma_{\oppInv{i}}}$. On the one hand we obtain from $\Cref{thm:DecomposeXSigma}$ the following Gauss decomposition:
    \begin{equation*}
        \overline{\varsigma_{\oppInv{i}}}y_{\oppInv{i}}(v')=y_{\oppInv{i}}(-v')\check{\beta}_{\oppInv{i}}\big(\iota(-v')\big)x_{\oppInv{i}}(v'^{-1})\,.
    \end{equation*}
    On the other hand, we can write $w_0=w\varsigma_i=\varsigma_{\oppInv{i}}=\varsigma_iw^{-1}$ for some $w\in W(R)$\,. Then
    \begin{align*}
        \overline{\varsigma_{\oppInv{i}}}\,\doverline{w_0}x_i(v)\overline{w_0}
        &= \doverline{w}x_i(v)\overline{\varsigma_iw^{-1}}\\
        &= \doverline{w}\Big(y_i(v^{-1})\check{\beta}_i\big(\iota(v)\big)x_i(-v^{-1})\Big)\doverline{w}^{-1}\,.
    \end{align*}
    Since $w(\beta_i)$ is a positive root, this is a Gauss decomposition. As this is unique, we have
    \begin{equation*}
        \check{\beta}_{\oppInv{i}}\big(\iota(-v')\big) = w\Big(\check{\beta}_i\big(\iota(v)\big)\Big) = w_0\Big(\check{\beta}_i\big(\iota(v)\big)^{-1}\Big) = \check{\beta}_{\oppInv{i}}\big(\iota(\oppInv{v})\big)
    \end{equation*}
    and thus since $\iota$ is injective, $\oppInv{v}=-v'$.
\end{proof}

\subsection{\texorpdfstring{$L$}{L}-distance revisited}
We are now prepared to take a closer look at the $L$-distance which we defined in \Cref{sec:pairsFlagsTheta}. First, we establish its basic properties:
\begin{proposition}\label{thm:HDistPropertiesTheta}
    Let $w\in W(R)$, and $(A_1,A_2)\in\Conf_2(\decFlagTheta)$ which are $w$-transverse.
    \begin{enumerate}
        \item The $L$-distance is not generally symmetric. In particular, we have if $w=\varsigma_i$
            \begin{align*}
                l(A_2,A_1) &= \check{\beta}_i(-1)\varsigma_i\big(l(A_1,A_2)^{-1}\big)\,.
            \end{align*}
            If $A_1\pitchfork A_2$, i.e. if $w=w_0$, then
            \begin{align*}
                l(A_2,A_1)=s_G w_0\big(l(A_1,A_2)^{-1}\big) = s_G \oppInv{l(A_1,A_2)}\,.
            \end{align*}
        \item For $l_1,l_2\in L_\Theta$
            \begin{equation*}
                l(A_1 l_1,A_2 l_2)=w^{-1}\big(l_1^{-1}\big) l(A_1,A_2)\,l_2\,.
            \end{equation*}
        \item Let $(A_2,A_3)$ be a pair of $w'$-transverse flags. If $\WeylLength(ww')=\WeylLength(w)+\WeylLength(w')$, then
            \begin{equation*}
                l(A_1,A_3)=(w')^{-1}\big(l(A_1,A_2)\big)l(A_2,A_3)\,.
            \end{equation*}
        
    \end{enumerate}
\end{proposition}

The proof uses the following:
\begin{lemma}\label{thm:WDistAdditiveTheta}
    Let $A_1,A_2,A_3\in\decFlagTheta$ such that $(A_1,A_2)$ are $w$-transverse and $(A_2,A_3)$ are $w'$-transverse for some $w,w'\in W(R)$. If $\WeylLength(ww')=\WeylLength(w)+\WeylLength(w')$, then $(A_1,A_3)$ are $ww'$-transverse.
\end{lemma}
\begin{proof}
    We may assume that $A_1=U_\Theta$, $A_2=gU_\Theta$ and $A_3=gg'U_\Theta$ for some $g,g'\in G$. Then the assumptions about the transversality can be expressed in terms of the Bruhat decomposition with respect to $P_\Theta$ as follows: $g\in P_\Theta w P_\Theta$, $g'\in P_\Theta w' P_\Theta$. By \Cref{thm:BruhatCellsMultiply} since $\WeylLength(ww') = \WeylLength(w) + \WeylLength(w')$ we have $gg' \in P_\Theta w P_\Theta w' P_\Theta = P_\Theta ww' P_\Theta$. Translating the cell back to transversality implies that $(A_1,A_3)$ is $ww'$-transverse as claimed.  
\end{proof}

\begin{proof}[Proof (of \Cref{thm:HDistPropertiesTheta}).]\leavevmode
    \begin{enumerate}
        \item By simple calculation we see if $g(A_1,A_2) = (U_\Theta, \overline{w}U_\Theta l)$ for some $g\in G,\,l\in L_\Theta$, then
            \begin{equation*}
                 l^{-1} \overline{w}^{-1} g(A_2,A_1) = \left( U_\Theta,  l^{-1}\overline{ w}^{-1} U_\Theta\right) = \left( U_\Theta,  \doverline{w^{-1}} U_\Theta\,w(l^{-1})\right)\,.
            \end{equation*}
            For $w=\varsigma_i$, we have
            \begin{equation*}
                \doverline{\varsigma_i^{-1}}=\doverline{\varsigma_i}=\doverline{\varsigma_i}^2\overline{\varsigma_i}=\check{\beta}_i(-1)\overline{\varsigma_i}\,.
            \end{equation*}
            Similarly, we have $\doverline{w_0}=s_G\overline{w_0}$\,.
        \item By definition, there exists $g\in G$ such that
            \begin{equation*}
                g(A_1l_1,A_2l_2)=\big(l_1U_\Theta,\overline{w}U_\Theta l(A_1,A_2) l_2\big)\,.
            \end{equation*}
        Thus,
            \begin{equation*}
                l_1^{-1}g(A_1l_1,A_2l_2)=\Big(U_\Theta,\overline{w}U_\Theta w^{-1}(l_1^{-1})l(A_1,A_2) l_2\Big)
            \end{equation*}
        which means $l(A_1l_1,A_2l_2)=w^{-1}(l_1)^{-1}l(A_1,A_2)l_2$ by definition.
        \item From \Cref{thm:WDistAdditiveTheta} we know $(A_1,A_3)$ are $ww'$-transverse. So there exists $g\in G$ such that
            \begin{equation*}
                g(A_1,A_3)=\Big(U_\Theta,\overline{ww'}U_\Theta l(A_1,A_3) \Big)\,.
            \end{equation*}
        Setting $A_2':=gA_2$, there exist $g_1,g_2\in G$ such that
            \begin{align*}
                &g_1\big(U_\Theta,A_2'\big) = \big(U_\Theta,\overline{w}U_\Theta l(A_1,A_2)\big)\\
                &g_2\big(A_2',\overline{ww'}U_\Theta l(A_1,A_3)) = \big(U_\Theta,\overline{w'}U_\Theta l(A_2,A_3)\big)\,.
            \end{align*}
        From the second equation, we deduce $A_2'=\overline{w}U_\Theta w'(l(A_1,A_3)l(A_2,A_3)^{-1})$, and thus the first equation implies
            \begin{equation*}
                l(A_1,A_2)=w'(l(A_1,A_3)l(A_2,A_3)^{-1})
            \end{equation*}
            which implies the lemma. 
    \end{enumerate}
\end{proof}

We now compute some explicit $L$-distances:
\begin{lemma}\label{thm:LdistSigma_i}
    The decorated flags $x_i(v)\overline{w_0}U_\Theta$ and $\overline{w_0}U_\Theta$ are $\varsigma_{\oppInv{i}}$-transverse if $v\in V_i^\times$. In this case the $L$-distance is given by:
    \begin{equation*}
        l\big(x_i(v)\overline{w_0}U_\Theta,\overline{w_0}U_\Theta\big) = \check{\beta}_{\oppInv{i}}\big(\iota(\oppInv{v})\big)\,.
    \end{equation*} 
\end{lemma}
\begin{proof}
    Let $w\in W(R)$ such that $w_0=\varsigma_i w$. Then $w_0=w\varsigma_{\oppInv{i}}$ We compute in $\Conf_2(\decFlagTheta)$ using \Cref{thm:DecomposeXSigma}:
    \begin{align*}
        \big(x_i(v)\overline{w}_0U_\Theta,\overline{w}_0U_\Theta\big)
        &= \big(x_i(v)\overline{\varsigma_iw}U_\Theta,\overline{w_0}U_\Theta\big)
        = \big(y_i(v^{-1})\check{\beta}_i(\iota(v))x_i(-v^{-1})\overline{w}U_\Theta,\overline{w_0}U_\Theta\big)\\
        &= \Big(x_i(-v^{-1})\overline{w}U_\Theta, \overline{w_0}U_\Theta\big(\check{\beta}_i(\iota(v))\big)^*\Big)
        = \Big(\overline{w}U_\Theta,\overline{w\varsigma_{\oppInv{i}}}U_\Theta\check{\beta}_{\oppInv{i}}\big(\iota(\oppInv{v})\big)\Big)\\
        &= \Big(U_\Theta,\overline{\varsigma_{\oppInv{i}}}U_\Theta\check{\beta}_{\oppInv{i}}\big(\iota(\oppInv{v})\big)\Big)\,.
    \end{align*}
    We have used that $\overline{w}^{-1}x_i(v')\overline{w}\in U_\Theta$ for any $v'\in V_i$, since $w(\beta_i)$ is a positive root.
\end{proof}

\begin{corollary}\label{thm:LdistanceW0}
    Let $D$ be a reduced expression of $w\in W(R)$ and $w_k = \varsigma_{i_1}\cdots \varsigma_{i_{k}}$ be the length $k$ prefixes of $D$. Let $\vec{v}=(v_1,\dots,v_m)$ with $v_j\in V_{i_j}^\times$ for $1\leq j\leq m=\WeylLength(w)$. Then: 
    \begin{align*}
        l\big(F_D(\vec{v})\overline{w_0}U_\Theta,\overline{w_0}U_\Theta\big) &=
        w_{m-1}^*(\check{\beta}_{\oppInv{i}_m}(\iota(\oppInv{v}_m))\cdot w_{m-2}^*(\check{\beta}_{\oppInv{i}_{m-1}}(\iota(\oppInv{v}_{m-1}))\cdots \check{\beta}_{\oppInv{i}_1}(\iota(\oppInv{v}_1))\\
        &= \prod_{j=1}^m w_{m-j}^* (\check{\beta}_{\oppInv{i}_{m-j+1}}(\iota(\oppInv{v}_{m-j+1}))),
    \end{align*}
    where $F^k_D$ is the prefix $R$-Lusztig map (see \Cref{def:LuszitigMap}).
\end{corollary}
\begin{proof}
    Let $F^k_D = \prod_{j=1}^{k}x_{i_j}(v_j)$ be a prefix $R$-Lusztig map  as in \Cref{def:LuszitigMap}.
    Since all $v_j$ are invertible, we have by \Cref{thm:LdistSigma_i} that the flags $(F^{k} \overline{w_0}U_\Theta, F^{k-1}\overline{w_0}U_\Theta)$ are $\varsigma_{\oppInv{i}_k}$-transverse with $L$-distance $l(F^k \overline{w_0}U\Theta, F^{k-1}\overline{w_0}U_\Theta) = \check{\beta}_{\oppInv{i}_k}(\iota(\oppInv{v}_k))$. Application of \Cref{thm:HDistPropertiesTheta}(3) concludes the proof.
\end{proof}
The previous statement allows for a simple computation for a very tricky decomposition if we combine it with the following observation:
\begin{proposition}
    Let $D = \varsigma_{i_1}\cdots \varsigma_{i_m}$ be a reduced expression in $W(R)$. For generic $\vec{v}=(v_1,\dots,v_m)$ with $v_j\in V_{i_j}^\times$, the element $F_D(\vec{v})\overline{\varsigma_{i_m}\cdots\varsigma_{i_1}}$ is Gauss decomposable with Levi part
    \begin{equation*}
        l\big(F_D(\vec{v})\overline{w_0}U_\Theta,\overline{w_0}U_\Theta\big)=[F_D(\vec{v})\overline{\varsigma_{i_m}\cdots\varsigma_{i_1}}]_0^*\,.
    \end{equation*}
\end{proposition}
\begin{proof}
    Let $w\in W(R)$ such that $w_0=\varsigma_{i_m}\cdots\varsigma_{i_1}w$, so that $w_0=w\varsigma_{\oppInv{i}_m}\cdots\varsigma_{\oppInv{i}_1}$. We Gauss decompose
    \begin{equation*}
        F_D(\vec{v})\overline{\varsigma_{i_m}\cdots\varsigma_{i_1}}=u_-lu_+
    \end{equation*}
    and compute in $\Conf_2(\decFlagTheta)$:
    \begin{align*}
        \big(F_D(\vec{v})\overline{w_0}U_\Theta,\overline{w_0}U_\Theta\big)
        &= \big(F_D(\vec{v})\overline{\varsigma_{i_m}\cdots\varsigma_{i_1}w}U_\Theta,\overline{w_0}U_\Theta\big)=\big(u_+\overline{w}U_\Theta,\overline{w_0}U_\Theta \oppInv{l}\big)\\
        &= \big(\overline{w}^{-1}u_+\overline{w}U_\Theta,\overline{w}^{-1}\overline{w_0}U_\Theta\oppInv{l}\big) = \big(U_\Theta,\overline{\varsigma_{\oppInv{i}_m}\cdots\varsigma_{\oppInv{i}_1}}U_\Theta\oppInv{l}\big)\,.
    \end{align*}
    The final equality requires a proof that $\overline{w}^{-1}u_+\overline{w}\in U_\Theta$, which follows from the root system calculus (\Cref{sec:RootSystemCalculus}). Explicitly, write
    \begin{align*}
        u_+=[F_D(\vec{v})\overline{\varsigma_{i_m}\cdots\varsigma_{i_1}}]_+
        &= x_{i_1}(v_1')\cdot\doverline{\varsigma_{i_1}}x_{i_2}(v_2')\overline{\varsigma_{i_1}}\cdots\doverline{\varsigma_{i_1}\cdots \varsigma_{i_{m-1}}}x_{i_m}(v_m')\overline{\varsigma_{i_{m-1}}\cdots \varsigma_{i_1}}\\
        &=\prod_{j=1}^m\doverline{\varsigma_{i_1}\cdots\varsigma_{i_{j-1}}}x_{i_j}(v_j')\overline{\varsigma_{i_{j-1}}\cdots\varsigma_{i_1}}
    \end{align*}
    for some $v_1'\in V_{i_1},\dots,v_m'\in V_{i_m}$. This is generically possible as $V_{i_j}^\times\subset V_{i_j}$ is dense. Since $\varsigma_{i_j}\cdots\varsigma_{i_1}w$ is reduced for any $j$, $w^{-1}\varsigma_{i_1}\cdots\varsigma_{i_{j-1}}(\beta_{i_j})$ is a positive root, which implies the claim.
\end{proof}
Now we can use \Cref{thm:LdistSigma_i} to compute the Levi part of elements of the form $u_w \overline{w^{-1}}$ where $u$ is an element of $U_\Theta$ given by a Lusztig type map for a reduced expression of $w$.
\begin{corollary}\label{thm:DecomposeUW}
    Let $D = \varsigma_{i_1}\cdots \varsigma_{i_m}$ be a reduced expression $w\in W(R)$ and $w_k$ the word given by the length $k$ prefix of $D$. Then 
    \begin{equation*}
        [F_D(\vec{v})\overline{\varsigma_{i_m}\cdots\varsigma_{i_1}}]_0^* = \prod_{j=1}^m w_{m-j}^* (\check{\beta}_{i_{m-j+1}}(\iota(\oppInv{v_{m-j+1}})))\,.
    \end{equation*}
\end{corollary}

Given a system of Jordan weights $\{\Lambda_i\}$, we can apply the weights to the $L$-distance to obtain coordinates, valued in $\Gamma(\jordan{J}_i)$. To lighten the notation we make the following definition with a small abuse of notation:
\begin{definition}
    Given a system of Jordan weights $\{\Lambda_i\} $ for $\Theta$, we also write $\Lambda_i$ for the $G$-invariant function on $\decFlagTheta^2$ defined by:
    \begin{equation*}
       \Lambda_i(A_1,A_2)=\Lambda_i(l(A_1,A_2))\,.
    \end{equation*}
    Also, for Gauss-decomposable $g\in G$ we define $\Lambda_i(g) = \Lambda_i(U_\Theta,\overline{w_0} g U_\Theta ) = \Lambda_i([g]_0)$.
\end{definition}

In the noncommutative rank 1 case, there is a unique nontrivial Jordan algebra associated to the $p^{th}$ root. This Jordan algebra has a pair of maps: the structure involution $\sigma$, and the adjugate $\tau$ (\Cref{def:structureInvolution,def:adjugateJordanAlgebra}).
\begin{corollary}\label{thm:ThetaWeightsOfReversedFlags}
    For $i< p$, $\Lambda_i(s_G A_2,A_1) = \Lambda_i(A_1,A_2)$. Moreover when $i=p$, we have
    \begin{equation*}
        \Lambda_p(s_G A_2, A_1) = \begin{cases}
            \sigma(\Lambda_p(A_1,A_2)) & \text{Jordan split type $A_1$}\\
            \tau(\Lambda_p(A_1,A_2)) & \text{Jordan split type $B_p$ with $p$ even}\\
            \sigma(\Lambda_p(A_1,A_2)) & \text{Jordan split type $B_p$ with $p$ odd}\\
            \sigma(\Lambda_p(A_1,A_2)) & \text{Jordan split type $G_2$}
        \end{cases}
    \end{equation*}
\end{corollary}
\begin{proof}
    By \Cref{thm:HDistPropertiesTheta}(1) this follows from computing the action of $w_0$ on each pinning map $\check{\beta}_i$. For $i< p$ this is independent of Jordan split type and  $w_0(\check{\beta}_i(t)) = \check{\beta}_i(t^{-1})$.
    For the special Jordan weight, we must compute each case separately. In type $A_1$, $w_0 = \varsigma_p$ and the action follows by \Cref{thm:ActionOfSigmaPonSpecialRoot}. Type $B_p$ is computed in \Cref{thm:W0ActionBp} and type $G_2$ in \Cref{thm:w0ActionG2}.
\end{proof}

\section{From \texorpdfstring{$\Confx_3(\decFlagTheta)$}{Conf3x(A)} to elementary configurations}\label{sec:triplesToElementaryConfigs}
In \Cref{sec:coordinateStrategy} we outlined the general strategy of constructing cluster-like coordinates on $\Confx_n(\decFlagTheta)$. We will now explain in detail how to decompose a configuration $[(A_1,A_2,A_3)]\in\Confx_3(\decFlagTheta)$ into configurations associated to simple reflections $\varsigma_i\in W(R)$.

\subsection{Interpolating sequence}
Choose a reduced expression of the longest word, $D(w_0) = \varsigma_{i_1}\cdots \varsigma_{i_n}$ and a system of Jordan weights. First, we use these choices to compute an \keyword{interpolating sequence} of decorated flags between $A_2$ and $A_3$.

\begin{proposition}\label{thm:flagSequenceTheta}
	Assume that $(A_2,A_3)\in\decFlagTheta^2$ is a generic pair of decorated flags, and that we have a reduced expression $D(w_0)=\varsigma_{i_1}\cdots \varsigma_{i_n}$ with $i_k\in\{1,\dots,p\}$ of the longest word in $W(R)$. Let $w_k = \varsigma_{i_1}\cdots \varsigma_{i_{k}}$. Then there exists a unique chain
	\begin{equation*}
		A_2=A_2^0,A_2^1,\dots,A_2^n=A_3
	\end{equation*}
	of decorated flags such that
	\begin{enumerate}[(i)]
	   \item $(A_2^k,A_2^{k-1})$ are $\varsigma_{i_k}$-transverse,
		\item If $\check{\gamma}_k:=w_{k-1}(\check{\beta}_{i_k})$ is simple (equal to $\check{\beta}_{j}$ for some $j$), then
            \begin{equation*}
                l\left(A_2^k,A_2^{k-1}\right) = \check{\beta}_{i_k}\left(\Lambda_{j}\left(l\left(A_3,A_2\right)\right)\right) \in \check{\beta}_{i_k}(\Gamma(\jordan{J}_{i_k}))
            \end{equation*}
        otherwise $l\left(A_2^k,A_2^{k-1}\right)=1$.
	\end{enumerate}
    This chain is called the \keyword{interpolating sequence} between $A_2$ and $A_3$ with respect to $D(w_0)$. 
\end{proposition}

\begin{proof}
    First we note that there is nothing to prove in the case when $p=1$ ($\Delta_\Theta$ is $A_1$).
    Since $(A_2,A_3)$ is generic, there exists $g\in G$ such that
    \begin{equation*}
        g(B_2,B_3)=\left(P_\Theta,\doverline{w_0}P_\Theta\right)\in\flagTheta^2
    \end{equation*}
    where $B_i:=\pi(A_i)$. Define
    \begin{equation*}
        B_2^k:=g^{-1}\doverline{\varsigma_{i_1}}\cdots \doverline{\varsigma_{i_k}}P_\Theta\,.
    \end{equation*}
    Then $B_2^0=B_2, B_2^n=B_3$, and $\left(B_2^k,B_2^{k-1}\right)$ is $\varsigma_{i_k}$-transverse. Condition (ii) along with our choice of Jordan weights uniquely determines the points $A_2^k$ in the fiber of $\pi$ over $B_2^k$: we have $l(A_3,A_2) =\prod_i\check{\beta}_i(\Lambda_i(l(A_3,A_2)))$ so we choose $A_2^k$ so that 
            \begin{equation*}
                l\left(A_2^k,A_2^{k-1}\right) = \check{\beta}_{i_k}(\Lambda_j(l(A_3,A_2)))
            \end{equation*}
        or $l\left(A_2^k,A_2^{k-1}\right)=1$ appropriately.

    This proves existence of an interpolating sequence. For uniqueness, we need to prove that the sequence $\left\{B_2^k\right\}$ is unique. This is done inductively. It suffices to prove the following: Suppose $B',B\in\flagTheta$ are $w\varsigma_i$-transverse where $w\varsigma_i$ is a reduced word in $W(R)$ for some $w\in W(R)$ and $i\in\{1,\dots,p\}$. Then the flag $B''\in\flagTheta$, that is $w$-transverse from $B'$ and $\varsigma_i$-transverse to $B$, is unique.

    Now, suppose that $B=gP_\Theta,\, B'=g'P_\Theta$, and that $B_1''=g_1P_\Theta$ and $B_2''=g_2P_\Theta$ are two flags which fulfill the condition above. This means that $g_1^{-1}g,\, g_2^{-1}g\in P_\Theta \varsigma_i P_\Theta$. Since $\varsigma_i^{-1}=\varsigma_i$, we have $g^{-1}g_2\in P_\Theta \varsigma_i P_\Theta$ and thus
    \begin{align*}
        g_1^{-1}g_2\in P_\Theta \varsigma_i P_\Theta \varsigma_i P_\Theta.
    \end{align*}
    $P_\Theta \varsigma_i P_\Theta \varsigma_i P_\Theta$ decomposes as a union of cells $\coprod_j P_\Theta \varsigma_i^{(j)} P_\Theta$ by \Cref{prop:double_cosets_A1}. We wish to show that $g_1^{-1}g_2$ is actually in the identity cell. Assume for contradiction that this is not the case. 
    If $g_1^{-1}g_2\in P_\Theta \varsigma_i^{(j)} P_\Theta$, $j\neq 0$, then
    \begin{align*}
        (g')^{-1}g_2=\left((g')^{-1}g_1\right)\cdot\left(g_1^{-1}g_2\right)\in P_\Theta w P_\Theta \varsigma_i^{(j)} P_\Theta = P_\Theta w\varsigma_i^{(j)} P_\Theta
    \end{align*}
     by \Cref{thm:bruhat_cells_mulitply_partial}. This contradicts the assumption that the $w$-distance between $B_2''$ and $B'$ is $w$.
\end{proof}
\begin{remark}\label{rem:LdistanceFromElementaryConfig}
    The $L$-distance between $A_3$ and $A_2$ is easily computed from the interpolating sequence by 
    \begin{equation*}
        l(A_3,A_2) = w_{k_1-1}(l(A_2^{k_1},A_2^{k_1-1}))\cdots w_{k_p-1}(l(A_2^{k_p},A_2^{k_p-1}))
    \end{equation*}
    where $\check{\gamma}_k:=w_{k-1}(\check{\beta}_{i_k})$ is equal to the simple coroot $\check{\beta}_j$.
\end{remark}

\begin{definition}\label{def:rescaling_action}
    We define an action of $(L_\Theta)^3$ on $\Conf_3(\decFlagTheta)$ by $(l_1,l_2,l_3)\cdot(A_1,A_2,A_3) = (A_1\oppInv{l}_1,A_2\oppInv{l}_2,A_3\oppInv{l}_3)$ called the \keyword{rescaling action}.
\end{definition}
The coinvariants of the rescaling action are configurations of undecorated flags.
\begin{lemma}\label{lem:rescaling_action_interpolating}
    The interpolating flags $B_2^k$ from the above construction are invariant under the rescaling action.
\end{lemma}
\begin{proof}
    In the construction $B_2^k$ only depends on the undecorated flags $B_2,B_3$, which are independent of the rescaling action.  
\end{proof}

\subsection{Elementary configurations and their partial potential}
For any reduced expression of $w_0\in W(R)$, the configurations represented by $(A_1,A_2^{k-1},A_2^k)$ are generically elementary configurations for $\varsigma_{i_k}$ which we define now.
\begin{definition}\label{def:ElementaryConfiguration}
    An \keyword{elementary configuration} of decorated flags for $\varsigma_i\in W(R)$ is an element of $\Conf_3(\decFlagTheta)$ which can be represented as $(A,A_l,A_r)\in\decFlagTheta^3$ such that:
    \begin{itemize}
        \item $(A,A_l)$ and $(A,A_r)$ are $w_0$-transverse,
        \item $(A_r,A_l)$ is $\varsigma_{i}$-transverse,
        \item $l(A_r,A_l) \in \check{\beta}_{i}(\Gamma(\jordan{J}_{i}))$\,.
    \end{itemize}
    We denote the space of elementary configurations for $\varsigma_i$ by $\Conf_3^i(\decFlagTheta)$.
\end{definition}

\input{Figures/figureElementaryConfigurations}

The space of elementary configurations admits a natural parametrization.
\begin{lemma}\label{thm:elementaryConfigPotentialTheta}
    Let $(A,A_l,A_r)\in\decFlagTheta^3$ represent an element of $\Conf_3^i(\decFlagTheta)$. Then there exist unique $g\in G$, $l_l,l_r\in L_\Theta$ and $a\in V_{i}^\times$ such that
    \begin{equation*}
        g(A,A_l,A_r)=\left(U_\Theta,\overline{w_0}U_\Theta l_l,x_{\oppInv{i}}(\oppInv{a})\overline{w_0}U_\Theta l_r\right)\,.
    \end{equation*}
\end{lemma}
\begin{definition}\label{def:partialPotentialTheta}
    The element $a\in V_{i}^\times$ only depends on $(A,\pi(A_l),\pi(A_r))\in\decFlagTheta\times\flagTheta^2$ and is called the \keyword{partial potential} of the elementary configuration, which is denoted by $\pot_i(A,A_l,A_r)$.
\end{definition}
\begin{proof}[Proof (of \Cref{thm:elementaryConfigPotentialTheta}).]
     Since $(A,A_l)$ are transverse we have unique $g\in G$ and $l_l\in L_\Theta$ such that $g(A,A_l)=(U_\Theta,\overline{w_0}U_\Theta l_l)$. Because we also have $(A,A_r)$ are transverse, $gA_r=u\overline{w_0}U_\Theta l_r$ for some $u\in U_\Theta,\, l_r\in L_\Theta$. But the $\varsigma_{i}$-transversality of $(A_r,A_l)$ implies $u\in U_\Theta\cap P_\Theta^\opp\varsigma_{i^*} P_\Theta^\opp$. This intersection is exactly $x_{i^{*}}(V_{i^*}^\times)$ by \Cref{thm:BruhatSigmaIntersection}. Thus there is an $a \in V_{i}^\times$ such that $u = x_{i^*}(\oppInv{a})$ as needed. Moreover changing $l_l$ and $l_r$ leaves $u$ fixed and so $a$ only depends on $\pi(A_l)$ and $\pi(A_r)$.
\end{proof}

\subsubsection{Computing partial potentials}
Let $(A,A_l,A_r)\in\decFlagTheta^3$ represent an elementary configuration in $\Conf_3^i(\decFlagTheta)$. The $L_\Theta$-distances $l_l= l(A,A_l),\,l_r = l(A,A_r)$ together with the partial potential induce an embedding
\begin{equation*}
    \Conf_3^i(\decFlagTheta)\hookrightarrow L_\Theta^2\times V_{i}^\times\,.
\end{equation*}
There is a relation which determines the image of the embedding which is given in \Cref{thm:elementaryConfigLeftRightFunctionsTheta} for lower noncommutative rank and in \Cref{rem:ElementaryNetworkEmbedding} for higher noncommutative rank.\medskip

Set $l_b:=l(A_r,A_l)\in\check{\beta}_i(\Gamma(\jordan{J}_i))$. We can express $l_b$ in terms of $l_l,\,l_r$, the partial potential $a:=\delta_i(A,A_l,A_r)\in V_{i}$.

\begin{proposition}\label{thm:elementaryConfigBottomLDistTheta}
    For $l_l,l_r,l_b$ and $a$ as above we have:
    \begin{equation*}
        l_b=\varsigma_{i}(l_r^{-1})\check{\beta}_{i}\big(\iota(a)\big)l_l\,. 
    \end{equation*}
\end{proposition}
\begin{proof}
    We see that $l_b  = \varsigma_{i}(l_r^{-1}) l(x_{\oppInv{i}}(\oppInv{a})\overline{w}_0U_\Theta,\overline{w}_0U_\Theta) l_l $ by \Cref{thm:HDistPropertiesTheta}(2). Then \Cref{thm:LdistSigma_i} tells us that the inner $L$-distance is given by $\check{\beta}_{i}(\iota(a))$ as needed. 
\end{proof}

We can also understand the structure of the rescaling action on the potentials:

\begin{proposition}\label{prop:rescaling_potentials}
    The rescaling action restricted to the top flag $A$ in an elementary configuration transforms the potential by $\delta_i(A\ell,A_l,A_r) = \beta_i(\varsigma_i(\ell))(\delta_i(A,A_l,A_r))$
\end{proposition}
\begin{proof}
The potential does not depend on on the decorations of $A_2,A_3$ so 
    $$ (A\oppInv{\ell},B_l,B_r) = \left(U_\Theta\oppInv{\ell},\overline{w_0}P_\Theta,x_{\oppInv{i}}(\oppInv{a})\overline{w_0}P_\Theta\right)=\left(U_\Theta,\overline{w_0}P_\Theta,{\oppInv{\ell}}^{-1}x_{\oppInv{i}}(\oppInv{a})\oppInv{\ell} \,\overline{w_0}P_\Theta\right) $$ 
    Now noting that $\doverline{w}_0 {\oppInv{\ell}}^{-1}x_{\oppInv{i}}(\oppInv{a})\oppInv{\ell}\overline{w}_0 = \ell  y_i(-a)\ell^{-1} = y_i(\beta_i(\varsigma_i(\ell)(-a))$ computes the action on the potential $a$. 
\end{proof}

\subsubsection{Opposite elementary configurations}

We now consider another parameterization of elementary configurations. These will be useful when computing \textit{flips} in \Cref{sec:FlipTriangulationGenericTheta} and actions of the cluster modular group in \Cref{sec:ClusterModularGroupAction}. 
\begin{definition}\label{def:OppositeElementaryConfiguration}
    The set of \keyword{opposite elementary configurations}, $\Conf_3^{\bar{i}}(\decFlagTheta)$, is the set of triples $(A^l,A^r,A)\in \decFlagTheta^3$ for which 
    \begin{itemize}
        \item $(A^l,A)$ and $(A^r,A)$ are $w_0$-transverse,
        \item $(A^l,A^r)$ is $\varsigma_{\oppInv{i}}$-transverse,
        \item $l(A^l,A^r) \in \check{\beta}_{\oppInv{i}}(\Gamma(\jordan{J}_{\oppInv{i}}))$\,.
    \end{itemize}
\end{definition}
\begin{lemma}
    There is an isomorphism
    \begin{equation*}
        \mathrm{flip}: \Conf_3^{i}(\decFlagTheta) \to \Conf_3^{\bar{\oppInv{i}}}(\decFlagTheta) \quad (A,A_l,A_r) \mapsto (A^l:=A_r,A^r:=A_l,s_GA)\,.
        \end{equation*}
\end{lemma}
We can write opposite elementary configurations in a convenient standard from: 
\begin{equation*}
    (A^l,A^r,A) = (U_\Theta (l^l)^*, y_{\oppInv{i}}(\oppInv{b})U_\Theta (l^r)^*, \overline{w_0}U_\Theta)
\end{equation*}
with $l^l = l(A^l,A),l^r = l(A^r,A) \in L_\Theta$ and $b\in V_i$.

\begin{definition}
    The element $b =: \delta_{\bar{i}}(A^l,A^r,A) \in V_i$ is the \keyword{partial potential of} $(A^l,A^r,A)$.
\end{definition}

The following lemma, proves this potential is same as the potential for an ordinary elementary configuration. 
\begin{lemma}
    Let $(s_GA,A_l,A_r) \in \Conf_3^{\oppInv{i}}(\decFlagTheta)$ be the flip of $(A^l,A^r,A)\in \Conf_3^{\bar{i}}(\decFlagTheta)$. Then $$\oppInv{\delta_{\bar{i}}(A^l,A^r,A)} = \delta_{\oppInv{i}}(s_GA,A_l,A_r).$$
\end{lemma}
\begin{proof}
    This follows by direct computation. Let $b = \delta_{\bar{i}}(A^l,A^r,A)$,  $l^l = l(A^l,A)$ and $l^r = l(A^r,A)$.
    \begin{align*}
        (s_G A, A_l, A_r) = (s_G A, A^r, A_l) \sim (s_g \overline{w}_0 U_\Theta, y_{\oppInv{i}}(\oppInv{b})U_\Theta (l^r)^*, U_\Theta (l^l)^*) \,.
    \end{align*}
    By \Cref{thm:HDistPropertiesTheta}, $\oppInv{(l^r)} = l(s_GA,A^r) = l(s_GA,A_l) =: l_l$ and $\oppInv{(l^l)} = l(s_G A,A_r)=: l_r$. Applying $\overline{w}_0$ to the configuration results in  
    \begin{align*}
    (s_G A, A_l, A_r) \sim ( U_\Theta, \overline{w}_0y_{\oppInv{i}}(\oppInv{b})\doverline{w}_0\overline{w}_0U_\Theta l_l, \overline{w}_0U_\Theta l_r)\,.
    \end{align*}
    By \Cref{thm:w0conjugationPinning}, $\overline{w}_0y_{\oppInv{i}}(\oppInv{b})\doverline{w}_0 = x_{i}(-b)$. Thus
     \begin{align*}
    (s_G A, A_l, A_r) \sim ( U_\Theta, x_i(-b)\overline{w}_0U_\Theta l_l, \overline{w}_0U_\Theta l_r) \sim ( U_\Theta, \overline{w}_0U_\Theta l_l, x_i(b)\overline{w}_0U_\Theta l_r)\,.
    \end{align*}
    From this standard form we see $\delta_{\oppInv{i}}(s_G A, A_l,A_r) = b$ as claimed.
\end{proof}

\begin{corollary}\label{thm:elementaryConfigBottomLDistTheta_opposite}
    For $l^l=l(A^l,A) ,l^r=l(A^r,A),l^t = l(A^l,A^r)$ and $b=\delta_{\bar{i}}(A^l,A^r,A)\in V_i$ as above we have:
    \begin{equation*}
        l^t=\oppInv{\left(l^r\check{\beta}_{i}\big(\iota(b)\big)\varsigma_{i}(l^l)^{-1}\right)}\,. 
    \end{equation*}
    Thus the potential of an opposite elementary configuration can be computed as follows,
    \begin{equation*}
        \check{\beta}_{i}\big(\iota(b)\big) = (l^r)^{-1}(l^t)^*\varsigma_{i}(l^l)\,.
    \end{equation*}
\end{corollary}
\begin{proof}
    As before take $(A^l,A^r,A)$ to be the flip of the elementary configuration $(s_GA,A_l,A_r)\in \Conf^{\oppInv{i}}_3$. By \Cref{thm:elementaryConfigBottomLDistTheta}, $l_b= l(A_r,A_l) = \varsigma_{\oppInv{i}}(l_r^{-1})\check{\beta}_{\oppInv{i}}\big(\iota(\oppInv{b})\big)l_l$. Then since  $l^r = \oppInv{l_l}$ and $l^l=\oppInv{l_r}$ we have, 
    \begin{equation*}
        l^t=l_b= \varsigma_{\oppInv{i}}(l_r^{-1})\check{\beta}_{\oppInv{i}}\big(\iota(\oppInv{b})\big)l_l =
        \oppInv{\left(l^r\check{\beta}_{i}\big(\iota(b)\big)\varsigma_{i}(l^l)^{-1}\right)} \,.
    \end{equation*}
\end{proof}

\subsubsection{Reflection maps}

\begin{definition}
    The \keyword{left reflection map} $r:\Conf_3^{i}(\decFlagTheta) \to \Conf_3^{\bar{i}}(\decFlagTheta)$ is given by 
\begin{equation*}
    r(A,A_l,A_r) = (A',A,A_r)
\end{equation*}
where $A'$ is the unique decorated flag so that 
\begin{equation*}
    (A',A) \text{ are $\varsigma_{\oppInv{i}}$-transverse},\quad (A',A_l) \text{ are $\varsigma_{\oppInv{i}}w_0$-transverse}, \quad l(A',A) = \oppInv{l(A_r,A_l)}\,.
\end{equation*}
    There is an analogous \keyword{right reflection map} which creates an opposite configuration by replacing $A_r$ with a new flag $A'$.
\end{definition}

The left reflection map is visualized in \Cref{fig:reflectionMap}.

\input{Figures/figureReflectionMap}

\begin{proposition} The standard form for the left reflection is 
    $$r(U_\Theta,\overline{w}_0U_\Theta l_l,x_{\oppInv{i}}(\oppInv{a})\overline{w}_0U_\Theta l_r) 
    =(U_\Theta \oppInv{(l^l)}, y_{\oppInv{i}}(\oppInv{b})U_\Theta \oppInv{(l^r)}, \overline{w_0}U_\Theta)$$
    where $b= \beta_i(l_r^{-1})(a^{-1})$,  $l^l=\varsigma_{i}(l_b^{-1}l_ll_b^{-1})$ and $l^r=l_r$.
\end{proposition}
\begin{proof}
    First we compute that with  $l' = \varsigma_{\oppInv{i}}(\oppInv{l_b})^{-1}$, 
    $$r(U_\Theta,\overline{w}_0U_\Theta l_l,x_{\oppInv{i}}(\oppInv{a})\overline{w}_0U_\Theta l_r) 
    =(\doverline{\varsigma}_{\oppInv{i}}U_\Theta l',U_\Theta ,x_{\oppInv{i}}(\oppInv{a})\overline{w}_0U_\Theta l_r)\,.$$
    Next, we use \Cref{thm:DecomposeXSigma} to decompose $x_{\oppInv{i}}(-\oppInv{a})\doverline{\varsigma_i}=y_{\oppInv{i}}(-(\oppInv{a})^{-1})\check{\beta}_{\oppInv{i}}(\iota(\oppInv{a}))x_{\oppInv{i}}((\oppInv{a})^{-1})$
    so that this configuration is the same as 
    $$
    (U_\Theta \check{\beta}_{\oppInv{i}}(\oppInv{a})l',y_{\oppInv{i}}((\oppInv{a})^{-1})U_\Theta ,\overline{w}_0U_\Theta{l}_r) $$
    which is the same as 
    $$ (U_\Theta \oppInv{l}_r \check{\beta}_{\oppInv{i}}(\oppInv{a})l',\,\oppInv{l}_ry_{\oppInv{i}}((\oppInv{a})^{-1})){\oppInv{l}_r}^{-1}U_\Theta \oppInv{l}_r,\,\overline{w}_0U_\Theta).$$
    By applying the opposition involution to 
    $\oppInv{l}_ry_{\oppInv{i}}((\oppInv{a})^{-1})){\oppInv{l}_r}^{-1}$, we can compute 
    $$ \left(\oppInv{l}_ry_{\oppInv{i}}((\oppInv{a})^{-1})){\oppInv{l}_r}^{-1}\right)^* = l_r^{-1}x_i(a^{-1})l_r = x_i(\beta_i(l_r^{-1})(a^{-1})).$$
    Therefore we have $b= \beta_i(l_r^{-1})(a^{-1})$ and $l^r=l_r$. By rearranging \Cref{thm:elementaryConfigBottomLDistTheta} we have 
    $$\check{\beta}_{i}\big(\iota(a)\big) = \varsigma_{i}(l_r)l_bl_l^{-1}.$$
    Since the structure involution fixes the image of $\iota$, we have
    $$\check{\beta}_{i}\big(\iota(a)\big) = \varsigma_i\big(\check{\beta}_{i}(\iota(a))\big) =\varsigma_{i}(l_ll_b^{-1})l_r^{-1} \,. $$
    Now we compute 
    \begin{equation*}
        \oppInv{(l')} \check{\beta}_{i}(\iota(a))l_r =  \varsigma_{i}(l_b)^{-1}\check{\beta}_{i}(\iota(a))l_r =  \varsigma_{i}(l_b)^{-1}\varsigma_{i}(l_ll_b^{-1})l_r^{-1}l_r = \varsigma_{i}(l_b^{-1}l_ll_b^{-1})\,, 
    \end{equation*}
    which finally implies $l^l=\varsigma_{i}(l_b^{-1}l_ll_b^{-1})$.
\end{proof}

\begin{lemma}\label{lem:reflected_potential}
    The reflected potential $b=\beta_i(l_r^{-1})(a^{-1})$ is given in terms of $L$-distances in the original configuration by
    \begin{equation*}
        l_r^{-1}l_ll_b^{-1} = \check{\beta}_i(\iota(b))\,.
    \end{equation*}
\end{lemma}
\begin{proof}
    Starting with  $l_b=\varsigma_{i}(l_r^{-1})\check{\beta}_{i}\big(\iota(a)\big)l_l$ we rearrange to find 
    $$ l_r^{-1}l_ll_b^{-1} = l_r^{-1}\check{\beta}_{i}\big(\iota(a^{-1})\big)\varsigma_{i}(l_r) = \check{\beta}_{i}\big(\iota(\beta_i(l_r^{-1})(a^{-1}))\big)= \check{\beta}_{i}(b).$$
\end{proof}

\section{Cluster Coordinates on elementary configurations of \texorpdfstring{$\Theta$-}{Theta }flags}\label{sec:CoordsOnElementaryCofigurations}
Now that we understand how to reduce a configuration of three decorated flags to a series of elementary configurations, we can start the process of constructing coordinates by first endowing the elementary configuration spaces with appropriate coordinates which are suitable for amalgamation in the next step.

\subsection{Partial potentials in noncommutative rank 1}
We now \keyword{specialize to noncommutative rank 1}, i.e. the Jordan split type is either $A_1,\,B_p$ or $G_2$. Observe that in these cases the opposition involution on the (co)root system is trivial, while it may be nontrivial on $L_\Theta$. Recall that we index the Dynkin diagram such that the noncommutative root has index $p$. We denote the corresponding Jordan algebra $\jordan{J}_p=:\jordan{J}$. All roots $\beta_i$ with $i<p$ are commutative.
\begin{corollary}\label{thm:potentialFormulaTheta}
    The partial potential $a = \delta(A,A_l,A_r)$ for a $\varsigma_p$-elementary configuration fulfills
    \begin{align*}
        \sigma\big(\Lambda_p(l_r)\big)^{-1}\Lambda_p(l_b)\Lambda_p(l_l)^{-1}\cdot \Lambda_{p-1}(l_r)=\iota(a)
    \end{align*}
    with $l_l = l(A,A_l) ,l_r=l(A,A_r),l_b=l(A_r,A_l)$.
\end{corollary}
\begin{proof}
    Rearrange the equation from \Cref{thm:elementaryConfigBottomLDistTheta} and apply $\Lambda_p$ to both sides. By definition $\Lambda_p$ only sees the image of $\check{\beta}_p$. So by \Cref{thm:ActionOfSigmaPonOtherRoots,thm:ActionOfSigmaPonSpecialRoot} we see that $\Lambda_p(\varsigma_p(l_r)) = \sigma(\Lambda_p(l_r))^{-1}\Lambda_{p-1}(l_r)$.
\end{proof}

\begin{proposition}\label{thm:potentialFormulaSplit}
    Let $(A,A_l,A_r)$ represent an element in $\Conf_3^{i}(\decFlagTheta)$ for $i < p$. The partial potential is computed by the following formula:
    \begin{align*}
        \pot_{i}(A,A_l,A_r)=~&\frac{\Lambda_{i}(l_b)}{\Lambda_{i}(l_l)\Lambda_{i}(l_r)}\cdot \Lambda_{i-1}(l_l)\Lambda_{i+1}(l_l)\\
        \pot_{p-1}(A,A_l,A_r)=~&\frac{\Lambda_{p-1}(l_b)}{\Lambda_{p-1}(l_l)\Lambda_{p-1}(l_r)}\cdot \Lambda_{p-2}(l_l) N(\Lambda_p(l_l)) \,.
    \end{align*}
\end{proposition}
\begin{proof}
This is \cite[Lemma 11.12]{goncharov2019quantum}. It also follows by analogous computation to \Cref{thm:potentialFormulaTheta}. 
\end{proof}

\begin{corollary}\label{thm:elementaryConfigLeftRightFunctionsTheta}
    Let $(A,A_l,A_r)$ represent a configuration in $\Conf_3^p(\decFlagTheta)$. Then for any $i<p$
    \begin{equation*}
        \Lambda_i(A,A_l)=\Lambda_i(A,A_r)\,.
    \end{equation*}
\end{corollary}

\begin{proof}
    This is essentially the same as \cite[Lemma 11.13]{goncharov2019quantum}. 
    We know that $l_b\in\check{\beta}_p(\Gamma(\jordan{J}))$ by definition of the elementary configuration. Applying $\Lambda_j$ for $j\neq i$ to the equation for $l_b$, we find $1=\Lambda_j(\varsigma_i(l_l))^{-1}\Lambda_j(l_r)$ and thus
    \begin{equation*}
        \Lambda_j(l_r)=\Lambda_j(\varsigma_i(l_l))=\Lambda_j(l_l)
    \end{equation*}
    where the last equality holds because $j\neq i$.
\end{proof}

\subsection{Elementary seeds in lower noncommutative rank}\label{sec:elementarySeedsLowerNCrank}
To each elementary configuration we associate an \keyword{elementary seed} of a polygonal cluster algebra. The cluster coordinates will be given by Jordan weights applied to $l(A,A_l)$, $l(A,A_r)$, $l(A_r,A_l)$.

The quivers of the elementary seed $\elemSeed_\circ(i)$ for $1\leq i\leq p$ will be based on the quivers defined in \cite{goncharov2019quantum}: We think of the plane as having different levels $1,\dots,p$. $\elemSeed_\circ(i)$ has one node labeled $j$ on every level $j\neq i$, and two nodes labeled $i_l,i_r$ on level $i$. These are connected by an arrow from $i_l$ to $i_r$, and there are half-arrows from $i_r$ to $i+1$ and $i-1$, and from $i-1$ and $i+1$ to $i_l$. The nodes on level $p$ have weight $r$, on all other levels they have weight 1. The elementary quiver $\elemSeed(i)$ is obtained from $\elemSeed_\circ (i)$ by adding an additional frozen node labeled $i_b$ of the same weight as the nodes on level $i$ attached in an oriented cycle of full edges. See \Cref{fig:ElementaryQuivers} for examples.

\input{Figures/figureElementaryQuiver.tex}

The cluster coordinates are given as follows:
\begin{equation}
    \begin{aligned}
        a_{i_l} =~& \Lambda_i(A,A_l) \hspace{5pc} a_{i_r} = \Lambda_i(A,A_r) \hspace{5pc }a_{i_b} = \Lambda_i(A_r, A_l) \\
        a_{j} =~& \Lambda_j(A,A_l) = \Lambda_j(A,A_r) \hspace{1pc} \text{for } j \neq i\\   
    \end{aligned}
\end{equation}

For $i < p$, the only noncommutative node can be weaved or switched without affecting the mutation relations at any other node. In fact, the potential can be recovered by the same formula as in \cite{goncharov2019quantum} when we remember that small nodes only ``see'' $N(a_{p})$ in this case.

The quiver for $i=p$ requires additional decoration to specify a quiver in a polygonal cluster algebra. In \Cref{fig:elementaryQuiverTheta} we see this both as an ordered quiver and as a polygonal tile representing the pruned quiver. These decorations are chosen so that the potential is recovered (up to commutative multiplication) by the angle based near the top vertex.

\begin{lemma}\label{thm:partialPotentialFromCoordsTheta}
     The potential $\pot_p(A,A_l,A_r)$ is computed as
    \begin{equation}
        \iota(\pot_p(A,A_l,A_r)) = \Delta \frac{1}{N(a_{p_l})N(a_{p_r})} a_{p-1}
    \end{equation}
    where the angle $\Delta$ is defined by the angle path going around the polygon which starts near the end of the angle indicator, as shown in green in \Cref{fig:elementaryQuiverTheta}.
\end{lemma}
\begin{proof}
    This is a simple computation. According to \Cref{def:anglePCA} the angle is given by:
    \begin{equation*}
          \Delta= \sau(a_{p_r}) (a_{p_b}) \tau(a_{p_l})
    \end{equation*}
    with $a_{p_b}=1$ for the quiver $\elemSeed_\circ(p)$. Scaling by the norms converts $\tau(a_{p_l})$ and $\sau(a_{p_r})$ into $a_{p_l}^{-1}$ and $\sigma(a_{p_r}^{-1})$ as needed to match the formula in \Cref{thm:potentialFormulaTheta}. 
\end{proof}

\input{Figures/figureElementaryQuiverTheta}

\subsection{Elementary networks in higher noncommutative rank}\label{sec:elementaryNetworks}
Now let $G\simeq \GL_{p+1}(\mathcal{R})$ for $\mathcal{R}$ a finite dimensional algebra over $\K$. Recall that in this case we have $L_\Theta\simeq \GL_1(\mathcal{R})^{p+1}$.

Let $(A,A_l,A_r)$ represent a configuration in $\Conf_3^i(\decFlagTheta)$.  This determines three $L$-distances: 
\begin{align*}
    l_l &= l(A,A_l) = (L_1,\cdots,L_{p+1}), \\\
    l_r &= l(A,A_r) = (R_1,\cdots,R_{p+1}),\\
     l_b &= l(A_r,A_l) = (1,\cdots,B_i,B_{i+1},\cdots,1).
\end{align*}

We construct an elementary network associated with this configuration modeled on the plabic graphs associated to wiring diagrams \cite{fwz_IntroToClusters_PlabicGraphs}. The outer edges of the network $\elemNetwork(i)$ are assigned elements of $\GL_1(\mathcal{R})$ from the three $L$-distances as in \Cref{fig:ElementaryNetworksWeighted}.  All the internal edge weights are determined by this assignment and the rule that the cyclic product around each vertex is 1.

\input{Figures/figureElementaryNetworks}

\begin{remark}\label{rem:ElementaryNetworkEmbedding}
    The vertex relations for $\elemNetwork(i)$ immediately imply that for $j \notin \{i,i+1\}$, $R_{j} = L_{j}$. The final vertex relation around the central $\bullet$ vertex has the form $ B_{i+1}L_{i+1}^{-1}\cdot R_{i+1}B_{i} \cdot L_i^{-1}R_i = 1$.
\end{remark}
This is equivalent to \Cref{thm:elementaryConfigBottomLDistTheta} which relates the partial potential $D =\delta_i(A,A_l,A_r) \in V_{\oppInv{i}}$ to the $L$-distances. In this case, \Cref{thm:elementaryConfigBottomLDistTheta} takes the form
\begin{equation*}
    \begin{bmatrix} \vdots\\ R_{i+1}^{-1}\\ R_{i}^{-1}\\ \vdots \end{bmatrix} \begin{bmatrix} \vdots\\D\\D^{-1}\\ \vdots \end{bmatrix} \begin{bmatrix} \vdots\\ L_i\\ L_{i+1}\\ \vdots \end{bmatrix} = \begin{bmatrix} \vdots\\ R_{i+1}^{-1}DL_{i}\\ R_{i}^{-1}D^{-1}L_{i+1}\\ \vdots \end{bmatrix} = \begin{bmatrix} \vdots\\ B_i\\ B_{i+1}\\ \vdots \end{bmatrix}\,.
\end{equation*}

\begin{proposition}\label{thm:potential_calculation_higherNCRank}
    Choose $\Theta$ weights $\Lambda_i$  which are faithful on the root $\beta_i$. Then the potential $D=\delta_i(A,A_l,A_r)$ is given by 
    $$\iota(D^*) = \sigma_i(\Lambda_i(l_r))\Lambda_i(l_b)\Lambda_i(l_l)^{-1}.$$
\end{proposition}

Thus, we can use the partial potential to label the networks with a free collection of elements in $\GL_1(\mathcal{R})$, $L_1,\cdots,L_{p+1}, B_i,B_{i+1}, D$ as in \Cref{fig:ElementaryNetworksWeightedByPotential}. The analogous weighting for $\elemNetwork_\circ(i)$ is given by setting $B_i = B_{i+1} = 1$ and removing the legs.
\input{Figures/figureElementaryNetworksPotentialWeights}

We will also use these networks to describe elementary configurations for groups $G$ with a $C_p$-grading. In this case, $\mathcal{R}$ has an anti-involution $\star$ and $G\simeq \SP_{2p}(\mathcal{R},\star)$ by the recognition theorem (\Cref{thm:RecognitionCp}). Then $G$ is a subgroup of $\GL_{2p}$ and will use elementary networks of $\GL_{2p}$ to describe those for $G$. In particular, for $i< p$ the elementary configuration $\Conf_3^i(\decFlagTheta)$ is represented by the amalgamation of $\elemNetwork(i)$ and $\elemNetwork(\oppInv{i})$ for $\GL_{2p}$ with the additional constraint that elements of $l\in L_\Theta, l = (L_1,\dots,L_{2p}) \subset \GL_{2p}(\mathcal{R})$ now satisfy that $L_i= \star(L_{2p-i+1}^{-1})$. 
We will assign elementary networks with $2p$ strands which are symmetric by reflecting along the line between strands $p,p+1$ exchanging labels under this action. For $\Conf_3^p(\decFlagTheta)$ we use the network $\elemNetwork(p)$ for $\GL_{2p}$ with the extra symmetry on strand weights. In \Cref{fig:C2elemNetworks} we see the two networks needed to parameterize type $C_2$.
\input{Figures/figureC2ElementaryNetworks}

In order to recover the $L$-distances from the network, we consider a path from infinity to the ground which passes through $p+1$ strands. If the strand crossed at height $i$ is labeled $X$ with $\bullet$ on the right, the $i^{th}$ component of the $L$-distance is $X$. If instead $\circ$ is on the right, the $i^{th}$ component is $X^{-1}$. When this path is homomotopic to a side of the network without crossing any bridges or legs it recovers the $L$-distance of that side. 

Similarly, the path parallel to the ground traveling from the right side to the left side recovers $l_b$. The index of the legs is given by the strand it lands on after sliding it to the left past all bridges. By construction this path always has $\circ$ vertices on the right, so all the entries of the $L$-distance will be the inverse of the weight in the network.

\section{Triples of \texorpdfstring{$\Theta$}{Theta}-flags}\label{sec:flagTriplesTheta}

The results of the previous section imply the following lemma: 

\begin{lemma}\label{lem:jordanpoints_elemenseed}
    The $\{\jordan{J}_i\}$ points of the noncommutative cluster variety $\A_{\elemSeed(i)}$ parameterize the elementary configuration space $\Conf_3^i(\decFlagTheta)$. 
\end{lemma}
Note that there are no mutable vertices or possible square moves at this point, so this cluster variety is really just a noncommutative torus. 

Now we will amalgamate these cluster varieties into one which describes generic triples of flags.

Consider a triple $(A_1,A_2,A_3)$ of pairwise transverse flags in $\decFlagTheta$. Then there exist $g\in G$, $l_2,l_3\in L_\Theta$ and $u_3\in U_\Theta$ such that
\begin{equation*}
    g(A_1,A_2,A_3)=(U_\Theta,\overline{w_0}U_\Theta l_2,u_3\overline{w_0}U_\Theta l_3)\,.
\end{equation*}
\begin{proposition}\label{thm:ParameterizationOfConf3Theta}
    The map
    \begin{align*}
        \Confx_3(\decFlagTheta) & \to L_\Theta \times L_\Theta \times U_\Theta \\
        (A_1,A_2,A_3) & \mapsto (l_2,l_3,u_3)
    \end{align*}
    is an open dense embedding.
\end{proposition}
\begin{proof}
  The image of the map on the $U_\Theta$-factor is $U_\Theta\cap P_\Theta^\opp w_0 P_\Theta^\opp$, which is an open dense subset of $U_\Theta$.
\end{proof}

We now describe how to assemble the cluster charts for elementary configurations into a cluster chart for a triple of flags: We fix a reduced expression $D(w_0)=\varsigma_{i_1}\cdots \varsigma_{i_n}$ for the longest word $w_0$ in $W(R)$. Given a triple $(A_1,A_2,A_3)\in\decFlagTheta$ which represents an element of $\Confx_3(\decFlagTheta)$, we may use this reduced expression to compute an interpolating sequence $\{A_2^k\}$ between $A_2$ and $A_3$. As we have seen we obtain elementary configurations $(A_1,A_2^k,A_2^{k+1})$ in this way, and we constructed elementary seeds/networks for each elementary configuration space. We can now amalgamate these elementary seeds according to the reduced expression $D(w_0)$ to obtain a cluster chart for $\decFlagTheta$.

\subsection{Amalgamating}
Depending on the noncommutative rank, we need to amalgamate either polygonal quivers or noncommutative networks.

\subsubsection{Noncommutative rank 1}
The seeds for each elementary quiver are planar, and so it makes sense to amalgamate adjacent frozen nodes at the same level. For each $k$, let $k^+$ be the next index for which $i_k = i_{k^+}$. At the end of the process every node at level $i_k$ in the elementary seeds between $i_k$ and $i_{k^+}$ will be amalgamated. Due to \Cref{thm:elementaryConfigLeftRightFunctionsTheta} these nodes are assigned equal functions, and thus every node has an associated function on $\Confx_3(\decFlag)$. Observe that the extra decoration which turned $\elemSeed(p)$ into an ordered quiver is chosen such that it makes sense to amalgamate on level $p$. The half-arrows either combine to give a full arrow if they are oriented the same way or cancel if they have opposite orientations. Internal nodes become unfrozen.

The extra nodes on the bottom are connected with half edges oriented from left to right according to the adjacency of the coroots $\check{\beta}_k$ in the associated Dynkin diagram. For more details see \cite{goncharov2019quantum}.

The resulting quiver has three components of frozen nodes connected by half arrows: the left nodes of the initial elementary seed, the right nodes of the final elementary seed, and the bottom nodes consisting of the extra frozen nodes. By \Cref{rem:LdistanceFromElementaryConfig}, $\Delta_p(A_3,A_2) = w_{k-1}(\check{\beta}_p(A))$ where $w_{k-1}$ is the prefix of the reduced expression which returns $\check{\beta}_p$ to itself and $l(A_2^{k},A_2^{k-1}) = \check{\beta}_p(A)$. As we saw in \Cref{thm:ActionOfSigmaPonSpecialRoot}, the action of $w_{k-1}$ might not be trivial on the special coroot space. This is not a problem, as the difference can only be application of the involution $\sau$. The final step in this construction is then to weave at the frozen weight $r$ node so the function agrees with $\Delta_p(A_3,A_2)$. We compute the action of $w_{k-1}$ specifically for type $B_p$ in \Cref{thm:InterpolatingSeqInvolutionBp} and for type $G_2$ in \Cref{thm:potentialInvolutionG2}. 

\subsubsection{Higher noncommutative rank}
Our choice of Jordan weights factors into our construction of the interpolating sequence. In higher noncommutative rank we amalgamate our networks simply by gluing their adjacent boundary black vertices on the same strands together. This produces a grounded wiring network with edge weights as the gluing respects the cyclic product condition. However, the amalgamated network has $2p$ legs and we would like a network with exactly $p+1$ legs so that the ground side matches the left and right. 

Sliding (\Cref{def:LegSlide}) each leg to the left past all the brides identifies a unique strand for each leg. All legs with the same label can be replaced by a single leg whose weight is the product of the weights. By the construction of the system of Jordan weights, the product is independent of the order even though the ring itself is noncommutative. 

\begin{definition}
    A grounded wiring network obtained by sliding all the legs to the left of all the bridges and reducing to at most one leg per strand is called \keyword{left slid}.
\end{definition}

\begin{proposition}
        The left slid grounded wiring network obtained is independent of the choice of system of Jordan weights. 
\end{proposition}

\begin{proof}
    Let $\tilde{A_2}$ be such that its undecorated flag is the same as $A_2$ but $l(A_3,\tilde{A_2}) =1$. Now, every choice of Jordan weights gives the same interpolating sequence between $\tilde{A}_2$ and $A_3$. The grounded wiring network obtained from amalgamating elementary networks for the interpolating sequence between $\tilde{A}_2$ and $A_3$ is the same as the left slid network other than the leftmost labels and the legs. Clearly this is true for any choice of Jordan weights. Furthermore, different systems Jordan weights change the decomposition $l_b = \prod \check{\beta}_i(\Lambda_i(l_b))$ but not the final projection to the $i^{th}$ component. This projection is the weight of the leg of the $i^{th}$ strand and thus is independent of the choice of Jordan weights.  Therefore the left slid network is independent as claimed.
\end{proof}

\begin{corollary}
    The amalgamated elementary networks associated to two different choices of Jordan weights are equivalent by leg slides. 
\end{corollary}

\subsection{Cluster charts}
The result of the construction is a set of functions on $\Confx_3(\decFlagTheta)$ and an ordered quiver or network which indexes the functions.
\begin{definition}
    A set of functions on $\Confx_3(\decFlagTheta)$ arising from this construction is called a \keyword{cluster chart} for $\Confx_3(\decFlagTheta)$. These functions live on the nodes of an ordered quiver or edges of a network we denote by $Q_{D,1}$, where $D$ is the reduced expression of $w_0$. We write $\A_{D,1}$ for the noncommutative torus for which $Q_{D,1}$ is the group of monomial functions. 
\end{definition}
\begin{remark}
    The superscript $1$ is used to distinguish these charts from the `rotated' ones that we will introduce in \Cref{sec:rotationCharts}.
\end{remark}
We note that unlike the commutative charts, these functions satisfy multiplicative relations. Nevertheless they replace the ordinary cluster charts of \cite{goncharov2019quantum} and can be mutated like ordinary cluster charts.

\begin{definition}
    A \keyword{generalized mutation} is any of the following operations from in \Cref{part:NoncomClusterVarieties}.
    \begin{itemize}
    \item Mutation of polygonal quivers (\Cref{def:JordanPointMutation})
    \item Square moves of weighted networks (\Cref{def:SquareMove})
    \item Switching/ Weaving of noncommutative nodes (\Cref{def:weavingQuiver})
    \item Leg slides (\Cref{def:LegSlide})
    \item Regrounding a planar network (\Cref{def:ChangeOfGround})
\end{itemize}
\end{definition}

The first two operations are birational transformations and replace the standard mutation of commutative cluster algebras. Switching or weaving at a noncommutative node $k$ results in essentially the same set of functions with a different representative of $\{X_k, \sigma(X_k), \tau(X_k), \sau(X_k)\}$ chosen for the node. Similarly changing the ground, changes the underlying combinatorial representation but leaves all the functions the same.

The final operation of a leg slide is slightly more complicated. It results in ``the same chart'', but with the functions changed by a noncommutative monomial transformation.
\begin{definition}
    Let $H$ be a group, $h_1,h_2 \in H$. A \keyword{noncommutative monomial transformation} is a map $H\to H$ expressed by $m_{h_1,h_2}(g):= h_1gh_2^{-1}$.
\end{definition}

There is another important monomial transformation on cluster charts for $\decFlagTheta$ given by the rescaling action (\Cref{def:rescaling_action}).
\begin{proposition}\label{prop:rescaling_is_monomial}
    The rescaling action acts on the functions in any cluster chart via pull back along a noncommutative monomial transformation. 
\end{proposition}
\begin{proof}
    This follows from \Cref{thm:HDistPropertiesTheta} and \Cref{lem:rescaling_action_interpolating}; the rescaling action rescales all of the flags in the interpolating sequence used to make our cluster charts and the second property of the $L-$distance is clearly a monomial transformation. 
\end{proof}

\noindent Our main goal will be to prove:
\begin{theorem}\label{thm:clusterChartsTriplesRelatedByMutationTheta}
    The noncommutative cluster variety $\A_{|Q_{D,1}|}$ contains a torus for each cluster chart associated to all reduced expressions $D'$ of $w_0$ and the rotated cluster charts. Its $\{\jordan{J}_i\}$ points parameterize open dense subsets of $\Confx_3(\decFlagTheta)$.
\end{theorem}
The cluster charts we have seen so far depend on the choice of a reduced expression for $w_0\in W(R)$ in Steps 2 and 3. We will prove \Cref{thm:clusterChartsTriplesRelatedByMutationTheta} for charts obtained by changing the reduced expression in \Cref{sec:changingReducedExpressionTheta}.

First, we establish that these cluster charts deserve to be called charts for $\Confx_3(\decFlagTheta)$. For this we need:
\begin{lemma}
\label{thm:tripleFlagsPolyCAEvaluation}
    The functions from a cluster chart evaluated on any element of $\Confx_3(\decFlagTheta)$ define a $\{\jordan{J}_i\}$-point of $\A_{D,1}$.
\end{lemma}
\begin{proof}
    In the quiver case  we check that all angles are mapped to the image of $\iota$. This holds by the construction of the elementary seeds, the angles at the weight $r$ nodes recover $\iota(\oppInv{v})$ for $v\in V$ the potential of the elementary quiver. In the network case, the Jordan algebra is special and the edge weights live in the specializing ring $\mathcal{R}$.
\end{proof}

\noindent Next we need the converse, that each Jordan point corresponds to a triple of decorated flags.
\begin{theorem}\label{thm:coordinatesTriplesTheta}
    The $\jordan{J}$-points of $\A_{D,1}$ parametrize an open dense subset of $\Confx_3(\decFlagTheta)$.
\end{theorem}
\begin{proof}
    We will show that the $\jordan{J}$-points parameterize the set of elements of $\Confx_3(\decFlagTheta)$ which can be represented as
    \begin{equation*}
        (A_1,A_2,A_3)=(U_\Theta,\overline{w_0}U_\Theta l_2,F_D(t_{i_1},\dots,t_{i_n})\overline{w_0}U_\Theta l_3)
    \end{equation*}
    with $t_{i_k}\in V_i^\times$. Since the $R$-Lusztig map (\Cref{def:LuszitigMap}) has open dense image in $U_\Theta$, this set of triples is dense in $\Confx_3(\decFlagTheta)$.
    
    Notice that we use the $R$-Lusztig map for the same reduced expression $D$ for the chart. The image of the interpolating sequence $\{A_2^k\}$ under the projection $\pi:\decFlagTheta\to\flagTheta$,
    \begin{equation*}
        B_2^k=\left(\prod_{j=1}^{k-1} x_{i_j}(t_j) \right)\overline{w_0}P_\Theta\,,
    \end{equation*}
    recovers the Lusztig map factor by factor. Thus the partial potential of the elementary configuration $(A_1,A_2^{k-1},A_2^k)$ recovers the parameter $t_k$. By \Cref{thm:potentialFormulaTheta}, the partial potential is a simple function of the $\jordan{J}$-points. Since $l(A_1,A_2)=l_2$ and $l(A_1,A_3)=l_3$, we can compute the remaining parameters from the points associated to the left/right edges of the triangle. Therefore, any such configuration can be reconstructed from the $\jordan{J}$-coordinates as needed.
\end{proof}

\subsection{Changing reduced expression}\label{sec:changingReducedExpressionTheta}
    We will now explain how cluster charts coming from different reduced expressions of $w_0\in W(R)$ are related. Reduced expressions can be changed from one to another by the following moves \cite{bjorner2005combinatorics}:
    \begin{proposition}
        Let $\{\varsigma_i\,|\, 1\leq i \leq p\}$ denote the standard generators of the Weyl group $W(R)$ corresponding to the root grading. Any two reduced expressions of $w\in W$ can be related by a finite sequence of the following \keyword{braid moves}:\medskip
        
        \begin{tabular}{rll}
            \textbf{Two Move} &  $\varsigma_i\varsigma_j=\varsigma_j\varsigma_i$ &for $|i-j|>1$.\\
            \textbf{Three Move} & $\varsigma_i\varsigma_{i+1}\varsigma_i=\varsigma_{i+1}\varsigma_i\varsigma_{i+1}$ &for $1\leq i\leq p-2$\\
            \textbf{Four Move} &  $\varsigma_{p-1}\varsigma_p\varsigma_{p-1}\varsigma_p=\varsigma_{p-1}\varsigma_p\varsigma_{p-1}\varsigma_p$   & ($B_p/C_p$ case)\\
            \textbf{Six Move} & $\varsigma_{1}\varsigma_2\varsigma_{1}\varsigma_2\varsigma_{1}\varsigma_2=\varsigma_2\varsigma_{1}\varsigma_2\varsigma_{1}\varsigma_2\varsigma_1$ & ($G_2$ case).
        \end{tabular}
    \end{proposition}

    Thus, to prove \Cref{thm:clusterChartsTriplesRelatedByMutationTheta} it suffices to show that each braid move is realized by a sequence of mutations. 
    \begin{lemma}[Two Move]\label{thm:twoMoveTheta}
        The seeds given by two words which differ by a two move only differ in the indexing of the functions. 
    \end{lemma}
        \begin{proof}
            Let $1\leq i,j\leq p$ with $|i-j|>1$. Clearly, the quiver obtained by amalgamating $\elemSeed(i)$ and $\elemSeed(j)$ is the same as the one obtained by amalgamating $\elemSeed(j)$ with $\elemSeed(i)$. According to \Cref{thm:elementaryConfigLeftRightFunctionsTheta} all the functions can be computed in terms of the two decorated flags with $W$-distance $\varsigma_i\varsigma_j=\varsigma_j\varsigma_i$, and are thus independent of whether we use $\varsigma_i\varsigma_j$ or $\varsigma_j\varsigma_i$ to obtain the interpolating flag.
        \end{proof}

        \input{Figures/figureThreeMove}

    The three move in noncommutative rank 1 is proved essentially the same way as in \cite{goncharov2019quantum}:
    \begin{lemma}[Three move, NC Rank 1]\label{thm:threeMoveTheta}
        The seeds given by two words which differ by a three move $\varsigma_i \varsigma_{i+1} \varsigma_i \rightarrow \varsigma_{i+1} \varsigma_i \varsigma_{i+1}$ are related by a single mutation at the single mutable node of the subquiver $\elemSeed(i,i+1, i)$.
    \end{lemma}
        \begin{proof}
            We see that the subquivers $\elemSeed(i,i+1,i)$ and $\elemSeed(i+1,i,i+1)$ in \Cref{fig:ThreeMove}. Label the functions as indicated in the figure. By \Cref{thm:elementaryConfigLeftRightFunctionsTheta} we observe that all the functions at the frozen nodes can be calculated in terms of the decorated flags with $W$-distance $\varsigma_i\varsigma_{i+1}\varsigma_i$ which are independent of this choice of interpolating sequence. This justifies labeling the frozen nodes the same in both quivers in \Cref{fig:ThreeMove}. Thus, it suffices to verify the following identity relating the functions at the mutable node: 
            \begin{equation*}
                XX'=X_i^l X_{i+1}^r X_i^e + X_{i+1}^l X_i^r X_{i+1}^e\,,
            \end{equation*}
            which can equivalently be written as
            \begin{equation*}
                \frac{1}{X_i^l X_i^r}\cdot X'X_{i-1} = \frac{X_i^e}{XX_i^r}\cdot X_{i+1}^rX_{i-1} + \frac{X_{i+1}^e}{XX_i^l}\cdot X_{i+1}^lX_{i-1}\,.
            \end{equation*}
            The terms in this equation are by \Cref{thm:potentialFormulaSplit} the partial potentials for the elementary configurations for $\varsigma_i$ associated to the two interpolating sequences. Thus the equation holds by the additivity of the potential at height $i$. Essentially, this boils down to the equation \cite{lusztig1994total}
            \begin{equation*}
                x_i(a)x_{i+1}(b)x_i(c)=x_{i+1}\left(\frac{bc}{a+c}\right)x_i(a+c)x_{i+1}\left(\frac{ab}{a+c}\right)\,.
            \end{equation*}
        \end{proof}
    We leave the remaining calculations to later sections, the higher noncommutative rank three move in \Cref{sec:ThreeMoveHigherNoncom}, the four move for $B_p$ in in \Cref{sec:fourmove}, the four move for $C_p$ in in \Cref{sec:FourMoveHigherNoncom},  and the six move for $G_2$ in \Cref{sec:sixmove} respectively. 

    \subsection{Rotation}\label{sec:rotationCharts}
    The notion of positivity in flag varieties is inspired by tuples on the circle which follow an orientation. In this case any cyclic rearrangement preserves this notion of positivity. This is also true for generic configurations in $\flagTheta$ but not in decorated flag varieties. Instead, we need to use the central element $s_G\in G$:

    \begin{definition}
        The \keyword{twisted cyclic shift} map $\tw:\Confx_3(\decFlagTheta)\to\Confx_3(\decFlagTheta)$ is defined by
        \begin{align*}
            \tw:(A_1,A_2,A_3)\mapsto (s_GA_2,A_3,A_1)\,.
        \end{align*}
    \end{definition}

    We can obtain additional cluster charts of $\Confx_3(\decFlagTheta)$ by precomposing the above construction with the twisted cyclic shift:
    \begin{definition}
        The \keyword{rotated cluster charts} are obtained by applying the previous construction to the triples $\tw(A_1,A_2,A_3)$, respectively $\tw^2(A_1,A_2,A_3)$. The coordinate functions are indexed by quivers/networks $Q_{D,2}$, respectively $Q_{D,3}$, where $D$ is a reduced expression of $w_0\in W(R)$.
    \end{definition}
    Clearly, all the results about the original cluster charts also apply to the rotated ones. For the discussion of decorated twisted local systems in \Cref{sec:localSystems} the following result is crucial:

    \begin{proposition}\label{thm:coordinateSystemsS3MutationTheta}
        The twisted cyclic shift is realized by a sequence of generalized mutations. 
    \end{proposition}

    \begin{proof}
        This is trivial to prove in the $A_1$ case as our cluster charts are already rotationally symmetric up to the orientation of the edges of the triangle. But the choices were made exactly so that the difference between the ordered quivers is a switch at the corresponding node, while the functions are related by the opposition involution, which is exactly $\sigma$ in the $A_1$ case.

        The $A_p$ case is proved by changing the choice of ground in the honeycomb $GL_{p+1}$ network (\Cref{fig:GlnNetwork}). The $GL_{p+1}$ network is the network constructed using the word $\varsigma_1\cdot \varsigma_2\varsigma_1 \cdot \varsigma_3\varsigma_2\varsigma_1\cdot \ldots \cdot \varsigma_p\varsigma_{p-1}\ldots \varsigma_1$. So we first use braid transformations to move to this chart, then we change the choice of ground, and finally we undo the braid transformations. The $C_p$ case follows from the $A_{2p-1}$ case since in this case our charts are really just the same as charts for groups of type $A_{2p-1}$ and we can rotate using the same procedure. 
        
        For the $B_p$ and $G_2$ cases we prove this in \Cref{sec:RotationThetaBp,sec:sixmove} via an explicit calculation. This is slightly more involved even for the edge functions since the opposition involution is not always $\sigma$ but can also be $\tau$ on $\Gamma(\jordan{J})$, in which case one has to weave and switch at the corresponding node.
    \end{proof}

    \begin{corollary}\label{thm:twPositiveTriples}
        If $c\in\Confx_3(\decFlagTheta)$ is a positive configuration, so is $\tw(c)$.
    \end{corollary}

\subsection{Changing orientations}\label{sec:orientations}
Given a triple $(A_1,A_2,A_3)\in\decFlagTheta^3$, we have assigned a total order given by $A_1\to A_3\to A_2$ in the construction of cluster charts, and we have computed the functions associated to frozen edges accordingly: If $A_i\to A_j$, the functions on the corresponding edge are given by $\Lambda_i(A_i,A_j)$. We have seen that for the rotated cluster charts the edge functions are related to the original ones via the opposition involution. 

In this spirit we could upgrade the construction of cluster charts even further by allowing arbitrary orientations for the edges of the triangle and computing the edge functions as follows: If the orientation of an edge $A_i\to A_j$ is compatible with the total order, we assign the functions $\Lambda_i(A_i,A_j)$, otherwise we assign $\oppInv{\Lambda_i(A_j,A_i)}$.

\begin{example}
    Consider the triangle with the side edges oriented ``up'' from $A_2 \rightarrow A_1$ and $A_3 \rightarrow A_1$ when $G$ is $B_2$-graded. In \Cref{fig:orientationChangeB2}, we see the polygonal quiver for the standard case and with the edge orientation changed. In this case the opposition involution fixes the root system, but acts via the involution $\tau$ on the $p=2$ root space (\Cref{thm:potentialInvolutionBp}). We observe that in order to amalgamate the elementary quiver to the sides, we first weave and switch at the boundary. Note these operations preserve the angle and the partial potential.
    
    \input{Figures/figureOrientationChangeB2}
\end{example}
\begin{example}
    We now consider a similar example when $G$ is $A_2$-graded. In this group, the opposition involution acts by reversing the order of components and applying the inverse: $l = [A,B,C] \mapsto [C^{-1},B^{-1},A^{-1}]$). This action preserves the edge weights assigned to the edges and the network is unchanged.
\end{example}

\section{Tuples of \texorpdfstring{$\Theta$}{Theta}-flags}\label{sec:TuplesOfThetaFlags}

We are now prepared to conclude the detailed explanation of the construction of cluster coordinates on $\Confx_n(\decFlagTheta)$ which we outlined in \Cref{sec:coordinateStrategy}. Thus let $(A_1,\dots,A_n)$ represent an element of $\Confx_n(\decFlagTheta)$. Recall that we start by associating the flags to the vertices of an $n$-gon in a counterclockwise manner. Next, fix a triangulation of the $n$-gon. We will use this data to obtain an element of $\Confx_3(\decFlagTheta)$ for every triangle as follows:

First, orient all arcs of the triangulation, including the sides of the $n$-gon via the total order given by reading clockwise around the $n$-gon:
\begin{equation*}
    A_1 \rightarrow A_n \rightarrow A_{n-1} \rightarrow \dots \rightarrow A_{3} \rightarrow A_2 \,.
\end{equation*}
Given a triangle $T$, let $i,j,k$ be the vertices of $T$ ordered such that $A_i \rightarrow A_k \rightarrow A_j$ in the total order we just defined. Then the configuration associated to $T$ is represented by the triple $(A_i, A_j, A_k)\in\decFlagTheta^3$.

For every triangle we choose one side and a reduced expression of the longest element $w_0\in W(R)$. Using the construction explained in \Cref{sec:flagTriplesTheta} we obtain a (possibly rotated) cluster chart on $\Confx_3(\decFlagTheta)$ which we evaluate on the configuration associated to the triangle.

For every triangle we have a seed which we view as inscribed in it. We then amalgamate these seeds along the internal arc of the triangulation. As explained in \Cref{sec:orientations}, in the $B_p$ and $G_2$ cases we may need to weave/switch so that the edge functions for the individual triangles become compatible with the orientations of the sides.

In the $A_p$ and $C_p$ cases, we need to amalgamate our seeds by gluing the networks along the boundary $\bullet$ vertices with the same index on the edges of the triangles. Recall that by construction, the number of legs is equal to the number of strands and each leg has a unique index given by performing a full slide to the left. 

In this way we obtain an ordered quiver or network inscribed into the $n$-gon whose nodes index a set of functions on $\Confx_n(\decFlagTheta)$.

\begin{definition}
    A set of functions on $\Confx_n(\decFlagTheta)$ arising from the above construction is called a \keyword{cluster chart}. We will write $\ASpace{R}{n}$ for the noncommutative cluster variety associated to (any) cluster chart.
\end{definition}
Recall that $R$ here refers to the Jordan split type of the group we are considering.

\input{Figures/figureHeptagonFanTriangulationTheta}

\begin{theorem}
\label{thm:coordinatesDetectTuplePositivityTheta}
    The noncommutative cluster variety $\ASpace{R}{n}$ contains seeds for all cluster charts. Moreover the $\{\jordan{J}_i\}$ points of any cluster chart parametrize an open dense subset of $\Confx_n(\decFlagTheta)$.
\end{theorem}

\begin{proof}
    Very nice charts are given by the \keyword{fan triangulation}, the oriented triangulation for which all internal arcs start at $A_1$, see \Cref{fig:HeptagonFlagsFanTriangulationTheta}. 
    Any point in $\Confx_n(\decFlagTheta)$ has a representative of the form
    \begin{equation*}
        (U_\Theta,\overline{w_0}U_\Theta l_2, u_3\overline{w_0}U_\Theta l_3, u_3 u_4 \overline{w_0}U_\Theta l_4,\dots,u_3\cdots u_n \overline{w_0}U_\Theta l_n)\,.
    \end{equation*}
    This provides an open dense embedding of $\Confx_n(\decFlagTheta)\to L_\Theta^{n-1}\times U_\Theta^{n-2}$. Now we simply extend the proof of \Cref{thm:coordinatesTriplesTheta}: If $u_3,\dots,u_n$ lie in the image of the Lusztig map, they can be reconstructed from the coordinates one by one, just like the elements $l_2,\dots,l_n$. Then \Cref{thm:CoordinatesDetectTriplePosTheta} shows that a positive $J$-point reconstructs a positive $u_i$ and $l_i$. This is the definition of a positive tuple of decorated flags. 
 
    In order to extend this proof to charts for arbitrary triangulations, we will see that any two such charts are connected by series of cluster mutations (\Cref{thm:clusterChartsTuplesRelatedByMutationTheta}). The resulting map from the set parametrized by one cluster chart to the set parametrized by the other is defined on an open dense subset. Moreover such mutations map positive $\jordan{J}$-points to positive $\jordan{J}$-points.
\end{proof}

\subsection{Flip}\label{sec:FlipTriangulationGenericTheta}
In order to extend the proof of \Cref{thm:coordinatesDetectTuplePositivityTheta} for the fan triangulation to any cluster chart, we will prove:
\begin{theorem}\label{thm:clusterChartsTuplesRelatedByMutationTheta}
    Any two cluster charts for $\Confx_n(\decFlagTheta)$ are related by general mutations. 
\end{theorem}
 We already proved \Cref{thm:clusterChartsTriplesRelatedByMutationTheta}, which accounts for the different choices within each individual triangles. Since flips relate any two triangulations, it suffices to prove that flips in the triangulation can be realized by generalized mutations.
 To understand the flip, we adapt Goncharov-Shen's approach by considering four flags on a triangulated square, where we take interpolating sequences across the top and bottom edges (\cite[Section 11]{goncharov2019quantum}). In this case, the elementary configurations on the top can be identified with flipped opposite elementary configurations (\Cref{def:OppositeElementaryConfiguration}). These amalgamate nicely with any of the standard elementary configurations allowing for many mixed charts. These charts are parametrized by words in $W(R)\times W(R)$, where we denote the generators of the bottom edge by $\varsigma_i:=(\varsigma_i,1)$ and the generators of the top edge by $\varsigma_{\overline{i}}:=(1,\varsigma_{\oppInv{i}})$ for $1\leq i\leq p$. 
 
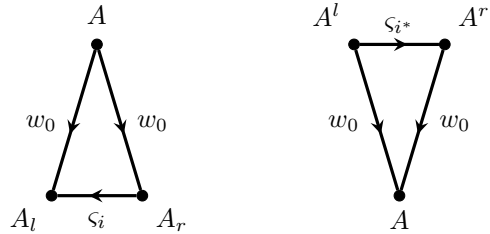
 By amalgamating the corresponding elementary and opposite elementary charts, we associate a cluster chart to every reduced expression of $(w_0,w_0)\in W(R)\times W(R)$, similar to \Cref{fig:flipSetup}. In the initial setup, all elements of the form $\varsigma_i$ come before elements of the form $\varsigma_{\overline{i}}$ and after the flip the order is reversed. Thus, to prove \Cref{thm:clusterChartsTuplesRelatedByMutationTheta} it suffices to prove that the following moves on these cluster charts can be realized by mutations:\medskip

\begin{tabular}{ll}
    \textbf{Two Move} &  $\varsigma_i\varsigma_{\overline{j}}=\varsigma_{\overline{j}}\varsigma_i$ for $i\neq j$.\\
    \textbf{Twisted Two Move} & $\varsigma_i\varsigma_{\overline{i}}=\varsigma_{\overline{i}}\varsigma_i$.
\end{tabular}\medskip

\noindent The two move in this extended version, can be treated the same as before (\Cref{thm:twoMoveTheta}).

\begin{lemma}\label{thm:twistedTwoMoveTheta}
    The seeds given by two words which differ by a twisted two move $i\overline{i}\to\overline{i}i$ are related by a single mutation  on the cluster chart $\elemSeed(i\overline{i})$. 
\end{lemma}

    Let us make some preparatory observations: Consider a triple $(A,A_l,A_r)$ which represents an element of $\Conf_3^i(\decFlagTheta)$. There is a unique (undecorated) flag $B'\in\flagTheta$ such that $(A_r,B')$ are $\varsigma_i$-transverse and $(A,B')$ are $w_0\varsigma_i$-transverse. We can interpret the triple $(A_r,B',\pi(A_l))\in\decFlagTheta\times\flagTheta^2$ as a triple in the flag variety for $\SLJ$, and thus we can assign a (total) potential $W(A_r,B',\pi(A_l))$ to the triple in analogy to the partial potentials defined in \Cref{def:partialPotentialTheta}. Concretely this means:
    \begin{equation*}
        W(U_\Theta,\varsigma_iP_\Theta,x_i(v)\varsigma_iP_\Theta)=v\,.
    \end{equation*}
    Clearly, this potential is additive: Since $x_i(v)x_i(w)=x_i(v+w)$ we have
    \begin{equation*}
        \begin{split}
            W(U_\Theta,\varsigma_iP_\Theta,x_i(v) & x_i(w)\varsigma_iP_\Theta) = v+w\\
            &= W(U_\Theta,\varsigma_iP_\Theta,x_i(v)\varsigma_iP_\Theta)+W(U_\Theta,x_i(v) \varsigma_iP_\Theta,x_i(v)x_i(w)\varsigma_iP_\Theta)\,.
        \end{split}
    \end{equation*}
    We stress that this is really a potential on the $\SLJi$-flag variety, which replaces the $\SL_2$-flag varieties appearing in \cite{goncharov2019quantum}. Now we compute this potential using the parameters $l_l,l_r\in L_\Theta$, and $a\in V_i^\times$ for the elementary configuration $(A,A_l,A_r) = \left(U_\Theta,\overline{w}_0U_\Theta l_l,x_{i^*}(a^*)\overline{w}_0U_\Theta l_r\right)\,$. 

    First we choose an alternate representative of the triple where we can see $B' = \varsigma_iP_\Theta$:
    \begin{equation*}
        (A,A_l,A_r)=(l_r^{-1}\doverline{w}_0U_\Theta,l_r^{-1}\doverline{w}_0 x_{\oppInv{i}}(-\oppInv{a})\overline{w}_0 U_\Theta l_l,U_\Theta ) \,.
    \end{equation*}
     Using the opposition involution (\Cref{def:OppositionInvolution}) and the decompositions from \Cref{thm:DecomposeXSigma}, we compute
    \begin{align}\label{eqn:Sl2likePotentialInTermsOfUsualPotential}
        \begin{aligned}
            W(A_r,B',\pi(A_l))&=W(U_\Theta,\varsigma_i P_\Theta,l_r^{-1}\doverline{w}_0 x_{i^*}(-\oppInv{a})\overline{w}_0 P_\Theta)\\
        &= W(U_\Theta,\varsigma_i P_\Theta,l_r^{-1} y_{i}(a) P_\Theta)\\
        &= W(U_\Theta,\varsigma_i P_\Theta,l_r^{-1}\doverline{\varsigma}_i x_{i}(-a)\overline{\varsigma}_i P_\Theta)\\
        &= W\left(U_\Theta,\varsigma_i P_\Theta,l_r^{-1}\doverline{\varsigma}_i y_{i}\left(-{a}^{-1}\right)\check{\beta}_i(\iota(-a))x_i\left({a}^{-1}\right) P_\Theta\right)\\
        &= W\left(U_\Theta,\varsigma_i P_\Theta,l_r^{-1}x_i\left({a}^{-1}\right) \varsigma_i P_\Theta\right)\\
        &= \beta_i(l_r^{-1})\left({a}^{-1}\right)\,.
        \end{aligned}
    \end{align}
    Notice that this potential is exactly the same as the potential of the reflected configuration. Therefore we can use \Cref{lem:reflected_potential} to compute 
    \begin{equation*}
         \check{\beta}_i(\iota(W(A_r,B',\pi(A_l)))) = l_r^{-1}l_ll_b^{-1}
    \end{equation*}

    \begin{proof}
    Now, consider a configuration of four decorated flags represented by $(A^l,A_l,A_r,A^r)\in\decFlagTheta^4$ such that $(A^l,A_l,A_r)$ and $(A^r,A^l,A_r)$ represent elementary configurations in $\Conf_3^i(\decFlagTheta)$ and $\Conf_3^{\bar{i}}(\decFlagTheta)$ respectively, as in \Cref{fig:twistedTwoMoveSetupTheta}. Let $a,b\in V_i$ be the potentials of these configurations. Notice that $(A^r,A_l,A_r)$ is also a configuration in $\Conf_3^i(\decFlagTheta)$, denote by $c$ its potential. Let $l_l^l = l(A^l,A_l)$, $l_r^l = l(A^l,A_r)$, $l_r^r = l(A^r,A_r) \in L_\Theta$.

    Let $B_1,B_2\in\flagTheta$ such that $(A_r,B_1)$ and $(A_r,B_2)$ are $\varsigma_i$-transverse and $(B_1,A^l)$ and $(B_2,A^r)$ are $\varsigma_i w_0$-transverse. By the previous calculations we have
    \begin{align*}
         W(A_r,B_2,\pi(A_l)) &= W(A_r,B_1,B_2)+W(A_r,B_1,\pi(A_l))\\
          \beta_i\left((l^r_r)^{-1}\right)\left({c}^{-1}\right) &= \beta_i\left((l^l_r)^{-1}\right)\left({a}^{-1}\right) + b \,.
    \end{align*}
    By applying $\check{\beta}_i\circ\iota$, each of the three terms can be computed as a product of $L$-distances,
    \begin{equation*}
        (l^r_r)^{-1}l^r_ll_b^{-1} \qquad (l^l_r)^{-1}l^l_ll_b^{-1} \qquad(l^r_r)^{-1}(l^l_r)^*\varsigma_{i}(l^l_r)\,.
    \end{equation*}
   
    \input{Figures/figureTwistedTwoMoveSetupTheta}
    
    To finish the proof we apply the appropriate Jordan weights to obtain the mutation formula in the chart $E(i\bar{i})$. In noncommutative rank 1, when $i=p$, the coordinate are given by the following functions assigned to nodes as in \Cref{fig:twistedTwoMoveSetupTheta}:
    \begin{equation*}
        \begin{split}
            A &:= \Lambda_p(A^l,A_l)\\
            C &:= \Lambda_p(A^r,A_r)\\
            E &:= \Lambda_p(A^l,A_r)
        \end{split}\hspace{2cm}
        \begin{split}
            B &:= \Lambda_p(A_r,A_l)\\
            D &:= \Lambda_p(A^l,A^r)\\
            x &:= \Lambda_{p-1}(A^l,A_r)
        \end{split}
    \end{equation*}
    where $x$ is the function assigned to the additional small node on the level $p-1$. Explicitly, the function $D^*$ depends on the type, if we are in type $B_{2p}$ then $\oppInv{D} =\sau(D)$, otherwise $\oppInv{D} = D$.
    
    Now, we wish to prove that the set of functions where we replace $E$ with $E':=\Lambda_p(A^r,A_l)$ is compatible with the mutation in the polygonal cluster algebra, i.e. we need to show
    \begin{equation*}
        E'=CE^{-1}A+x\sigma(D^*)\sigma(E)^{-1}B\,,
    \end{equation*}
    or equivalently
    \begin{equation*}
        C^{-1}E' B^{-1} = E^{-1}AB^{-1} + xC^{-1}\sigma(D^*)\sigma(E)^{-1}\,.
    \end{equation*}
    This is exactly what we get if we apply $\Lambda_p$ to each term in earlier, which completes the proof in the noncommutative rank 1 case.

    In higher noncommutative rank, we would have to check the formula for adding potentials agrees with performing a square move. This was done in \cite[Theorem 6.5]{goncharov2021spectral}. 
\end{proof}

\begin{proof}[Proof (of \Cref{thm:clusterChartsTuplesRelatedByMutationTheta}).]
    For any single flip we first rotate so the quivers are in standard form as in \Cref{fig:flipSetup}.
    
    \input{Figures/figureFlipSetup}
    
    Then the flip is realized by transforming the word $r_1r_2$ to $r_2r_1$ (where $r_2$ is written in the generators $s_{\overline{i}}$). By \Cref{thm:twistedTwoMoveTheta} this is realized by a sequence of cluster mutations, see \Cref{fig:flip} for the case $B_2$.
\end{proof}

\input{Figures/figureFlip}

We extend the twisted cyclic shift map to $n$-tuples as follows:
\begin{definition}
    The \keyword{twisted cyclic shift} map $\tw:\Confx_n(\decFlag)\to\Confx_n(\decFlag)$ for $n\geq 3$ is given by
    \begin{equation*}
        \tw:(A_1,A_2,A_3,\dots,A_n)\mapsto (s_GA_2,A_3,\dots,A_n,A_1)\,.
    \end{equation*}
\end{definition}
\begin{corollary}\label{thm:TupleTwPositiveMapTheta}
    Let $c\in\Confx_n(\decFlagTheta)$. The cluster charts for $c$ and $\tw(c)$ are related by generalized mutations.
\end{corollary}
\begin{proof}
    A chart for the twisted cyclic shift can be obtained by rotating a triangulation for the original tuple. This will result in the edges incident to $A_2$ having the opposite orientation from the standard total order previously used in the construction. However since $A_2$ has been replaced with $s_G A_2$, we know the functions on these edges differ by switches/weaves of the polygonal cluster chart. Applying these switches and weaves results in a chart for the original tuple.
\end{proof}

\section{Noncommutative rank 1 detailed examples}\label{sec:NoncomRank1Examples}
We will now explain the detailed calculations missing in the preceding proofs case by case.
\subsection{\texorpdfstring{Type $A_1$}{Type A1}}\label{sec:typeA1}
The abstract construction becomes much simpler in Jordan split type $A_1$. In fact, we recover the coordinates on the space of positive configurations of decorated flags which were constructed for $G\simeq \SP(2n,\R)$ in \cite{alessandrini2019noncommutative}. The polygonal cluster algebras reduce to the noncommutative surface cluster algebra of Berenstein and Retakh \cite{berenstein2018noncommutative}.

\subsubsection{\texorpdfstring{Jordan}{Jordan}-pinning}\label{sec:A1pinningExamples}
Let $G$ be a Jordan-split group over $A_1$. In \Cref{sec:A1gradedGroups} we defined a Jordan pinning of $G$ to consist of the following three maps:
\begin{equation*}
    x : V \to G, \hspace{2pc}  y : V \to G, \hspace{2pc} \check{\beta}: \Gamma(\jordan{J}) \to G \,.
\end{equation*} 
We assume now that the Jordan pinning provides an isomorphism of groups $G\cong\SLJ$, so that the weight $\Lambda_\beta:G_0 \to \mathcal{R}$ is well-defined as the dual to $\check{\beta}$.

\begin{example}
Consider the Jordan algebra $\jordan{J}=H(\mathcal{R},\sigma)$ for $\mathcal{R}$ a ring with involution $\sigma$. Then $V = \mathcal{R}^\sigma$ and  $\tilde{\Gamma}(\jordan{J}) = \mathcal{R}^\times$. As both $\mathcal{R}^\sigma$ and $\mathcal{R}^\times$ naturally sit inside $\mathcal{R}$, we can represent $G$ as a subgroup $\SP_2(\mathcal{R},\sigma)$ inside $\Mat_2(\mathcal{R})$ as defined in \cite{alessandrini2019noncommutative}. The pinning maps become for $B,C \in \mathcal{R}^\sigma$ and $A \in \mathcal{R}^\times$:
    \[ x(B) = \begin{bmatrix}
        1 & B \\
        & 1
    \end{bmatrix} \hspace{2pc} y(C) = \begin{bmatrix}
        1 & \\
        C & 1
    \end{bmatrix} \hspace{2pc} \check{\beta}(A) = \begin{bmatrix}
        A & \\ &\sigma(A)^{-1}
    \end{bmatrix}.\]
    The weight $\Lambda_\beta$ is nothing more than the upper left entry of the matrix. 
    This mirrors the split-real $A_1$ story exactly, except the matrix entries can belong to a noncommutative ring. 
In this case each $\Theta$-flag can be identified with a decorated Lagrangian as follows:
\begin{equation*}
    \begin{bmatrix}
        A & B \\ C & D
    \end{bmatrix}U_\Theta \mapsto \begin{bmatrix}
        A \\ C
    \end{bmatrix}.
\end{equation*} 
    
\end{example}

\subsubsection{Pairs of flags}
Recall that the invariant of a pair of decorated $\Theta$-flags is called the $L$-distance. We now unwind the properties of the $L$-distance from \Cref{sec:pairsFlagsTheta} in this case. 
\begin{lemma}\label{thm:LDistPropA1}
    We have for any two transverse flags $A_1, A_2$:
    \begin{align*}
        l(A_2,A_1) =& -\sigma(l(A_1,A_2))\\
        l(s_G A_2, A_1) =& \sigma(l(A_1,A_2)) \,.
    \end{align*}
\end{lemma}

\begin{example}
    We return to our running example of $\jordan{J} = H(\mathcal{R},\sigma)$. In this case the $L$-distance is closely related to the symplectic $\Lambda$-length defined in  \cite{alessandrini2019noncommutative}. The  \keyword{symplectic $\Lambda$-length} is a $G$-invariant function of a pair of transverse decorated Lagrangians given by the formula $\Lambda([A,C],[B,D]) = \sigma(A)D-\sigma(C)B$. 
\end{example}

\begin{lemma}
    Let $(A_1,A_2)$ be a pair of transverse decorated flags. Then $\Lambda_\beta(l(A_1,A_2))$ is equal to the symplectic $\Lambda$-length.  
\end{lemma}
\begin{proof}
    Let $l = \check{\beta}(A)$ be a generic element of the Levi subgroup. Then since both functions are $G$-invariant it suffices to check the relation on any representative $(U_\Theta, \overline{w}_0 U_\Theta l)$:
    \[\left(\begin{bmatrix}
        1 & 0\\
        0 & 1
    \end{bmatrix}, \begin{bmatrix} 0 & -1\\1 & 0\end{bmatrix}\begin{bmatrix}1 & 0 \\ 0 & 1\end{bmatrix}\begin{bmatrix}
        A & 0\\ 0 & \sigma(A)^{-1}
    \end{bmatrix}\right) = \left(\begin{bmatrix}
        1 & 0\\
        0 & 1
    \end{bmatrix}, \begin{bmatrix}
        0 & -\sigma(A)^{-1}\\
        A & 0
    \end{bmatrix} 
    \right) \mapsto \left(\begin{bmatrix}
        1 \\
        0 
    \end{bmatrix}, \begin{bmatrix}
        0 \\
        A 
    \end{bmatrix} 
    \right) \,.\]
    By construction $\Lambda_\beta(l(A_1,A_2)) = \Lambda_\beta(\check{\beta}(A)) = A$. The symplectic $\Lambda$-length is $\sigma(1)\cdot A - \sigma(0)\cdot 0$, which agrees. 
\end{proof}

\begin{example}
    Let us consider another explicit example, $G=\Spin(2,n)$. We explicitly fix
    \begin{equation*}
        Q = \begin{bmatrix}
            & & -1 \\
            & J & \\
            -1 & & 
        \end{bmatrix}\hspace{1cm}\text{with }
        J= \begin{bmatrix}
            & & 1 \\
            & -\Id_{n-2} & \\
            1 & & 
        \end{bmatrix}
    \end{equation*}
    with the induced bilinear form
    \begin{equation*}
        b(v,w) = v^T Q w
    \end{equation*}
    and quadratic form $q(v):=b(v,v)$ on $V=\R^{n+2}$. This form has signature $(2,n)$ so that $V\cong\R^{2,n}$. Let $e_i$ denote the standard basis vectors and let $V':=\langle e_2,\dots,e_{n+1}\rangle$ be the subspace on which the restriction of $b$ is induced by $J$. Now $\Spin(2,n)$ can be realized in the Clifford algebra $\CL(V)$ over $V$.
    
    $G$ can be understood as a group of type $A_1$ over the Jordan algebra $\jordan{J}=\CL(1,n-1)$ of Clifford type. We have a Jordan pinning given by
    \begin{equation}\label{eq:thetaPinningSO2n}
        \begin{array}{rrl}
            x: & V'\to G\,,& v\mapsto 1+e_1v/\sqrt{2}\\
            y: & V'\to G\,, & v\mapsto 1+e_{n+2}Jv/\sqrt{2}\\
            \check{\beta}: & \Gamma^0(V')\to G\,, & t\mapsto\frac{1}{N(t)}\big(1+(1-N(t))e_1e_{n+2}/2\big)t\,.
        \end{array}
    \end{equation}
    Since $J$ has signature $(1,n-1)$, the even Clifford group $\Gamma^0(V')\cong\Gamma^0(1,n-1)$ can be identified with the structure group $\Gamma(\jordan{J})$ and embeds naturally into $\CL(V)$.
    
    $G$ admits a 2:1 projection map
    \begin{equation*}
        \rho: \Spin(2,n)\twoheadrightarrow\SO_0(2,n)=:G'
    \end{equation*}
    given by the twisted adjoint representation \cite{lawson2016spin}. If we compose the Jordan pinning maps with $\rho$, we obtain
    \begin{equation*}
        \begin{array}{rrl}
            x'=\rho\circ x: & V'\to G'\,,& v\mapsto \Bigg[
            \begin{smallmatrix}
                1 & \sqrt{2}v^TJ & q(v)\\
                & \Id & \sqrt{2}v\\
                & & 1
            \end{smallmatrix}\Bigg]\\[15pt]
            y'=\rho\circ y: & V'\to G'\,, & v\mapsto \Bigg[
            \begin{smallmatrix}
                1 & & \\
                \sqrt{2}Jv & \Id & \\
                q(v) & \sqrt{2}v^T & 1
            \end{smallmatrix}\Bigg]\\[15pt]
            \check{\beta}'=\rho\circ\check{\beta}: & \Gamma^0(V')\to G'\,, & t\mapsto\Bigg[
            \begin{smallmatrix}
                N(t) & & \\
                & \tilde{\mathrm{Ad}}_t & \\
                & & N(t)^{-1}
            \end{smallmatrix}\Bigg]\,,
        \end{array}
    \end{equation*}
    where $\tilde{\mathrm{Ad}}:\Gamma^0(V')\to\SO(V')$ denotes the usual projection. Observe that $\check{\beta}'$ is a 2:1 map, while $\check{\beta}$ is injective.

    We can use the composite maps to understand that the maps in \Cref{eq:thetaPinningSO2n} actually provide a Jordan pinning, i.e. an isomorphism $G\cong\SLJ$. For this we consider the Lie algebra $\mathfrak{g}$ of $G$, which is isomorphic to the Lie algebra of $G'$ and prove that the maps induce an isomorphism $\mathfrak{g}\cong\slj$. First, we note that
    \begin{equation*}
        x'(v)=\exp\Bigg[\begin{smallmatrix}
            0 & \sqrt{2}v^TJ & 0\\
            & 0 & \sqrt{2}v \\
            & & 0
        \end{smallmatrix}\Bigg]\,,\hspace{1cm}
        y'(w)=\exp\Bigg[\begin{smallmatrix}
            0 & & \\
            \sqrt{2}Jw & 0 & \\
            0 & \sqrt{2}w^T & 0 
        \end{smallmatrix}\Bigg]\,.
    \end{equation*}
    Now the Lie bracket in $\mathfrak{g}$ is readily computed to be
    \begin{equation*}
        \left[\rule{0cm}{7mm}\Bigg[\begin{smallmatrix}
            0 & \sqrt{2}v^TJ & 0\\
            & 0 & \sqrt{2}v \\
            & & 0
        \end{smallmatrix}\Bigg],
        \Bigg[\begin{smallmatrix}
            0 & & \\
            \sqrt{2}Jw & 0 & \\
            0 & \sqrt{2}w^T & 0 
        \end{smallmatrix}\Bigg]\right]
        = \Bigg[\begin{smallmatrix}
            2b(v,Jw) & & \\
            & 2(vw^TJ-Jwv^T)J & \\
            & & -2b(v,w) 
        \end{smallmatrix}\Bigg]\,.
    \end{equation*}
    We need to prove that this corresponds to $\iota(v,w)$ under the identification of block diagonal elements in $\mathfrak{g}$ with $\mathfrak{g}_\jordan{J}$. Recall that the Clifford algebra $\CL(1,n)$ over $V'$ becomes a quadratic Jordan algebra by fixing an element $\Id=e\in V'$ with $q(e)=1$ and setting
    \begin{equation*}
        \iota(v)(w)=vewev\,.
    \end{equation*}
    Here we choose $e=(1,0,\dots,0,1)/\sqrt{2}$. Now we compute in the Clifford algebra:
    \begin{align*}
        \iota(v,w)(u) &= \big(\iota(x+u)-\iota(x)-\iota(u)\big)(w) = v\overbrace{ewe}^{w^*}u + uewev\\
        &= vw^*u + \big(2b(u,w^*)-w^*u\big)v = 2b(u,w^*)v + vw^*u-w^*\big(2b(v,u)-vu\big)\\
        &= 2\big(v b(w^*,u)-w^*b(v,u)+b(v,w^*)u\big) = 2\big(b(v,w*)+v(w^*)^TJ-w^*v^TJ\big)u\,.
    \end{align*}
    Observe that $w^*=ewe=2b(e,w)e-w=(2ee^TJ-\Id)w=Jw$. Thus, $\iota(v,w)\in\mathfrak{g}_\jordan{J}$ is in fact identified with the Lie bracket we computed before. The other relations are easy to check.

    \begin{remark}
        In \cite{rogozinnikov2025, rogozinnikov2020symplectic}, the group $\Spin(2,n)$ is identified with $\SP_2(\CL(1,n-1),\sigma)$ with the anti-involution $\sigma$ given by composing Clifford algebra transposition with conjugation by $e$.
    \end{remark}

    We now turn our attention to the decorated flag variety $\decFlagTheta=G/U_\Theta$ for $U_\Theta=x(V')$. Of course our general theory applies to it, but to make everything as explicit as possible, we consider the decorated flag variety $\decFlagTheta'=G'/U_\Theta'$ for $\SO_0(2,n)$ instead. It is double covered by $\decFlagTheta$ via the map induced by $\rho$ and can be identified with a subset of
    \begin{equation*}
        \left\{(v,B)\st v\in V \text{ with } q(v)=0\,, B \text{ basis of } \langle v\rangle^\perp/\langle v\rangle \text{ adapted to } J\right\}\,
    \end{equation*}
    via
    \begin{equation*}
        \SO_0(2,n)\ni x=\Big[v\; \tilde{B}\; v'\Big]\mapsto (v,B)\,,
    \end{equation*}
    i.e. $v$ is the first column of the matrix $x$, $\langle v\rangle^\perp/\langle v\rangle$ can be identified with the column span of the $(n+2)\times n$-matrix $\tilde{B}$, which then provides the basis $B$ through its columns. This is adapted to $J$, i.e. in this basis the restriction of $b$ to this space is induced by the matrix $J$.
    \begin{remark}
        One can explicitly describe the subset of $\{(v,B)\}$ that is identified with the decorated flag variety using conditions on the (time) orientation.
    \end{remark}

    Now let $F_1=(v_1,B_1)$ and $F_2=(v_2,B_2)$ be two decorated flags. We want to describe the configuration $[F_1,F_2]\in\Conf_2(\decFlagTheta')$ geometrically. First, we observe that we have
    \begin{itemize}
        \item $b(v_1,v_2)\neq 0$ and
        \item $\dim(\langle v_1\rangle^\perp \cap \langle v_2\rangle^\perp)=n$\,.
    \end{itemize}
    To see this, we first see that this is true for the standard decorated flags $[\Id]=\big(e_1,(e_2,\dots,e_{n+1})\big)$ and $[\overline{w_0}]=\big(e_{n+2},(-e_{n+1},e_3,\dots,e_n,-e_2)\big)$ (from our Jordan pinning we have $\overline{w_0}=-Q$).

    Thus, our first invariant of the configuration $[F_1,F_2]$ is simply $b(v_1,v_2)$. To use the rest of the data, set
    \begin{equation*}
        W:=\langle v_1\rangle^\perp\cap\langle v_2\rangle^\perp
    \end{equation*}
    and observe that the projections $W\to \langle v_i\rangle^\perp/\langle v_i\rangle$ are isomorphisms if $F_1$ and $F_2$ are transverse: Let $x\in\langle v_1\rangle^\perp$. Then $x-\frac{b(x,v_2)}{b(v_1,v_2)}v_1\in\langle v_2\rangle^\perp$ represents the same element in $\langle v_1\rangle^\perp/\langle v_1\rangle$, showing that the projection $W\to\langle v_1\rangle^\perp/\langle v_1\rangle$ is surjective. Because the spaces have the same dimension by transversality, it must in fact be an isomorphism. The same argument works for the other projection.

    This allows the construction of a second, noncommutative invariant of the configuration: The bases $B_i$ induce bases $\bar{B}_i$ of $W$, both of which are adapted to $J$. We can pair the basis vectors using $b$, to obtain an $n\times n$-matrix in $\SO(J)=\SO(1,n-1)$. If we view the bases $\bar{B}_i$ as $(n+2)\times n$-matrices (columns are the basis vectors), the invariant can be computed as $\bar{B}_1^TQ\bar{B}_2$.

    If $F_1=[\Id]$ is the standard decorated flag and
    \begin{equation*}
        F_2=\left[\rule{0cm}{5mm}\overline{w_0}\begin{bsmallmatrix}
            a & & \\
            & M & \\
            & & a_{-1}
        \end{bsmallmatrix}
        \right]
        =\left[\rule{0cm}{5mm}\begin{bsmallmatrix}
            & & a^{-1}\\
            & -JM & \\
            a & &
        \end{bsmallmatrix}
        \right]
    \end{equation*}
    we get
    \begin{equation*}
        \langle v_1\rangle^\perp \cap \langle v_2\rangle^\perp = \langle e_2,\dots,e_{n+1}\rangle = V'
    \end{equation*}
    and the bases are represented as
    \begin{equation*}
        \bar{B}_1=\Id_n\,,\qquad \bar{B}_2=-JM
    \end{equation*}
    so that the invariants are $b(v_1,v_2)=b(e_1,ae_{n+1})=-a$ and $b(\bar{B}_1,\bar{B}_2)=-M$.

    For two arbitrary transverse flags $F_1,F_2$ as before, assume that they are represented by the matrices $g_i=\Big[v_i\;\tilde{B}_i\;w_i\Big]$, where $\tilde{B}_i$ is an arbitrary lift of $B_i$ to $V$. To compute the second invariant from this form, we need to compute the lifts $\bar{B}_i$ of $B_i$ to $\langle v_1\rangle^\perp\cap\langle v_2\rangle^\perp$. Explicitly, these are given by
    \begin{equation*}
        \bar{B}_1 = \tilde{B}_1 - v_1\cdot\frac{v_2^TQ\tilde{B}_1}{b(v_1,v_2)}\,,\qquad\bar{B}_2 = \tilde{B}_2 - v_2\cdot\frac{v_1^TQ\tilde{B}_2}{b(v_1,v_2)}\,.
    \end{equation*}
    Thus the invariant is
    \begin{equation*}
        \bar{B}_1^TQ\bar{B}_2 = \tilde{B}_1^TQ\left(\Id - \frac{v_2v_1^TQ}{b(v_1,v_2)}\right) \tilde{B}_2\,.
    \end{equation*}

    On the other hand, the $L$-distance between the two decorated flags is given by
    \begin{equation*}
        l(g_1U_\Theta,g_2U_\Theta) = l\left(U_\Theta,\overline{w_0}\doverline{w_0}g_1^{-1}g_2U_\Theta\right) = \left[\doverline{w_0}g_1^{-1}g_2\right]_0
    \end{equation*}
    using the Gauss decomposition. Since $g_1^{-1}=Qg_1^TQ$ and $\doverline{w_0}=\overline{w_0}=-Q$ in $\SO_0(2,n)$ (this is not true in $\Spin(2,n)$!) we get
    \begin{equation*}
        \doverline{w_0}g_1^{-1}g_2 = -g_1^TQ g_2 =-
        \begin{bsmallmatrix}
            v_1^T\\ \tilde{B}_1^T\\ w_1^T
        \end{bsmallmatrix}
        \begin{bsmallmatrix}
            Qv_2 & Q\tilde{B}_2  & Qw_2
        \end{bsmallmatrix}=-
        \begin{bsmallmatrix}
            b(v_1,v_2) & v_1^TQ\tilde{B}_2 & b(v_1,w_2)\\
            \tilde{B}_1^TQv_2 & \tilde{B}_1^TQ\tilde{B}_2 & \tilde{B}_1^TQw_2\\
            b(w_1,v_2) & w_1^TQ\tilde{B}_2 & b(w_1,w_2)
        \end{bsmallmatrix}\,.
    \end{equation*}
    By computing the Gauss decomposition, we obtain
    \begin{equation*}
        l(g_1U_\Theta,g_2U_\Theta) = \begin{bsmallmatrix}
            -b(v_1,v_2) & &\\
            & \tilde{B}_1^TQv_2v_1^TQ\tilde{B}_2/b(v_1,v_2)-\tilde{B}_1^TQ\tilde{B}_2 & \\
            & & -1/b(v_1,v_2)
        \end{bsmallmatrix}\,.
    \end{equation*}
    Thus the geometric invariants are equivalent to the $L$-distance.
\end{example}

\subsubsection{Triples of flags}\label{sec:TriplesFlagsThetaA1}
    The theory for triples of transverse flags $(A_1,A_2,A_3)$, also dramatically simplifies in this case. 
    In Jordan split type $A_1$ there is a unique reduced expression of the longest word in the Weyl group, $w_0 = \varsigma_\beta$. Moreover, the interpolating sequence of flags is trivial as $(A_1,A_2,A_3)$ is already an elementary configuration of type $\varsigma_\beta$. 
    Therefore the cluster chart for a triple of decorated $\Theta$-flags is given by 
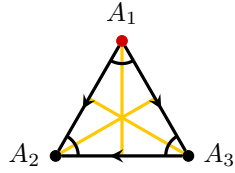
\begin{center}
    \begin{tikzpicture}
        \node[regular polygon,draw,regular polygon sides=3, minimum size=2cm] (p) at (0,0) {};
		
        \draw[color=Yellow, very thick] (p.side 2) to (p.corner 1);
        \draw[color=Yellow, very thick] (p.side 3) to (p.corner 2);
        \draw[color=Yellow, very thick] (p.side 1) to (p.corner 3);

        \centerarc[straightangle, very thick](p.corner 1)(240:300:.3)
        \centerarc[straightangle,very thick](p.corner 2)(0:60:.3)
        \centerarc[straightangle, very thick](p.corner 3)(120:180:.3)
            
        \node[flag, color = Red]  (F1)    at (p.corner 1)	[label = above: $A_1$]	{};
		\node[flag]               (F2)    at (p.corner 2)	[label = left: $A_2$]	{};
		\node[flag]               (F3)    at (p.corner 3)	[label = right: $A_3$]	{};

        \draw[qarrow, postaction = decorate] (F1) to (F2);
		\draw[qarrow, postaction = decorate] (F1) to (F3);
        \draw[qarrow, postaction = decorate] (F3) to (F2);
    \end{tikzpicture}
\end{center}

\noindent Let $M_{ij} = \Lambda_\beta(l(A_i,A_j))$ be the cluster coordinate associated to the directed side $(i,j)$. 
\begin{theorem}
    The potential of the configuration $(A_1, A_2,A_3)$ is given by the angle 
    \begin{equation}
        \Delta_g(M_{23}) = \sigma(M_{12})^{-1}M_{23} M_{13}^{-1}\,.
    \end{equation}
\end{theorem}
\begin{proof}
    The configuration $(A_1,A_2,A_3)$ is equivalent to $(U_\Theta,\overline{\omega_0}U_\Theta M_{12},x(v)\overline{\omega_0}U_\Theta M_{13})$ where $v$ is the potential.
    We then compute $M_{23} = \Lambda_\beta(l(A_2,A_3)) = \sigma(M_{12})vM_{13}$ and the formula follows.
\end{proof}
\begin{remark}
    In this case, the angle in the polygonal cluster algebra is the noncommutative angle of a noncommutative surface cluster algebra. The fact that the angle belongs to the vector space of the Jordan algebra is equivalent to the angle being fixed by $\sigma$ for $\jordan{J} = H(\mathcal{R},\sigma)$.
\end{remark}

    The final step to understanding a triple of flags is understanding twisted cyclic rotation. This case is very simple, as we observe that the only change to this seed under rotation is the orientation of the arcs bounding the triangular tile. As seeds of a polygonal cluster algebra this corresponds to performing a switch at these arcs. This changes the cluster coordinate by applying $\sigma$, which has, by \Cref{thm:LDistPropA1}, exactly the same effect as computing $\Lambda_\beta(l(s_G A_j, A_i))$ instead of $\Lambda_\beta(l(A_i,A_j))$ as needed. 

\subsubsection{Tuples of flags}
    In \Cref{sec:TuplesOfThetaFlags}, we saw that flipping an arc of the triangulation amounted to a series of moves on reduced expressions of the longest word in $W(R) \times W(R)$. The new move for $\Theta$-flags was the twisted two move at the noncommutative node. In Jordan split type $A_1$, this move suffices to flip the arc. 

    The proof of the twisted two move in \Cref{thm:twistedTwoMoveTheta} relied on the additivity of the potential, which we saw in \Cref{sec:TriplesFlagsThetaA1} is equivalent to the noncommutative surface angles. Consider a four tuple $(A_1,A_2,A_3,A_4)$ of pairwise transverse flags. The additivity of the potentials implies that 
    \begin{equation}\label{eq:sumofangles}
        \sigma(M_{14}^{-1})M_{42}M_{12}^{-1}=\sigma(M_{13}^{-1})M_{32}M_{12}^{-1}+\sigma(M_{14}^{-1})M_{43}M_{13}^{-1} \,.
    \end{equation}
    This in turn implies the following non commutative exchange relation: 
    \begin{equation}
        M_{42} = \sigma(M_{14})\sigma(M_{13}^{-1})M_{32}+ M_{43}M_{13}^{-1}M_{12} \,.
    \end{equation}
    \begin{remark}
        This is exactly the exchange relation on a noncommutative surface given in \cite{berenstein2018noncommutative}.
    \end{remark}

Our next step is to study the case of Jordan split type $B_p$, where more of the subtleties of the general construction appear.

\subsection{\texorpdfstring{Type $B_p$}{Type Bp}}

We now consider $G$ a group of Jordan split type $B_p$.

A root system of type $B_p$ contains $2p$ short roots and $2p(p-1)$ long roots. We chose a set of simple roots for the  $\Theta$-root system, $\{\beta_1, \cdots, \beta_{p-1},\beta_{p}\}$ with standard indexing coming from the Dynkin diagram so that $\beta_p$ is the unique short root. 

The long roots correspond to commutative Jordan algebras, and are therefore associated with a finite extension of $\K$ we will denote by $\K'$.  The short root spaces, $\{\pm\beta_p,\pm(\beta_p+\beta_{p-1}), \dots, \pm\sum_i\beta_i\}$ can be associated to a degree 2 Jordan algebra of Clifford type, $\jordan{J}(V,b,e)$ for $V$ a $\K'$ vector space, $b$ a bilinear form, and $e$ an element of norm $1$. 

The structure group is the even Clifford group $\Gamma^0(V,b)$. We recall that the map embedding the vector space into the structure group is given by $\iota(v) = ve$. 

The vector space $V$ comes with its adjugate involution which we denote by $\tau$. Since our Jordan algebra is of degree 2, this map is a linear map. Similarly, the generalized structure group comes with the two anti-involutions $\sigma, \tau$. These are given explicitly by the standard transposition map on the Clifford algebra
\begin{equation*}
    \tau(x) = x^T\,,\quad\mathrm{and}\quad \sigma(x)= ex^Te\,.
\end{equation*}
Since $\iota(\tau(v)) = \tau (\iota(v))$ we are justified in using the same name for both anti-involutions. Moreover we find that 
$L_\Theta$ is a central extension of $\Gamma^0(V,b)$
and we extend $\sigma$ and $\tau$ to anti-involutions on $L_\Theta$ by having them act trivially on the center. By further abuse of notation, we denote the resulting maps by $\sigma$ and $\tau$ again.

\subsubsection{\texorpdfstring{$\Theta$}{Theta}-pinning}\label{sec:BpThetaPinning}
    Choose a Jordan pinning for $G$ as discussed in \Cref{sec:ThetaPinning}. Since the Jordan algebra is always of Clifford type, the pinning maps have the following form: For $1\leq i<p$
    \begin{equation*}
        x_i \colon \K'\rightarrow G \hspace{2pc }y_i \colon \K'\rightarrow G \hspace{2pc} \check{\beta}_i \colon \K'^\times\rightarrow G \,,
    \end{equation*}
    and
    \begin{equation*}
        x_p \colon V\rightarrow G \hspace{2pc }y_p\colon V\rightarrow G \hspace{2pc} \check{\beta}_p \colon \Gamma^0(V,b)\rightarrow G \,.
    \end{equation*}
    We recall that the pinning provides lifts of the grading Weyl group to $G$ via:
    \begin{equation*}
        \overline{\varsigma}_i = y_i(1)x_i(-1)y_i(1) \hspace{2pc} \overline{\varsigma}_p = y_p(\Id)x_p(-\Id)y_p(\Id)\,.
    \end{equation*}
    We will use these to compute some identities for the Jordan pinning in this concrete case. 
    
\subsubsection{Weyl group}
    There are two nice reduced expressions of the longest word we will use in our calculations. The first is given by $p$ repetitions of the Coxeter element $\varsigma_1 \varsigma_2 \cdots \varsigma_p$. The second word is inductive and consists of multiplying the $p$ distinct \keyword{mountain blocks}, 
    \begin{equation*}
        M_k = \varsigma_k \varsigma_{k+1}\cdots \varsigma_{p-1}\varsigma_p\varsigma_{p-1}\cdots \varsigma_{k+1}\varsigma_k \,.
    \end{equation*}

    \noindent The mountain blocks provide a nontrivial action on the special root space $\mathfrak{u}_p$.
    \begin{lemma}\label{thm:MountainAction}
        For $k< p$, the action of $M_k$ on $\mathfrak{u}_p$ is given by 
        \begin{equation}
            M_k \cdot x_p(v) = x_p(\tau(v)) \,.
        \end{equation}
        Similarly $M_k$ acts on the short coroot by 
        \begin{equation}
            M_k \cdot \check{\beta}_p(A) = \check{\beta}_p(\sau(A)) \,.
        \end{equation}
    \end{lemma}
    \begin{proof}
        We begin by calculating the action of $M_{p-1}$ on $\check{\beta}_p$ using \Cref{thm:ActionOfSigmaPonSpecialRoot,thm:ActionOfSigmaPonOtherRoots}.
        \begin{align*}
            M_{p-1}\cdot \check{\beta}_p(A) =~& \varsigma_{p-1}\varsigma_p\varsigma_{p-1} \cdot \check{\beta}_p(A) = \varsigma_{p-1}\varsigma_p \cdot \check{\beta}_p(A)\check{\beta}_{p-1}(N(A)) \\
            =~& \varsigma_{p-1} \check{\beta}_p(\sigma(A)^{-1})\check{\beta}_{p-1}(N(A))\check{\beta}_p(N(A)1) \\
            =~& \check{\beta}_p(\sigma(A)^{-1})\check{\beta}_{p-1}(N(A)^{-1})\check{\beta}_{p-1}(N(A)^{-1})\check{\beta}_p(N(A)1)\check{\beta}_{p-1}(N(A)^2) \\
            =~& \check{\beta}_p(\sigma(A^{-1})N(a)) = \check{\beta}_p(\sau(A))\,.
        \end{align*}
        For $k< p-1$,  $M_k = wM_{p-1}w$ for $w$ a word consisting of reflections over roots with index less than or equal to $p-2$. As such $w$ acts trivially on the image of $\check{\beta}_p$ and thus arbitrary $M_k$ act as claimed.

        We now consider the action of $M_k$ on $x_p(v)$. Let $w = \varsigma_{p-1}\varsigma_{p-2}\cdots\varsigma_k$. Then $M_k = w^{-1}\varsigma_p w$. Inductively we observe $w$ moves $\mathfrak{u}_{\beta_p}$ to another short $\Theta$-root space which is perpendicular to $\mathfrak{u}_{\beta_p}$. The action of $\varsigma_p$ does not change the image root space and thus the action $w^{-1}$ returns us to $\mathfrak{u}_{\beta_p}$. As the pinning map is surjective, there is some $w \in V$ such that $x_p(w) = M_k\cdot x_p(v)$. We then compute the action on $x_p(v) \overline{\varsigma}_p$ using \Cref{thm:DecomposeXSigma}:
        \begin{align*}
         M_k \cdot (x_p(v) \overline{\varsigma}_p) =~& M_k\cdot(y(v^{-1}) \check{\beta}_p(\iota(v)) x_p(v)) = (M_k \cdot y(v^{-1})) (M_k \cdot \check{\beta}_p(\iota(v))) (M_k \cdot x_p(v)) \\
         =~& y(w_1^{-1}) \check{\beta}_p(\sau\iota(v)) x_p(w_2) \,.
        \end{align*}
        On the other hand $M_k \cdot \overline{\varsigma}_p = \overline{\varsigma}_p$. So 
        \begin{equation}
             M_k \cdot (x_p(v) \overline{\varsigma}_p) = x(w)\overline{\varsigma}_p = y(w^{-1}) \check{\beta}_p(\iota(w)) x_p(w) \,.
        \end{equation}
        The Gauss decomposition of elements is unique, so we equate the argument of the coroot. This implies $\iota(w) = \sau(\iota(v)) = \iota(\tau(v))$. Since $\iota$ is injective we then have $w = \tau(v)$ as claimed.
    \end{proof}

    We now calculate the opposition involution, as defined in \Cref{def:OppositionInvolution}.
    \begin{corollary}\label{thm:W0ActionBp}
        The action of $w_0$ on $\check{\beta}_p$ is given by
        \begin{equation*}
            w_0\cdot \check{\beta}_p(A) = \begin{cases}
                \check{\beta}_p\big(\tau(A^{-1})\big) & p \text{ even} \\
                \check{\beta}_p\big(\sigma(A^{-1})\big) & p \text{ odd}  \,.
            \end{cases}
        \end{equation*}
        Consequently, the opposition involution on $L_\Theta$ is
        \begin{equation*}
            \begin{array}{rrcl}
                \oppInv{(.)}: & \check{\beta}_i(x) &\mapsto &\check{\beta}_i(x)\\
                & \check{\beta}_p(A) &\mapsto &\begin{cases}
                    \check{\beta}_p\big(\tau(A)\big) & p \text{ even} \\
                    \check{\beta}_p\big(\sigma(A)\big) & p \text{ odd}  \,.
                \end{cases}
                 
            \end{array}
        \end{equation*}
    \end{corollary}
    \begin{proof}
        Write $w_0 = M_1\cdots M_p$ using the inductive word. Then since $M_p = \varsigma_p$, the result follows from applying the action of each subword computed in \Cref{thm:ActionOfSigmaPonSpecialRoot,thm:MountainAction}. 
    \end{proof}

    The induced involution on the $\Theta$-root space $\mathfrak{u}_{\beta_p}$ is given by:
    \begin{corollary}\label{thm:potentialInvolutionBp}
        The involution on the potential is given by $\oppInv{v} = \tau^{p-1}(v) = \begin{cases}
            \tau(v) & p \text{ even}\\
            v & p \text{ odd} \,.
        \end{cases}$
    \end{corollary}
    \begin{proof}
        We simply calculate:
        \begin{align*}
            \check{\beta}_p\big(\iota(\oppInv{v})\big) = \oppInv{\Big(\check{\beta}_p\big(\iota(v)\big)\Big)}
            =\oppInv{\Big(\check{\beta}_p(ve)\Big)} = \begin{cases}
                \check{\beta}_p(ev) = \check{\beta}_p(evee) = \check{\beta}_p\Big(\iota\big(\tau(v)\big)\Big) & p \text{ even}\\
                \check{\beta}_p(ve) = \check{\beta}_p\big(\iota(v)\big) & p \text{ odd}\,.
            \end{cases}
        \end{align*}
    \end{proof}

    \subsubsection{\texorpdfstring{$\Spin(p,q)$}{Spin(p,q)} with its positive pinning}

    Now we consider the real group $G\simeq\Spin(p,q)$. We can describe an explicit positive Jordan pinning using the projection $\rho:G\to G'\simeq\SO(p,q)^\circ$. Since the groups are locally isomorphic, the differential $\mathrm{d}_1\rho:\mathfrak{g}\to\mathfrak{g}'$ provides an isomorphism of Lie algebras. Moreover, $\rho$ identifies the unipotent subgroups $U_\Theta\cong U_\Theta'$, $U_\Theta^\opp\cong (U_\Theta')^\opp$ of $G$ and $G'$, and provides a 2:1-map $\rho:L_\Theta\twoheadrightarrow L_\Theta'$ of the Levi subgroups.

        To describe the subgroups of $G'$ explicitly, we choose a matrix representation of the group: We take a $(p+1,q+p+1)$-form on $\R^{2p+q+2}$ defined by $b(v,w)=v^TQw$ with
        \begin{equation*}
		  Q:=
		  \begin{bmatrix}
			0 & 0 & K\\
			0 & J & 0\\
			K^T & 0 & 0
		  \end{bmatrix}\,,
	    \end{equation*}
	    where $K\in\GL(p,\R)$, and $J\in\GL(q+2,\R)$ are given by 
	    \begin{equation} \label{eqn:R1qForms}
		  K=\left[\begin{smallmatrix}
                & & & (-1)^p\\
                & & \iddots & \\
                & 1 & & \\
                -1 & & & \\
            \end{smallmatrix}\right]
		  \hspace{5mm}\mathrm{and}\hspace{5mm}
		  J=\begin{bmatrix}
			0 & 0 & 1\\
			0 & -\Id_{q} & 0\\
			1 & 0 & 0
		  \end{bmatrix}\,,
	   \end{equation}
	   so that $G'\cong\mathrm{SO}(Q)^\circ$. It will be convenient to consider the quadratic form of signature $(1,q+1)$ on $V=\mathbb{R}^{q+2}$ defined as follows: For $v,w\in V$
	   \begin{align*}
		  b_J(v,w):=v^TJw\,,\\
		  N(v):=b_J(v,v)\,.
	   \end{align*}
       Of course this is just the restriction of $b$ to the subspace $V=\langle e_{p+1},\dots,e_{p+q+2}\rangle\cong\R^{1,q+1}$.

       Write for $a_1,\dots,a_p\in\R^\times$ and $M\in\SO(J)$
       \begin{equation*}
           l(a_1,\dots,a_p,M)\coloneq \left[
                \begin{smallmatrix}
                    a_1 & & & & & & \\
                    & \ddots & & & & & \\
                    & & a_p & & & & \\
                    & & & M & & & \\
                    & & & & a_p^{-1} & & \\
                    & & & & & \ddots & \\
                    & & & & & & a_1^{-1}
                \end{smallmatrix}\right]\in\SO(Q)\,.
       \end{equation*}
       Then the Levi subgroup of $G'$ is given explicitly by
       \begin{equation*}
           L_\Theta'=\scalebox{1.2}{\Bigg\{}l(a_1,\dots,a_p,M)\,\Bigg|\, a_1,\dots, a_p\in\R^\times,\, M\in\begin{cases}
                    \SO(J)^\circ & \prod a_i > 0\\
                    \SO(J)\setminus\SO(J)^\circ & \prod a_i<0
            \end{cases}\scalebox{1.2}{\Bigg\}}
       \end{equation*}
       where the relation between $\{a_i\}$ and $M$ comes from restricting to the identity component of $\SO(Q)$. The subgroup $U_\Theta'<G'$ consists of all matrices in $G'$ with ones on the diagonal and entries `above' $L_\Theta'$. The opposite unipotent subgroup $(U_\Theta')^\opp$ is obtained from $U_\Theta'$ by transposing the matrices.

       Consider the following maps: For $1\leq i<p$
       \begin{equation*}
       \begin{array}{rll}
            x_i':&\R\to U_\Theta'\,,  &t\mapsto \Id+t(E_{i,i+1}+E_{2p+q+2-i,2p+q+3-i})\\
            y_i':&\R\to (U_\Theta')^\opp\,, &t\mapsto \big(x_i'(t)\big)^T\\
            \check{\beta}_i':&\R^\times\to L_\Theta'\,, & t\mapsto l\Big(1,\dots,1,\underset{\color{Green}{i}}{t},\underset{\color{Green}{i+1}}{t^{-1}},1,\dots,1,\Id\Big)
       \end{array}
       \end{equation*}
       and
       \begin{equation*}
       \begin{array}{rll}
           x_p':& V\to U_\Theta'\,, &v\mapsto\left[\begin{smallmatrix}
               1 & & & & \\
               & \ddots & \sqrt{2}v^TJ & q_J(v) & \\
               & & \Id & \sqrt{2}v & \\
               & & & \ddots & \\
               & & & & 1
           \end{smallmatrix}\right]\\
           y_p':& V\to (U_\Theta')^\opp\,, &v\mapsto\big(x_p'(Jv)\big)^T\\
           \check{\beta}_p':& \Gamma^0(1,q+1)\to L_\Theta'\,, &x\mapsto l\Big(1,\dots,1,N(x),\tilde{\mathrm{Ad}}(x)\Big)
       \end{array}
       \end{equation*}
       where $\tilde{\mathrm{Ad}}:\Gamma^0(1,q+1)\to\SO(J)$ is the usual projection via the twisted adjoint representation, see e.g. \cite{lawson2016spin}. Observe that $\check{\beta}_p$ is a 2:1-map as $\check{\beta}_p(x)=\check{\beta}_p(-x)$, but that all the other maps are injective. Consequently, $L_\Theta$ decomposes as a direct product
       \begin{equation*}
           L_\Theta' = \bigoplus_{i=1}^p\mathrm{Im}\big(\check{\beta}_i\big)\overset{2:1}{\twoheadleftarrow} L_\Theta
       \end{equation*}
       where the double covering simply arises as the double covering of $\mathrm{Im}\big(\check{\beta}_p\big)$ by $\Gamma^0(1,q+1)$ which is a factor of $L_\Theta$. The maps above define a Jordan pinning of $G'$ which is obtained by postcomposing Jordan pinning maps of $G$ with the projection $\rho:G\to G'$.

       The Jordan pinning endows $V$ with the structure of a quadratic Jordan algebra as explained in \Cref{sec:jordan_algebras}, with the map $\iota:V\to\End(V)$ given by $v\mapsto ve$ where we can choose $e=(1,0,\dots,0,1)^T$ in the standard basis of $V$.

       As in the split case, the $\Theta$-fundamental weights are only defined for $G$, not for $G'$ where the map $\check{\beta}_p'$ is not injective. However, all the commutative fundamental weights still arise as pullbacks of characters of $L_\Theta'$ under $\rho:G\to G'$. Namely, $\Lambda_i^\Theta$ for $1\leq i<p$ is the pullback of
       \begin{equation*}
           L_\Theta'\to\R^\times\,,\quad l(a_1,\dots,a_p,M)\mapsto\prod_{k=1}^i a_k\,.
       \end{equation*}
       Accordingly, the extension of $\Lambda_i^\Theta$ to Gauss-decomposable elements in $G$ is simply the pullback of the $i$-th principal minor on $G'$.
       
       While the weight $\Lambda_p^\Theta:L_\Theta\to\Gamma^0(1,q+1)$ does not arise as a pullback, the following composition of maps does:
       \begin{equation*}
           L_\Theta\to\R^\times\times\Gamma^0(1,q+1)\,,\quad l\mapsto\Big(N\big(\Lambda_p^\Theta(l)\big),\tilde{\mathrm{Ad}}\big(\Lambda_p^\Theta(l)\big)\Big)\,.
       \end{equation*}
       This can be seen as a noncommutative version of the square of a character (as $N$ is a quadratic map) and it is the pullback of the map
       \begin{equation*}
           L_\Theta\to\R^\times\times\SO(J)\,,\quad l(a_1,\dots,a_p,M)\mapsto (a_1\cdots a_p,M)\,.
       \end{equation*}
       The corresponding map on the Gauss decomposable elements is the pullback of the map $\Phi_p$ on $G'$ defined as follows: Write $g\in G$ with respect to the quadratic form $Q$ defined above as 
        \begin{equation*}
            g=
		    \scalebox{0.8}{$
            \left[\begin{array}{c|c|c}
			     A_{11}  & A_{12} & A_{13}\\
                \hline
                A_{21}  & A_{22} & A_{23}\\
                \hline
                A_{31} & A_{32} & A_{33}\\
            \end{array}\right]$}.
        \end{equation*}
        For this we have:
        \begin{equation*}
            \Phi_p(g) = \big(\det(A_{11}),A_{22}-A_{21}A_{11}^{-1}A_{12}\big)\,.
        \end{equation*}

       To illustrate the lifts of the grading Weyl group, we restrict to the case $p=2$, i.e. $G'=\SO(3,q+3)$. The grading Weyl group is $W(R)=\langle s_1=:\varsigma_1,s_2s_3s_2=\varsigma_2\rangle<W$ which is of type $B_2$. Its standard generators and longest element lift to $G'$ using the Jordan pinning:
        \begin{align*}
            \overline{\varsigma_1}=\left[\begin{smallmatrix}
                0 & -1 & & & \\
                1 & 0 & & & \\
                & & \Id_{q+2} & & \\
                & & & 0 & -1\\
                & & & 1 & 0\\
            \end{smallmatrix}\right]\,,\hspace{1pc}
            \overline{\varsigma_2}=\left[\begin{smallmatrix}
                1 & & & & \\
                & 0 & 0 & 1 & \\
                & 0 & -J & 0 & \\
                & 1 & 0 & 0 & \\
                & & & & 1
            \end{smallmatrix}\right]\,,\hspace{1pc}
            \overline{w_0}=\overline{\varsigma_1\varsigma_2\varsigma_1\varsigma_2}= \left[\begin{smallmatrix}
                & & & & 1 \\
                & & & -1 & \\
                & & \Id_{q+2} & & \\
                & -1 & & & \\
                1 & & & &
            \end{smallmatrix}\right]\,.
        \end{align*}

	   As $W(R)$ is of type $B_2$, we can fix the reduced expression $w_0=\varsigma_1\varsigma_2\varsigma_1\varsigma_2$ to obtain the positive semigroup $U_\Theta^{>0}$ in $U_\Theta$:
	   \begin{equation*}
		  U_\Theta^{>0}=\left\{x_1(x)x_2(v)x_1(y)x_2(w)=:F(x,v,y,w)|x,y\in\R^{>0}, v,w\in \overset{\circ}{c}_{\beta}\right\}\,.
	   \end{equation*}
	   If we write this parametrization in terms of matrices by applying $\rho$, we find
	   \begin{align*}
		  \rho\big(F(x,v,y,w)\big)
            =\left[\begin{smallmatrix}
			     1 & x+y & \sqrt{2}(xv^TJ+(x+y)w^TJ) & xN(v+w)+yN(w) & xyN(v)\\
			     0 & 1 & \sqrt{2}(v+w)^TJ & N(v+w) & yN(v)\\
			     0 & 0 & Id_{q+2} & \sqrt{2}(v+w) & \sqrt{2}yv\\
			     0 & 0 & 0 & 1 & x+y\\
			     0 & 0 & 0 & 0 & 1
		  \end{smallmatrix}\right]\,.
	   \end{align*}
        This is (up to some factors of $\sqrt{2}$) the parametrization of $U_\Theta^{>0}$ described in \cite{guichard2022generalizing}.

\subsubsection{Pairs of flags} 

As in the generic case, the invariants of a pair of flags are the $\Theta$-fundamental weights applied to the $L$-distance between the flags. We now calculate explicitly the change in $L$-distance when we reverse a pair of transverse flags.

\begin{lemma}
    Suppose $A_1,A_2\in\decFlagTheta$ are transverse. Then
    \begin{equation*}
        l(s_G A_2,A_1)=\begin{cases}
            \tau\big(l(A_1,A_2)\big) & p \text{ even}\\
            \sigma\big(l(A_1,A_2)\big) & p \text{ odd} \,.
        \end{cases} 
    \end{equation*}
     Moreover, if $l_1,l_2\in L_\Theta$
    \begin{equation*}
        l(A_1l_1,A_2l_2)=\begin{cases}
            \tau(l_1)l(A_1,A_2)l_2 & p \text{ even}\\
            \sigma(l_1)l(A_1,A_2)l_2 & p \text{ odd}\,.
        \end{cases}
    \end{equation*}
\end{lemma}
\begin{proof}
    From \Cref{thm:HDistPropertiesTheta}, we know $l(s_G A_2,A_1) = \oppInv{l(A_1,A_2)}$. In \Cref{thm:W0ActionBp}, we computed the opposition involution depending on the parity of $p$ which completes the proof. 
    Similarly the scaling on $A_1$ factors as the action of $w_0(l_1^{-1})$ depends in the same way on the parity of $p$. 
\end{proof}

\begin{remark}
    Regardless of the parity of $p$, the difference between $l(A_1,A_2)$ and $l(s_GA_2,A_1)$ can be accounted for by weaving and/or switching in the polygonal cluster algebra.  
\end{remark}

\subsubsection{Triples of flags}
We now account for the last details in \Cref{sec:flagTriplesTheta} to build a cluster chart for a triple of flags $(A_1,A_2,A_3)$ in Jordan split type $B_p$. As in \Cref{sec:flagTriplesTheta} we will chose a reduced expression for $w_0$ and use it to compute an interpolating sequence of flags $A_2=A_2^{0}, A_2^{1},\cdots A_2^{p^2} = A_3$. This sequence allows us to obtain a chart by amalgamating charts for elementary configurations $(A_1,A_2^{k-1},A_2^{k})$. In \Cref{fig:QuiverTripleFlags}, we see the resulting quiver for the reduced expression $D(w_0) = (\varsigma_1\cdots \varsigma_p)^p$. The corresponding polygonal tiling for the pruned quiver is given in \Cref{fig:tripleFlagsDecPoly}.

It remains to show that all the choices made in the construction above result in cluster charts which are related by a sequence of mutations. In particular, we prove that cluster charts using two words related by the 4-braid move are connected by mutations. We will also give an explicit mutation sequence relating the charts obtained by twisted cyclic shift. 
\medskip
\input{Figures/figureQuiverTripleFlags}

\input{Figures/figureTripleFlagsDecPoly}

There is one subtlety in the amalgamation which we must discuss before studying the four-move or twisted cyclic shift.  In \Cref{fig:tripleFlagsDecPoly}, we see that the coloring of the bottom angles in the triangle depend on the parity of $p$. We must compute the difference between $\Lambda_p(A_2^{k},A_2^{k-1})$ and $\Lambda_p(A_3,A_2)$ which by \Cref{thm:HDistPropertiesTheta} is the action of $w_k = \varsigma_{i_1} \cdots \varsigma_{i_{k-1}}$ on $\check{\beta}_p$. 
\begin{lemma}\label{thm:InterpolatingSeqInvolutionBp}
    Fix a reduced expression $D=\varsigma_{i_1}\cdots \varsigma_{i_{p^2}}$ of $w_0$ and let $\{A_2^i\}$ be the corresponding interpolating sequence. Let $k$ be the index such that $l(A_2^k,A_2^{k-1}) = \check{\beta}_p(t) \neq\Id$. Let $n_p:=|\{j\in[1,k-1]\,|\, i_j=p\}|$ be the number of times $\varsigma_p$ occurs in $w_k$. Then
    \begin{equation*}
        \Lambda_p(A_3,A_2) = \sau^{n_p}\big(\Lambda_p(A_2^k,A_2^{k-1})\big) 
    \end{equation*}
\end{lemma}
\begin{proof}
From \Cref{thm:ActionOfSigmaPonSpecialRoot} we know that $\varsigma_p\big(\check{\beta}_p(A)\big)=\check{\beta}_p\big(\sau(A)/N(A)\big)$. By construction of the interpolating sequence and \Cref{thm:HDistPropertiesTheta} we know that
\begin{equation*}
    \Lambda_p(A_3,A_2) = (w_k)^{-1}\big(\Lambda_p(A_2^k,A_2^{k-1})\big)\in\mathrm{Im}(\check{\beta}_p)\,.
\end{equation*}
Since all the other reflections $\varsigma_i$ for $i\neq p$ only change the argument of $\check{\beta}_p$ by multiplying by central terms and the calculation is exactly as in the split case for these, $\Lambda_p(A_3,A_2)$ and $\Lambda_p(A_2^k,A_2^{k-1})$ only differ because of the action of $\varsigma_p$, which applies $\sau$ unlike in the commutative case. Since we do this $n_p$ times, this proves the statement.
\end{proof}

Therefore to amalgamate the elementary quiver with the bottom edge of the triangle, we must weave the bottom frozen node of the elementary quiver when $n_p$ is odd.
\begin{example}
    For $w_0 = (\varsigma_1 \cdots \varsigma_p)^p$, the relevant scaling occurs at $k= n$ and so $n_p = p-1$. This results in the difference of angle coloring in \Cref{fig:tripleFlagsDecPoly}.\\
    For $w_0 = (\varsigma_p \cdots \varsigma_1)^p$, the relevant scaling occurs at $k= 1$ and so $n_p = 0$. Therefore for this word, no weaving occurs at the bottom node regardless of the parity of $p$.\\
    For $w_0 = M_{p-1}M_p M_{p-2}M_{p-3}\cdots M_{1}$, the relevant scaling occurs at $k = 4$ and $w_k = M_{p-1}$. Thus $n_p=1$ and weaving occurs at the bottom node regardless of the parity of $p$. 
\end{example}

\begin{example}\label{ex:b2_explicit}
    Let us consider the explicit example $p=2$. Take an element of $\Confx_3(\decFlagTheta)$ represented by $(A_1,A_2,A_3)=(U_\Theta,\overline{w_0}U_\Theta l_2,F(x,v,y,w)\overline{w_0}U_\Theta l_3)$ with $l_2=\check{\beta}_1(a)\check{\beta}_2(B)$, $l_3=\check{\beta}_1(c)\check{\beta}_2(D)\in L_\Theta$ where $a,c\in\K'$ and $B,D\in\Gamma^0(V,b)$. Using \Cref{thm:HDistPropertiesTheta}(3) and \Cref{thm:flagSequenceTheta}, we compute
    \begin{align*}
        l(A_3,A_2) &= w_0(l_3)^{-1}\Big(\varsigma_1\varsigma_2\varsigma_1\big(\check{\beta}_2(\iota(\oppInv{w})\big) \varsigma_1\varsigma_2\big(\check{\beta}_1(y)\big) \varsigma_1\big(\check{\beta}_2(\iota(\oppInv{v}))\big) \check{\beta}_1(x)\Big)l_2\\
        &= w_0(l_3)^{-1}\Big(\check{\beta}_2(\iota(w))\cdot \check{\beta}_1(y)\check{\beta}_2(y) \cdot \check{\beta}_2(\iota(\oppInv{v}))\check{\beta}_1(N(v))\cdot \check{\beta}_1(x)\Big) l_2\\
        &= w_0(l_3)^{-1}\check{\beta}_1(xyN(v))\check{\beta}_2(y\iota(w)\iota(\tau(v))) l_2\\
        &= \check{\beta}_1\big(acxyN(v)\big)\check{\beta}_2\big(y\tau(D)wvB\big)
    \end{align*}
    where we used \Cref{thm:ActionOfSigmaPonOtherRoots,thm:ActionOfSigmaPonSpecialRoot} to compute the action of the Weyl group on the second line. We observe the first and last terms resulted in simple coroots. Therefore, the interpolating sequence for the reduced expression $w_0=\varsigma_1\varsigma_2\varsigma_1\varsigma_2$ is given by
    \begin{align*}
        A_2^0 &= A_2 = \overline{w_0}U_\Theta l_2\\
        A_2^1 &= x_1(x)\overline{w_0}U_\Theta\check{\beta}_1\big(cyN(v)N(B)\big)\check{\beta}_2\big(B\big)\\
        A_2^2 &= x_1(x)x_2(v)\overline{w_0}U_\Theta\check{\beta}_1\big(cyN(v)N(B)\big)\check{\beta}_2\big(cyve\sau(B)\big)\\
        A_2^3 &= x_1(x)x_2(v)x_1(y)\overline{w_0}U_\Theta\check{\beta}_1\big(c\big)\check{\beta}_2\big(cyve\sau(B)\big)\\
        A_2^4 &= A_3 = F(x,v,y,w)\overline{w_0}U_\Theta l_3
    \end{align*}
    The coordinate functions are easy to obtain from this: For instance, to the internal big node we assign
    \begin{align*}
        \Lambda_2(A_1,A_2^2) &= \Lambda_2\Big(U_\Theta,x_1(x)x_2(v)\overline{w_0}U_\Theta\check{\beta}_1\big(cyN(v)N(B)\big)\check{\beta}_2\big(cyve\sau(B)\big)\Big)\\
        &= \Lambda_2\Big(\check{\beta}_1\big(cyN(v)N(B)\big)\check{\beta}_2\big(cyve\sau(B)\big)\Big)=cyve\sau(B)\,.
    \end{align*}
    All the functions together with the quiver are shown in \Cref{fig:SO3qquiverFunctions}. 

    \input{Figures/figureSO3qquiverFunctions}
\end{example}

\subsubsection{Four move}\label{sec:fourmove}

We now prove that the four move, $\varsigma_{p-1}\varsigma_p\varsigma_{p-1}\varsigma_{p} \mapsto \varsigma_p\varsigma_{p-1}\varsigma_{p}\varsigma_{p-1}$ can be realized as mutation which is the last remaining part of the proof of \Cref{thm:clusterChartsTriplesRelatedByMutationTheta} in this case. We will see that cluster charts associated to two reduced expression which differ by a four move are related by a sequence of cluster mutations and a switch in the polygonal cluster algebra as defined in \Cref{def:weavingQuiver} and \Cref{def:weavingGroup}.

\begin{proposition}\label{thm:fourMoveTheta}
    The seeds given by two words which differ by the four move are related by a sequence of three mutations and a switch. If the mutable big node at height $p$ is labeled $2$ and the mutable small node is labeled $1$, the mutation sequence is $1,2,1$, then perform a switch at $2$.
\end{proposition}
It suffices to consider Jordan split type $B_2$, since the mutation sequence only sees levels $p$ and $p-1$. In this case there are only two different decompositions of $w_0$, 
\begin{equation*}
    D_1:=\varsigma_1\varsigma_2\varsigma_1\varsigma_2 \hspace{2pc} D_2:=\varsigma_2\varsigma_1\varsigma_2\varsigma_1.
\end{equation*}
Let $u=F_{D_1}(x,v,y,w)$. Now consider a configuration of flags $(A_1,A_2,A_3)=(U_\Theta,\overline{w_0}U_\Theta ,u\overline{w_0}U_\Theta)$ as in \Cref{ex:b2_explicit}. Notice that we set $l(A_1,A_3)=l(A_1,A_2)= \Id$ for simplicity. One can check by \Cref{thm:HDistPropertiesTheta} that these $L$-distances factor out of every computation.

We begin by computing the coordinates in \Cref{fig:FourMoveOrderedQuiverStart}. These are the coordinates we computed in \Cref{ex:b2_explicit}, explicitly:
\begin{center}
    \begin{tabular}{l l}
        $ X_2 = yve$ \hspace*{3pc} &$X_8 = ywv $\\
        $x_5 = yN(v)$  &$x_7 = xyN(v)$\,.
    \end{tabular}
\end{center}

\begin{lemma}\label{thm:B2FlipSeq1}
    The chart $Q_{D_1,1}$ is connected to $Q_{D_2,2}$ by mutation at node 1 and a switch and weave at the frozen nodes $X_1,X_8$. 
\end{lemma}
\begin{proof}
    This mutation is easy to calculate since it happens at a commutative node. 
    We have 
    \begin{equation}
        x_5' = \frac{xyN(v)+y^2N(v)}{yN(v)} = x+y
    \end{equation}

    Next we need to calculate the cluster chart associated to $Q_{D_2,2}$. When we apply the twisted cyclic shift to our triple we obtain: 
    \begin{align*}
        (s_GA_2,A_3,A_1)= &(s_G\overline{w_0}U_\Theta,u\overline{w_0}U_\Theta,U_\Theta)\\
        \intertext{which represents the same configuration as}
        &(u^{-1}s_G\overline{w_0}U_\Theta,\overline{w_0}U_\Theta,U_\Theta) =(u^{-1}\doverline{w_0}U_\Theta,\overline{w_0}U_\Theta,U_\Theta)\,.
    \end{align*}
    In this case the interpolating sequence of flags is easy to compute. We have 
    \begin{equation*}
        A_3^0=A_3,\hspace{2pc} A_3^1= \overline{\varsigma_1\varsigma_2\varsigma_1}U_\Theta,\hspace{2pc} A_3^2 =\overline{\varsigma_1\varsigma_2}U_\Theta,\hspace{2pc} A_3^3=\overline{\varsigma_1}U_\Theta,\hspace{2pc} A_3^4 = A_1\,.
    \end{equation*} 
    It suffices to show 
    \begin{align*}
        x_5' = &\Lambda_1(u^{-1}\doverline{w_0}U_\Theta, A_3^2) \\
        X_2 = & \Lambda_2(u^{-1}\doverline{w_0}U_\Theta, A_3^1) = \Lambda_2(u^{-1}\doverline{w_0}U_\Theta, A_3^2) \,.
    \end{align*}
    The last equality is by \Cref{thm:elementaryConfigLeftRightFunctionsTheta} and shows we can compute both values from the same $L$-distance. We observe $l(u^{-1}\doverline{w_0}U_\Theta, \overline{\varsigma_1\varsigma_2}U_\Theta) = l(U_\Theta,\overline{w_0}u\overline{\varsigma_1\varsigma_2}U_\Theta)$, so computing the $L$-distance reduces to computing $[u\overline{\varsigma_1\varsigma_2}]_0$.
    \begin{align*}
        [u\overline{\varsigma_1\varsigma_2}]_0 &= [x_1(x)x_2(v)x_1(y)\textcolor{red}{x_2(w)}\overline{\varsigma_1\varsigma_2}]_0 \\
        \intertext{Since $x_2(w)\overline{\varsigma_1\varsigma_2} = \overline{\varsigma_1\varsigma_2}\exp(w')$ for $w'\in \mathfrak{u}_{\beta_1+\beta_2}$, $x_2(w)$ escapes to the right.}
        &= [x_1(x)x_2(v)\textcolor{blue}{x_1(y)\overline{\varsigma_1}}\overline{\varsigma_2}]_0\\
        \intertext{Next, we decompose $x_1(y)\overline{\varsigma_1} = y_1(y^{-1})\check{\beta}_1(y)x_1(-y^{-1})$ using \Cref{thm:DecomposeXSigma}.}
        &= [x_1(x)x_2(v)\textcolor{blue}{y_1(y^{-1})\check{\beta}_1(y)x_1(-y^{-1})}\overline{\varsigma_2}]_0 \\
        \intertext{Now $x_1(-y^{-1})$ escapes past $\overline{\varsigma_2}$ as another element in $U_\Theta$ and $y_1(y^{-1})$ commutes with $x_2(v)$. }
        &= [\textcolor{orange}{x_1(x)y_1(y^{-1})}x_2(v)\check{\beta}_1(y)\overline{\varsigma_2}]_0 \\
        \intertext{We use \Cref{thm:DecomposeXY} to decompose $x_1(x)y_1(y^{-1})$ and apply $\varsigma_2$ to $\check{\beta}_1(y)$ and obtain:}
        &= [\textcolor{orange}{y_1\big((x+y)^{-1}\big)\check{\beta}_1\big((x+y)y^{-1}\big)x_1\big((x^{-1}+y^{-1})^{-1}\big)}\textcolor{blue}{x_2(v)\overline{\varsigma_2}}\check{\beta}_1(y)\check{\beta}_2(y)]_0 \\
        \intertext{Finally, we decompose $x_2(v)\overline{\varsigma_2}$ allowing all the remaining $x_i/y_i$ to escape on the right/left.}
        &= \check{\beta}_1(x+y)\check{\beta}_2(y\textcolor{blue}{ve}) \,.
    \end{align*}
    Then applying the Jordan weights, we see 
    \begin{equation*}
        \Lambda_1(u^{-1}\doverline{w_0}U_\Theta, A_3^2)= x+y=x_5'\quad\text{and}\quad \Lambda_2(u^{-1}\doverline{w_0}U_\Theta, A_3^1) = yve = X_2
    \end{equation*} as needed.
    The switch and weave at the frozen nodes happen since $Q_{D_2,2}$ calculate the functions on the edge $(A_1,A_2)$ and $(A_2,A_3)$ in the opposite direction from $Q_{D_1,1}$.
\end{proof}

\begin{lemma}
    The chart $Q_{D_1,1}$ is connected to the chart $Q_{D_2,3}$ by mutation and weave at node 2 and a switch and weave at the frozen nodes $X_1,X_3$.
\end{lemma}
\begin{proof}
    This is almost the same calculation as last time, except we have to use the noncommutative polygonal cluster algebra mutation rule. We have 
    \[X_2' = \sau(X_8)X_2^{-1}+\sigma(X_2)^{-1}x_5 = yewve \cdot y^{-1}evN(v)^{-1}+ y^{-1}evN(v)^{-1}\cdot yN(v) = e(w+v)\]

    To compute $Q_{D_2,3}$ we look at the configuration of flags $(A_3,s_GA_1,A_2)$ which is the same configuration as $(u\doverline{w_0}U_\Theta,U_\Theta,\doverline{w_0}U_\Theta)$. The interpolating sequence of flags is given by \[A_1^0 = U_\Theta \hspace{2pc} A_1^1=\doverline{\varsigma_2}U_\Theta, \hspace{2pc} A_1^2=\doverline{\varsigma_2\varsigma_1}U_\Theta, \hspace{2pc} A_1^3=\doverline{\varsigma_2\varsigma_1\varsigma_2}U_\Theta,\hspace{2pc} A_1^4=A_2\]
    As before, both interior functions can be computed at the second flag. So we want to show
    \begin{align*}
        &\Lambda_1(u\doverline{w_0}U_\Theta, A_1^2)= x_5 \\
        &\Lambda_2(u\doverline{w_0}U_\Theta, A_1^2)= \sau(X_2') 
    \end{align*}

    We compute this time $l(u\doverline{w_0}U_\Theta, \doverline{\varsigma_2\varsigma_1}U_\Theta) = [u^{-1}\doverline{\varsigma_2\varsigma_1}]_0$ and so 
    \begin{align*}
        [u^{-1}\doverline{\varsigma_2\varsigma_1}]_0&= [x_2(-w)x_1(-y)x_2(-v)x_1(-x)\doverline{\varsigma_2\varsigma_1}]_0 = [x_2(-w)x_1(-y)x_2(-v)\doverline{\varsigma_2\varsigma_1}]_0\\
        &= [x_2(-w)y_2(-v^{-1})x_1(-y)\doverline{\varsigma_1}\check{\beta}_2(\iota(v))\check{\beta}_1(N(v))]_0 = \check{\beta}_2(\iota(v+w))\check{\beta}_1(yN(v)).
    \end{align*}
    Recalling that $\iota(v+w)=(v+w)e$ we see that we needed to perform a weave at $X_2'$. The switch and weave at the frozen nodes occurs since $Q_{D_2,3}$ computes two of the edges in the opposite order from $Q_{D_1,1}$. 
\end{proof}

To complete the proof of the \Cref{thm:fourMoveTheta}, we simply observe that mutating at nodes $1,2,1$ and weaving at node $2$ moves us from $Q_{D_1,1}$ to $Q_{D_2,2}$ to $Q_{D_1,3}$ and finally to $Q_{D_2,1}$. Thus the two charts related by a four move are related by cluster mutation as needed.

\input{Figures/figureFourMoveOrderedQuivers}

\input{Figures/figureB2mutationTheta}

\subsubsection{Rotation}\label{sec:RotationThetaBp}

This section is dedicated to the proof of \Cref{thm:coordinateSystemsS3MutationTheta} for the $B_p$ case. The case $p=2$ has incidentally already been proven in the previous section. 

We will do this by considering a particular reduced expressions of $w_0$. Let $M= M_1M_2\dots M_p$ be the reduced expression given by the product of the mountain blocks in order and let $C= (\varsigma_1\cdots\varsigma_p)^p$ be the Coxeter word. We denote by $M'$, $C'$ the reversed words. 

\begin{proposition}
    The cluster charts $Q^i_{M}$ and $Q^{i+1}_{M'}$ are related by a sequence of cluster mutations at commutative nodes and switches and weaves at noncommutative nodes
\end{proposition}

Proving this proposition will allow us to change to any rotated cluster chart since we can use two moves, three moves and four moves to change $Q^2_{M'}$ into $Q^2_{M}$ and then apply the same sequence to change to $Q^3_{M'}$.

\begin{proof}    

 In the split real case ($G=\Spin(p,p+1)$) Le \cite{le2019cluster} constructs an explicit sequence of mutations which relate the cluster charts $Q^1_C$ and $Q^2_{C'}$. We will denote this sequence by $L$. There are two important facts to know about this sequence:
    \begin{enumerate}
        \item It only mutates at small nodes.
        \item It flips all the arrows between big nodes which correspond to digons in \Cref{fig:tripleFlagsDecPoly}, i.e. reverses their direction and changes their color.
    \end{enumerate}

The following lemma is a simple exercise in rearranging reduced expressions:
\begin{lemma}
    The reduced expression $M$ can be changed into $C$ only by using two moves and three moves.
\end{lemma}

    Let $R$ be the sequence of mutations realizing the change from $Q^1_M$ to $Q^1_C$ and $R'$ the sequence which moves $Q^2_{C'}$ to $Q^2_{M'}$. Our candidate sequence is given by $RLR'$ (read as a sequence of mutations from left to right). Since the mutation corresponding to the three move does not change any arrows between the noncommutative nodes, this composite sequence also has the effect of flipping all of the arrows which are in the digons. 

    We now want to compare directly the functions of our mutated chart $RLR'(Q^1_M)$ and those of $Q^2_{M'}$. Due to Le's work we know that the functions on the small nodes transform correctly, and moreover we know that the norms of the noncommutative coordinates have not changed (and hence have transformed correctly). Therefore, we can make our lives much easier by considering a seed for which all of the commutative nodes and every norm is equal to 1  and explicitly computing the transformation.

    Consider the triple of flags $(A_1,A_2,A_3)=(U_\Theta,\overline{w_0}U_\Theta l_2,u\overline{w_0}U_\Theta l_3)$ such that $l_2=\check{\beta}_p(A), l_3=\check{\beta}_p(B)$ and $u=F_M(\boldsymbol{t})$ such that $N(A)=N(B)=1$ and all the real entries of $\boldsymbol{t}$ are 1 and the vector entries given in order by $v_1,\dots,v_p$ with $N(v_i)=1$. We will denote for convenience $V_i:=\iota(v_i)$ and we recall that $\sigma(V_i)=V_i$. We will also denote for convenience $\oppInv{l}:=\sau^{p-1}(l)$ for elements of $L_\Theta$.
    
    Then we can compute all of the coordinates of $Q^1_M$. We first see that $\Lambda_i(A_3,A_2)=1$ for $i<p$ since all the norms and commutative factorization parameters are 1. $\Lambda_p(A_3,A_2)$ is given by
    \begin{align*}
        X_{p}= \Lambda_p(A_3,A_2) &= \sigma(\oppInv{B})\cdot (M_1\cdots M_{p-1}(\oppInv{V_p}))\cdot (M_1\cdots M_{p-2}(\oppInv{V_{p-1}}))\cdots \oppInv{V_1} \cdot A \\
        &= \sigma(\oppInv{B})V_p\tau(V_{p-1})V_{p-3}\cdots \oppInv{V_1}A\,.
    \end{align*}

    The other noncommutative functions can be computed in order using the formulas for the potentials. We have 
    \begin{align*}
        &X_0= A,\quad  X_1= \tau(\oppInv{V_1})\sau(A),\quad X_2=\tau(\oppInv{V_2})(\oppInv{V_1})A,\quad X_3=\tau(\oppInv{V_3})(\oppInv{V_2})\tau(\oppInv{V_1})\sau(A),\dots \\
        & X_i = \tau(\oppInv{V_{i}})(\oppInv{V_{i-1}})\cdots \sau^i(A) \quad X_{p+1}= B\,.
    \end{align*}

    Now we compute the functions in the rotated chart, $Q^2_{M'}$ which amounts to considering the configuration $(u^{-1}\doverline{w_0}U_\Theta l_2,\overline{w_0}U_\Theta l_3,U_\Theta)$. An interpolating sequence between $\overline{w_0}U_\Theta l_3$ and $U_\Theta$ using the word $M'$ begins with 
    $A_3^1 = \overline{M_1\cdots M_{p-1}}U_\Theta$. This is where the scaling corresponding to $\check{\beta}_p$ should happen and so $l(\overline{M_1\cdots M_{p-1}}U_\Theta,A_3)= B= l(A_1,A_3)$. Therefore writing $M=\varsigma_{i_1}\cdots\varsigma_{i_{p^2}}$, we have the interpolating sequence is just given by $A_3^k= \overline{\varsigma_{i_1}\cdots\varsigma_{i_{n-k}}}U_\Theta$. This makes it easy to compute the corresponding functions. 

    The internal coordinates we are interested in are $X_i'= \Lambda_p(u^{-1}\doverline{w_0}U_\Theta l_2,\overline{M_1\cdots M_i}U_\Theta)$. To compute $X_i'$, we write $u= u_1(v_1) u_2(v_2) \dots u_p(v_p)$ where each $u_i(v_i)$ is the part of the product of $F_M(\boldsymbol{t)}$  which corresponds to the mountain block $M_i$ e.g $u_2(v_2) = x_2(1)x_3(1)\cdots x_p(v_2) \cdots x_2(1)$.

    Now we manipulate $$\Lambda_p(u^{-1}\doverline{w_0}U_\Theta l_2,\overline{M_1\cdots M_i}U_\Theta) = w_0(l_2^{-1})\Lambda_p(U_\Theta,\overline{w_0}u\overline{M_1\cdots M_i}U_\Theta) = \sigma(\oppInv{A})[u\overline{M_1\cdots M_i}]_0 .$$

    Next, we use the fact that $u_k(v_k)\overline{M_i}=\overline{M_i}u_k(\tau(v_k)) $ for $i\neq k $ to write
    \begin{align*}
        [u\overline{M_1\cdots M_i}]_0 &= [u_1(v_1)\overline{M_1}u_2(\tau(v_2))\overline{M_2}\cdots u_i(\tau^{i-1}(v_i)\overline{M_i}]_0
    \end{align*}
    Looking at just $\Lambda_p[u_i(v_i)\overline{M_i}]_0$ we find that the only contribution happens when we decompose $x_p(v_i)\overline{\varsigma_p\varsigma_{p-1}\dots\varsigma_i}$ giving us a $\check{\beta}_p(v_i)$. We still have to worry about the remaining upper triangular part $x_p(-v_i^{-1})\overline{\varsigma_{p-1}\dots\varsigma_i}$. Luckily this reflected into a root space which commutes with all the other $u_k(t)$ and $\overline{M_k}$ for $k>i$ and we can remove it. 

    Therefore, we have 
    \begin{align*}
        [u_1(v_1)\overline{M_1}u_2(\tau(v_2))\overline{M_2}\cdots u_i(\tau^{i-1}(v_i)\overline{M_i}]_0 = \check{\beta}_p(V_1)\check{\beta}_p(\tau(V_2))\dots\check{\beta}_p(\tau^{i-1}(V_i))\,.
    \end{align*}

    Putting it all together, we have 
    \begin{equation*}
        X_i'= \oppInv{\sigma(A)}V_1\tau(V_2)\dots\tau^{i-1}(V_i) \,,
    \end{equation*}
    and the frozen variables $X_0' = \oppInv{\sigma(A)}$, $X_{p+1}'= B$, $X_{p}=\oppInv{\sigma(A)}V_1\tau(V_2)\dots\oppInv{V_p}B$.

    Therefore, we have to switch at every variable other than $X_p$. For $p$ even we also weave at each $X_{2i}$, but for $p$ odd we weave at each $X_{2i+1}$  except for the frozen variable $X_{p}$. 

    Comparing the decorated quivers after weaving, we see that they agree, proving the theorem.

\end{proof}

\subsection{\texorpdfstring{Type $G_2$}{Type G2}}

We now consider Jordan split type $G_2$. As in Jordan split type $B_p$ we only need to prove two lemmas in order to use the general theory to construct cluster charts for arbitrary configurations of flags. 
\begin{description}
    \item[\Cref{thm:clusterChartsTriplesRelatedByMutationTheta}] The two charts on triples of flags related by the six move are connected by a sequence of cluster mutations.
    \item[\Cref{thm:coordinateSystemsS3MutationTheta}] Two charts related by twisted cyclic shift are connected by a sequence of cluster mutations. 
\end{description}

In this case, $G$ is associated to a commutative Jordan algebra given by a field extension $\K'$ associated to the root $\beta_1$ and a degree $3$ Jordan algebra over $\K'$ which we denote simply by $\jordan{J} = (V,\iota,1) $ which is associated to the short root $\beta_2$.  The grading Weyl group, $W(R)$, is generated by $\varsigma_1,\varsigma_2$. For ease of reading in this section,  we will denote $\beta_1 =\alpha$ and $\beta_2= \beta$. 

For any Jordan algebra of degree 3, the norm map is degree $3$, and the adjugate map is degree 2. We can define a trilinear form $3N(u,v,w):=N(u+v+w)-N(u+v)-N(u+w)-N(v+w)+N(u)+N(v)+N(w)$, and a bilinear product $2 (v\times w) :=\tau(v+w)-\tau(v)-\tau(w)$. This bilinear product is called the Freudenthal product.

A choice of Jordan pinning consists of the following compatible collection of maps:
\begin{equation*}
    \begin{aligned}
        x_1 &:\K' \to G\\
        x_2 &:V \to G
    \end{aligned}\hspace{2pc}
    \begin{aligned}
        \check{\alpha} &:\K'^\times \to G\\
        \check{\beta} &:\Gamma(\jordan{J}) \to G
    \end{aligned}\hspace{2pc}
    \begin{aligned}
        y_1 &:\K' \to G\\
        y_2 &:V \to G\,.
    \end{aligned}
\end{equation*}

There are two different Lusztig maps coming from the two reduced expressions of $w_0$, $D_1,D_2$: 
\begin{align*}
    F_1(a,u,b,v,c,w) &= x_1(a)x_2(u)x_1(b)x_2(v)x_1(c)x_2(w) \\
    F_2(u,a,v,b,w,c) &= x_2(u)x_1(a)x_2(v)x_1(b)x_2(w)x_1(c)
\end{align*}
with $a,b,c\in\K'$ and $u,v,w\in V$.

\begin{lemma}\label{thm:potentialInvolutionG2}
    The opposition involution $v\to \oppInv{v}$ on $\mathfrak{u}_\beta$  is trivial. 
\end{lemma}
\begin{proof}
    It suffices to compute the action of $\varsigma_2 w_0 = \varsigma_1\varsigma_2\varsigma_1\varsigma_2\varsigma_1$ on $\check{\beta}(v)$  using \Cref{thm:ActionOfSigmaPonOtherRoots,thm:ActionOfSigmaPonSpecialRoot}.
    \begin{align*}
        \varsigma_1\varsigma_2\varsigma_1\varsigma_2\varsigma_1 \cdot \check{\beta}(v) &= \varsigma_1\varsigma_2\varsigma_1\varsigma_2 \cdot \check{\alpha}(N(v))\check{\beta}(v) = \varsigma_1\varsigma_2\varsigma_1 \cdot \check{\alpha}(N(v))\check{\beta}(\tau(v)) \\
        &= \varsigma_1\varsigma_2\cdot \check{\alpha}(N(v))\check{\beta}(\tau(v)) = \varsigma_1 \cdot \check{\alpha}(N(v))\check{\beta}(v) = \check{\beta}(v) \,.
    \end{align*}
    As in \Cref{thm:potentialInvolutionBp}, the action on the coroot induces the action on the root space. 
\end{proof}
\begin{corollary}\label{thm:w0ActionG2}
    The action of $w_0$ on $\check{\beta}_2$ is given by $w_0\cdot \check{\beta}_2(v) = \check{\beta}_2(\sigma(v^{-1}))$.
\end{corollary}

\subsubsection{Rotation and six move}\label{sec:sixmove}

We now wish to prove \Cref{thm:clusterChartsTriplesRelatedByMutationTheta} for $G_2$. Our proof strategy will incidentally also construct a rotation sequence and prove \Cref{thm:coordinateSystemsS3MutationTheta}. 

Let $u= F_1(a,u,b,v,c,w)$. For simplicity, we consider the triple of flags $(A_1,A_2,A_3)=(U_\Theta,\overline{w_0}U_\Theta,u\overline{w_0}U_\Theta)$ with two trivial $L$-distances. We now compute the functions associated to the cluster chart $Q_{D_1,1}$ (\Cref{fig:G2InitialSeed}).

First we compute the $L$-distance $l(u\overline{w_0}U_\Theta,\overline{w_0}U_\Theta)$ as in \Cref{rem:LdistanceFromElementaryConfig}. We find $l(A_3,A_2) = \check{\alpha}\left(aN(u)b^2N(v)^2c\right) ~\check{\beta}\left(\iota(u)b\iota(\tau(v))c\iota(w)\right)$ and so 
$$\Delta_1(A_3,A_2) = ab^2cN(u)N(v)\hspace{3pc} \Delta_2(A_3,A_2) = bc\iota(w)\iota(\tau(v))\iota(u).$$

To keep the notation light, we will suppress the function $\iota$ and identify $u = \iota(u)$. This only affects the exceptional case, as we can take $\iota$ to be the identity map in all other cases.

Using \Cref{thm:partialPotentialFromCoordsTheta}, the formulas for the potentials in terms of the coordinates, we obtain a series of equations each containing exactly one new interior node.   Following the order of the terms in our reduced expression, we compute first a commutative variable we call $x$, then the non-commutative variable $X$, followed by $y$ and finally $Y$. We label the nodes corresponding to these variables 1,2,3,4 respectively: 
\begin{align*}
\begin{aligned}
    a_1 = x &= b^2cN(u)N(v)\\
    a_3 = y &= b^3c^2N(v)^2N(u)
\end{aligned}\hspace{3pc}
\begin{aligned}
     a_2 = X &= b^2cN(v)\tau(u) \\
     a_4 = Y &= bc\tau(v)u .
\end{aligned}
\end{align*}

\input{Figures/figureG2initial}

We now consider the mutation sequence, read left to right, $S= \mu_3\mu_1\mu_2\mu_3$. 
\begin{lemma}
    The mutation sequence $S$ followed by a switch at $4$ changes $Q_{D_1,1}$ to $Q_{D_2,2}$.
\end{lemma}
\begin{proof}
    This follows from a very similar calculation to \Cref{thm:B2FlipSeq1}. First we calculate the mutation sequence following our polygonal cluster algebra mutation rules (\Cref{fig:G2FlipSeq1}). 
    The first mutation happens at the variable $y$ so we calculate $$y'= \frac{b^5c^4N(v)^3N(u)^2+b^6c^3N(v)^3N(u)^2}{b^3c^2N(v)^2N(u)} = b^2cN(v)N(u)(c+b).$$

    The next two mutations commute with each other, and so in turn we find 
    $$x' = \frac{ab^2cN(u)N(v)+b^2cN(v)N(u)(b+c)}{b^2cN(u)N(v)} = a+b+c$$ and 
    $$ X' = \sigma(X)^{-1}\tau(Y)+ X^{-1}y' = vc+u(c+b)$$
    and finally 
    $$ y''= \frac{ab^2cN(u)N(v)N(X')+ (a+b+c)(b^3c^3N(v)^2N(u))}{b^2cN(v)N(u)(c+b)} = \frac{(aN(X')+bc^2N(v)(a+b+c))}{b+c} .$$

    Next we need to calculate the cluster chart associated to $Q_{D_2,2}$. Applying the twisted cyclic shift, we obtain $(s_GA_2,A_3,A_1)= (s_G\overline{w_0}U_\Theta,u\overline{w_0}U_\Theta,U_\Theta)$ which is the same configuration as $(u^{-1}s_G\overline{w_0}U_\Theta,\overline{w_0}U_\Theta,U_\Theta)$. As in \Cref{thm:B2FlipSeq1}, the interpolating sequence is given by successively removing generators from the end of $D_2$:
    \begin{equation*}
       A_3^0=A_3,\hspace{1pc} A_3^1= \overline{\varsigma_1\varsigma_2\varsigma_1\varsigma_2\varsigma_1}U_\Theta, \hspace{1pc} A_3^2 =\overline{\varsigma_1\varsigma_2\varsigma_1\varsigma_2}U_\Theta,\hspace{1pc}\dots\hspace{1pc}, A_3^5=\overline{\varsigma_1}U_\Theta, \hspace{1pc} A_3^6 = A_1 .
    \end{equation*}
    To see this is the same chart as the result of the mutation sequence we will show:
    \begin{align}
        &\Delta_1(u^{-1}s_GA_2, A_3^4)= x' \label{eq:G2FlipSeq1Eq_x'}\\
        &\Delta_2(u^{-1}s_GA_2, A_3^4)= X' \label{eq:G2FlipSeq1Eq_X'}\\ 
        &\Delta_1(u^{-1}s_GA_2, A_3^3)= y'' \label{eq:G2FlipSeq1Eq_y''}\\ 
        &\Delta_2(u^{-1}s_GA_2, A_3^1)= \sigma(Y) .\label{eq:G2FlipSeq1Eq_Y}
    \end{align}
    We begin by computing $l(u^{-1}s_GA_2, A_3^4) = l(U_\Theta, \overline{w_0} u\overline{\varsigma_1\varsigma_2}U_\Theta)$. So we decompose $u\varsigma_1\varsigma_2$ using properties of the Jordan pinning as follows:
    \begin{align*}
        [u\varsigma_1\varsigma_2&]_0 = [x_1(a)x_2(u)x_1(b)x_2(v)x_1(c)\textcolor{red}{x_2(w)}\overline{\varsigma_1\varsigma_2}]_0 \\
        &= [x_1(a)x_2(u)x_1(b)x_2(v)\textcolor{blue}{x_1(c)\overline{\varsigma_1}}\overline{\varsigma_2}]_0\\
        &= [x_1(a)x_2(u)x_1(b)\textcolor{red}{x_2(v)}\textcolor{blue}{y_1(c^{-1})\check{\alpha}(c)x_1(-c^{-1})}\textcolor{red}{\varsigma_2}]_0\\
        &= [x_1(a)x_2(u)\textcolor{orange}{x_1(b)y_1(c^{-1})}\textcolor{red}{y_2(v^{-1})\check{\beta}(v)x_2(-v^{-1})}\check{\alpha}(c)\check{\beta}(c)]_0\\
        &= [x_1(a)x_2(u)\textcolor{orange}{y_1((b+c)^{-1})\check{\alpha}(1+bc^{-1})x_1((b^{-1}+c^{-1})^{-1})}y_2(v^{-1})\check{\beta}(v)\check{\alpha}(c)\check{\beta}(c)]_0\\
        &= [\textcolor{blue}{x_1(a)y_1((b+c)^{-1})}x_2(u)y_2((1+bc^{-1})v^{-1})\check{\alpha}(1+bc^{-1})\textcolor{red}{x_1((b^{-1}+c^{-1})^{-1})}\check{\beta}(v)\check{\alpha}(c)\check{\beta}(c)]_0\\
        &= [\textcolor{blue}{\check{\alpha}(1+a(b+c)^{-1})x_1((a^{-1}+(b+c)^{-1})^{-1})}\textcolor{orange}{x_2(u)y_2((1+bc^{-1})v^{-1})}\check{\alpha}(c+b)\check{\beta}(vc)]_0\\       
        &= [\check{\alpha}(1+a(b+c)^{-1})x_1((a^{-1}+(b+c)^{-1})^{-1})\textcolor{orange}{\check{\beta}(1+u(1+bc^{-1})v^{-1})}\check{\alpha}(c+b)\check{\beta}(vc)]_0\\
        &= \check{\alpha}(a+b+c)\check{\beta}(vc+u(c+b))\,.
    \end{align*}
    Thus $\Delta_1(u^{-1}s_GA_2,A_3^4) = a+b+c = x'$ and $\Delta_2(u^{-1}s_GA_2,A_3^4) = vc+u(c+b) = X'$, which proves \Cref{eq:G2FlipSeq1Eq_x',eq:G2FlipSeq1Eq_X'}. 
    Next we decompose $u\varsigma_1\varsigma_2\varsigma_1$ to compute $l(u^{-1}s_GA_2,A_3^3)$. We can use the previous calculation if we keep track of the $x_i$ which don't escape past the new $\overline{\varsigma_1}$: 
    \begin{align*}
        [u\varsigma_1\varsigma_2\varsigma_1&]_0 = [x_1(a)x_2(u)x_1(b)x_2(v)x_1(c)x_2(w)\varsigma_1\varsigma_2\varsigma_1]_0 \\
        &= [\check{\alpha}(a+b+c)\check{\beta}(X')x_1\left(aN(X')(b+c)^{-1}(a+b+c)^{-1}\right)x_2(z)x_1(N(v)c^2b(b+c)^{-1})\varsigma_1]_0\\
        &= [\check{\alpha}(a+b+c)\check{\beta}(X')\check{\alpha}(N(v)c^2b(b+c)^{-1}+aN(X')(a+b+c)^{-1})]_0\\
        &= [\check{\alpha}(y'')\check{\beta}(X')]_0\,.
    \end{align*}
    Applying $\Delta_1$ we verify \Cref{eq:G2FlipSeq1Eq_y''}.
    For \Cref{eq:G2FlipSeq1Eq_Y} we note the final $x_2(w)$ escapes leaving $[x_1(a)x_2(u)x_1(b)x_2(v)x_1(c)\overline{\varsigma_1\varsigma_2\varsigma_1\varsigma_2\varsigma_1}]_0$ a product of root spaces against the opposite word.   Thus we compute the $L$-distance step by step as in \Cref{thm:DecomposeUW} to see: 
    \begin{align*}
         [u\varsigma_1\varsigma_2\varsigma_1\varsigma_2\varsigma_1]_0 =& [x_1(a)x_2(u)x_1(b)x_2(v)x_1(c)\varsigma_1\varsigma_2\varsigma_1\varsigma_2\varsigma_1]\\
         =~& \varsigma_1\varsigma_2\varsigma_1\varsigma_2(\check{\alpha}(a))\cdot \varsigma_1\varsigma_2\varsigma_1(\check{\beta}(u))\cdot\varsigma_1\varsigma_2(\check{\alpha}(b))\cdot\varsigma_1(\check{\beta}(v)) \check{\alpha}(c) \\
         =~&\check{\alpha}\left(aN(u)b^2N(v)c\right)~ \check{\beta}(ub\tau(v)c) \,.
    \end{align*}
    Thus $\Delta_2(u^{-1}s_GA_2,A_3^1) = bc u\tau(v) = \sigma(bc \tau(v)u) = \sigma(Y)$ as needed.
\end{proof}

\input{Figures/figureG2flipSeq1}

Now we consider another mutation sequence $T = \mu_2\mu_3\mu_4\mu_2$. This sequence is obtained from $S$ by switching the roles of the big and small nodes. 
\begin{lemma}
    The mutation sequence $T$ changes our cluster chart to $Q_{D_2,3}$
\end{lemma}
\begin{proof}
    This calculation is almost identical to the previous one, only the roles of the two roots are reversed. We first calculate the effect of the cluster mutations. We have 
    \begin{align*}
        X'&= \sigma(X)^{-1}y+ xX^{-1}Y\\
         &= (b^2cN(v)\tau(u))^{-1}b^3c^2N(v)^2N(u)+ b^2cN(u)N(v)(b^2cN(v)\tau(u))^{-1}bc\tau(v)u \\
         &= bcN(v)u+bcu\tau(v)u\,.
    \end{align*}
    
    The next two mutations commute and we find
    \begin{align*}
        Y' =& \sigma(Y)^{-1}\sigma(X')+ bcw\tau(v)uY^{-1}=v+u+w  \,,\\
         y'=& \frac{N(bcN(v)u +  bcu\tau(v)u)+b^4c^2N(u)^2N(v)^2}{b^3c^2N(v)^2N(u)}= cN(v + u)+bN(u) \,.
    \end{align*}
    The final mutation at $2$ produces 
    $$X''= x\sigma(Y')\sigma(X')^{-1}+y'bcw\tau(v)uX'^{-1} =\left(bN(u)(u+v+w)+ y'wv^{-1}u\right)(u+uv^{-1}u)^{-1}\,.$$

    Now we want to compare these functions to those of the cluster chart $Q_{D_2,3}$. This chart uses the configuration of flags $(A_3,s_GA_1,A_2)$ which is the same configuration as $(u\doverline{w_0}U_\Theta,U_\Theta,\doverline{w_0}U_\Theta)$. The interpolating sequence of flags is given by adding one generator in $D_2$ at a time:
    \begin{equation*}
        A_1^0 = A_1^1, \hspace{1pc} A_1^1=\doverline{\varsigma_2}U_\Theta,\hspace{1pc} A_1^2=\doverline{\varsigma_2\varsigma_1}U_\Theta,\hspace{1pc}\dots\hspace{1pc},A_1^5=\doverline{\varsigma_2\varsigma_1\varsigma_2\varsigma_1\varsigma_2}U_\Theta,\hspace{1pc} A_1^6=A_2.
    \end{equation*}

    Thus we want to prove that:
    \begin{align*}
        &\Delta_1(u\doverline{w_0}U_\Theta, A_1^2)= y' \\
        &\Delta_2(u\doverline{w_0}U_\Theta, A_1^2)= X' \\ 
        &\Delta_2(u\doverline{w_0}U_\Theta, A_1^3)= X'' \\ 
        &\Delta_1(u\doverline{w_0}U_\Theta, A_1^5)= x \,.
    \end{align*}
    First we compute $l(u\doverline{w_0}U_\theta, A_1^2)$ by decomposing $[u^{-1}\doverline{\varsigma_2\varsigma_1}]_0$:
    \begin{align*}
        &[x_2(-w)x_1(-c)x_2(-v)x_1(-b)\textcolor{blue}{x_2(-u)}{x_1(-a)}\textcolor{blue}{\doverline{\varsigma_2}}\doverline{\varsigma_1}~]_0\\
        &=[x_2(-w)x_1(-c)x_2(-v)\textcolor{red}{x_1(-b)}\textcolor{blue}{y_2(-u^{-1})\check{\beta}(u)}\textcolor{red}{\doverline{\varsigma_1}}]_0\\
        &=[x_2(-w)x_1(-c)\textcolor{orange}{x_2(-v)y_2(-u^{-1})}\textcolor{red}{y_1(-b^{-1})\check{\alpha}(b)}\check{\alpha}(N(u))\check{\beta}(u)]\\
         &=[x_2(-w)x_1(-c)\textcolor{orange}{y_2(-(v+u)^{-1})\check{\beta}(1+vu^{-1})x_2(-(v^{-1}+u^{-1})^{-1})}y_1(-b^{-1})\check{\alpha}(b)\check{\alpha}(N(u))\check{\beta}(u)]_0\\
         &=[\textcolor{red}{x_2(-w)y_2(-(v+u)^{-1})}\textcolor{blue}{x_1(-c)y_1(-b^{-1}N(1+vu^{-1}))}\check{\beta}(1+vu^{-1})\check{\alpha}(b)\check{\alpha}(N(u))\check{\beta}(u)]_0  \\
         &=[\textcolor{red}{\check{\beta}(1+w(v+u)^{-1})}\textcolor{blue}{\check{\alpha}(1+cb^{-1}N(1+vu^{-1}))}\check{\beta}(1+vu^{-1})\check{\alpha}(b)\check{\alpha}(N(u))\check{\beta}(u)]_0  \\
         &=\check{\beta}(u+v+w)\check{\alpha}(bN(u)+cN(u+v))\,.
    \end{align*}
   This proves $\Delta_1(u\doverline{w_0}U_\Theta, A_1^2) = bN(u)+cN(u+v) = y'$ and $\Delta_2(u\doverline{w_0}U_\Theta, A_1^2)=u+v+w = X'$ as needed.
   To compare with $A_1^3$, we add another $\doverline{\varsigma_2}$ on the right and find 
    \begin{align*}
    [u^{-1}&\doverline{\varsigma_2\varsigma_1\varsigma_2}]_0= \\
    &[\check{\beta}(u+v+w)\check{\alpha}(y')\check{\beta}\left((u+v)w^{-1}(u+v)+(v+u))^{-1}(y')+bN(u)(uv^{-1}u+u)^{-1}\right)]_0\,.
    \end{align*}
    By refactoring we see that $$((u+v)w^{-1}(u+v)+(v+u))^{-1} = (w+u+v)^{-1}(w)(u+v)^{-1}= (w+u+v)^{-1}(w)v^{-1}u(uv^{-1}u+u)^{-1}.$$
    Thus 
     $$[u^{-1}\doverline{\varsigma_2\varsigma_1\varsigma_2}]_0=\check{\beta}(X'')\check{\alpha}(y') $$ as desired. 
    For the final equation use \Cref{thm:DecomposeUW} as before to compute $l(u\doverline{w_0}U_\Theta, A_1^5)$:
    \[ [u^{-1}\doverline{\varsigma_2\varsigma_1\varsigma_2\varsigma_1}]_0 = \check{\alpha}(cN(v)b^2N(u))\check{\beta}(wc\tau(v)bu) \]
    Thus $\Delta_1(u\doverline{w_0}U_\Theta, A_1^5)= b^2cN(u)N(v) = x$ as needed.
\end{proof}

\input{Figures/figureG2flipSeq2}

\begin{proposition}[\Cref{thm:coordinateSystemsS3MutationTheta} for $G_2$]
    The twisted cyclic shift and the 6-move are realized by polygonal cluster mutations in Jordan split type $G_2$.
\end{proposition}
\begin{proof}
    Each of the mutation sequences $S,T$ are involutions. Starting at $Q_{D_1}^1$ and applying $S$ we go to $Q_{D_2}^2$. Next applying $T$ we move to $Q_{D_1}^3$. This is the rotation. Applying $S$ again moves us to $Q_{D_2}^1$, which is the 6-move. 
\end{proof}

\begin{remark}
    The final mutated cluster variables appear to be rational functions in the factorization parameters. 
    One can use the trilinear form and Freudenthal product to rewrite them as as \textit{polynomials} in the factorization parameters, compare to \cite{guichard2026algebraic}
\end{remark}

\section{Higher noncommutative rank}\label{sec:HigherRankExamples}
The key difference between groups of higher noncommutative rank and the previous groups is the lack of a unique system of Jordan weights. Previously the system of Jordan weights were used to define the coordinates associated to the nodes of the cluster chart. Instead, we use a Jordan pinning of the group to construct a noncommutative grounded wiring network (\Cref{sec:NoncomNetworks}) which will parameterize configurations of flags, elements of $G$, and surface representations. This description will match the noncommutative cluster structure of \cite{goncharov2021spectral}.

\subsection{Type \texorpdfstring{$A_p$}{Ap}}\label{sec:NoncomAp}
Recall that the recognition theorem (\Cref{thm:RecognitionAp}) states that for $p>2$ any $A_p$ graded Lie algebra is a central extension of $\sl_{p+1}(\divisionAlg)$ for an associative unital algebra $\divisionAlg$ over a field $\mathbb{K}$. The corresponding $\Theta$-simply connected group is $\SL_{p+1}(\divisionAlg)$. However it will be convenient to consider the central extension $\GL_{p+1}(\divisionAlg)$ of invertible $(p+1)\times (p+1)$ matrices with entries in $\divisionAlg$ as well. 

We will assume for now that if $p=2$ then we are in the classical case, so that $G$ is not a form of the group $E_6$.

\subsubsection{\texorpdfstring{$\Theta$-}{Theta }pinning}
We now construct a Jordan pinning as discussed in \Cref{sec:ThetaPinning}. In this case, each Jordan algebra $\jordan{J}_i$ is of associative type and we identify the structure group $\Gamma(\jordan{J})$ with a subgroup of $\divisionAlg^\times \times \divisionAlg^\times$. Then a choice of pinning maps is the following:
\begin{align*}
    x_i(M) = {\begin{bsmallmatrix}
        \Id & & & & \\
          & \ddots & & & \\
          & & \Id & M & &\\
          & & &  \Id & &\\
          & & &  & \ddots &\\
          & & &  & & \Id
    \end{bsmallmatrix}}\qquad
    \check{\beta}_i(A,B) = {\begin{bsmallmatrix}
        \Id & & & &\\
         & \ddots & & & \\
          & & A &  & &\\
          & &  & B & &\\
          & & &  & \ddots &\\
          & & &  & & \Id
    \end{bsmallmatrix}}\qquad
    y_i(M) = {\begin{bsmallmatrix}
        \Id & & & & \\
          & \ddots & & & \\
          & & \Id &  & &\\
          & & M & \Id & &\\
          & & &  & \ddots &\\
          & & &  & & \Id
    \end{bsmallmatrix}}
\end{align*}
In order to construct the corresponding network, we must make a choice of system of Jordan weights. One such system is to take $\Lambda_1(L_1,\cdots,L_{p+1}) = (L_1,L_2)$ and $\Lambda_i(L_1,\cdots,L_{p+1}) = (\Id_k,L_{i+1})$ for all $i> 1$.

\begin{example}
    Let $G=\GL_{3k}(\R)$ and choose $\Theta$ as in $A_{3k-1}(2,k)$ to realize $G$ as a central extension of the group of Jordan split type $A_2(\jordan{M}_k(\R))$. The Jordan pinning then has the form 
        \begin{align*} 
    x_1(M) = 
        \begin{bmatrix}
        \Id_k & M & 0\\ 0 &\Id_k & 0\\ 0& 0& \Id_k
    \end{bmatrix} \quad  \check{\beta}_1(A,B) = \begin{bmatrix}
        A  & 0 & 0\\ 0 & B & 0\\ 0& 0& \Id_k
    \end{bmatrix} \quad  y_1(M) = \begin{bmatrix} \Id_k & 0 & 0 \\ M & \Id_k & 0 \\ 0 & 0 & \Id_k\end{bmatrix}\\
    x_2(M) = \begin{bmatrix}
        \Id_k & 0 & 0\\ 0 &\Id_k & M\\ 0& 0& \Id_k
    \end{bmatrix} \quad \check{\beta}_2(A,B) = \begin{bmatrix}
        \Id_k  & 0 & 0\\ 0 &A & 0\\ 0& 0&  B
    \end{bmatrix} \quad y_2(M) = \begin{bmatrix} \Id_k & 0 & 0 \\ 0 & \Id_k & 0 \\ 0 & M & \Id_k\end{bmatrix}.
    \end{align*}

    The following are two possible systems of Jordan weights: 
    \begin{align*}
        \Lambda_1 \left(\begin{bmatrix} A& & \\& B & \\ & & C \end{bmatrix}\right) &= (A,B) \qquad \Lambda_2 \left(\begin{bmatrix} A& & \\& B & \\ & & C \end{bmatrix}\right) = (\Id_k,C) \qquad \text{or}\\
        \Lambda_1 \left(\begin{bmatrix} A& & \\& B & \\ & & C \end{bmatrix}\right) &= (A,\Id_k) \qquad \Lambda_2 \left(\begin{bmatrix} A& & \\& B & \\ & & C \end{bmatrix}\right) = (B,C) \,.
    \end{align*}
\end{example}

Consider the triple of flags $(U_\Theta,\overline{w}_0U_\Theta l_2,u\overline{w}_0U_\Theta l_3)$ with $u=x_1(M_1)x_2(M_2)x_1(M_3)$. We write $l_2=\begin{bmatrix}   A& B& C  \end{bmatrix}$ and $l_3=\begin{bmatrix}   D& E& F  \end{bmatrix}$. Using the first system of Jordan weights we calculate 
\begin{align*}
    l(u\overline{w}_0lU_\Theta l_3,\overline{w}_0lU_\Theta l_2) &= \overline{w_0}(l_3^{-1}) [\oppInv{(\varsigma_1\varsigma_2)}\check{\beta}_{\oppInv{1}}(\iota(M_3)) \oppInv{(\varsigma_1)}\check{\beta}_{\oppInv{2}}(\iota(M_2)) \check{\beta}_{\oppInv{1}}(\iota(M_1))] l_2 \\
    & = \begin{bmatrix}F^{-1}\\E^{-1}\\D^{-1}\end{bmatrix} \begin{bmatrix} M_3\\ M_3^{-1}\\\Id_k \end{bmatrix}\begin{bmatrix}M_2\\\Id_k \\ M_2^{-1}\end{bmatrix}\begin{bmatrix} \Id_k \\ M_1 \\ M_1^{-1}\end{bmatrix}\begin{bmatrix} A\\B\\C\end{bmatrix} = \begin{bmatrix}
        F^{-1}M_3M_2A\\  E^{-1}M_3^{-1}M_1B\\ D^{-1}M_2^{-1}M_1^{-1}C 
    \end{bmatrix}\\
    & = \check{\beta}_1(F^{-1}M_3M_2A,E^{-1}M_3^{-1}M_1B) \check{\beta}_2(\Id_k,D^{-1}M_2^{-1}M_1^{-1}C) .
\end{align*}
If we choose the word $\varsigma_1\varsigma_2\varsigma_1$, we obtain the following interpolating sequence of flags:
\begin{align*}
    A_3^0 &=\overline{w}_0U_\Theta \cdot \check{\beta}_1(A,B)\check{\beta}_2(\Id_k,C) \\
    A_3^1 &= x_1(M_1)\overline{w}_0U_\Theta \cdot \check{\beta}_1(A,M_2D)\check{\beta}_2(\Id_k,M_1 B)  \\
    A_3^2 &= x_1(M_1)x_2(M_2)\overline{w}_0U_\Theta \cdot \check{\beta}_1(D,M_2A)\check{\beta}_2(\Id_k,M_1 B)\\
    A_3^3 &= x_1(M_1)x_2(M_2)x_1(M_3)\overline{w}_0U_\Theta \cdot\check{\beta}_1(D,E)\check{\beta}_2(\Id_k,F).
\end{align*}
The amalgamated network is shown in \Cref{fig:A2ElementaryNetworkA}.

Using the second system of Jordan weights, we obtain a slightly different interpolating sequence: 
\begin{align*}
    A_3^0 &=\overline{w}_0U_\Theta \cdot \check{\beta}_1(A,\Id_k)\check{\beta}_2(B,C) \\
    A_3^1 &= x_1(M_1)\overline{w}_0U_\Theta \cdot \check{\beta}_1(A,\Id_k)\check{\beta}_2(M_2D,\textcolor{orange}{M_3 E})  \\
    A_3^2 &= x_1(M_1)x_2(M_2)\overline{w}_0U_\Theta \check{\beta}_1(D,\Id_k)\check{\beta}_2(M_2A,\textcolor{orange}{M_3 E})\\
    A_3^3 &= x_1(M_1)x_2(M_2)x_1(M_3)\overline{w}_0U_\Theta \check{\beta}_1(D,\Id_k)\check{\beta}_2(E,F).
\end{align*}
The new amalgamated network is shown in \Cref{fig:A2ElementaryNetworkB}. The difference between the two networks is only on the strands at level 3 between the two bridges and the legs labeled 2. They are clearly equivalent under leg slide and also both equivalent to the left slid network in \Cref{fig:A2LeftSlidNetwork}. 
\input{Figures/figureA2ElementaryNetworks}
\input{Figures/figureA2leftSlidNetwork}

\subsubsection{Three move}\label{sec:ThreeMoveHigherNoncom}
\begin{lemma}\label{thm:ThreeMoveHigherNoncomRank}
    The networks given by two seeds which differ by a three move $\varsigma_i \varsigma_{i+1}\varsigma_{i} \rightarrow \varsigma_{i+1}\varsigma_{i}\varsigma_{i+1}$ are related by leg slides and a single square move at the single internal face of $\elemNetwork(j+1,j,j+1)$ where $j=\oppInv{(i+1)}$.
\end{lemma}
\begin{proof}
    \input{Figures/figureThreeMoveHigherRank}
    We begin by fully sliding all the legs in each network to the right. Showing these left slid networks are related by a square move proves the result. In \Cref{fig:ThreeMoveHigherRank}, we see the portion of each network corresponding to the three move. 

    In order to use the relations for the square move given in \Cref{eqn:SquareMove}, we identify the edges of the network as follows :
    \begin{align*}
    \begin{aligned}
        a_1&=M_2^{-1}M_3^{-1}M_1\\
        b_1&=N_2
    \end{aligned}\qquad
    \begin{aligned}
        a_2&=M_1^{-1}\\
        b_2&= N_1
    \end{aligned}\qquad
    \begin{aligned}
        a_3&= M_1\\
        b_3&=N_1^{-1}
    \end{aligned}\qquad
    \begin{aligned}
        a_4&=M_2\\
        b_4&=N_1N_3^{-1}N_2^{-1} \,.
    \end{aligned}
    \end{align*}
    We then have
    \begin{align*}
        N_1 &= b_2 = (1+a_4a_1a_2a_3)^{-1}a_4 = (1+M_3^{-1}M_1)^{-1}M_2 = (M_3+M_1)^{-1}M_3 M_2\\
        N_2 &= b_1 = (1+a_3a_4a_1a_2)a_2^{-1}a_1^{-1}a_4^{-1} = (1+M_1M_3^{-1})M_3 = (M_3+M_1)\\
        N_3 &= N_2^{-1} b_4^{-1} N_1 = (M_3+M_1)^{-1}[(1+M_3M_1^{-1})^{-1}M_1^{-1}]^{-1} (M_3+M_1)^{-1}M_3M_2 \\
            &= (M_3+M_1)^{-1}M_1M_2\,.
    \end{align*}
    
    Next we compare these relations to the information from the flags. The first network encodes the configuration of flags $(U_\Theta, \overline{w_0}U_\Theta, u\overline{w_0}U_\Theta)$ where $u$ is decomposed as $x_i(M_1)x_{i+1}(M_2)x_i(M_3)$. The second network instead decomposes $u$ as $x_{i+1}(N_1)x_i(N_2)x_{i+1}(N_3)$. We compute the portion of the matrix with rows and columns $i,i+1,i+2$ as follows:
    \begin{align*}
        x_i(M_1)x_{i+1}(M_2)x_i(M_3) = {\setlength\arraycolsep{0.8 pt}\begin{bmatrix}
            1 & M_1+M_3 & M_1 M_2\\
              & 1 & M_2 \\
              &   & 1
        \end{bmatrix}}
        =
         {\setlength\arraycolsep{0.8 pt}\begin{bmatrix}
            1 & N_2 & N_2 N_3\\
              & 1 & N_1+N_3\\
              &   & 1
        \end{bmatrix}}=x_{i+1}(N_1)x_{i}(N_2)x_{i+1}(N_3) .
    \end{align*}
    For these decompositions to be equal we must have
    \begin{equation*}
        N_1 = (M_1+M_3)^{-1}M_3M_2, \hspace{2pc} N_2 = M_1+M_3, \hspace{2pc} N_3=(M_1+M_3)^{-1}M_1M_2\,,
    \end{equation*}
    which agrees with the square move calculation above.
\end{proof}

\subsection{\texorpdfstring{Type $C_p$}{Type Cp}} 
Now consider a group $G$ with a $C_p$ grading. By the recognition theorem (\Cref{thm:RecognitionCp}), for $p \geq 4$ the Lie algebra of $G$ is $\slj$ for $\jordan{J}$ a Hermitian type Jordan algebra sitting inside an associative algebra $\mathcal{R}$ with anti-involution $\star$. Once again, we will assume that we are not in the exceptional case, i.e. $G$ is not a form of the group $E_7$. 

When we use the symplectic form 
$$\begin{bsmallmatrix}
    & & & & -\Id\\
    & & & \Id & \\
    & & \iddots & & \\
    & -\Id & & & \\
    \Id & & & &
\end{bsmallmatrix},$$
$\mathfrak{g}$ sits inside a Lie algebra $\tilde{\mathfrak{g}}$ with an $A_{2p-1}$ grading. 

\subsubsection{\texorpdfstring{$\Theta$-}{Theta }pinning}
We can use a pinning and the networks from \Cref{sec:NoncomAp} to pin and parameterize $G$ as discussed in \Cref{sec:elementaryNetworks}. Let $\tilde{x}_i,\tilde{y}_i,\check{\gamma}_i$ be a $A_{2p-1}$ pinning of $\tilde{G}$. Then a pinning for $G$ is given by:
\begin{align*}
    \begin{aligned}
        x_i(M) &= \tilde{x}_i(M)\tilde{x}_{2p-i}(\star M)\\
        x_p(M) & = \tilde{x}_p(M)
    \end{aligned} \hspace{1.4pc}
     \begin{aligned}
        \check{\beta}_i(A,B) &= \check{\gamma}_i(A,B)\check{\gamma}_{2p-i}(\star B,\star A)\\
        \check{\beta}_p(A) & = \check{\beta}_{p}(A,\star A^{-1})
    \end{aligned} \hspace{1.4pc}
    \begin{aligned}
        y_i(M) &= \tilde{y}_i(M )\tilde{y}_{2p-i}(\star M )\\
        y_p(M) & = \tilde{y}_p(M)\,.
    \end{aligned}
\end{align*}

\begin{example}
    Concretely consider the group $\SP_{2*6}(\R)$ which has a $C_3$-grading with $\jordan{J} = \jordan{H}_2(\R)$. We can choose the following pinning,
    \begin{align*}
        \begin{aligned}
            x_1 &\colon \jordan{M}_2(\R) \rightarrow \SP_{2*6}(\R)\\
            M &\mapsto \begin{bsmallmatrix}\Id_k & M & & &&\\ & \Id_k & 0 & & &\\ & & \Id_k & 0 & & \\ & & & \Id_k &0 & \\& & & & \Id_k &M^{T} \\ & & & & & \Id_k \end{bsmallmatrix}\\
            x_2 &\colon \jordan{M}_2(\R) \rightarrow \SP_{2*6}(\R)\\
            M &\mapsto \begin{bsmallmatrix}\Id_k & 0 & & && \\ & \Id_k &M & & & \\ & & \Id_k &0 & &\\ & & & \Id_k & M^{T} & \\& & & & \Id_k &0 \\& & & & & \Id_k \end{bsmallmatrix}\\
            x_1 &\colon \jordan{H}_2(\R) \rightarrow \SP_{2*6}(\R)\\
            M &\mapsto \begin{bsmallmatrix}\Id_k & 0 & & && \\ & \Id_k & 0& & & \\ & & \Id_k &M & &\\ & & & \Id_k & 0 & \\& & & & \Id_k &0 \\& & & & & \Id_k \end{bsmallmatrix}
        \end{aligned}\quad
        \begin{aligned}
            \check{\beta}_1 \colon M_2(&\R)\times M_2(\R) \rightarrow \SP_{2*6}(\R)\\
            (A,B) &\mapsto \begin{bsmallmatrix}A &  & & && \\ & B &  & &&\\& & \Id_k & &&\\ & & & \Id_k & &\\ & & & & B^{-T} &  \\ & & & & & A^{-T}\end{bsmallmatrix} \\
            \check{\beta}_2 \colon M_2(&\R)\times M_2(\R) \rightarrow \SP_{2*6}(\R)\\
            (A,B) &\mapsto \begin{bsmallmatrix}\Id_k &  & & &&\\ & A &  & && \\& & B& &&\\ & & & B^{-T} & &\\& & & & A^{-T}& \\& & & & & \Id_k \end{bsmallmatrix}\\
            \check{\beta}_3 \colon M_2(&\R) \rightarrow \SP_{2*6}(\R)\\
            A &\mapsto \begin{bsmallmatrix}\Id_k &  & & &&\\ & \Id_k &  & && \\& & A& &&\\ & & & A^{-T} & &\\& & & & \Id_k & \\& & & & & \Id_k \end{bsmallmatrix}\\
        \end{aligned}\quad
        \begin{aligned}
            y_1 &\colon \jordan{M}_2(\R) \rightarrow \SP_{2*6}(\R)\\
            M &\mapsto \begin{bsmallmatrix}\Id_k &  & & &&\\ M & \Id_k &  & & &\\ & 0 & \Id_k & & & \\ & & 0& \Id_k & & \\& & & 0& \Id_k & \\ & & & & M^{T}& \Id_k \end{bsmallmatrix}\\
            y_2 &\colon \jordan{M}_2(\R) \rightarrow \SP_{2*6}(\R)\\
            M &\mapsto \begin{bsmallmatrix}\Id_k &  & & && \\ 0 & \Id_k & & & & \\ & M & \Id_k & & &\\ & & 0& \Id_k & & \\& & & M^{T}& \Id_k & \\& & & & 0& \Id_k \end{bsmallmatrix}\\
            y_3 &\colon \jordan{H}_2(\R) \rightarrow \SP_{2*6}(\R)\\
            M &\mapsto \begin{bsmallmatrix}\Id_k & & & && \\0 & \Id_k & & & & \\ & 0& \Id_k & & &\\ & & M & \Id_k & & \\& & & 0& \Id_k & \\& & & & 0 & \Id_k \end{bsmallmatrix}.
        \end{aligned}
    \end{align*}
    If we chose the system of Jordan weights $\Lambda_1(A,B,C) = (A,1)$, $\Lambda_2(A,B,C) = (B,1)$, and $\Lambda_3(A,B,C) = C$ and reduced expression $\varsigma_1\varsigma_2\varsigma_3 \varsigma_1\varsigma_2\varsigma_3 \varsigma_1\varsigma_2\varsigma_3$ we obtain the  network in \Cref{fig:C3network} parameterizing $\Confx_3(\decFlagTheta)$.
    \input{Figures/figureC3network}
\end{example}
\subsubsection{Four move}\label{sec:FourMoveHigherNoncom}
We now wish to prove \Cref{thm:clusterChartsTriplesRelatedByMutationTheta} for $C_p$. For $i<p$, the three move $\varsigma_{i-1}\varsigma_{i}\varsigma_{i-1} \mapsto \varsigma_{i}\varsigma_{i-1}\varsigma_{i}$ is realized by two and three moves in the unfolded $A_{2p-1}$ network. Therefore the calculation in \Cref{sec:ThreeMoveHigherNoncom} covers these cases. The remaining case is the four move $\varsigma_{p-1}\varsigma_{p}\varsigma_{p-1}\varsigma_{p} \mapsto \varsigma_{p}\varsigma_{p-1}\varsigma_{p}\varsigma_{p-1}$. This move lives inside a left slid network of type $C_2$ and so it suffices to verify the following:
\begin{lemma}
    The left slid grounded wiring networks of type $C_2$ for $D_1 = \varsigma_{p-1}\varsigma_{p}\varsigma_{p-1}\varsigma_{p} $ and $D_2 = \varsigma_{p}\varsigma_{p-1}\varsigma_{p}\varsigma_{p-1}$ are related by a sequence of four square moves.
\end{lemma}
\begin{proof}
    We begin by lifting both words to words in the $A_3$ Weyl group. These words are related by a sequence of two moves and three moves and thus can be realized by a sequence of square moves. In \Cref{fig:FourMoveNetwork}, we show an explicit sequence of square moves in faces $3,2,1,3$ which realizes the transformation from $D_1$ to $D_2$. Since the weighted network associated to a given word and choice of Jordan weights is unique it must be the case that the weights of the final network agree with those on the chart associated to $D_2$ even though we passed though asymmetric charts to get there. 
    \input{Figures/figureFourMoveNetwork}
\end{proof}

\subsubsection{Rotation and \texorpdfstring{$C_p$}{Cp} networks}

We already proved that rotation by mutation for groups of type $C_p$ is possible by simply rotating the associated chart in for the group of type $A_{2p-1}$ which contains it. However, we make here an interesting observation. 

Groups of type $B_2$ can also be considered as type $C_2$ where the Jordan algebra on the long root is commutative. For type $B_2$ we rotate in an entirely different way; a single mutation at either the small or large node has the effect of flipping to a rotated chart with the other reduced word, recall \Cref{fig:B2mutationTheta}. However in the $C_2 $ case we are entirely unable to do either of these options. In \Cref{fig:FourMoveNetwork} we see that we cannot mutate at the square labeled 2 (which would be the small node for $B_2$) and we cannot mutate at square 1 immediately after mutating at 3 and obtain a network corresponding to a reduced expression in $A_3$ (which would be like mutating at the large node). 

We propose the following remedy: define mutation at the square 2 by the composition of square moves and regrounding which rotates the $C_2$ chart $Q_{D_1,1}$ to $Q_{D_2,2}$ and mutation at squares 1 and 3 by the composition moving $Q_{D_1,1}$ to $Q_{D_2,3}$. 
The composition of these two mutations clearly has order 3 in analogy to the $B_2$ case. 

We could try to define a whole theory of type $C_p$ networks using these rules as mutations. We note that one must work with grounded networks to make sense of this since the new composite mutations involve regrounding the network.

\subsection{Exceptional Groups}
The exceptional groups with $A_2$, $C_3$ and $F_4$ grading cannot be described using noncommutative networks or polygonal cluster algebras. They light their own path, following their own exceptional rules. 

However, our rules for constructing cluster charts still gives us a collection of $L$-distances which parametrize a given configuration of flags. However we are not aware of a simple set of rules for changing between charts that we could call cluster transformations. It would be interesting to try and isolate such a set of transformations. They would include
\begin{itemize}
    \item Octonionic $3$ move for $A_2$,
    \item Octonionic rotation for $A_2$,
    \item Octonionic rotation for $C_3$,
    \item Noncommutative rotation for $F_4$.
\end{itemize}
The flip and 4-move actually are already understood as they live in $A_1$ and $C_2$ types respectively. 

\begin{remark}
    The quaternionic $F_4$ case is interesting to consider since it both contains a $B_3$ which is given by a noncommutative polygonal cluster algebra and a $C_3$ which is given by a noncommutative network. Somehow the correct structure here would marry these two worlds together. 
\end{remark}

\part{Applications}\label{part:Applications}
We will now discuss some concrete applications of the theory developed in the preceding parts. First, we will use the coordinates from \Cref{part:ClusterlikeCoordinates} to obtain coordinates on subsets of $U_\Theta$ and $G$. We use the combinatorics to compute the Gauss decomposition of group elements, their opposite decomposition, as well as the decomposition of the product. This has useful applications to positivity which we will discuss in the subsequent section. There we will consider the special cases in which $(G,\Theta)$ is such that all Jordan algebras admit a positive structure. Then a positive Jordan pinning induces a positive structure on $G$ and its decorated flag variety, which is actually the same as described in \cite{guichard2022generalizing}. We will discuss how the cluster coordinates from \Cref{part:ClusterlikeCoordinates} detect this positivity on the decorated flag variety. Using the previous discussion we will apply this to detect positivity of elements in $G$, allowing us to obtain alternative proofs for some statements about the positive semigroups from \cite{guichard2022generalizing, guichard2026geometric, guichard2026algebraic}. On the other hand, we will explain how to extend the coordinates to moduli spaces of decorated local systems.

\section{Cluster-like coordinates on groups}\label{sec:coordinatesGroup}
 We will provide coordinates on an open dense subset of $G$ using our previous construction of coordinates on $\Confx_4(\decFlagTheta)$. This will eventually also give us positivity tests when $G$ has a $\Theta$-positive structure in the sense of Guichard-Wienhard \cite{guichard2022generalizing}: given an element $g\in G$, we can decide whether $g$ is in the positive semigroup by evaluating certain functions on  $G$ in $g$.

\subsection{Group elements as configurations of flags}
Fix $G$, a Jordan split group of type $R$, with associated parabolic $P_\Theta$.
For two elements $w_1,w_2\in W(R)$ of the Weyl group of $R$, the \keyword{double $R$-Bruhat cell} is given by
\begin{equation*}
    C^{w_1,w_2}:=P_\Theta w_1 P_\Theta\cap P_\Theta^\opp w_2 P_\Theta^\opp\,.
\end{equation*}
In this case we note that $G$ does \textit{not} decompose as a disjoint union of these sets, which may not even be cells. 

If $w_0\in W$ is the longest element, the set
\begin{equation*}
    C^{e,w_0}= P_\Theta\cap P_\Theta^\opp w_0 P_\Theta^\opp
\end{equation*}
is open and dense in $P_\Theta$.

The cluster structure we have given on $\Confx_3(\decFlag)$ is actually and extension of a cluster structure on $P_\Theta$ which we describe in detail in the next section. In analogy to \Cref{thm:ParameterizationOfConf3Theta}, we have
\begin{equation*}
    \Confx_3(\decFlag)\cong C^{e,w_0} \times L_\Theta.
\end{equation*}
Under this isomorphism, the cell $C^{e,w_0}$ corresponds to configurations for which one $L$-distance is $\Id$. 

Of course, one can carry out the same construction for the cell $C^{w_0,e}$, which is open dense in $P_\Theta^\opp$. We will depict these correspondences as in \Cref{fig:BorelTriangles}, where we think of three flags representing a configuration which corresponds to an element of $P_\Theta$ or $P_\Theta^\opp$ as sitting at the vertices of the triangles. As $P_\Theta$ is a semidirect product of $L_\Theta$ and $U_\Theta$, we may choose this correspondence such that one edge of the triangle essentially parametrizes $L_\Theta$ as indicated.

\input{Figures/figureBorelTriangles}

All of this will be made more precise in the subsequent sections but for now it suffices to know that open dense subsets of the two opposite parabolic subgroups are isomorphic to subsets of $\Confx_3(\decFlag)$. Consequently, we have cluster charts on these subsets. Using cluster amalgamation, we obtain charts on an open dense subset of $G$ as indicated in \Cref{fig:BorelTriangles}: Namely, there is a map
\begin{equation*}
    G_0 = U_\Theta^\opp L_\Theta U_\Theta \supset C^{w_0,e}C^{e,w_0} \to\Conf_4(\decFlag)
\end{equation*}
which comes from identifying $C^{w_0,e}$ and $C^{e,w_0}$ 
with $\Confx_3(\decFlag)$.
\begin{remark}
    Generic elements of $G_0$ lie in the domain and map to $\Confx_4(\decFlag)$ under this map. More precisely, we will see that this is the case for elements which lie in $C^{w_0,e}C^{e,w_0}\cap U_\Theta L_\Theta U_\Theta^\opp$, i.e. also admit an opposite Gauss decomposition.
\end{remark}
Using this perspective and the cluster coordinates developed in the preceding section, one can concretely compute Gauss decompositions: Given $g\in C^{w_0,e}C^{e,w_0}\subset G_0$, we can consider the corresponding configuration $c_g\in\Conf_4(\decFlag)$ and compute the cluster coordinates using the triangulation in \Cref{fig:BorelTriangles}. Now a flip in the triangulation may be realized by cluster mutations, and the corresponding coordinates may be used to compute an opposite Gauss decomposition $g=u_+lu_-$ whenever this exists.

\input{Figures/figureGroupCoordinatesFlip}

Similarly, for two elements $g_1,g_2\in C^{w_0,e}C^{e,w_0}$, we can compute the Gauss decomposition of $g_1g_2$ whenever it exists: Use cluster amalgamation again and consider a configuration of 6 decorated flags represented by gluing the two squares to obtain a hexagon as indicated in \Cref{fig:GroupProduct}. By retriangulating this, we may obtain a 4-gon corresponding to $g_1g_2$, which allows us to compute the Gauss decomposition of $g_1g_2$.

\input{Figures/figureGroupProduct}

We will now give more precise explanations of all the above.

\subsection{Parabolic subgroups}
We define the following map
\begin{align*}
    \iota_\Theta: P_\Theta &\to \Conf_3(\decFlagTheta)\\
    b &\mapsto \big(b^{-1}\overline{w_0}U_\Theta,\doverline{w_0}U_\Theta,s_GU_\Theta\big) =: c_b\,.
\end{align*}
\begin{lemma}
    The map $\iota_\Theta$ is an embedding, which sends the double $R$-Bruhat cell $P_\Theta \cap P_\Theta^\opp w_0 P_\Theta^\opp$ to $\Confx_3(\decFlagTheta)$.
\end{lemma}
\begin{proof}
    To see that $\iota_\Theta$ is injective, observe that $c_b = c_{b'}$ if and only if $b'b^{-1}\overline{w_0}U_\Theta = \overline{w_0}U_\Theta$. This is equivalent to $b'b^{-1}\in U_\Theta^\opp$. But by assumption $b'b^{-1}\in P_\Theta$ and thus $b = b'$.

    For given $b\in P_\Theta$, $c_b\in\Confx_3(\decFlagTheta)$ if and only if $b^{-1}\overline{w_0}U_\Theta\pitchfork\doverline{w_0}U_\Theta$. Equivalently, $b\doverline{w_0}$ is Gauss-decomposable, i.e. $b\doverline{w_0}\in U_\Theta^\opp L_\Theta U_\Theta = P_\Theta^\opp w_0 U_\Theta^\opp\overline{w_0}$. This, in turn, is equivalent to $b\in P_\Theta^\opp w_0 P_\Theta^\opp$.
\end{proof}

Now, since we have cluster coordinates on $\Confx_3(\decFlagTheta)$, we obtain cluster coordinates on $P_\Theta$ via the map $\iota_\Theta$. Let us describe the coordinate functions a bit more explicitly from the point of view of $U_\Theta$:
\begin{definition}\label{def:generalizedMinorsTheta}
    Let $w_1,w_2\in W(R), g\in\doverline{w_1}G_0\,\overline{w_2}^{-1}$. For $1\leq i\leq p$ the \keyword{generalized $R$-minor} $\Delta_i^{w_1,w_2}$ is defined as
    \begin{equation*}
        \Delta_i^{w_1,w_2}(g):=\Lambda_i\left(\left[\doverline{w_1^{-1}}g\overline{w_2}\right]_0 \right)=\Lambda_i\left(\doverline{w_1^{-1}}g\overline{w_2}\right)\,.
    \end{equation*}
\end{definition}
This uses the Gauss decomposition with respect to $P_\Theta$. 
Note that these minors depend on the choice of Jordan weights. 
This is a straightforward generalization of the generalized minors introduced in \cite{fomin1999double} and the quasideterminants of \cite{GelfandRetakhDeterminantsMatricesNoncommutative1991}.

\begin{proposition}\label{thm:functionsAsMinorsTheta}
    Let $D(w_0)=\varsigma_{i_1}\cdots \varsigma_{i_n}$ be a reduced expression of $w_0$. The cluster coordinates of $\iota_\Theta(b)$ in the cluster chart associated to this reduced expression are given by $\Delta_{i_k}^{e,w_k}(b)$ where $w_k=\varsigma_{i_n}\cdots\varsigma_{i_{k+1}}$.
\end{proposition}

\begin{remark}
    These generalized minors correspond to the flag minors which are often used to characterize positivity in $\SL_n(\R)$.
\end{remark}

\begin{proof}
    Let us represent $\iota_\Theta(b)$ as $(A_1,A_2,A_3)=(s_Gb^{-1}\overline{w_0}U_\Theta,\overline{w_0}U_\Theta,U_\Theta)$. Then the interpolating sequence $A_2= A_2^0,A_2^1,\dots,A_2^n=A_3$ is given by $A_2^k = \overline{w_k}U_\Theta$. Consequently, the coordinates are
    \begin{equation*}
        \Lambda_{i_k}(s_Gb^{-1}\overline{w_0}U_\Theta,\overline{w_k}U_\Theta)=\Lambda_{i_k}\big(U_\Theta,\overline{w_0}(b\overline{w_k})U_\Theta\big)=\Lambda_{i_k}(b\overline{w_k})
    \end{equation*}
    as claimed.
\end{proof}

Analogously, one can define the map
\begin{align*}
    \iota_\Theta^\opp: P_\Theta^\opp &\to \Conf_3(\decFlagTheta)\\
    b &\mapsto \big(\doverline{w_0}U_\Theta,b U_\Theta,U_\Theta\big)\,.
\end{align*}
Using the opposition involution, we may also write
\begin{equation*}
    \iota_\Theta^\opp(b)=\tw^2\big(\iota_\Theta(\oppInv{b})\big).
\end{equation*}
Thus $\iota_\Theta^\opp$ has the same properties as $\iota_\Theta$ and we can carry out the analogous discussion. In this case, we can express the cluster coordinates as generalized $R$-minors by considering the configuration $\tw\big(\iota_\Theta^\opp(b)\big)$. For a reduced expression as above, the coordinates are given by $\Delta_{i_k}^{e,w_k}(\oppInv{b})$ using the opposition involution.

\subsection{Full group}
Let us now combine the two maps $\iota_\Theta$ and $\iota_\Theta^\opp$ in order two associate a configuration of 4 decorated flags to a generic group element $g\in G$: Namely, we define the composite map
\begin{align*}
    I_\Theta: G &\to \Conf_4(\decFlagTheta)\\
    g &\mapsto \big(\doverline{w_0}U_\Theta,g\overline{w_0}U_\Theta,gU_\Theta,U_\Theta).
\end{align*}

\begin{lemma}
    The map $I_\Theta$ is an embedding which sends the generic subset $C^{e,w_0}C^{w_0,e}\cap C^{w_0,e}C^{e,w_0}$ to $\Confx_4(\decFlagTheta)$.
\end{lemma}

\begin{proof}
    Let $g\in G$.
    If $I_\Theta(g)=I_\Theta(g')$, then $gU_\Theta=g'U_\Theta$ as well as $g\overline{w_0}U_\Theta=g'\overline{w_0}U_\Theta$. Equivalently, $g^{-1}g'\in U_\Theta\cap U_\Theta^\opp=\{\Id\}$.

    Suppose $I_\Theta(g)\in\Confx_4(\decFlagTheta)$, i.e. the flags in the configuration are transverse. The nontrivial cases are:
    \begin{enumerate}
        \item $U_\Theta\pitchfork g\overline{w_0}U_\Theta$: This is equivalent to $g\in P_\Theta P_\Theta^\opp$.
        \item $\doverline{w_0}U_\Theta\pitchfork gU_\Theta$: This is equivalent to $g\in P_\Theta^\opp P_\Theta$.
        \item $U_\Theta\pitchfork gU_\Theta$: This is equivalent to $g\in P_\Theta w_0 P_\Theta$.
        \item $\doverline{w_0}U_\Theta\pitchfork g\overline{w_0}U_\Theta$: This is equivalent to $g\in P_\Theta^\opp w_0 P_\Theta^\opp$.
    \end{enumerate}
    Consequently, $g$ can be written as a product $g=u_-lu_+$ for $u_-\in U_\Theta^\opp$, $l\in L_\Theta$, $u_+\in U_\Theta$. Since $g\in P_\Theta w_0 P_\Theta$, we have $u_-\in C^{w_0,e}$. Similarly, we find $u_+\in C^{e,w_0}$ and thus $g\in C^{w_0,e}C^{e,w_0} $. Analogously, one sees that $g\in C^{e,w_0}C^{w_0,e}$.
\end{proof}

The proof shows that we have in fact carried out the construction shown in \Cref{fig:BorelTriangles}.

\begin{remark}
    Choose a cluster chart as in \Cref{fig:flipSetup} and evaluate it on $I_\Theta(g)$ for $g= u_-lu_+\in G$ as above. The functions can then be computed as in \Cref{thm:functionsAsMinorsTheta} to be the flag minors
    \begin{align*}
        \Delta_{i_k}^{e,w_k}(g) &= \Lambda_{i_k}(l)\Delta_{i_k}^{e,w_k}(u_+)\,,\\
        \Delta_{j_k}^{e,w_k'}(\oppInv{g}) &= \Lambda_{j_k}(\oppInv{l})\Delta_i^{e,w_k'}(\oppInv{u}_-)\,.
    \end{align*}
    where $D_1(w_0) = \varsigma_{i_1}\cdots\varsigma_{i_n}$ and $D_2(w_0) = \varsigma_{j_1}\cdots\varsigma_{j_n}$ are the two reduced expressions defining the cluster chart, and $w_k$ and $w_k'$ are the reversed suffixes of these. 
\end{remark}

\subsection{Gauss decompositions}
The construction we presented above is adapted to the Gauss decomposition $P_\Theta^\opp P_\Theta$. One can carry out an analogous discussion for the opposite Gauss decomposition $P_\Theta P_\Theta^\opp$ (essentially exchanging the two opposite parabolics). In this case, one starts with the embeddings
\begin{align*}
    \bar{\iota}_\Theta: P_\Theta &\to \Conf_3(\decFlagTheta)\\
    b &\mapsto (b\doverline{w_0}U_\Theta,s_G U_\Theta,\overline{w_0}U_\Theta)
\end{align*}
and $\bar{\iota}_\Theta^\opp(b)=\tw^2\big(\bar{\iota}_\Theta(\oppInv{b})\big)$. Gluing the resulting configurations in an analogous fashion, we have a map
\begin{align*}
    \bar{I}_\Theta: G &\mapsto\Conf_4(\decFlagTheta)\\
    g &\mapsto (gU_\Theta,s_G U_\Theta,\overline{w_0}U_\Theta, g\overline{w_0}U_\Theta)\,.
\end{align*}
We immediately see that
\begin{equation*}
    I_\Theta(g) = \bar{I}_\Theta(\oppInv{g})\,.
\end{equation*}
By realizing the flip of the diagonal of the square as indicated in \Cref{fig:GroupCoordinatesFlip}, we can thus explicitly compute the Gauss decomposition $\oppInv{g}=u_+'l'u_-'$ from the Gauss decomposition $g=u_-lu_+$. Similarly, we can realize the flips indicated in \Cref{fig:GroupProduct} by mutations. This allows us to compute the Gauss decomposition of a product $g_1 g_2$ from the Gauss decompositions of $g_1$ and $g_2$.

\subsection{Generalized double \texorpdfstring{$R$}{R}-Bruhat cells}
We can generalize the cluster structure given on the full group to any double $R$-Bruhat cell, $C^{u,v} = P_\Theta uP_\Theta \cap P_\Theta^\opp vP_\Theta^\opp$ for $u,v\in W(R)$. 

\begin{definition}
     Let $\Conf_4^{u,v}(\decFlagTheta)$ be the configurations of decorated flags  $(A_1,A_2,A_3,A_4) \in \Conf_4(\decFlagTheta)$ which satisfy 
        \begin{enumerate}
            \item $(A_1,A_2)$ and $(A_4,A_3)$ are $w_0$-transverse, 
            \item $(A_1,A_4)$ are $u$-transverse and $(A_2,A_3)$ are $v$-transverse.
        \end{enumerate}
\end{definition}

We write generators of $W(R)\times W(R)$ as $i = (e,\varsigma_{i})$ and $\overline{i} =  (\varsigma_{\oppInv{i}},e)$.
Given a reduced expression for $(u,v)$ we can make a cluster chart by replacing each generator  $i$ (resp. $\overline{i}$) with $E(i)$ (resp. $E(\bar{i})$) and amalgamating. Denote by $\A_{R}^{u,v}$ the noncommutative cluster variety which has this as an initial seed. 

\begin{theorem}\leavevmode
    \begin{enumerate}
        \item $\A_{R}^{u,v}$ is independent of the choice of reduced expression of $(u,v)$.
        \item The $\{\jordan{J}_i\}$ points of each cluster chart of $\A_{R}^{u,v}$ parameterize a dense open subset of  $\Conf_4^{u,v}(\decFlagTheta)$.
        \item  There is an injective map $I_\Theta: C^{u,\oppInv{v}}\to \Conf_4^{u,v}(\decFlagTheta)$  given by  $$I_\Theta(g) = (U_\Theta,\overline{w_0}U_\Theta,g\overline{w_0}U_\Theta,gU_\Theta).$$ The configurations in the image satisfy  $l(A_1,A_2)=l(A_4,A_3)=\Id$.
        \item The $\{\jordan{J}_i\}$ points of each cluster chart which land in the image of $I_\Theta$ parameterize an open subset of $C^{u,\oppInv{v}}$.
    \end{enumerate}
\end{theorem}

\begin{proof}
    The cluster variety is independent of the choice of reduced expression for $(u,v)$ by \Cref{thm:clusterChartsTriplesRelatedByMutationTheta} and \Cref{thm:clusterChartsTuplesRelatedByMutationTheta} which showed any two reduced expressions are related by cluster mutation. 
    
    The second statement follows from our construction of the cluster charts; each cluster chart is built using elementary pieces which parameterize elementary configurations. Gluing these configurations gives the appropriate configuration space. 
    
    For the third statement, first observe that the pair $(P_\Theta,gP_\Theta)$ is $u$ transverse if and only if $g\in P_\Theta u P_\Theta$. Similarly, the pair $(P_\Theta^\opp,gP_\Theta^\opp)$ is $v$-transverse if and only if $g\in P_\Theta^\opp \oppInv{v} P_\Theta^\opp$.
     So both transversality conditions occur exactly when $g\in C^{u,\oppInv{v}}$ and $I_\Theta(g)$ lies in $\Conf_4^{u,v}(\decFlagTheta)$. Moreover, the elements in the image of $I_\Theta$ satisfy the condition on $L$-distances. To see that the map is injective, we observe that a configuration of flags in $\Conf_4^{u,v}(\decFlagTheta)$ can be represented as a quadruple in $(G/U_\Theta)^4$ in a unique way so that $A_1 = U_\Theta$ and $A_2=\doverline{w_0}U_\Theta$ with $A_3$ and $A_4$ fixed. Since the pairs $(A_1,A_2)$ and $(A_3,A_4)$ have the same transversality and $L$ distance there must be a unique element of $G$ which sends the pair $(U_\Theta, \doverline{w_0}U_\Theta)$ to $(A_4,A_3)$. For each point in the image of $I_\Theta$, this determines a unique an element $g\in C^{u,\oppInv{v}}$.

    The fourth statement now follows from the previous three statements. The $\{\jordan{J}_i\}$ points describing configurations which land in the image of $I_\Theta$ can be pulled back by $I_\Theta$ to parameterize an open dense subset of $C^{u,\oppInv{v}}$.
\end{proof}

The noncommutative cluster variety $\A_R^{u,v}$ is a noncommutative deformation of the cluster structures on double Bruhat cells for split groups described in \cite{berenstein2005cluster}. The $R=A_p$ case should recover much of the theory of noncommutative double Bruhat decompositions described in \cite{BR_noncommutative_double_bruhat}. Furthermore, we expect that one can continue this entire story to a theory of ``$R$-Braid varieties'' as a generalization of \cite{casals2025cluster}.

\section{Positivity}\label{sec:ThetaPositivity}
For a split algebraic group $G$ over a field $\K$ which contains a positive subset $\K^{>0}$, Lusztig introduced a positive semigroup $G^{>0}$, as well as positive semigroups in the unipotent subgroup $U_\Delta$ and $U_\Delta^\opp$ \cite{lusztig1994total}. In particular, this provides positive semigroups for any split-real Lie group using $\R^{>0}$. The construction relies on a pinning of the group $G$. We will now explain how a Jordan pinning can provide more general positive semigroups. This includes Lusztig's construction and recover Guichard-Wienhard's positive structures with respect to $\Theta$ \cite{guichard2022generalizing}.

\subsection{Positive \texorpdfstring{$\Theta$}{Theta}-pinnings}
Recall from \Cref{def:JordanAlgebraPositive} that a \keyword{positive structure} on a quadratic Jordan algebra $\jordan{J}$ over a formally real field $\K$ is given by a choice of \keyword{positive cone} $V^{>0}\subset V^\times$, which is closed under addition, scalar multiplication by positive elements, and inverse, and moreover contains $\Id$ and is preserved by the action of the structure group $\Gamma(\jordan{J})$.

If each isotopy class of Jordan algebra appearing contains a Jordan algebra with positive structure, we have a distinguished class of Jordan pinnings.
\begin{definition}\label{def:positivePinning}
    A \keyword{positive Jordan pinning} of $G$ with respect to a root system $R$ is a Jordan pinning such that each $\jordan{J}_i$ admits a positive structure.
\end{definition}
A choice of positive structure on each $\jordan{J}_i$ identifies a positive cone in each simple $R$-root space by $x(V_i^{>0})$. This is reminiscent of Guichard-Wienhard's notion of positivity with respect to $\Theta$: Let $\Theta$ be a subset of simple restricted roots of a real Lie group $G$ and $Z_\Theta$ be the center of the corresponding Levi-subgroup. Let $\mathfrak{u}_\beta$ be a weight space of the adjoint action of $Z_\Theta$ on $\mathfrak{g}$.
\begin{definition}
    $G$ has a \keyword{$\Theta$-positive structure} if each $\mathfrak{u}_\beta$ contains a proper convex cone $\poscone{}\subset \mathfrak{u}_\beta$ which is invariant under the action of $L_\Theta^\circ$, the connected component of the identity in the Levi subgroup. 
\end{definition}

It is not immediately obvious that the pairs $(G,\Theta)$ for which the grading induced by $P_\Theta$ on $G$ admits a positive Jordan pinning are the same as those for which $G$ admits a $\Theta$-positive structure. However, we can obtain the same classification: Simply combine the classification of Jordan-compatible $\Theta$ (\Cref{fig:possibleWeylTypedynks_all}) with that of Jordan algebras with a positive structure (\Cref{thm:posJordanAlgsClassification}) to obtain:

\begin{theorem}[\cite{guichard2022generalizing},Theorem 1.1]
    A simple Lie group, $G$, with root system grading induced by $P_\Theta$ admits a positive Jordan pinning  (respectively a $\Theta$-positive structure) if and only if $(G,\Theta)$ belongs to the following list:
    \begin{enumerate}
        \item (Split Real) $G$ is split real and $\Theta=\Delta$.
        \item (Type $A_1 $) $G$ is a Hermitian Lie group of tube type (with restricted root system $C_r$) and $\Theta=\{\alpha_r\}$ the unique long simple root.
        \item (Type $B_p$) $G$ is locally isomorphic to $\SO(p+1,q+p+1)$ (with restricted root system $B_{p+1}$) and $\Theta=\{\alpha_1,\dots,\alpha_p\}$ the collection of all long simple roots. 
        \item (Type $G_2$) $G$ is one of the real forms $F_4^{(4)},E_6^{(2)},E_7^{(-5)},E_8^{(-24)}$ whose restricted root system is type $F_4$ and $\Theta=\{\alpha_1,\alpha_2\}$ is the set of long simple roots. 
    \end{enumerate}
\end{theorem}
\begin{remark}
    Observe that all pairs $(G,\Theta)$ in the classification are of NC rank 0 (split real) or NC rank 1 but never of higher noncommutative rank.
\end{remark}

\subsection{Positive semigroups}
Let $\{x_i,y_i,\check{\beta}_i\}$ be a positive Jordan pinning of $G$ of type $R$, and fix a choice of positive structure $V_i^{>0}$ on each Jordan algebra $\jordan{J}_i$. Following Lusztig \cite{lusztig1994total, guichard2022generalizing}, we may define a positive subset of the unipotent subgroup $U_\Theta$.
\begin{definition}
    Let $D = \varsigma_{i_1}\varsigma_{i_2}\cdots \varsigma_{i_k}$ be a reduced expression of $w_0$, the longest element of the Weyl group $W(R)$. Recall the maps
    \begin{align*}
        F_D:~&V_{i_1}\times V_{i_2}\times\dots\times V_{i_k} \to U_\Theta\,, \hspace{6mm} (t_1,\dots,t_k) \to \prod_{j=1}^k x_{i_j}(t_j)\\
        F_D^\opp:~&V_{i_1}\times V_{i_2}\times\dots\times V_{i_k} \to U_\Theta^\opp\,, \quad (t_1,\dots,t_k) \to \prod_{j=1}^k y_{i_j}(t_j)
    \end{align*}
    which are called the \keyword{$R$-Lusztig maps}. The images of $V_{i_1}^{>0}\times\dots\times V_{i_k}^{>0}$ are called the \keyword{totally positive unipotent semigroups} $U_\Theta^{>0}$ and $U_\Theta^{\opp,>0}$. The \keyword{positive semigroup $G^{>0}$} in $G$ is defined as
    \begin{equation*}
        G^{>0} := U_\Theta^{\opp,>0}L_\Theta^\circ U_\Theta^{>0}
    \end{equation*}
    where $L_\Theta^\circ$ denotes the identity component of $L_\Theta$.
\end{definition}
\begin{remark}
    These are the same semigroups as the ones defined in \cite{guichard2022generalizing} and investigated in \cite{guichard2026geometric, guichard2026algebraic}. In those works, they show that over $\mathbb{R}$ one can give a more geometric definition of the semigroup and address the question of a parametrization later.
\end{remark}

We can use our framework from \Cref{sec:coordinatesGroup} to give alternative proofs of some statements from \cite{guichard2022generalizing, guichard2026geometric, guichard2026algebraic}.

\begin{theorem}
    Let $(G,\Theta)$ be such that the grading induced by $P_\Theta$ on $G$ admits a positive Jordan pinning of $G$. Let $D$ be a reduced expression of the longest element $w_0\in W(R)$. Then
    \begin{enumerate}
        \item The maps $F_D$ and $F_D^\opp$ are injective on $V_{i_1}^{>0}\times\dots\times V_{i_k}^{>0}$.
        \item The subsets $U_\Theta^{>0}\subset U_\Theta$ and $U_\Theta^{\opp,>0}\subset U_\Theta^\opp$ are independent of $D$.
        \item $U_\Theta^{>0}$ and $U_\Theta^{\opp,>0}$ are semigroups.
    \end{enumerate}
\end{theorem}
\begin{proof}
    We will give all proofs for the unipotent subgroup $U_\Theta$. Of course they are completely analogous for $U_\Theta^\opp$.

    For statement (1), consider the configuration of decorated flags represented by $(A_1,A_2,A_3)=\big(U_\Theta,\overline{w_0}U_\Theta, F_D(t_1,\dots,t_k)\overline{w_0}U_\Theta\big)$ with $t_{j}\in V_{i_j}^{>0}$. Using the reduced expression $D$ to compute the interpolating sequence between $A_2$ and $A_3$, we get $A_2^l=\prod_{j=1}^l x_{i_j}(t_j)$. Thus the partial potential for the elementary configuration $(A_1,A_2^{k-1},A_2^k)$ is $t_k$. We know that this can be computed (as an angle) from the coordinates of the configuration.

    For statement (2), consider the same configuration of 3 decorated flags with coordinates computed with respect to the reduced expression $D$. We will see in the next subsection that the coordinates provide a positive $\jordan{J}$-point of the seed group (here we use that the noncommutative rank is at most 1, thus there is at most one noncommutative Jordan algebra $\jordan{J}$ appearing). The coordinates for another reduced expression $D'$ can be computed through a sequence of mutations as explained in \Cref{sec:changingReducedExpressionTheta}. Since mutations sends positive $\jordan{J}$-points to positive $\jordan{J}$-points, this means that the positive subset defined by the map $F_{D'}$ is the same.

    For statement (3), take $u_1,u_2\in U_\Theta^{>0}$ and consider the configuration represented by
    \begin{equation*}
        (A_1, A_2, A_3, A_4)=(U_\Theta,\overline{w_0}U_\Theta,u_1\overline{w_0}U_\Theta,u_1 u_2\overline{w_0}U_\Theta)\,.
    \end{equation*}
    If we compute coordinates as indicated in \Cref{fig:unipotentFlip}, we obtain a positive $\jordan{J}$-point of the corresponding seed group, just like above.
    
    \input{Figures/figureUnipotentFlip}
    
    Using \Cref{sec:FlipTriangulationGenericTheta} we can realize the flip in the triangulation by mutations, meaning that in particular the triple $(A_1,A_2,A_4)$ provides a positive point. Thus $u_1u_2\in U_\Theta^{>0}$.
\end{proof}
This proof hints already at a notion of positive configurations of decorated flags, which can be detected by our coordinates. We will explore this next.

\subsection{Configurations of Flags}
Throughout this section assume that $(G,\Theta)$ is such that the grading on $G$ induced by $P_\Theta$ admits a positive Jordan pinning. Then the positive structure on $U_\Theta$ induces a notion of positivity on the (decorated) flag variety: Any flag $P'$ transverse to the fixed flag $P_\Theta\in G/P_\Theta=\flagTheta$ can be written as $u \overline{w_0}P_\Theta$ with some $u\in U_\Theta$. We call a tuple $(P_1,...,P_n)\in\flagTheta^n$ of pairwise transverse flags \keyword{positive} if there exist $g\in G$ and $u_2,\dots,u_{n-1}\in U_\Theta^{>0}$ such that
\begin{align*}
    g(P_1,\dots,P_n)=(P_\Theta,\overline{w_0}P_\Theta,u_2 \overline{w_0}P_\Theta,u_2u_3 \overline{w_0}P_\Theta,\dots,u_2\cdots u_{n-1} \overline{w_0}P_\Theta)\,.
\end{align*}

Since the fibers of $\pi$ are $L_\Theta$-torsors, we need elements of $L_\Theta$ into account when extending this definition to the decorated flag variety $\decFlagTheta$:

\begin{definition}\label{def:positiveTupleDecFlags}
    A tuple $(A_1,\dots,A_n)\in\decFlagTheta^n$ of pairwise transverse decorated flags is called \keyword{positive} if there exist $g\in G$, $l_2,\dots,l_{n}\in L_\Theta^\circ$ and $u_2,\dots,u_{n-1}\in U_\Theta^{>0}$ such that
    \begin{align*}
        g(A_1,\dots,A_n)=(U_\Theta,\overline{w_0}U_\Theta l_2, u_3\overline{w_0}U_\Theta l_3, u_3 u_4 \overline{w_0}U_\Theta l_4,\dots,u_3\cdots u_n \overline{w_0}U_\Theta l_n)\,.
    \end{align*}
\end{definition}

Positivity is by definition invariant under the diagonal $G$-action on $\decFlagTheta^n$. Thus, positivity descends to $\Confx_n(\decFlagTheta)$, and we can talk about \keyword{positive configurations} of decorated flags.

\begin{remark}
    This setup also covers the case where $G$ is split real and $\Theta=\Delta$. Then, all the definitions agree with those given in \cite{fock2006moduli}. More concretely, take $P_\Theta=B$, a Borel subgroup with unipotent subgroup $U_\Theta=U$, and use the Weyl group $W=W(R)$ to define transverse flags in this case. The fibers of the projection $\pi:\decFlag\to\flag$ are $H$-torsors, where $H=B\cap B^\opp$ is a Cartan subgroup. 
\end{remark}

We can obtain a characterization of $\Theta$-positive configurations of 3 decorated flags in terms of the cluster charts constructed in \Cref{part:ClusterlikeCoordinates}:
\begin{proposition}\label{thm:CoordinatesDetectTriplePosTheta}
    A configuration in $\Confx_3(\decFlagTheta)$ is positive if and only if it is parameterized by a positive $\jordan{J}$-point. Explicitly, all the functions from any cluster chart take values in $\R^{>0}$, respectively the identity component of $\Gamma^0(\jordan{J})$ and all the angles take values in $V^{>0}$.
\end{proposition}
\begin{proof}
    Suppose a configuration in $\Confx_3(\decFlagTheta)$ is parameterized by a positive $\jordan{J}$-point. Such a configuration can be represented as $(U_\Theta,\overline{w_0}U_\Theta l_2,u\overline{w_0}U_\Theta l_3)$ with $u\in U_\Theta$, $l_2,l_3\in L_\Theta$. The elements $l_2$ and $l_3$ are computed as the product of $\check{\beta}_i(s_i)$ where $s_i$ is in the identity component by assumption. Thus $l_2,l_3\in L_\Theta^\circ$. Furthermore $u = F_D(t_1,\cdots, t_n)$ where $t_i$ are the partial potentials which are positive by assumption. Thus $u\in U_\Theta^{>0}$ as needed.\medskip

    Conversely, suppose $(A_1,A_2,A_3)=(U_\Theta,\overline{w_0}U_\Theta l_2,u\overline{w_0}U_\Theta l_3)$ represents a positive configuration, i.e. $l_2,l_3\in L_\Theta^\circ$ and $u=F(t_1,\dots,t_{n})$ with $t_i\in \K^{>0}$ or $ V^{>0}$ appropriately. Then the functions $\Lambda_i(A_1,A_2)$, $\Lambda_i(A_1,A_3)$ and $\Lambda_i(A_3,A_2)$ take values in $\K^{>0}$ or $\Gamma^\circ(\jordan{J})$. To see that all the other functions fulfill the positivity conditions, we proceed inductively along the interpolation sequence. At each elementary configuration we have the that all coordinates are positive except for the rightmost node. Moreover the potential of this configuration is $t_{i_k}$ is assumed to be positive. Thus the equation in \Cref{thm:partialPotentialFromCoordsTheta} can be rearranged to recover the rightmost function as a product of positive functions. It remains to verify every other angle in each polygonal tile is positive. By construction, one angle in each tile recovers the potential and thus is positive. All other angles are obtained by operations which preserve the positive cone, since all the coordinates have been shown to lie in $\R^{>0}$, respectively the identity component of $\Gamma^0(\jordan{J})$. 
\end{proof}

\begin{remark}
    The positivity criterion from \Cref{thm:CoordinatesDetectTriplePosTheta} is not minimal in the sense that there are too many conditions to check. From \Cref{thm:ParameterizationOfConf3Theta} we expect that it should be sufficient to check $2p + n$ conditions where $n$ is the length of the longest word.
\end{remark}
Let $k$ be the number of times $p$ occurs in the reduced expression of the longest word. This is an invariant of the (non-split) Jordan split type (independent of the reduced expression) and is equal to $1$ in type $A_1$, $p$ in type $B_p$ and $3$ in type $G_2$. This number is also the number of polygonal tiles in any cluster chart for the triple of flags. 
\begin{corollary}\label{thm:minimalTestTriplePositivityTheta}
    A configuration in $\Confx_3(\decFlag)$ is positive if and only if it fulfills the following $2p+n$ conditions:
    \begin{enumerate}
        \item All the $n-p+1-k$ functions associated to small nodes in the interior of the triangle take values in $\K^{>0}$. 
        \item The $p-1$ function associated to small nodes on the bottom edge of the triangle take values in $\K^{>0}$.
        \item The $p$ functions on the left and right edges of the triangle take values in $\K^{>0}$ or $\Gamma(\jordan{J})^\circ$ depending on the weight.
        \item One angle in each of the $k$ polygonal tiles is equal to $\iota(v)$ for $v \in V^{>0}$.
    \end{enumerate}
\end{corollary}
\begin{proof}
    The forward implication is clear as the $2p+n$ conditions are a subset of the conditions in \Cref{thm:CoordinatesDetectTriplePosTheta}. To show the reverse implication, suppose that a configuration represented by $(A_1,A_2,A_3)$ fulfills the smaller set of conditions. We prove that this implies the full set of conditions from \Cref{thm:CoordinatesDetectTriplePosTheta}. First, we prove that all functions in the interior associated to weight $r$ nodes take values in $\Gamma(\jordan{J})^\circ$. This is proved inductively from each edge for each elementary seed of the form $\elemSeed_\circ(p)$. One angle in this tile is assumed to be positive. Moreover the angle can be expressed as a monomial where only one coordinate is not assumed to be in $\Gamma(\jordan{J})^\circ$ by inductive assumption. As $\Gamma(\jordan{J})^\circ$ is a subgroup, this proves the missing coordinate lies in $\Gamma(\jordan{J})^\circ$ as well. The elementary seed $\elemSeed(p)$ occurs only once, and so the inductive proof above is enough to guarantee the left and right nodes are in $\Gamma(\jordan{J})^\circ$. This leaves one missing node, the weight $r$ node on the bottom edge of the triangle, which we can conclude is positive since the angle is positive. 
    
    Every other angle can be constructed from the original angles by conjugation by weight $r$ nodes, which preserves positivity and these coordinates are in $\Gamma(\jordan{J})^\circ$. The configuration then satisfies the condition of \Cref{thm:CoordinatesDetectTriplePosTheta} and so is positive as needed.
\end{proof}

Of course this result extends to configurations of $n$ decorated flags:
\begin{theorem}
    A configuration in $\Confx_n(\decFlagTheta)$ is positive if and only if it is parametrized by a positive $\jordan{J}$-point for some cluster chart. In this case, it provides a positive $\jordan{J}$-point for any cluster chart.
\end{theorem}

\begin{proof}
    From the definition it is clear that a configuration of decorated flags is positive if and only if all the configurations of three decorated flags defined by the fan triangulation (\Cref{fig:HeptagonFlagsFanTriangulationTheta}) are positive. Now apply \Cref{thm:CoordinatesDetectTriplePosTheta} to see that a configuration is positive if and only if it is parametrized by a positive $\jordan{J}$-point for any cluster chart coming from the fan triangulation. Combining \Cref{thm:clusterChartsTuplesRelatedByMutationTheta} and \Cref{thm:admissibleMutationPositivity} we conclude that this is the case for any cluster chart.
\end{proof}

\subsection{Positive elements as configurations of flags}
Recall the connection between group elements and configurations of flags described in \Cref{sec:coordinatesGroup}. Let us briefly revisit this, now through the lens of positivity. Thus assume that $G$ has a positive Jordan pinning.
\begin{lemma}
    The embedding $\iota_\Theta:P_\Theta\to\Conf_3(\decFlagTheta)$ preserves $\Theta$-positivity, i.e. if $b\in L_\Theta^\circ U_\Theta^{>0}$, then $c_b$ is a positive configuration.
\end{lemma}
\begin{proof}
    If $b = lu \in L_\Theta^\circ U_\Theta^{>0}$, then by definition the configuration
    \begin{equation*}
        \big(U_\Theta,l^{-1}\overline{w_0}U_\Theta,u\overline{w_0}U_\Theta\big)
    \end{equation*}
    is positive. By applying the twisted cyclic shift, we obtain another positive configuration by \Cref{thm:twPositiveTriples}, which is
    \begin{align*}
        \Big[\big(l^{-1}\doverline{w_0}U_\Theta, u\overline{w_0}U_\Theta,U_\Theta\big)\Big] = \Big[\big(u^{-1}l^{-1}\overline{w_0}U_\Theta, \doverline{w_0}U_\Theta, s_G U_\Theta\big)\Big] = c_b\in\Confx_3(\decFlagTheta)\,.
    \end{align*}
\end{proof}

This allows us to provide positivity test in terms of generalized $R$-minors, see \Cref{def:generalizedMinorsTheta}:

\begin{corollary}[Positivity test]
    Let $D(w_0)=\varsigma_{i_1}\cdots \varsigma_{i_n}$ be a reduced expression of $w_0$. An element $b\in P_\Theta$ is positive if and only if
    \begin{enumerate}
        \item all the functions $\Delta_{i_k}^{e,w_k}(b)$ take values in $\K^{>0}$ whenever $i_k\neq p$ or $\Gamma(\jordan{J})^\circ$ for $k=1,\dots,n$,
        \item the functions $\Delta_{j}^{e,w_0}(b)$ take values in $\K^{>0}$ for $j=1,\dots,p-1$, and $\Delta_p^{e,w_0}(b)\in\Gamma(\jordan{J})^\circ$,
        \item certain combinations of the noncommutative functions land in $\iota(V^{>0})$.
    \end{enumerate}
\end{corollary}
This is a direct consequence of \Cref{thm:CoordinatesDetectTriplePosTheta}, and the combinations in (3) correspond to the angles. Their explicit form depends on the choice of reduced expression.\medskip

Moving on to the full group $G$ we have similarly:
\begin{lemma}
    The embedding $I_\Theta:G\to\Conf_4(\decFlagTheta)$ preserves positivity, i.e. sends $\Theta$-positive elements to positive configurations.
\end{lemma}
\begin{proof}
    Write $g=u_-lu_+$ again. Using the approach developed in the previous sections, we know that $I_\Theta(g)$ is positive iff the following configurations of 3 flags are positive:
    \begin{align*}
        &\Big[\big(\doverline{w_0}U_\Theta, g\overline{w_0}U_\Theta, gU_\Theta\big)\Big] = \Big[\big((lu_+)^{-1}\doverline{w_0}U_\Theta, \overline{w_0}U_\Theta, U_\Theta\big)\Big]=\iota_\Theta(lu_+)\\
        &\Big[\big(\doverline{w_0}U_\Theta, gU_\Theta, U_\Theta\big)\Big] = \Big[\big(\doverline{w_0}U_\Theta, u_-lU_\Theta, U_\Theta\big)\Big]=\iota_\Theta^\opp(u_-l)\,.
    \end{align*}
    By our previous discussion this is the case exactly when $u_-\in U_\Theta^{\opp,>0}$, $l\in L_\Theta^\circ$, and $u_+\in U_\Theta^{>0}$, so when $g$ is $\Theta$-positive.
\end{proof}

Finally recall that we have related different Gauss decompositions to flips in triangulations. Since mutations and the opposition involution preserve positivity, we have the following immediate corollaries:
\begin{corollary}
    The $\Theta$-positive part of $G$ is
    \begin{equation*}
        G_\Theta^{>0} = U_\Theta^{\opp,>0}L_\Theta^\circ U_\Theta^{>0} = U_\Theta^{>0}L_\Theta^\circ U_\Theta^{\opp,>0}\,.
    \end{equation*}
\end{corollary}
\begin{corollary}
    $G_\Theta^{>0}$ is a semigroup.
\end{corollary}

\section{Decorated twisted local systems}\label{sec:localSystems}
    In their pioneering work, Fock and Goncharov connected cluster algebras and positivity to local systems and higher Teichm\"uller theory \cite{fock2006moduli}. For split algebraic groups $G$ and a surface $S$ with boundary, they considered two moduli spaces of $G$-local systems with additional `decoration' associated to the boundary components. They proved that these spaces are positive varieties and that their positive parts may be considered as analogues of Teichm\"uller space. We will now explain how to adapt their discussion of one of these moduli spaces to any algebraic group which admits a unique system of Jordan weights. First, let us recall the basic setup.

    \subsection{Marked surfaces}
    Let $S=\bar{S}\setminus(D_1\cup\dots\cup D_k)$ be a compact, connected, oriented surface obtained by removing finitely many non-intersecting open disks from the closed surface $\bar{S}$.
    \begin{definition}
        A \keyword{marked surface} $\hat{S}$ consists of a surface $S$ as above, together with a set $M=\{x_1,\dots,x_r\}\subset\partial S$ of \keyword{marked points}. The \keyword{punctured boundary} is $\partial\hat{S}:=\partial S\setminus M$.
    \end{definition}
    Moreover, we will always assume that $\partial S\neq\emptyset$. A marked surface $\hat{S}$ of genus $g$ admits a finite area hyperbolic structure with nonempty geodesic boundary which takes each marked point $x_i$ to a boundary cusp if and only if
    \begin{equation}\label{eq:markedSurfHypCondition}
        4g-4+2k+r>0 \,.
    \end{equation}
    Alternatively, we may consider the surface $S'$ obtained by shrinking all unmarked boundary components to punctures. Topologically it is the same as $S$, but the hyperbolic structure on $S'$ will have a cusp at each puncture.\medskip

    Now, let $n$ be the number of components of $\partial\hat{S}$, and pick a set of \keyword{distinguished points} $\{y_1,\dots,y_n\}=:N$, one per component. Let $M'\subset S'$ be the set of marked points after shrinking unmarked holes. Similarly let $N'$ be the set of distinguished points after shrinking. The distinguished points on unmarked boundary components simply become the punctures of $S'$.

    \begin{definition}
        An \keyword{ideal triangulation} of $\hat{S}$, is a maximal collection of arcs in $S'\setminus M'$, all of which start and end at points in $N'$, considered up to homotopy with fixed endpoints.
    \end{definition}

    \begin{remark}
        When \Cref{eq:markedSurfHypCondition} holds, we can also choose a finite area hyperbolic metric on $S'\setminus N'$. Then an ideal triangulation of this hyperbolic surface in the usual sense is equivalent to an ideal triangulation of the marked surface.
    \end{remark}

    Using this alternative approach, endow $S'\setminus N'$ with a hyperbolic metric. This provides us with a developing map from the universal cover $\tilde{S'}$ to the hyperbolic plane $\Hyp$. The preimage of the set $N'$ of distinguished points under the covering map $\tilde{S'}\to S'$, is a set $\farey(\hat{S})\subset\partial\Hyp$, called the \keyword{Farey set}. It is a cyclic set on which $\pi_1(S)$ acts.
    
    \subsection{Decorated twisted local systems}
    Let $G$ be a reductive algebraic group over $\K$ and $S$ a surface. Recall that a \keyword{$G$-local system} on $S$ is a flat principal $G$-bundle with base space $S$.
    \begin{remark}
        The moduli space of $G$-local systems up to isomorphism can be identified with the space of representations $\pi_1(S)\to G$ up to conjugation via the holonomy representation.
    \end{remark}
    
    Fock and Goncharov connected this to configurations of decorated flags by considering a slightly more complicated moduli space. We first need to handle the element $s_G$. Instead of $S$, consider the \keyword{punctured tangent bundle}, denoted $\puncTB$, which is obtained from the tangent bundle $TS$ by removing the zero section. Its fundamental group provides a central extension of $\pi_1(S)$ described by the short exact sequence
    \begin{equation*}
        \Z\overset{\iota}{\to}\pi_1(\puncTB)\to\pi_1(S)\,.
    \end{equation*}
    The subgroup $\iota(\Z)<\pi_1(\puncTB)$ is generated by $s$, a loop around the fiber of $\puncTB\to S$. 
    
    \begin{definition}\label{def:twistedLS}
        A \keyword{twisted $G$-local system} on $S$ is a $G$-local system on $\puncTB$ such that the holonomy around any fiber $\puncTB\to S$ is $s_G$.
    \end{definition}
    This is well-defined since $s_G$ is at most of order 2.

    \begin{remark}
        The moduli space of twisted $G$-local systems is isomorphic to the space of twisted representations up to conjugation: A \keyword{twisted representation} is a representation $\varrho:\pi_1(\puncTB)\to G$ for which $\varrho(s)=s_G$.
    \end{remark}

    \begin{remark}
        If $G$ is a group where $s_G=1$, the twisting is trivial and we can safely consider local systems on $S$ instead.
    \end{remark}

   Now assume that $(G,\Theta)$ is Jordan compatible. We will add some additional data to twisted local systems, which morally consists of a choice of decorated flag in $\decFlagTheta$ at each distinguished point in $N$. This follows Le's description for the split case \cite{le2019cluster}, adapted to our setting. For a twisted $G$-local system $L$ on a marked surface $S$, consider the associated flag bundle $L_\decFlag$ on $\puncTB$, i.e. the associated bundle with fiber $\decFlagTheta$. We may lift the boundary of $S$ to $\puncTB$ by considering the positively oriented unit tangent vector (after endowing $S$ with an auxiliary metric). If $\hat{S}$ is a marked surface with underlying surface $S$, we may lift every component of the punctured boundary $\partial\hat{S}$ to $\puncTB$ in this way to obtain a set of arcs and loops, which we call the \keyword{lifted boundary}. Furthermore, we can restrict the flag bundle $L_\decFlag$ to these.
    
    \begin{definition}\label{def:decTwistedLS}
        A \keyword{decorated twisted $G$-local system} on a marked surface $\hat{S}$ is a twisted local system $L$ on $S$ together with flat sections of $L_\decFlag$ restricted to any component of the lifted boundary. We call this choice the \keyword{decoration}. $\Loc{G}{\hat{S}}$ denotes the moduli space of decorated twisted $G$-local systems on $\hat{S}$ up to isomorphism.
    \end{definition}

    Suppose that $\hat{S}$ has a boundary component with no marked point on it. The lifted boundary component represents an element of $\gamma\in\pi_1(\puncTB)$. For any twisted $G$-local system $L$ the associated holonomy $\varrho(\gamma)\in G$ is obtained by considering a point $p$ in the fiber of the local system over a base point $x$ on $\gamma$. Parallel transport takes this to another point $p.\varrho(\gamma)$ in the same fiber, using the right action of $G$. The decoration on this boundary component is a flat section of $L_\decFlag$ along the lifted boundary $\gamma$. From the construction of the associated bundle $L_\decFlag$, it is clear that the flat section is fully determined by its value at $x$ which can be represented as $(p,A)$ for some $A\in\decFlagTheta$. Since it is a flat section, we can follow $\gamma$ and arrive at the point $\big(p.\varrho(\gamma),A\big)$ in the fiber over $x$, which is equivalently represented as $\big(p,\varrho(\gamma)A\big)$. Since the section has to close up, the decorated flag $A$ has to be fixed by $\varrho(\gamma)$.\bigskip

     We will now connect this to our previous discussion of configurations of decorated flags:
    \begin{example}\label{ex:ASpaceDisk}
        Consider the marked surface $\hat{D}_n=(D,\{x_1,\dots,x_n\})$, i.e. a disk with $n$ marked points on the boundary. Since $\pi_1(D)$ is trivial, there is a unique twisted local system $L$ on $D$. Thus a point of $\ASpace{R}{\hat{D}_n}$ is completely determined by the decoration. Denote by $C_1,\dots,C_n$ the boundary intervals as shown in \Cref{fig:ASpaceDisk}. The decoration consists of a flat section of $L_\decFlag$ along any lifted $C_i$, or equivalently a point $a_i$ in each fiber $L_\decFlag|_{v_i}$ with fixed distinguished points $y_i\in\partial D\cap C_i$, and $v_i\in T'_{y_i}D$ the counterclockwise oriented tangent vectors of $\partial D$ at $y_i$. 
         
         In order to identify this data with a configuration of flags, we must move the points in the different fibers to a chosen common base point. Let $\gamma_i$ be a path in $D$ from $y_{i}$ to $y_{i+1}$ as in \Cref{fig:ASpaceDisk}. These lift to paths in $\puncTB$ using their unit tangent vectors. Let $P_i$ denote parallel transport along the lift of $\gamma_i$ to $\puncTB$. In this way we obtain $n$ elements in $L_\decFlag|_{v_1}$:
        \begin{equation*}
            (a_1,\,P_n\cdots P_2a_2,\,P_n\cdots P_3 a_3,\,\dots,\,P_n a_n) \;.
        \end{equation*}

        \input{Figures/figureASpaceDisk}
        
        Since $L_\decFlag|_{v_1}\cong\decFlagTheta$ with the isomorphism determined up to the right $G$-action, the $n$-tuple of elements obtained in this way identifies with an element $x\in\Conf_n(\decFlagTheta)$. If we used $y_2$ as our base point instead, we would obtain the tuple
        \begin{equation*}
            (a_2,\,P_1P_n\cdots P_3a_3,\,\dots,\,P_1P_na_n,\,P_1a_1)
        \end{equation*}
        in $L_\decFlag|_{v_2}$. To compare with our original choice, we use parallel transport to obtain the tuple 
        \begin{equation*}
            (P_n\cdots P_2a_2,\,s_GP_n\cdots P_3a_3,\,\dots,\,s_GP_na_n,\,s_G a_1) = \tw(x)
        \end{equation*}
        in $L_\decFlag|_{v_1}$.
        This shows that the choice of base point only changes the tuple by the twists and thus
        \begin{equation*}
            \Loc{G}{\hat{D}_n}\cong\Conf_n(\decFlagTheta)/\tw\,.
        \end{equation*}
    \end{example}

    \begin{remark}
        We see here that the space $\Loc{G}{\hat{D}_n}$ is a more natural object to study then the configuration space $\Conf_n(\decFlagTheta)$. This explains the complications related to choosing a vertex of the $n$-gon in our discussion of cluster coordinates for $\Conf_n(\decFlagTheta)$.
    \end{remark}

    \subsection{Cluster coordinates}
    We will now reduce the study of decorated twisted local systems to the study of (twisted) configurations of decorated flags, which we already discussed extensively in \Cref{part:ClusterlikeCoordinates}. This discussion follows \cite{fock2006moduli}.

    Consider a marked surface $\hat{S}$ as above, and fix an ideal triangulation $T$ of $\hat{S}$. Denote the set of triangles by $\tr(T)$ and the set of edges by $\ed(T)$.
    For every $t\in\tr(T)$, we obtain an element of $\Loc{G}{\hat{D}_3}$ as follows. Given a decorated twisted local system on $\hat{S}$, we obtain a twisted local system on $t$ via the pullback along the embedding given by the triangulation $T$. The decoration on $S$ induces flat sections of the associated bundle near the vertices of $t$, which we consider distinguished points of $t$.  By choosing a marked point $x_e$ on each edge of $t$, we identify $\hat{D}_3$ with $t$. Similarly, we can thicken every edge $e\in\ed(T)$ to obtain an embedded digon in $S$. This defines an element of $\Loc{G}{\hat{D}_2}$ by pullback along the embedding. 

    Therefore the triangulation $T$ provides a map
    \begin{equation}\label{eq:triangulationMapSplit}
        \phi_{T,v}:\Loc{G}{\hat{S}}\to \prod_{t\in\tr(T)}\Loc{G}{\hat{D}_3}\,\times\prod_{e\in\ed(T)}\Loc{G}{\hat{D}_2}\cong \prod_{t\in\tr(T)}\tilde{\Conf_3}(\decFlagTheta)\,\times\prod_{e\in\ed(T)}\tilde{\Conf_2}(\decFlagTheta)
    \end{equation}
    with $\tilde{\Conf_k}(\decFlag):=\Conf_k(\decFlag)/\tw$. The last identification was described in \Cref{ex:ASpaceDisk}. Every point in the image clearly has the property that further restricting the decorated local system from a triangle to one of its edges yields the decorated local system associated to that edge. This can be seen as a gluing condition identifying the image of $\phi_T$.

    We will amalgamate the cluster charts we defined on triples of flags to obtain a chart on $\Loc{G}{\hat{S}}$ via the embedding in \Cref{eq:triangulationMapSplit}. However this chart is only defined on a generic subset of $\Loc{G}{\hat{S}}$ which is determined by the triangulation.
    \begin{definition}
        The space $\LocT{G}{\hat{S}}{T}$ of \keyword{$T$-transverse} decorated twisted local systems is the maximal subspace of $\Loc{G}{\hat{S}}$ such that
        \begin{equation*}
            \phi_T(\LocT{G}{\hat{S}}{T})\subset\prod_{t\in\tr(T)}\tilde{\Confx_3}(\decFlag)\,\times\prod_{e\in\ed(T)}\tilde{\Confx_2}(\decFlag)\,.
        \end{equation*}
    \end{definition}
    To obtain a cluster chart on the space $\LocT{G}{\hat{S}}{T}$, we choose
    \begin{enumerate}
        \item an orientation of every edge in $\ed(T)$,
        \item a vertex $v(t)$ for every triangle $t\in\tr(T)$,
        \item a reduced expression of the longest word $w_0$ for every $t\in\tr(T)$.
    \end{enumerate}
    We write $\mathcal{D}$ for this set of extra choices. 
    First, focus on a triangle $t\in\tr(T)$: The choice of $v(t)$ fixes an identification of the local system restricted to $t$ with $\Confx_3(\decFlagTheta)$. Using the reduced expression for $t$, we obtain a cluster chart for the restriction via the construction in \Cref{sec:flagTriplesTheta}. For such charts the sides are always oriented as in \Cref{fig:TriangleFlags}. Whenever this does not match the chosen orientations, we change the function as explained in \Cref{sec:orientations} and perform a switch (and weave) at the corresponding node (depending on the type).

    \input{Figures/figureTriangleFlags}
    
    Since the orientation of each edge was fixed the values of the cluster coordinates agree along each edge shared by two triangles after applying this construction to all triangles. Consequently, the charts amalgamate together.
    
    \begin{definition}
        A set of functions on $\LocT{G}{\hat{S}}{T}$ arising from the above construction is called a \keyword{cluster chart} adapted to the triangulation $T$. We write $\ASpace{T}{\mathcal{D}}$ for the noncommutative torus associated to this cluster chart. We write $\ASpace{R}{\hat{S}}$ for the noncommutative cluster variety which has this cluster chart as an initial seed. 
    \end{definition}

    These actually provide coordinates for $\LocT{G}{\hat{S}}{T}$:
    \begin{theorem}\label{thm:coordsASpaceSplit}\leavevmode
    \begin{enumerate}
        \item The $\{\jordan{J}_i\}$ points of $\ASpace{T}{\mathcal{D}}$ parametrize an open dense subset of $\LocT{G}{\hat{S}}{T}$.
        \item $\ASpace{R}{\hat{S}}$ contains all possible cluster charts related to different choices $T'$ and $\mathcal{D}'$. 
    \end{enumerate}
        
    \end{theorem}
    \begin{proof}
        We proved in \Cref{thm:coordinatesTriplesTheta} that an open dense subset of $\tilde{\Conf}_3(\decFlag)\cong\Loc{G}{\hat{D}_3}$ can be reconstructed from the cluster coordinates in each triangle. It remains to show that these restricted local systems determine the full local system. This follows identically to the split case \cite{fock2006moduli}, essentially by reconstructing the full local system from the restricted ones. Let us recall the argument for completeness:\medskip

        First, construct from the fixed triangulation $T$ a groupoid $\mathscr{G}_T$ which is equivalent to the fundamental groupoid of $\puncTB$: Fix one point $p(e)$ on every edge $e\in\ed(T)$ and two vectors $v_0(e),v_1(e)=-v_0(e)\in T'_{p(e)}S$, which we can think of as being tangent to $e$. The collection $\left\{v_i(e)\,|\,e\in\ed(T),i\in\Z/2\Z\right\}$ defines the objects of $\mathscr{G}_T$. The set of morphisms is defined as being generated by the following two types:
        \begin{enumerate}
            \item In $T'_{p(e)}S$ there are two preferred homotopy classed of paths connecting $v_1(e)$ to $v_2(e)$, the counterclockwise and clockwise rotation.
            \item Consider $t\in\tr(T)$. For two sides $e_1,e_2$ of $t$, such that $e_2$ is clockwise from $e_1$, suppose that $v_{i_1}(e_1)$ and $v_{i_2}(e_2)$ point at the same vertex of $t$. Define a morphism from $v_{i_1}(e_1)$ to $v_{i_2}(e_2)$ by picking the homotopy class of the path moving the minimal amount counterclockwise.
        \end{enumerate}
        The objects of $\mathscr{G}_T$ and these generators of $\mathrm{mor}(\mathscr{G}_T)$ for one triangle are shown in \Cref{fig:reconstructionHexagon}. Observe that the path around the digon on a side of the triangle corresponds to a rotation of the vector by $2\pi$ while the path around the hexagon inside the triangle corresponds to a rotation by $4\pi$.

        \input{Figures/figureReconstructionHexagon}

        The idea of the proof is now to identify decorated twisted local systems on $\hat{S}$ with certain functors $\mathscr{G}_T\to G$. We first construct the functor for a decorated twisted $G$-local system $L$. Focusing on an edge $e$, the decoration determines a point $a_i(e)$ in the fibers $L_\decFlag|_{v_i(e)}$ for $i\in\Z/2\Z$ of the associated bundle over the two points $v_i(e)$. The canonical projection $\decFlagTheta\to\flagTheta$ determines two points $b_i(e)$ in the corresponding fibers of $L_\flag$. Parallel transport of these points along either path from $v_i(e)$ to $v_{i+1}(e)$, gives us, with a slight abuse of notation, a transverse pair $\big(a_i(e),b_{i+1}(e)\big)\in(L_\decFlag\times L_\flag)|_{v_i(e)}$, which defines a unique point $l_i(e)\in L|_{v_i(e)}$ for $i\in\Z/2\Z$.
        
        Suppose that $\gamma:v_{i_1}(e_1)\leadsto v_{i_2}(e_2)$ is an oriented path in the graph on $S$ obtained by inscribing the hexagons as explained above. This is the same as a morphism in the category $\mathscr{G}_T$. To this we assign an element $g_\gamma\in G$ defined by
        \begin{equation*}
            l_{i_2}(e_2)=\big(P^\gamma l_{i_1}(e_1)\big).g_\gamma
        \end{equation*}
        where $P^\gamma$ denotes the parallel transport along $\gamma$ induced by the connection on $L$, and we have used the simply transitive $G$-action on the fibers of $L$. This clearly defines a functor $\mathscr{G}_T\to G$, i.e. we have:
        \begin{itemize}
            \item[(a)] Reversing a path gives the inverse element: $g_{\bar{\gamma}}=g_\gamma^{-1}$.
            \item[(b)] It respects composition: $g_{\gamma*\delta}=g_\gamma g_\delta$.
        \end{itemize}
        Since it comes from a twisted local system:
        \begin{itemize}
            \item[(c)] The loop around a digon is mapped to $s_G$.
            \item[(d)] The loop around the hexagon is mapped to $\Id$.
        \end{itemize}
        Moreover, it satisfies:
        \begin{itemize}
            \item[(e)] For any path $\varepsilon$ along a digon: $g_\varepsilon\in \overline{w_0}H$.
            \item[(f)] For any path $\nu$ connecting to sides of a triangle: $g_\nu\in U$.
        \end{itemize}
        The two types of edges are also indicated in \Cref{fig:reconstructionHexagon}.

        Conversely, any functor which satisfies the conditions (c) and (d) defines a twisted local system. If it also has (e) and (f) one obtains a decorated twisted local system. In particular, we can reconstruct a decorated twisted local system from the functor obtained above. Since the construction only depends on the restrictions of $L$ to the triangles, this concludes the proof of the first statement.

    For a fixed triangulation $T$, different cluster charts arise by choosing a different vertex or different reduced expressions of $w_0$. In this case the theorem is implied by \Cref{thm:clusterChartsTriplesRelatedByMutationTheta} and \Cref{thm:coordinateSystemsS3MutationTheta}. Relating cluster charts for two different triangulations can be done by realizing a sequence of flips connecting the two triangulations. This was proved in \Cref{thm:clusterChartsTuplesRelatedByMutationTheta}. This finishes the proof of the second statement.
    \end{proof}
    
    \subsection{Computing holonomies}\label{sec:holonomies}
    The reconstruction procedure as reviewed in the proof of \Cref{thm:coordsASpaceSplit} works by inscribing a graph into the surface which is determined by the ideal triangulation, and assigning elements in $G$ to every (oriented) edge. Then, the holonomy of the local system can be computed by homotoping a loop in $\puncTB$ onto the graph and taking the product of group elements along it. Here is a small example, which is related to our discussion in \Cref{sec:coordinatesGroup}.
    
    \input{Figures/figureLocalSystemSquareSplit}
    
    \begin{example}\label{ex:localSystemSquare}
        Consider the space $\Loc{G}{\hat{D}_4}$. It is (non-canonically) isomorphic to $\Conf_4(\decFlagTheta)$ as described in \Cref{ex:ASpaceDisk}. Thus we have an embedding
        \begin{equation*}
            G\to\Loc{G}{\hat{D}_4}\,,\quad g\mapsto (\doverline{w_0}U_\Theta,g\overline{w_0}U_\Theta,gU_\Theta,U_\Theta)
        \end{equation*}
        induced by the map $I_\Theta:G\to\Conf_4(\decFlagTheta)$ defined in \Cref{sec:coordinatesGroup}. We picture the disk with 4 marked points as a square with distinguished points as the vertices, and fix an ideal triangulation to obtain the local system shown in \Cref{fig:localSystemSquareSplit}. Here, we the local system is identifies with the configuration by moving every section to the upper left corner of the square.
        
        Now we factorize $g=u_1^*lu_2$ with $u_1,u_2\in U_\Theta$ and $l\in L_\Theta$. Then the group elements associated to the oriented arcs can be computed explicitly from the two elements in $\Confx_3(\decFlag)$ associated to the two triangles. This also shows that the group element for the highlighted path in $\puncTB$ is
        \begin{equation*}
            \overline{w_0}u_1^{-1}\doverline{w_0}lu_2=u_1^*lu_2=g\,.
        \end{equation*}

        Of course, this element is not a holonomy of the local system on the disk, where the only nontrivial loop is around the fiber of $\puncTB\to S$. However, since the upper and lower edge of the square correspond to the same configuration of two decorated flags, we may glue them to obtain a twisted decorated local system on the cylinder. This has nontrivial topology and the highlighted path closes up to represent the loop around the cylinder, lifted along its tangent vector to a loop in $\puncTB$. Thus, the holonomy for this loop is $g$ by our previous calculation.
    \end{example}

    Before we discuss some further examples of how to compute holonomies, let us simplify the situation by passing from the group $G$ to the centerless group $G'=G/Z_G$. A twisted $G$-local system $L$ on $\puncTB$, defines a $G'$-local system $L'$ on $S$. First, choose a nowhere vanishing vector field. This provides a pullback bundle on $S'$. Now, take the quotient by the $Z_G$ action on the fibers to obtain $L'$. Importantly, the isomorphism class of $L'$ is independent of the chosen vector field since the holonomy around fibers of $L$ is $s_G$, which is mapped to $\mathrm{Id}_{G'}$. This construction induces a projection
    \begin{equation*}
        p:\Loc{G}{\hat{S}}\to\mathscr{X}_{G',\hat{S}}
    \end{equation*}
    of the $\mathscr{A}$-space to the $\mathscr{X}$-space of framed $G'$-local systems on $S$. For full details of the notion of a framed local system see \cite{fock2006moduli}. We will only need the fact  that a decoration of the twisted local system $L$ determines a decoration of the local system $L'$.

    To compute the holonomy of the local system $L'=p(L)$, we consider the graph defined by the triangulation and project all the elements assigned to the edges to $G'$. Since the group elements for the two oriented edges of the graph along an edge of the triangulation only differ by $s_G$, they are projected to the same element in $G'$, and we may collapse the edges together. After that, the holonomy is computed in the same way as before. Let us consider two types of loops in $S$ as indicated in \Cref{fig:loopsOnSurface}.

    \input{Figures/figureLoopsOnSurface}

    If the loop can be thickened into an annulus with one distinguished point on each boundary component, we can change the triangulation to be in the situation of \Cref{ex:localSystemSquare} and can use the indicated path to compute the holonomy of the loop.

    Now suppose that the loop cuts off a surface which does not contain any marked points. Up to homotopy, we may assume that the boundary contains a distinguished point, so that the triangulation looks as in the pasting scheme shown in \Cref{fig:pastingScheme} for a genus 2 surface. The path following the boundary, can be homotoped onto the graph coming from the triangulation as shown in the figure, which provides a way to compute its holonomy.

    \input{Figures/figurePastingScheme}

    \subsection{\texorpdfstring{$L_\Theta$-}{L Theta }action at marked points and potentials}
    Given a decorated local system $\mathcal{L}$ we can alter its decoration at a particular marked point by rescaling the decorated flags there by an element of $L_\Theta$.  

    Specifically let $m$ be a marked point and $l\in L_\Theta$. We define $(l)_{m}\cdot \mathcal{L}$ by sending all flags $A$ which appear at the marked point $m$ to $A\oppInv{l}$. Note that monodromy operators send $A$ to $gA$ which commutes with this scaling, so that this action is well-defined. This is an extension of the rescaling action on configurations of flags to the moduli space of decorated local systems.
    
    \begin{lemma}\label{lem:rescaling_is_monomial_localsystems}
        This action is a monomial transformation of the cluster variables in every cluster chart. 
    \end{lemma}
    \begin{proof}
        This follows from \Cref{prop:rescaling_is_monomial}, as our cluster charts art built out of charts for triples of flags.
    \end{proof}

    We can choose a cluster chart for which the flag at the marked point $m$ is always chosen to be opposite to the side to which the interpolating sequence is calculated. This identities a number of partial potentials $\delta_i^{k}:\mathcal{L} \to  V_i$ which depend on the cluster chart chosen. 
    \begin{proposition}\leavevmode
    \begin{itemize}
        \item The total potential, $W_i(m) = \sum_k\delta_i^{k}$, is independent of the chosen cluster chart.
        \item The action of $L_\Theta$ on this marked point satisfies that $(l)_m\cdot\delta_i^{k} = \beta_i(\varsigma_i(l))(\delta_i^{k})$ for each partial potential. Since $\beta_i$ is a linear map this applies to the total potential as well.
    \end{itemize}
           
    \end{proposition}
    \begin{proof}
        Consider a single triangle $(A_m,A_2,A_3)$. The potentials do not depend on the decorations of $A_2,A_3$ so we can write the configuration on the triangle as $(U_\Theta,P_\Theta^{\opp},uP_\Theta^{\opp})$ for some $u$.
            
        In \Cref{thm:coordinatesTriplesTheta}, we saw the partial potentials within a single triangle recover the factorization parameters in the $\Theta$-Luztig map for $u$. Their sum is independent of the choice of decomposition of $w_0$, since this is the image under the quotient $U_\Theta \to U_\Theta/[U_\Theta,U_\Theta]\simeq \prod_i V_i$
        This proves that this part of the sum does not depend on the choice of $w_0$ in the cluster charts. When we change triangulation by flipping an arc which touches $m$ we create a singe triangle from two, see \Cref{fig:unipotentFlip}. In proving that this is realized by cluster mutations, we proved that the potential of the new triangle is the sum of the adjacent two potentials. Therefore the total sum is the same in all cluster charts. 

        Since each partial potential sits inside a triple, the second statement follows from \Cref{prop:rescaling_potentials}. 
    \end{proof}

    \subsubsection{Cluster \texorpdfstring{$\mathcal{X}$}{X} coordinates}
    In noncommutative rank 0, there is a family of coordinates called $\mathcal{X}$ coordinates which are independent of the decorations of the flags. Analogous coordinates when $G$ is Jordan split over $A_1$ were studied in \cite{alessandrini2019noncommutative}. 
    In general noncommutative rank, we propose that the set of all ratios
    \begin{equation*}
        X_m = \{\iota(\delta_i^{r})^{-1}\iota(\delta_i^{s})\}
    \end{equation*}  
    of partial potentials in a given cluster chart associated with the same marked point $m$ is a potential generalization. 
    \begin{lemma}
        The conjugacy classes of elements of $X_m$ are independent of the rescaling action at $m$. 
    \end{lemma}
    \begin{proof}
        We compute that $\iota(\beta_i(\varsigma_i(l))(v)) = \sigma_i(\beta_i(l)^{-1})\iota(v)(\beta_i(l)^{-1})$. So for any ration in $X_m$,  $(l)_m\cdot( \iota(\delta_i^{r})^{-1}\iota(\delta_i^{s})) = \beta_i(l)\left(\iota(\delta_i^r)^{-1}\iota(\delta_i^s)\right)\beta_i(l)^{-1}  $. This shows the rescaling action at $m$, conjugates the ratio as needed.
    \end{proof}
    We leave the noncommutative cluster structure of these ratios for future work.

    \subsection{Action of the cluster automorphism group}\label{sec:ClusterModularGroupAction}

    In commutative cluster algebra theory, a \keyword{cluster homomorphism} can be defined as a map between two cluster algebras which maps cluster variables to cluster variables. A \keyword{quasi-cluster homomorphism} the composition of a usual cluster homomorphism followed by a monomial transformation in the frozen cluster variables of the target algebra, see \cite{fraser_quasi-homomorphisms_2016}.

    A typical way to produce a quasi-cluster automorphism is to consider a sequence of cluster mutations which takes a starting seed to a new seed which has different cluster variables but with an isomorphic mutable portion of the quiver which defines it.

    For a split group $G$, Goncharov and Shen constructed a subgroup of the cluster automorphism group $\Gamma_{G,\hat{S}}$ which acts on $\Loc{G}{\hat{S}}$ by  \cite{goncharov2019quantum}. This group contains the mapping class group of $\hat{S}$, the outer automorphism group $\Out(G)$, a copy of the Weyl group, $W(G)$, for each puncture on $\hat{S}$, and a (possibly a quotient) of the Artin-Tits braid group of type $G$, $\textrm{Br}(G)$, for each boundary component of $\hat{S}$. It is natural to conjecture that this group is in fact the entire cluster automorphism group.

    Now let $G$ be Jordan split of type $R$ and let $\Gamma_{R,\hat{S}}$ be the subgroup of cluster automorphism group associated to the split group of type $R$.
    \begin{theorem}\label{thm:clustermodulargroup}
        The group $\Gamma_{R,\hat{S}}$ acts on $\ASpace{R}{\hat{S}}$ by noncommutative (quasi) cluster automorphisms.
    \end{theorem}

    The proof follows the split case, by ensuring that the action of each part can be realized noncommutatively as a sequence of mutations. 

    \subsubsection*{The mapping class group}
    The mapping class group acts on a triangulation of $\hat{S}$ by sending it to another triangulation. Any two triangulations are related by a sequence of flips and we proved in \Cref{sec:FlipTriangulationGenericTheta} that each flip is realized by a sequence of generalized cluster mutations. Furthermore the resulting seed is isomorphic (including frozen nodes/ boundary faces) and so gives an isomorphism of the cluster structure.

    \subsubsection*{Out(G)}
    Similarly $\Out(R)$ is realized by automorphisms of the Dynkin diagram $R$. This acts on reduced words for $w_0$ by permuting the letters via the automorphism of $R$. This action on $\ASpace{R}{\hat{S}}$ acts by changing the reduced expression for $w_0$ in each triangle in this way. Since we can realize this change by generalized mutations, we can realize this part of the action. 
    
   \begin{remark}
       The only noncommutative rank greater than 1 group with a nontrivial automorphism group is type $A_p$.
   \end{remark} 

    \subsubsection*{Weyl Group action}

    Consider a surface, $\hat{S}$, with a puncture $p$.  There is a birational action of the Weyl group, $W(R)$, on $\ASpace{R}{\hat{S}}$ associated to $p$ as follows: 

\begin{construction}
    Each decorated local system $\mathcal{L}\in \Loc{G}{\hat{S}}$ determines a decorated flag, $A$, associated to the puncture $p$. Without loss of generality, we can assume that $\pi(A)=P_\Theta$. Its decoration determines total potentials $w_i\in V_i$. For $\varsigma_i\in W(R)$, define $\varsigma_i\cdot \mathcal{L}$ by changing the decoration of $A$ to $A\,\oppInv{[\check{\beta}_i(\iota(w_i))]} $. 
\end{construction}

    \begin{lemma}
        The construction above is actually an action of the Weyl group.
    \end{lemma}
    \begin{proof}
        The fact that $\iota(a^{-1})(a)=a^{-1}$ implies that $\varsigma_i\cdot( \varsigma_i\cdot\mathcal{L})= \mathcal{L}$.
        The rest of the Weyl group relations can each be checked in the $A_2,B_2/C_2, G_2$ cases respectively.
    \end{proof}

    \begin{proposition}
        The action of $W(R)$ is realized by generalized mutations. 
    \end{proposition}     
    \begin{proof}
        Consider $\hat{S}$ a punctured disk with two marked points on its boundary and suppose that $G$ is Jordan split of type $A_1$. We consider the mutation sequence which switches the usual triangulation of this punctured bigon with a tagged triangulation. \cite[Example 5.15]{BerensteinEtAlNoncommutativeMarkedSurfaces2025-II},  computes the noncommutative cluster algebra in this situation, from which it is easy to see that the potential (called the total noncommutative angle there) is transformed to its inverse by this mutation sequence, which is the result of the change of decoration described above. 

        The $A_1$ case on any surface (with at least 2 marked points) can be reduced to this situation by moving to a triangulation which contains the desired marked point in a bigon. 

        Now consider $G$ Jordan split of general type. As before we can reduce to the punctured bigon case. Choose a cluster chart where the interpolating flags are placed on the sides away from the puncture. Let $m$ be the number of instances of the reflection $\varsigma_{\oppInv{i}}$ in these interpolating sequences. Consider the level $i$ subquiver or the subnetwork between strands $i$ and $ i+1$ which surrounds the puncture. This subquiver/subnetwork is exactly the one which is associated to a group of type $A_1$ on a punctured disk with $m$ marked points on the boundary. We use the mutation sequence applied only to this subquiver/network which changes the tagging at the central puncture. Since this only changes the cluster variables at level $i$, the $L$-distances must change as desired by the $A_1$ calculation.
    \end{proof}

    \subsubsection*{Braid Group action}
    
    Consider a surface, $\hat{S}$, with a boundary component $b$.  There is a birational action of the Artin-Tits braid group of type $R$, $B(R)$, on $\Loc{G}{\hat{S}}$ associated to $b$ as follows:

\begin{construction}
    Suppose that there are $d$ marked points on this boundary component and choose a particular one as a basepoint.
    A decorated local system $\mathcal{L}\in \Loc{G}{\hat{S}}$ determines a sequence of decorated flags $\{A_1,A_2,\dots,A_d,A_{d+1}=gA_1,\dots\}$ where $g\in G$ is the monodromy of $\mathcal{L}$ around the boundary. Write $\{B_i\}$ for the corresponding sequence of undecorated flags. Assume now that $d$ is even or the opposition involution is trivial.  For $\vartheta_i$ a braid generator we define the action $\vartheta_i\cdot\mathcal{L}$ in the following way:

    Pick $\boldsymbol{i}$ a reduced expression of $w_0$ which starts with $\varsigma_i$. We can find an interpolating sequence of flags from $B_i \to B_{i+1}$ following alternatively $\boldsymbol{i}$ and $\boldsymbol{i}^*$. Call the first new flags found along these interpolations $B_i'$ so that $(B_i,B_i')$ are $\varsigma_i$- or $\varsigma_{\oppInv{i}}$-transverse. Let $A_i'$ be the decorated flag over $B_i'$ with $l(A_i,A_i')=1$. We define the action of $\vartheta_i\cdot\mathcal{L}$ by transforming the flags $\{A_i\}$ to $\{A_i'\}$.

    If $d$ is odd and the opposition involution is nontrivial, we only define the action for the subgroup generated by positive braids $\omega\in B(\Theta)$ which are invariant under the opposition involution. In this case, the previous construction still works. 
\end{construction}

    \begin{proposition}
        The braid group action is realized by quasi-cluster isomorphisms.
    \end{proposition}

    \begin{proof}
        We follow the same proof as \cite[Section 14]{goncharov2019quantum}. The proof breaks into several cases depending on $\hat{S}$. In each case the mutation sequence realizing the braid symmetry is built using mutation sequences which change reduced expressions, flip edges, and the (left) reflection map, which is a monomial map. 
        
        To understand the quasi-cluster nature, we first see that if all of the $L$-distances around the boundary we are considering are $1$ then this action is a cluster automorphism on the nose. This is because the reflection map does not change any cluster variables in this case, and that the new flags $A_i'$ from the construction are actually the flags that would be found in an interpolating sequence between $A_i$ and $A_{i+1}$ which starts with $\varsigma_{i^{(*)}}$.

        Finally the situation when the boundary $L$-distances are not equal to 1 can be remedied by first rescaling the decorations so that they are 1 i.e. change $A_1$ to $A_1l(A_2,A_1)^{-1}$ and so on, then compute the mutation sequence, and put the scaling back. By \Cref{lem:rescaling_is_monomial_localsystems} this amounts to a monomial transformation in the frozen variables which come from the boundary $L$-distances. 
    \end{proof}

\subsection{Spectral description of types \texorpdfstring{$A_p$ or $C_p$}{Ap or Cp}}\label{sec:spectral_discription}

    Goncharov and Kontsevich give a \textit{spectral description} of  the space of twisted decorated local systems on a surface. 

    Each plabic network, $P$, encodes a topological surface called the spectral surface. 
\begin{definition}
    The \keyword{spectral surface} $\Sigma_P$ associated to $P$ is obtained as follows: we put a disk around each white and black vertices and connect them by a strip with a right-hand twist oriented from $\circ$ to $\bullet$. 
\end{definition}
\noindent The spectral surface is invariant under the local equivalence moves. 

\begin{definition}\label{conds:goncharov_konsevich}
    Goncharov and Kontsevich define a twisted version of the seed group $\quiverAlgebra_P$ (\Cref{def:NetworkGroup}) that we denote by $\bar{\quiverAlgebra}_P$.  The \keyword{twisted seed group} on a plabic graph $P$ is the quotient of the free group generated by edges of $P$ by the following relations:
\begin{itemize}
    \item The clockwise products around the the $\circ$ vertices is $-1$.
    \item The clockwise products around the the $\bullet$ vertices is $(-1)^{v(\bullet)-1}$, where $v(\bullet)$ is the valency of the vertex.
\end{itemize}
Denote by $\bar\A_P$ the noncommutative torus with group of functions given by $\bar\quiverAlgebra_P$.
\end{definition}
 The only difference between this torus and the torus $\A_P$ is in the sign of the product around each vertex, which is always $+1$ in the case of $\A_P$.

Now let $S$ be a marked surface and fix a triangulation $T$. The coordinates of Goncharov and Kontsevich are given by embedding the honeycomb network $P_m$ (\Cref{def:HoneycombNetwork}) in each triangle and amalgamating the twisted seed groups. 

\begin{definition}
    A local system is \keyword{fully twisted} if the monodromy around each contractable loop is $-1$.
\end{definition}

\begin{theorem}[\cite{goncharov2021spectral} Lemma 4.2,Theorem 1.6,Theorem 6.2]\leavevmode

\begin{enumerate}
    \item Let $\mathcal{R}$ be a noncommutative ring. The $\mathcal{R}^\times$ valued points of $\bar\quiverAlgebra_P$ parameterize fully twisted rank 1 $\mathcal{R}$-local systems on $\Sigma_P$. Equivalently, $\bar{\A}_P$ the moduli space of fully twisted rank 1 local systems on $\Sigma_P$.
    \item The space $\bar{\A}_{P_m}$ is birationally equivalent to the space of configurations of complete decorated flags in an $m$-dimensional vector space over $\mathcal{R}$.
    \item Let $\Sigma$ be a marked surface and $\Delta$ a triangulation of $\Sigma$. We get a plabic network on $\Sigma$ by inscribing $P_m$ into each triangle of $\Delta$, which we denote by $P_{\Delta,m}$. Then the space $\bar\A_{P_{\Delta, m}}$ is birationally equivalent to the moduli space of fully twisted rank $m$ local systems on $\Sigma$.
\end{enumerate}
    
\end{theorem}

For $G$ Jordan split of type $A_p$, we have $s_G=(-1)^{p}$. The group $\SL_m(\mathcal{R})$ is Jordan split of type $A_{m-1}$
When we turn our eyes back to our coordinates, we find 

\begin{corollary}\leavevmode\label{cor:spectraldiscription}
\begin{enumerate}
\item ${\A}_P$ is the moduli space of rank 1 (untwisted) local systems on $\Sigma_P$.
    \item The space ${\A}_{P_p}$ is birationally equivalent to the space of configurations of complete decorated flags in an $p+1$-dimensional vector space over $\mathcal{R}$.
    \item The space $\A_{P_{\Delta, p+1}}$ is birationally equivalent to the moduli space of ($s_G$) twisted rank $p+1$ local systems on $\Sigma$.
\end{enumerate}
\end{corollary}

These statements seem at odds with each other when $p$ itself is odd; how can the space of twisted and untwisted local systems on the spectral cover parameterize the same space of configurations of flags? The solution to this is that there is a (generally non-canonical) way to \textit{untwist} a local system. One isomorphism is given by taking a point of $\A_{P_{\Delta, p+1}}$ for $p$ odd and multiplying some of the edges by $-1$ so that the products around the black and white vertices satisfy \Cref{conds:goncharov_konsevich}. For example in \Cref{fig:GL4untwisting}, multiplying the dashed blue edges by $-1$ gives an untwisting of $p=3$. There are multiply ways to pick these edges so the untwisting is not canonical.

\input{Figures/figureGl4untwisting}

In Jordan split type $C_p$, our cluster charts are really charts from $A_{2p-1}$ which satisfy an extra symmetry. This symmetry corresponds to a symmetry of the spectral surface, and the local systems on the spectral surface which are invariant under this symmetry parametrize the appropriate spaces of configurations decorated flags and moduli of decorated local systems.

\subsection{Positivity}
    Now assume that $G$ has a $\Theta$-positive structure, i.e. each $\jordan{J}_i$ is formally real. 
    We have now endowed the subset $\LocT{G}{\hat{S}}{T}\subset\Loc{G}{\hat{S}}$ defined by an ideal triangulation with a coordinate chart. We use this to define a positive structure on $\Loc{G}{\hat{S}}$:
    \begin{definition}
        A decorated twisted local system $L\in\Loc{G}{\hat{S}}$ is \keyword{$\Theta$-positive} if it provides a positive $\{\jordan{J}_i\}$-point of $\ASpace{R}{\hat{S}}$. The set of these is denoted $\Loc{G}{\hat{S}}^{>0}$.
    \end{definition}
    \begin{theorem}
        The set $\Loc{G}{\hat{S}}^{>0}$ is well-defined, i.e. independent of all the choices made in the construction.
    \end{theorem}
    \begin{proof}
        This is a direct consequence of \Cref{thm:coordsASpaceSplit} and \Cref{thm:admissibleMutationPositivity}.
    \end{proof}

    $\Theta$-positive decorated twisted local systems may be understood in a more geometric fashion: Recall the Farey set for a surface $S$, $\farey(\hat{S})$ which can be realized as a subset of $\partial\Hyp\cong S^1$. The 2:1 cover $S^1\to S^1$ induces a double cover $\tilde{\farey}(\hat{S})\to\farey(\hat{S})$ with nontrivial automorphism $\alpha$. $\tilde{\farey}(\hat{S})$ is a cyclic set with an action of $\pi_1(\puncTB)$, defined by $s.x=\alpha(x)$ where $s\in\pi_1(\puncTB)$ denotes the generator of the central subgroup $\mathbb{Z}$ as before.
    \begin{definition}
        A \keyword{$\pi_1(\puncTB)$-equivariant configuration of flags} is a map $\eta \colon \tilde{\farey}(\hat{S}) \rightarrow \A$ such that for each $\gamma \in \pi_1(\puncTB)$ the configurations $\eta$ and $\eta \circ \gamma$ agree. Formally for each $\gamma$ there is a fixed $g_\gamma \in G$ such that for any $x \in \tilde{\farey}(\hat{S})$:
        \begin{equation*}
            \eta(\gamma. x) = g_\gamma. \eta(x) 
        \end{equation*}
    \end{definition}

        A $\pi_1(\puncTB)$-equivariant configuration of flags defines a twisted representation by $\varrho(\gamma) = g_\gamma$.

    \begin{theorem}\label{thm:boundaryMapSplit}
        The space $\Loc{G}{\hat{S}}$ parametrizes $\pi_1(\puncTB)$ equivariant configurations of flags.
    \end{theorem}
    We call the elements $\varrho$ and $\eta$ corresponding to $L\in\Loc{G}{\hat{S}}$ its \keyword{holonomy representation} and \keyword{boundary map}, respectively. The proof of this theorem is completely analogous to the split case, see \cite{fock2006moduli}.

    We can use this to give a different characterization of the space $\Loc{G}{\hat{S}}^{>0}$:
    \begin{theorem}[\cite{fock2006moduli}, Section 8.6]
        The elements of $\Loc{G}{\hat{S}}^{>0}$ are those with a positive boundary map.
    \end{theorem}
    Here, a map $\eta:\tilde{\mathcal{F}}(\hat{S})\to\decFlag$ is called \keyword{positive} if any \keyword{coherent} lift of a cyclically ordered tuple in $\mathcal{F}(\hat{S})$ to the double cover is mapped to a positive configuration of decorated flags. A coherent lift of a tuple $(x_1,\dots,x_k)$ is one obtained as follows: Take the path going through the points in $\mathcal{F}(\hat{S})$ in order, and lift it to a path in $\tilde{\mathcal{F}}(\hat{S})$ starting at a lift $\tilde{x}_1$ of $x_1$. This passes through exactly one preimage under the projection of every point in the tuple. These define the lift. The notion of positivity is independent of the chosen coherent lift.\medskip
    
    Finally, let us consider the holonomy representation of a positive local system. Recall the projection $p:\Loc{G}{\hat{S}}\to\mathcal{X}_{G',\hat{S}}$. The holonomy representation $\rho$ of the image of a decorated twisted local system under $p$ is a map $\rho:\pi_1(S)\to G'$.
    \begin{theorem}
        Let $\rho$ be the holonomy representation of a local system in $p\left(\Loc{G}{S}^{>0}\right)$. Let $\gamma$ be a simple closed curve on $S$ which is not homotopic to a puncture. Then $\rho(\gamma)$ is conjugated to an element of the $\Theta$-positive semigroup of $G'$. 
    \end{theorem}

    \begin{proof}
        Simply consider the two types of loops shown in \Cref{fig:loopsOnSurface}. We may compute their holonomy in the way described in \Cref{sec:holonomies}, i.e. by fixing any ideal triangulation and considering the corresponding graph.
        
        Since we assume that we are working with a positive local system, we see from \Cref{fig:localSystemSquareSplit} that following the graph counterclockwise around a vertex of the triangulation gives us an element of $U_\Theta^{>0}$. Following an edge of the triangulation on the other hand gives us an element of the form $w_0 l$ with $l\in L_\Theta^\circ$. Moreover we have for $u\in U_\Theta^{>0}$ that
        \begin{equation*}
            \oppInv{u}=w_0u^{-1}w_0\in U_\Theta^{\opp,>0}\,.
        \end{equation*}

        For the different kinds of loops, consider the setup as in \Cref{fig:localSystemSquareSplit} and \Cref{fig:pastingScheme}, respectively. By making all the choices as indicated, we see that $\rho(\gamma)$ lies in $U_\Theta^{\opp,>0}L_\Theta^\circ U_\Theta^{>0} L_\Theta^\circ = G^{>0}$ or $U_\Theta^{>0}L_\Theta^\circ U_\Theta^{\opp,>0} L_\Theta^\circ = G^{>0}$.
    \end{proof}

\bibliographystyle{alpha}
\bibliography{main.bib}

\end{document}